\documentclass[a4paper]{article}
\usepackage[frenchb]{babel}
\usepackage[latin1]{inputenc}
\usepackage{amsmath}
\usepackage{array}
\usepackage{amsthm}
\usepackage{amssymb}
\input xy
\xyoption{all}

\author{Aur\'elien DJAMENT\thanks{djament@math.univ-paris13.fr}\;\thanks{http://www.math.univ-paris13.fr/~djament/}\;\thanks{LAGA, Institut Galil\'ee,
universit\'e Paris 13,
99 avenue J.-B. Clément,
93430 VILLETANEUSE, FRANCE}}

\title{Cat\'egories de foncteurs en grassmanniennes}

\date{Octobre 2006 ; révisions mineures en novembre 2006}

\newcommand{\A}{{\mathcal{A}}}
\newcommand{\B}{{\mathcal{B}}}
\newcommand{\C}{{\mathcal{C}}}
\newcommand{\D}{{\mathcal{D}}}
\newcommand{\F}{{\mathcal{F}}}
\newcommand{\E}{{\mathcal{E}}}

\newcommand{\I}{{\mathcal{I}}}
\newcommand{\J}{{\mathcal{J}}}
\newcommand{\K}{{\mathcal{K}}}
\newcommand{\N}{{\mathcal{N}}}
\newcommand{\FF}{{\mathbb{F}_2}}
\newcommand{\col}{{\rm colim}\,}
\newcommand{\Gr}{{\mathcal{G}r}}
\newcommand{\Pl}{{\mathbf{Pl}}}
\newcommand{\T}{{\mathcal{T}}}
\newcommand{\agr}{{{\rm A}_\Gr}}
\newcommand{\agrr}{{\overline{{\rm A}}_\Gr}}
\newcommand{\kk}{{\Bbbk}}
\newcommand{\wt}{\widetilde}

\newcommand{\ptt}{{\,\widetilde{\otimes}\,}}

 \newcommand{\go}{{\guillemotleft~}}
\newcommand{\gf}{{~\guillemotright\,}}

\newtheorem{thm-intro}{Théorème}
\newtheorem{cor-intro}[thm-intro]{Corollaire}
\newtheorem{conj-intro}[thm-intro]{Conjecture}
\newtheorem{pr-intro}[thm-intro]{Proposition}

\newtheorem{theo}{Théorème}[section]
\newtheorem{pr}[theo]{Proposition}
\newtheorem{cor}[theo]{Corollaire}
\newtheorem{lm}[theo]{Lemme}

\newtheorem{prdef}[theo]{Proposition et définition}
\newtheorem{conj}[theo]{Conjecture}

\theoremstyle{definition}
\newtheorem{defi}[theo]{Définition}
\newtheorem{nota}[theo]{Notation}
\newtheorem{hyp}[theo]{Hypothèse}
\newtheorem{conv}[theo]{Convention}

\newtheorem{def-intro}[thm-intro]{Définition}

\theoremstyle{remark}
\newtheorem{rem}[theo]{Remarque}
\newtheorem{ex}[theo]{Exemple}

\begin{document}


\maketitle

\bigskip

\noindent
\textit{Résumé : }  Soit $\F$ la catégorie des foncteurs entre espaces
  vectoriels sur un corps fini. Les {\em catégories de foncteurs en
  grassmanniennes} sont obtenues en remplaçant la source de cette catégorie
  par la catégorie des couples formés d'un espace vectoriel et d'un
  sous-espace. Ces catégories possèdent une très riche structure algébrique
  ; nous étudions notamment leurs objets finis et leurs propriétés
  homologiques. Nous donnons des applications à la filtration de Krull
  de la catégorie $\F$ et à la $K$-théorie stable des corps finis.

\bigskip

\bigskip

\noindent
\textit{Abstract \textbf{(Grassmannian functor categories)} : } Let $\F$ be the category of functors between
  vector spaces over a finite field. The {\em grassmannian functor
    categories} are obtained by replacing the source of this category
  by the category of  pairs formed by a vector space and a
  subspace. These categories have a very rich algebraic structure ; we
  study in particular their finite objects and their homological
  properties. We give applications to the Krull filtration of the
  category $\F$ and to the stable $K$-theory of finite fields.

\bigskip

\bigskip

\noindent
\textit{Classification mathématique par sujets : } 18A25, 18G15. Secondaire : 16P60, 18A40, 18C15, 18D15,
  19D99, 20C33, 55S10.

\bigskip

\noindent
\textit{Mots clefs : } Catégories de foncteurs, algèbre homologique, groupes linéaires sur les corps
  finis, grassmanniennes, filtration de Krull, $K$-théorie stable,
  représen\-ta\-tions modulaires, (co)monades.

\pagebreak

\tableofcontents

\section*{Introduction}

Depuis une vingtaine d'années, de nombreux travaux
ont mis en évidence les liens féconds entre les catégories de
foncteurs, la topologie algébrique, les représentations modulaires des
groupes finis et plusieurs théories cohomologiques (cf. \cite{FFPS} et
le chapitre $13$ de \cite{Loday}).  Nous montrons
comment progresser dans cette voie à l'aide de nouvelles catégories de
foncteurs. Nous obtenons ainsi des résultats nouveaux sur des
catégories de foncteurs désormais classiques, notamment la catégorie
des foncteurs entre espaces vectoriels sur un corps fini. Nous
démontrons en particulier un théorème d'annulation
cohomologique très général, que l'on applique à la
$K$-théorie stable des corps finis. Les deux principes intuitifs
suivants guident notre travail : d'une part, la démonstration de propriétés
d'annulation cohomologique dans une catégorie de foncteurs est souvent
plus facile en transitant par une autre catégorie de foncteurs~\guillemotleft~plus grosse~\guillemotright~; d'autre part, la structure des objets d'une catégorie
de foncteurs peut souvent se ramener à celle
d'objets~\guillemotleft~plus petits~\guillemotright~(donc mieux compris) d'une autre catégorie
de foncteurs.


\medskip

Dans tout cet article, $\kk$ désigne un corps {\em fini} ; on note $\E_\kk$ la
catégorie des espaces vectoriels sur $\kk$, $\E^f_\kk$ la
sous-catégorie pleine des espaces de dimension finie et $\F(\kk)$ la
catégorie des foncteurs de $\E^f_\kk$ vers
$\E_\kk$ (la mention du corps $\kk$ sera souvent omise dans les notations, par la suite).  La catégorie $\F(\kk)$ est abélienne ; elle possède suffisamment
d'objets projectifs et injectifs. Des
calculs cohomologiques puissants ont été réalisés dans cette
catégorie, notamment avec les travaux de Franjou, Lannes et Schwartz
(\cite{FLS}) et de Franjou, Friedlander, Scorichenko et Suslin
(\cite{FFSS}). 

L'identification de la cohomologie de Mac Lane (définie dans
\cite{ML-coh}) comme un cas particulier de cohomologie fonctorielle
par Jibladze et Pirashvili (\cite{JP}) illustre  l'intérêt des
calculs cohomologiques dans la catégorie $\F(\kk)$. L'isomorphisme entre la $K$-théorie stable de $\kk$ et l'homologie
dans $\F(\kk)$ pour des systèmes de coefficients polynomiaux, établi
indépendamment par Betley (cf. \cite{Bet}) et Suslin (cf. appendice de
\cite{FFSS}), a fourni une autre motivation majeure à l'étude des
propriétés homologiques de $\F(\kk)$. Les articles \cite{PAQ},
\cite{PB} et \cite{PW} montrent comment aborder d'autres théories homologiques à l'aide de certaines catégories de
foncteurs.

\smallskip

L'étude systématique de la catégorie $\F$ est menée depuis le début
des années $1990$ (voir les articles de Kuhn \cite{K1}, \cite{K2} et \cite{K3}), à la suite des liens établis par Henn, Lannes et
Schwartz dans \cite{HLS} (voir aussi l'ouvrage \cite{LS} de Schwartz) entre les modules
instables sur l'algèbre de Steenrod et cette catégorie ; néanmoins, sa structure {\em globale} demeure fort mystérieuse. En effet, si la
compréhension de ses objets de
longueur finie peut se réduire à celle de  $\kk$-algèbres de dimension
finie, il en va bien différemment de ses objets de longueur infinie,
dont l'étude se heurte à des problèmes profonds de compréhension
{\em générique} des représentations modulaires (i.e. de compréhension
des liens, notamment cohomologiques, entre les représentations de plusieurs groupes) --- la catégorie
$\F(\kk)$ a d'ailleurs été nommée catégorie des  représentations
génériques des groupes linéaires sur $\kk$ par Kuhn, qui a justifié
cette terminologie dans les trois articles susmentionnés.

Ainsi, la conjecture suivante, proposée par Lannes et Schwartz, s'est révélée l'un des problèmes les
plus difficiles à résoudre dans cette catégorie.

\begin{conj-intro}[Conjecture artinienne]\label{canu} La catégorie $\F(\kk)$ est
  localement noethé\-rien\-ne.
\end{conj-intro}

Les raisons de la dénomination paradoxale de cette conjecture et une brève discussion de celle-ci sont données à
    la fin de la section~\ref{s-rf}.

\medskip

Explicitons les catégories de foncteurs qui permettent de progresser
dans l'étude de la conjecture artinienne. On note $\Gr(V)$ la grassmannienne des
sous-espaces vectoriels de $V$.  Nous regarderons $\Gr$ comme un
foncteur de $\E^f_\kk$ vers la catégorie des ensembles.

\begin{def-intro} Soit $\E^f_\Gr(\kk)$ la
catégorie des couples $(V,W)$, où $V$ est un $\kk$-espace vectoriel de
dimension finie et $W$ un élément de  $\Gr(V)$, et dont les morphismes
$(V,W)\to (V',W')$ sont les applications linéaires $f : V\to V'$
telles que $f(W)=W'$. La {\em catégorie de
  foncteurs en grassmanniennes} $\F_\Gr(\kk)$  est la catégorie des foncteurs de $\E^f_\Gr(\kk)$
vers $\E_\kk$.
\end{def-intro}

Il est également naturel de s'intéresser à la catégorie, notée
$\wt{\E}^f_\Gr(\kk)$, qui a les mêmes objets que $\E^f_\Gr(\kk)$, et dont les
morphismes $(V,W)\to (V',W')$ sont les applications linéaires $f : V\to V'$
telles que $f(W)\subset W'$. La catégorie des foncteurs de
$\wt{\E}^f_\Gr(\kk)$ vers $\E$, notée $\wt{\F}_\Gr(\kk)$, est
également considérée dans cet article ; elle constitue un
adjuvant pour l'étude de la catégorie $\F_\Gr(\kk)$.

Pour tout entier positif $n$, on note $\F_{\Gr,n}$ la
sous-catégorie pleine de $\F_\Gr$ (rappelons que la mention du corps
$\kk$ est désormais sous-entendue) des foncteurs $X$ tels que
$X(V,W)\neq 0$ si $\dim W\neq n$. Cette catégorie peut également se
décrire comme une catégorie de foncteurs. L'étude des catégories $\F_\Gr$
et $\F_{\Gr,n}$ constitue
le principal sujet de cet article ; le lien entre $\F_\Gr$ et les
$\F_{\Gr,n}$ provient de ce qu'il existe une stratification de
$\F_\Gr$ par des sous-catégories épaisses  $\F_{\Gr,\leq n}$, stratification dont les sous-quotients sont
équivalents aux catégories $\F_{\Gr,n}$.

On peut également ramener l'étude des catégories $\F_{\Gr,n}$ à celle
de catégories plus simples : le groupe linéaire $GL_n(\kk)$ intervient
naturellement dans la structure de $\F_{\Gr,n}$, qui constitue une
sorte de produit semi-direct tordu entre la catégorie des
$\kk[GL_n(\kk)]$-modules et la catégorie de foncteurs $\F_{\Pl,n}$
définie comme suit. On note $\E^f_{\Pl,n}$ la
catégorie des objets $V$ de $\E^f$ munis d'un monomorphisme
$\kk^n\hookrightarrow V$. La catégorie $\F_{\Pl,n}$ est la catégorie des foncteurs
de $\E^f_{\Pl,n}$ vers $\E$. Outre ses liens avec $\F_{\Gr,n}$,
cette catégorie possède un intérêt intrinsèque, car elle est
équivalente à la catégorie des comodules de la catégorie~$\F$ sur
le foncteur injectif standard associé à l'espace vectoriel~$\kk^n$.

Mentionnons une autre catégorie de foncteurs considérée dans cet
article. Soit $\E^f_{surj}$ la sous-catégorie de $\E^f$ ayant les
mêmes objets et dont les morphismes sont les épimorphismes de
$\E^f$. La catégorie, notée $\F_{surj}$, des foncteurs de
$\E^f_{surj}$ vers $\E$ joue un rôle important dans l'étude de la
catégorie $\F_\Gr$. En effet, si $f : (V,W)\to (V',W')$ est un
morphisme de $\E^f_\Gr$, alors $f$ induit un épimorphisme de $W$ sur
$W'$ : on obtient ainsi un foncteur $\E^f_\Gr\to\E^f_{surj}$, donné sur les objets par $(V,W)\mapsto W$,
puis par précomposition un foncteur fondamental $\F_{surj}\to\F_\Gr$.
Une catégorie analogue à $\F_{surj}$ s'obtient à partir des monomorphismes de $\E^f$
; on la note $\F_{inj}$. Elle possède des liens étroits avec les
systèmes de coefficients introduits par Dwyer dans \cite{Dw} et joue
un rôle essentiel dans la comparaison entre la $K$-théorie stable de
$\kk$ et les groupes d'extension dans $\F(\kk)$ (théorème de
Betley-Suslin susmentionné).

\smallskip

Revenons à la catégorie $\F_\Gr$. Le foncteur d'oubli $\E^f_\Gr\to\E^f$ donné par $(V,W)\mapsto
V$ induit par précomposition un foncteur $\iota : \F\to\F_\Gr$. Il est
adjoint à droite au foncteur  d'{\em intégrale en grassmanniennes} $\omega :
\F_\Gr\to\F$  donné par
$$\omega(X)(V)=\bigoplus_{W\in\Gr(V)}X(V,W).$$
Ce foncteur constitue l'outil le plus puissant pour relier les catégories
$\F_\Gr$ et $\F$. Son importance  est d'abord illustrée par le résultat formel
suivant, dans lequel le foncteur $\omega(\kk)=\kk[\Gr]$ de $\F$ est muni de la
structure comultiplicative déduite du fait que ce foncteur est la
linéarisation d'un foncteur ensembliste.

\begin{pr-intro} Le foncteur $\omega:
\F_\Gr\to\F$ induit une équivalence entre $\F_\Gr$ et la catégorie des
$\kk[\Gr]$-comodules de $\F$.
\end{pr-intro}

On peut donner une description similaire des catégories $\F_{\Gr,n}$ en
termes de comodules ; celle-ci montre, par dualité, que $\F_{\Gr,n}$ est étroitement
liée à la catégorie des $\bar{D}(n)$-modules de $\F$, où
$\bar{D}(n)$ est un foncteur qui a rapport à l'algèbre de Dickson ---
cf. l'article \cite{GP2} de Powell, qui a annoncé l'importance de
cette catégorie de modules.

\smallskip

L'étude de la structure élémentaire de la catégorie $\F_\Gr$ repose
principalement sur le {\em foncteur différence}, analogue à
l'endofoncteur du même nom dans $\F$. Ce foncteur, noté $\Delta^\Gr$, est
donné par le scindement canonique
$$X(V\oplus\kk,W)\simeq X(V,W)\oplus\Delta^\Gr(X)(V,W).$$

Rappelons que le foncteur différence $\Delta$ de $\F$ est quant à lui caractérisé
par le scindement canonique $F(V\oplus\kk)\simeq
F(V)\oplus\Delta(F)(V)$. Comme dans la catégorie $\F$, on introduit la
définition suivante :

\begin{def-intro} Un objet $X$ de $\F_\Gr$ est dit {\em polynomial} s'il
  existe un entier $n$ tel que $(\Delta^\Gr)^n(X)=0$.
\end{def-intro}

Dans ce qui suit, nous nommons {\em fini} un objet de longueur finie
d'une catégorie abélienne. Un objet {\em localement fini} est un objet
colimite d'objets finis.

\begin{pr-intro}\label{pri-ffgr} Les foncteurs finis de la catégorie
  $\F_\Gr$ sont polynomiaux.
\end{pr-intro}

Ce résultat permet de décrire explicitement les objets simples de la
catégorie $\F_\Gr$ à partir des objets simples de $\F$ et des
représentations simples des groupes linéaires.

Nombre des foncteurs entre la catégorie $\F_\Gr$ et la catégorie $\F$,
ou d'autres qui lui sont étroitement reliées, possèdent de bonnes
propriétés de commutation au foncteur différence. Par exemple, il existe un isomorphisme canonique
$\iota\circ\Delta\simeq\Delta^\Gr\circ\iota$ de foncteurs
$\F\to\F_\Gr$. En revanche, le
foncteur composé $\Delta\circ\omega : \F_\Gr\to\F$ diffère du
foncteur $\omega\circ\Delta^\Gr$, qu'il contient comme facteur
direct. Ainsi, l'image par le foncteur d'intégrale en grassmanniennes
d'un objet fini de la catégorie $\F_\Gr$ n'est généralement pas un
objet analytique (i.e. colimite d'objets polynomiaux) de la catégorie
$\F$. D'ailleurs, tous les objets projectifs de type fini de $\F$ sont
l'image par le foncteur $\omega$ de foncteurs finis de $\F_\Gr$. Le
principe général d'étude de la catégorie $\F$ à partir de la catégorie
$\F_\Gr$ consiste à ramener l'étude des objets de type fini de $\F$ à
celle des objets finis de $\F_\Gr$ grâce au foncteur $\omega$. La mise
en \oe uvre de ce principe s'avère ardue ; en effet, les objets de
type fini de la catégorie $\F$ ne sont pas tous isomorphes à l'image par $\omega$
d'un objet fini de $\F_\Gr$. Un énoncé précis  que tous les
foncteurs de type fini connus dans la catégorie $\F$ vérifient sera
discuté dans cet article ; la conjecture~\ref{canuf} présentée
ci-dessous en constitue une variante.

\smallskip

Venons-en maintenant à des propriétés  profondes des catégories de
foncteurs en grassmanniennes. Notre résultat d'annulation
cohomologique principal est le suivant, dans lequel $\I$ désigne
l'endofoncteur de $\F_\Gr$ donné par
$\I(X)(V,W)=\underset{B\in\Gr(W)}{\bigoplus}X(V,B)$. Ce foncteur
conserve les objets localement finis.

\begin{thm-intro}\label{thip} Soient $X$ et $Y$ des objets de
  $\F_\Gr$, $X$ étant supposé localement fini. Il existe un
  isomorphisme gradué naturel  ${\rm
    Ext}^*_\F(\omega(X),\omega(Y))\simeq {\rm Ext}^*_{\F_\Gr}(X,\I(Y))$.
\end{thm-intro}

Dans le corollaire suivant, le foncteur $\omega_n$ désigne la
restriction du foncteur $\omega$ à la sous-catégorie $\F_{\Gr,n}$ de $\F_\Gr$.

\begin{cor-intro}\label{fondcr-intr} Soient $k$ et $n$ deux entiers naturels, $X$
  un objet localement fini de $\F_{\Gr,k}$ et $Y$ un objet de $\F_{\Gr,n}$.
\begin{enumerate}\item Si $k<n$, alors ${\rm
    Ext}^*_{\mathcal{F}}(\omega_k (X),\omega_n(Y))=0$.
\item Si $k=n$, alors le morphisme naturel ${\rm
    Ext}^*_{\F_{\Gr,n}}(X,Y)\to {\rm
    Ext}^*_{\mathcal{F}}(\omega_n (X),\omega_n(Y))$ induit par $\omega_n$ est un isomorphisme.
\end{enumerate}
\end{cor-intro}

Le théorème~\ref{thip} et le corollaire~\ref{fondcr-intr} permettent
d'une part de généraliser le théorème de Betley-Suslin sur la
$K$-théorie stable de $\kk$, d'autre part de mener de nombreux calculs
cohomologiques dans la catégorie $\F$ utiles pour comprendre sa
structure, comme le théorème de l'appendice de \cite{GP5}. 

\smallskip

Un moyen efficace d'appréhender la structure globale d'une catégorie
abélienne $\A$ consiste à étudier sa filtration de Krull, notée
$(\K_n(\A))$ (cf. paragraphe~\ref{par-krull}). Le
corollaire~\ref{fondcr-intr}, qui indique une
hiérarchie dans la~\guillemotleft~taille~\guillemotright~des objets du
type $\omega_n(X)$ selon la valeur de~$n$, où
$X$ est un objet localement fini de $\F_{\Gr,n}$, suggère la description conjecturale suivante de la filtration de Krull de la
catégorie $\F$. Celle-ci précise considérablement la
conjecture~\ref{canu} et constitue une motivation essentielle à l'étude
des catégories de foncteurs en grassmanniennes. 

\begin{conj-intro}[Conjecture artinienne extrêmement
  forte]\label{canuf} Pour tout entier positif $n$, le foncteur
$\omega_n$ induit une équivalence entre la sous-catégorie pleine
$\F^{lf}_{\Gr,n}$ des objets localement finis de $\F_{\Gr,n}$ et la catégorie
quotient $\K_n(\F)/\K_{n-1}(\F)$.
\end{conj-intro}

Non seulement cet énoncé entraîne toutes les formes renforcées de la
conjecture artinienne émises jusqu'à présent, mais il implique aussi
un grand nombre d'autres résultats profonds sur la catégorie $\F$. Dans l'article
\cite{art3}, nous poursuivrons l'étude du foncteur $\omega$ et
montrerons notamment une forme affaiblie de la
conjecture~\ref{canuf}, avec des applications à de nouveaux cas de la
conjecture artinienne. La plupart des résultats du présent article et
de \cite{art3} sont déjà exposés dans la thèse de doctorat de l'auteur
(\cite{these}) lorsque $\kk$ est le corps $\FF$ à
deux éléments.

Pour établir la conjecture artinienne extrêmement
  forte, nous pensons qu'il sera nécessaire de combiner les
catégories de foncteurs en grassmanniennes avec de nouveaux outils
liés à la théorie des représentations.

\medskip

Les méthodes de cet article peuvent  se
généraliser : la définition,
les constructions et propriétés de base de la catégorie $\F_\Gr(\kk)$ se transposent sans changement
si l'on remplace la catégorie source $\E^f_\kk$ de $\F(\kk)$ par une
catégorie abélienne essentiellement petite $\A$ dans laquelle les
ensembles de morphismes et de sous-objets sont finis. En effet, dans
cette situation, l'adjonction entre les foncteurs $\iota :
\F\to\F_\Gr$ et $\omega$ s'étend aussitôt, et l'on dispose de
foncteurs différences. Des généralisations dans un cadre non abélien sont également envisageables.

\paragraph*{Organisation de l'article} La première
 section rappelle les rudiments nécessaires sur la catégorie
$\F$ et la conjecture artinienne. La deuxième traite des
catégories $\F_{surj}$ et $\F_{inj}$, à la fois d'un point de vue
intrinsèque et d'un point de vue préliminaire à l'étude des catégories de foncteurs en
grassmanniennes et de la $K$-théorie stable de~$\kk$.

Les sections \ref{sctccf} et \ref{sctcatb} s'attachent à des
constructions catégoriques nécessaires à la deuxième partie. Celle-ci
introduit les catégories de
foncteurs en grassmanniennes $\F_\Gr$, $\F_{\Gr,n}$, $\F_{\Pl,n}$ et
$\wt{\F}_\Gr$ et donne leurs propriétés de base : outre leur description fonctorielle,
les trois premières d'entre elles sont identifiées comme catégories de comodules, et 
traitées de manière monadique, à l'aide du théorème de Beck. Les
objets finis de la catégorie $\F_\Gr(\kk)$, notamment, sont étudiés en
détails. On établit ainsi la
proposition~\ref{pri-ffgr}, dont on tire les conséquences, ainsi que les propriétés d'adjonction élémentaires entre les
très nombreux foncteurs définis, qui font la richesse de la structure
de ces catégories.

La dernière partie expose des applications des constructions des
parties précédentes, à l'aide du foncteur d'intégrale en
grassmanniennes $\omega : \F_\Gr(\kk)\to\F(\kk)$, dont on démontre,
dans la section~\ref{sscep}, la
propriété d'annulation cohomologique fondamentale, le
théorème~\ref{thip}, ainsi que le corollaire~\ref{fondcr-intr}. Les applications à la catégorie $\F(\kk)$,
données dans la section~\ref{s-fkf}, seront approfondies dans
\cite{art3} ; la section~\ref{s-rhid}, qui tire des résultats de
la section~\ref{sscep} des conséquences
apparemment internes aux catégories de foncteurs en grassmanniennes, se révèle également utile dans $\F(\kk)$
puisqu'elle sous-tend la démarche de \cite{art1}. Enfin, la
section~\ref{parfij} est consacrée à une propriété d'annulation cohomologique des systèmes de
coefficients, déduite de celle du foncteur $\omega$, et qu'on applique à la $K$-théorie
stable.

Les appendices fixent nos notations et rappellent quelques résultats
connus sur trois leitmotive de cet article : les adjonctions (à la base
de plusieurs descriptions des catégories de foncteurs en
grassmanniennes comme du théorème d'annulation cohomologique
principal), les propriétés de finitude des catégories abéliennes (dans
l'optique de la conjecture artinienne) et les catégories de foncteurs.

\paragraph*{Remerciements} L'auteur témoigne sa chaleureuse
reconnaissance à Geoffrey Powell pour l'attention qu'il a portée à ce
travail durant toutes les étapes de sa réalisation. Il remercie
également Lionel Schwartz et Christine Vespa pour leurs conseils.

\subsection*{Notations et conventions} 
\begin{enumerate}\item {\bf fondamentales :}

\begin{enumerate}\item Nous noterons $\Gr(V)$ l'ensemble des sous-espaces vectoriels
  d'un espace vectoriel $V$. Si $I$ est une partie de $\mathbb{N}$,
  nous noterons $\Gr_I(V)$ le sous-ensemble de $\Gr(V)$ constitué des
  sous-espaces dont la dimension appartient à $I$.
%
\item Si $\C$ et $\D$ sont deux catégories, $\C$ étant {\em
    essentiellement petite}, nous noterons $\mathbf{Fct}(\C,\D)$ la
  catégorie des foncteurs de $\C$ vers $\D$.
\end{enumerate}
\item {\bf couramment utilisées :}

\begin{enumerate}
\item La caractéristique du corps fini $\kk$ sera notée $p$ ; $q=p^d$
  désignera son cardinal. Lorsqu'aucune confusion ne pourra en résulter, nous omettrons
  toute mention du corps $\kk$ dans les notations. 
\item Nous noterons $\bf{Ens}$ la catégorie des ensembles, et
  $\mathbf{Ens}^f$ la sous-catégorie pleine des ensembles finis.
\item  Soit $E$ un ensemble. Nous noterons $\kk[E]$
  le $\kk$-espace vectoriel  somme directe de copies de $\kk$ indexées par $E$. On peut voir
  l'association $E\mapsto\kk[E]$ comme un foncteur de $\bf{Ens}$ vers $\E_\kk$.

Nous noterons
  $[e]$ l'élément de la base canonique de $\kk[E]$
  associé à un élément $e$ de $E$.
\item Nous désignerons par $\mathbf{Mod}_A$ la catégorie des modules à
  droite sur un anneau $A$ et $_A\mathbf{Mod}$ la catégorie des
  $A$-modules à gauche. Nous adopterons plus généralement ces
  notations lorsque $A$ est un objet d'une catégorie monoïdale
  symétrique muni d'une structure d'algèbre.

Dans le cas d'un objet $C$ muni d'une structure de coalgèbre, nous
noterons $\mathbf{Comod}_C$ la catégorie des $C$-comodules à droite. 

Dans le cas d'une (co)algèbre (co)unitaire, les morphismes seront
toujours censés préserver la (co)unité.
\item Nous noterons ${\rm Ob}\,\C$ la classe des objets d'une
  catégorie $\C$. Si $X$ et $Y$ sont deux objets de $\C$, on
  note :
\begin{itemize}\item  ${\rm hom}_\C (X,Y)$ l'ensemble\,\footnote{Dans toutes les catégories que nous
    considérerons, la classe des morphismes entre deux objets sera un ensemble.} des
  morphismes de $X$ dans $Y$ ;
\item ${\rm End}_\C(X)$ le monoïde ${\rm
    hom}_{\C}(X,X)$ des endomorphismes de $X$ ;
\item  ${\rm Aut}_\C(X)$ le groupe des automophismes
  de $X$ ;
\item ${\rm Pl}_\C(X,Y)$ l'ensemble des monomorphismes de $X$ dans
  $Y$ ;
\item ${\rm Epi}_\C(X,Y)$ l'ensemble des épimorphismes de $X$ vers
  $Y$ ;
\item ${\rm Iso}_\C(X,Y)$ l'ensemble des isomorphismes de $X$ vers $Y$.
\end{itemize}
L'indice $\C$ sera omis quand aucune confusion ne peut
en résulter.

Enfin, $\C^{op}$ désignera la catégorie opposée de $\C$.
\item On désigne par $\mathbb{N}$ l'ensemble des entiers positifs ou
  nuls, et par $\mathbb{N}^*$ l'ensemble des entiers strictement positifs.

On adoptera également les abréviations suivantes  :
\begin{itemize}
\item $n=\{n\}$, 
\item $\leq n=\{i\in\mathbb{N}\,|\,i\leq n\}$, 
\item $\geq n=\{i\in\mathbb{N}\,|\,i\geq n\}$.
\end{itemize}
\end{enumerate}
\item {\bf plus secondaires :}

\begin{enumerate}\item Nous noterons $V^*$ le dual d'un espace vectoriel $V$, et
  $W^\perp$ l'orthogonal dans $V^*$ d'un sous-espace $W$ de~$V$.
\item\label{notmo-i} Soit $M$ un monoïde.
\begin{enumerate}\item Nous noterons $\kk[M]$ l'algèbre de $M$
  sur $\kk$.
\item Nous désignerons par $\underline{M}$
la catégorie à un seul objet de monoïde d'endomorphismes $M$.
\item  Pour toute catégorie $\C$, nous noterons $\C_M$ la catégorie
  $\mathbf{Fct}(\underline{M},\C)$. C'est la catégorie des {\em objets
    de $\C$ munis d'une action de $M$}. En effet, les objets de $\C_M$
  sont les objets $X$ de $\C$ munis d'un morphisme de monoïdes $M\to
  {\rm End}_\C (X)$. Nous noterons $O^\C_{M} : \C_{M}\to\C$ le
  foncteur d'oubli.
\end{enumerate}

La catégorie $(\E_\kk)_M$ est ainsi équivalente à $_{\kk[M]}\mathbf{Mod}$.
\item  Une catégorie {\em $\kk$-linéaire} est une catégorie $\I$ telle
  que, pour tous objets $A$ et $B$ de $\I$, les ensembles ${\rm hom}_\I(A,B)$ sont munis d'une
structure de $\kk$-espace vectoriel, de sorte que la composition des
morphismes soit $\kk$-bilinéaire. On ne suppose pas $\I$ additive.
\end{enumerate}
\end{enumerate}

\textit{D'autres notations utilisées dans tout l'article sont
  introduites dans les appendices.}



\part{Préliminaires}\label{p-prlm}

Après avoir rappelé les propriétés de base de la catégorie~$\F$, dont
l'étude constitue l'une des principales motivations de cet article,
nous traitons d'une catégorie analogue, notée $\F_{surj}$, avec
plus de détail. Cette catégorie s'avère essentielle autant pour la
considération de la catégorie $\F$ que pour celle de la catégorie de
foncteurs en grassmanniennes $\F_\Gr$ que nous introduirons dans la partie~\ref{p-dfg} ; elle possède en
outre un intérêt intrinsèque.

Les sections~\ref{sctccf} et~\ref{sctcatb} constituent le
soubassement des catégories de foncteurs en grassmanniennes,
notamment $\F_\Gr$ : la première sous-tend sa description comme catégorie de
comodules, la seconde sa description fonctorielle.

\section{Rappels sur la catégorie $\F$}\label{s-rf}

Les résultats rappelés dans l'appendice~\ref{apfct} s'appliquent à la
catégorie $\F(\kk)$ ; l'hypothèse~\ref{hypf3} est vérifiée car
le corps $\kk$ est \textit{fini}. Le foncteur projectif standard $P^{\E^f_\kk}_V$,
où $V$ est un $\kk$-espace vectoriel de dimension finie, sera
simplement noté $P_V$ ; de même, l'injectif standard
$I^{\E^f_\kk}_V$ sera noté $I_V$.

La plupart des résultats de cette section sont contenus dans \cite{HLS}, \cite{K1} ou~\cite{LS},
par exemple. 

Comme le foncteur de dualité $(\E^f_\kk)^{op}\to\E^f_\kk$ est une équivalence
de catégories, la proposition/définition~\ref{prdcff} procure un
foncteur $D : \F(\kk)^{op}\to\F(\kk)$, donné sur les objets par
$(DF)(V)=F(V^*)^*$. Il induit une équivalence de catégories entre
$(\F^{df}(\kk))^{op}$ et $\F^{df}(\kk)$, où l'on note $\F^{df}(\kk)$ la
sous-catégorie épaisse de $\F(\kk)$ constituée des foncteurs prenant
des valeurs de dimension finie.

\begin{defi}\label{dfad} Un objet $F$ de $\F$ est dit {\em auto-dual}
  s'il existe un isomorphisme $u : F\xrightarrow{\simeq}DF$ invariant par l'isomorphisme d'adjonction ${\rm hom}_\F
  (F,DF)\xrightarrow{\simeq}{\rm hom}_\F (F,DF)$. 
\end{defi}

\begin{ex} Les objets simples de $\F$ sont auto-duaux (cf.~\cite{K2}).
\end{ex}

\paragraph*{Décomposition scalaire et tors de Frobenius} Suivant la notation~\ref{not-ds}, nous poserons
$\F_i(\kk)=\{F\in {\rm Ob}\,\F\,|\,\forall\lambda\in\kk\quad
F(\lambda.id)=\lambda^i.id\}$. La catégorie $\F_0(\kk)$ est
la sous-catégorie des foncteurs constants de $\F(\kk)$, elle est
canoniquement équivalente à $\E_\kk$ ; nous {\em identifierons} par la suite ces
deux catégories. La décomposition scalaire d'un foncteur $F$ de
$\F(\kk)$ (proposition/définition~\ref{pra-ds}) s'écrit donc
$F\simeq F(0)\oplus F_1\oplus\dots\oplus F_{q-1}$. On notera
$\bar{F}=F_1\oplus\dots\oplus F_{q-1}$.

Si $V$ est un objet de $\E^f_\kk$, nous noterons $P_{V,i}$ pour
$(P_V)_i$. Pour $V=\kk$, on obtient la décomposition de $P_\kk$ en
somme directe de projectifs indécomposables. Explicitement,
$P_{\kk,0}$ est le foncteur constant $\kk$, et pour $1\leq i\leq q-1$,
on peut voir $P_{\kk,i}(V)$ comme le quotient de $P_\kk(V)$ par
l'espace vectoriel engendré par les $[\lambda v]-\lambda^i[v]$, pour
$\lambda\in\kk$ et $v\in V$.

Soit $\phi$ l'automorphisme de la catégorie
$\E^f_\kk$ obtenu en tordant l'action des scalaires par
l'automorphisme de Frobenius de $\kk$. Le foncteur de précomposition
$\phi^* : \F(\kk)\to\F(\kk)$ est appelé  {\em tors de Frobenius}. Il
induit, pour tout entier naturel $i$, une équivalence de catégories entre $\F_i(\kk)$ et $\F_{pi}(\kk)$.


\paragraph*{Changement de corps} Soit $K$ une extension finie de
$\kk$. Les foncteurs d'extension des scalaires $t : \E_\kk\to\E_K$ et
de restriction des scalaires $\tau : \E_K\to\E_\kk$ sont mutuellement
adjoints (i.e. $t$ est adjoint à droite et à gauche à $\tau$), de même que leurs
variantes entre $\E_\kk^f$ et $\E^f_K$, que nous noterons encore $t$
et $\tau$ par abus. On en déduit que les foncteurs $\tau^* t_* :
\F(\kk)\to\F(K)$ et $\tau_* t^* : \F(K)\to\F(\kk)$ sont mutuellement
adjoints. On rappelle que l'indice (resp. exposant) étoilé indique la
postcomposition (resp. précomposition) ---
cf. appendice~\ref{apfct}. Pour plus de détails à ce sujet, voir~\cite{FFSS}, §\,3.

Le foncteur $\tau_* t^* : \F(K)\to\F(\kk)$ est appelé encore {\em
  restriction des scalaires}, et $\tau^* t_* :\F(\kk)\to\F(K)$ {\em induction}.
%

\paragraph*{Le foncteur différence et les objets polynomiaux de
  $\F(\kk)$} Ces notions dépendent exclusivement du caractère additif
de la catégorie source $\E^f_\kk$. Si $V$ est un objet de cette
catégorie, le {\em foncteur de décalage} par $V$ est l'endofoncteur de $\F(\kk)$,
noté $\Delta_V$, de précomposition par $\cdot\oplus V :
\E^f_\kk\to\E^f_\kk$. Cette construction est fonctorielle en $V$ ;
ainsi, du fait que la composition $0\to\kk\to 0$ est l'identité, et
que $\Delta_0\simeq id_\F$, on obtient un scindement $\Delta_\kk\simeq
id\oplus\Delta$, où $\Delta$ est appelé {\em foncteur différence} de $\F(\kk)$.

Le résultat classique suivant (cf. \cite{GP4}, appendice) est une application des
propositions~\ref{prdffdec} et~\ref{prdffdec2}.

\begin{pr}\label{adj-fctd}\begin{enumerate}\item\begin{enumerate}\item Le foncteur de
    décalage $\Delta_V$ est adjoint à droite à $\cdot\otimes P_V$ ;
    cette adjonction est naturelle en $V$.
\item Le foncteur différence $\Delta$ est adjoint à droite à $\cdot\otimes\bar{P}_\kk$.
\end{enumerate}
\item \begin{enumerate}\item Le foncteur de
    décalage $\Delta_V$ est adjoint à gauche à $\cdot\otimes I_V$ ;
    cette adjonction est naturelle en $V$.
\item Le foncteur différence $\Delta$ est adjoint à gauche à $\cdot\otimes\bar{I}_\kk$.
\end{enumerate}
\end{enumerate}
\end{pr}

\begin{nota}\label{not-cfn} Nous désignerons par $\F^n(\kk)$ la
  sous-catégorie épaisse de $\F(\kk)$ noyau du foncteur exact
  $\Delta^{n+1}$, où $n\in\mathbb{N}$. Pour $n$ entier strictement négatif ou
  $-\infty$, on convient que $\F^n(\kk)$ est la
  sous-catégorie réduite à $0$.
%
\end{nota}

\begin{defi}\label{defi-cfn} Soit $F$ un foncteur de $\F(\kk)$.
\begin{enumerate}\item On dit que $F$ est {\em polynomial} s'il existe
  $n$ tel que $F\in {\rm Ob}\,\F^n(\kk)$. Le plus petit $n$ ayant
  cette propriété s'appelle alors le {\em degré} de $F$, et se note
  $\deg F$.
\item On dit que $F$ est {\em analytique} s'il est colimite de ses
  sous-foncteurs polynomiaux. On note $\F_\omega(\kk)$ la
  sous-catégorie pleine de $\F(\kk)$ des foncteurs analytiques.
\item On dit que $F$ est {\em coanalytique} s'il est limite de ses
  quotients polynomiaux.
\end{enumerate} 
\end{defi}

La proposition suivante se trouve démontrée dans \cite{K1} (§\,$4$), par
exemple. Nous détaillons la démonstration du corollaire qui s'en
déduit afin d'illustrer des raisonnements utiles dans d'autres
catégories de foncteurs.

\begin{pr}\label{prf-fpol}\begin{enumerate}\item Les foncteurs
    projectifs standard de $\F(\kk)$ sont
    coanalytiques.
\item Les foncteurs injectifs standard de
  $\F(\kk)$ sont analytiques.
\end{enumerate}
\end{pr}

On rappelle (cf. §\,\ref{par-anfd}) qu'un objet pf$_\infty$
(resp. co-pf$_\infty$) est un objet qui admet une résolution
projective (resp. injective) de type fini.  

\begin{cor}\label{crf-fpolf}\begin{enumerate}\item Un foncteur de
    $\F(\kk)$ est fini si et seulement s'il est polynomial et à
    valeurs de dimension finie.
\item La sous-catégorie pleine $\F^f(\kk)$ des objets finis de
  $\F(\kk)$ est préservée par le produit tensoriel et par le foncteur
  différence.
\item\label{aserls} Les objets finis de $\F(\kk)$ sont
  pf$_\infty$ et co-pf$_\infty$.
\item La sous-catégorie $\F_\omega(\kk)$ de $\F(\kk)$ est épaisse.
\item Un objet de $\F(\kk)$ est de co-type fini si et seulement s'il est
  analytique et de socle fini.
\end{enumerate}
\end{cor}

\begin{proof} Le foncteur $(\Delta,ev_0) : \F\to\F\times\E$ (on
  rappelle que $ev_0$ désigne le foncteur d'évaluation en $0$ ---
  cf. notation~\ref{postpre}) est exact et fidèle, car le noyau de
  $\Delta$ est constitué des foncteurs constants. La
  proposition~\ref{pr-ff} montre alors qu'un foncteur $F$ tel que
  $\Delta F$ est fini et que $F(0)$ est un espace vectoriel de
  dimension finie est lui-même fini. Par récurrence, on en déduit
  qu'un foncteur polynomial à valeurs de dimension finie est fini. La
  réciproque résulte de la proposition~\ref{prf-fpol}, puisqu'un objet
  simple se plonge dans un injectif standard.

La deuxième assertion découle de la première, en raison de
l'isomorphisme naturel $\Delta(F\otimes G)\simeq (\Delta F\otimes
G)\oplus (F\otimes\Delta G)\oplus (\Delta F\otimes\Delta G)$.

Comme le foncteur exact et fidèle $(\Delta,ev_0)$ commute aux
colimites, la proposition~\ref{crpf2} montre qu'un foncteur $F$ tel
que $F(0)$ est de dimension finie et que $\Delta F$ est pf$_\infty$
est lui-même pf$_\infty$. On en déduit le caractère pf$_\infty$ des
objets finis ; le cas co-pf$_\infty$ est dual.

La quatrième assertion se déduit de la précédente via la proposition~\ref{prevff}.

Un objet de co-type fini est toujours de socle fini ; il est aussi
analytique, dans~$\F$, parce que les injectifs standard le sont. La réciproque est donnée par la proposition~\ref{cotfs}.\end{proof}

L'assertion~\ref{aserls} de ce corollaire, particulièrement
importante, est due à L. Schwartz (cf.~\cite{LS}).

\paragraph*{Les objets simples de $\F(\kk)$} La propriété suivante est
une conséquence classique du caractère polynomial des foncteurs finis
de $\F(\kk)$. Dans~\cite{K4}, on trouvera une démonstration de cet
énoncé, sous une forme plus générale (qui ne suppose plus $\kk$ premier).

\begin{pr}\label{prec-f}  On
  suppose $d=1$, i.e. le corps $\kk$ premier. Soit $n$ un entier naturel.

Il existe un diagramme
  de recollement
$$\xymatrix{\F^{n-1}(\kk)\ar[r] & \F^n(\kk)\ar[r]\ar@/_/[l]\ar@/^/[l] & \mathbf{Mod}_{\kk[\Sigma_{n}]}\ar@/_/[l]\ar@/^/[l]
}.$$
\end{pr}

Par la proposition~\ref{crgrrec}, on en déduit le résultat suivant.

\begin{cor}\label{crf-sf2} Supposons $d=1$. La décomposition par le degré polynomial induit un isomorphisme de
    groupes de Grothendieck
$$G^f_0(\F(\kk))\simeq\bigoplus_{n\in\mathbb{N}}G^f_0(\mathbf{Mod}_{\mathbb{F}_p[\Sigma_{n}]}).$$
\end{cor}

\begin{pr}\label{crf-sf} Les automorphismes des objets simples de $\F(\kk)$ sont réduits
  à $\kk$ : $\kk$ est un {\em corps de décomposition} de la catégorie $\F(\kk)$.
\end{pr}

Comme le produit tensoriel est un bifoncteur exact en chaque variable,
que son unité $\kk$ est finie et qu'il préserve les objets finis, il induit une structure d'{\em anneau}
sur $G^f_0(\F(\kk))$.

\smallskip

La proposition suivante est établie dans \cite{K2}.

\begin{pr}[Kuhn]\label{comgr-ffs2} Si $S$ est un objet simple de $\F$, alors
  $S(\kk^{\oplus n})$ est un $GL_n(\kk)$-module nul ou simple pour
  tout $n\in\mathbb{N}$.
\end{pr}

\paragraph*{La conjecture artinienne} La terminologie employée dans ce
qui suit est rappelée dans l'appendice~\ref{apfca}.

 Nous donnons  différentes formulations
équivalentes de la conjecture artinienne (conjecture~\ref{canu}). Pour
la démonstration et plus de détails, voir \cite{GP1} ou~\cite{these}.

\begin{pr}\label{cael1} Les assertions suivantes sont équivalentes.
\begin{enumerate}\item\label{asca4} La catégorie $\F$ est localement
  noethérienne.
\item\label{asca5} Les projectifs standard $P_V$ sont des objets noethériens de $\F$.
\item\label{asca8} La catégorie $\F$ est co-localement artinienne.
\item\label{asca1} Les injectifs standard $I_V$ sont des objets artiniens de $\F$.
\item\label{asca3} Le produit tensoriel de deux objets artiniens de
  $\F$ est artinien.
\item\label{asca7} Le produit tensoriel de deux objets noethériens de
  $\F$ est noethérien.
\item\label{asca10} Le socle d'un objet de type fini de $\F$ est fini.
\item Les objets injectifs de $\F_\omega$ restent injectifs dans $\F$.
\item\label{ascae} Pour tout foncteur de type fini $F$ de $\F$,
  l'ensemble des classes d'isomorphisme de foncteurs simples $S$ de
  $\F$ tels que ${\rm Ext}^1_\F (F,S)\neq 0$ est fini.
\end{enumerate}
\end{pr}

L'assertion \ref{asca1} constitue la forme originelle de la conjecture
(l'intérêt particulier pour les injectifs provenant des liens avec les
modules instables sur l'algèbre de Steenrod --- cf.~\cite{HLS}), c'est
pourquoi elle est appelée {\em conjecture artinienne}. 

Les premières avancées significatives sur la conjecture artinienne remontent à Piriou
(\cite{Piriou}). Powell a accompli plusieurs progrès importants
(cf.~\cite{GP3}, \cite{GP2} et~\cite{GP5}). Dans \cite{art1}, nous démontrons de nouveaux cas
de cette conjecture ; les résultats de \cite{art3} (ou \cite{these})
déduits de l'étude des catégories de foncteurs en grassmanniennes surpassent ces résultats.

Il est naturel de poser la question du caractère
  localement noethérien ou co-localement artinien  d'autres
  catégories de foncteurs $\mathbf{Fct}(\I,\E)$, où $\I$ est une
  catégorie essentiellement petite vérifiant l'hypothèse~\ref{hypf3}. On ne saurait obtenir une réponse affirmative en
  général, comme le montre l'exemple élémentaire suivant. Cependant, toutes les catégories de foncteurs que nous rencontrerons par
  la suite semblent \go très probablement\gf localement
  noethériennes et co-localement artiniennes.

\begin{ex}\label{exnca} Soit $\I$ la petite catégorie telle que ${\rm
    Ob}\,\I=\mathbb{N}$ et ${\rm hom}_\I (n,m)=*$ (ensemble à un
  élément) si $n=0$ ou $n=m$, $\varnothing$ sinon. Le projectif standard
  de type fini $P^\I_0$ de $\mathbf{Fct}(\I,\E)$ est le foncteur constant $\kk$,
  puisque $0$ est objet initial de $\I$. Ce n'est pas un objet
  noethérien. En effet, pour tout entier $n>0$, l'image $A_i$ du
  monomorphisme $\underset{0<i\leq n}{\bigoplus}P^\I_i\to P^\I_0$ dont chaque
  composante $P^\I_i\to P^\I_0$ est induite par le morphisme $0\to i$
  associe à l'objet $i$ l'espace vectoriel $\kk$ si $0<i\leq
  n$, $0$ sinon. Aussi $(A_i)_{i>0}$ est-elle une suite strictement
  croissante de sous-objets de $P^\I_0$.
\end{ex}


Powell a émis la version renforcée suivante de la conjecture
artinienne.

\begin{conj}[Conjecture artinienne forte]\label{ca2} Pour tout $n\in\mathbb{N}$,
  le foncteur $P_{\kk^{\oplus n}}$ est noethérien de type $n$.
\end{conj}

\section{La catégorie $\F_{surj}$}\label{sct-surj}

\begin{nota} Nous noterons $\E^f_{surj}(\kk)$ et $\E^f_{inj}(\kk)$ les
sous-catégories de $\E^f_\kk$ ayant les mêmes objets et dont les
morphismes sont respectivement les épimorphismes et les monomorphismes
de $\E^f_\kk$. 

Si $I$ est une partie de $\mathbb{N}$, nous noterons $\E^I_{surj}(\kk)$ la
sous-catégorie pleine de $\E^f_{surj}(\kk)$ dont les objets sont les
espaces vectoriels dont la dimension appartient à $I$.

On introduit enfin les catégories de foncteurs
$$\F_{surj}(\kk)=\mathbf{Fct}(\E^f_{surj}(\kk),\E_\kk),$$
$$\F_{inj}(\kk)=\mathbf{Fct}(\E^f_{inj}(\kk),\E_\kk),$$
$$\F^I_{surj}(\kk)=\mathbf{Fct}(\E^I_{surj}(\kk),\E_\kk).$$
\end{nota}

On rappelle que, conformément aux conventions générales, on omettra
souvent la mention au corps $\kk$ dans la notation de ces catégories.

L'idée d'étudier une catégorie de foncteurs en utilisant
  la catégorie de foncteurs auxiliaire obtenue en ne conservant comme
  morphismes à la source que les épimorphismes, implicite dans les
  travaux de Suslin, apparaît chez Scorichenko. Elle est exposée  dans
  l'article {\em \go Stable K-theory is bifunctor homology\gf} (§\,4) de Franjou et Pirashvili du volume \cite{FFPS}.

Nombre des résultats présentés dans cette section sont connus, mais ils
ne semblent pas avoir fait l'objet d'une exposition systématique dans
la littérature.

\subsection{Généralités}\label{parsgs} Le foncteur de dualité
  $(\E^f)^{op}\to\E^f$ induit une équivalence de catégories entre
  $(\E^f_{surj})^{op}$ et $\E^f_{inj}$. On en déduit une  dualité
  entre les catégories $\F_{surj}$ et $\F_{inj}$ (cf. proposition/définition~\ref{prdcff}). Par conséquent, nous nous bornerons souvent à énoncer
  les résultats relatifs à l'une des deux catégories, privilégiant,
  dans cette section, la catégorie $\F_{surj}$. 

La considération d'une seule de ces catégories n'aurait néanmoins pas suffit
à nos investigations ultérieures : la catégorie $\F_{surj}$ intervient
dans l'étude des catégories de foncteurs en grassmanniennes ; en
retour, celles-ci fournissent des renseignements non triviaux sur la
catégorie $\F_{inj}$ --- cf. section~\ref{parfij}.

\begin{nota}\begin{enumerate}\item Si $V$ est un espace vectoriel de dimension finie, nous noterons
  simplement $P^{surj}_V$ le projectif
  standard $P^{\E^f_{surj}}_V$, et $I^{surj}_V$ l'injectif standard
  $I^{\E^f_{surj}}_V$. Nous noterons de même
  $P^{inj}_V$ et
  $I^{inj}_V$ pour $P^{\E^f_{inj}}_V$
  et $I^{\E^f_{inj}}_V$ respectivement.
\item L'espace vectoriel $\kk^{\oplus n}$ de dimension $n$ sera noté $E_n$. 
\item Nous noterons $ev_n=ev_{E_n} : \F_{surj}\to\E$ et  ${\rm ev}_n={\rm
  ev}_{E_n} :
  \F_{surj}\to\, _{\kk[GL_n(\kk)]}\mathbf{Mod}$ les foncteurs
  d'évaluation sur l'espace vectoriel $E_n$ (cf. notation~\ref{postpre}). Nous noterons aussi $ev_n$ et  ${\rm ev}_n$ les foncteurs analogues de source $\F_{inj}$.
\end{enumerate}
\end{nota}

Les considérations du paragraphe~\ref{p-prec} fournissent le résultat
suivant. La notation $\leq i$ en exposant désigne, selon nos
conventions générales, la partie de $\mathbb{N}$ évidente.

\begin{pr}\label{precsurj} Pour tout $n\in\mathbb{N}$, la
  sous-catégorie $\E^{\leq n-1}_{surj}$ est une sous-catégorie
  complète à droite de $\E^{\leq n}_{surj}$ (cf. définition~\ref{df-scc}). Par conséquent (cf. corollaire~\ref{crf-rec}), on a
  un diagramme de recollement
$$\xymatrix{\F^{\leq n-1}_{surj}\ar[r]|-{\mathcal{P}} &
\F^{\leq n}_{surj}  \ar[r]|-{{\rm ev}_n}\ar@/_/[l] \ar@/^/[l]^-{\mathcal{R}} &
 _{\kk[GL_n(\kk)]}\mathbf{Mod}\ar@/_/[l]\ar@/^/[l]^-{\mathcal{P}}}$$
où $\mathcal{R}$ est le foncteur de restriction et
$\mathcal{P}$ le prolongement par zéro.
\end{pr}

\begin{rem} Vespa a étudié (cf. \cite{CV}, chapitre 3 et \cite{CV2}, §\,$4$) une catégorie
  de foncteurs $\F_{iso}$ dont la source s'obtient à partir d'une
  catégorie d'espaces quadratiques dont toutes les flèches sont des
  monomorphismes, et de ce fait très analogue à $\E^f_{inj}$. Elle utilise une construction de foncteurs
  de Mackey qui permet d'éviter une partie des difficultés
  des catégories $\F_{inj}$ ou $\F_{surj}$ : le diagramme de
  recollement analogue à celui de la proposition~\ref{precsurj} est
  trivial dans $\F_{iso}$ (i.e. la catégorie qui apparaît au
  centre du diagramme se scinde en le produit des deux autres
  catégories --- cf. \cite{CV2}, théorème $4.2$), contrairement à ce qui advient dans $\F_{surj}$ ou $\F_{inj}$. Ce phénomène est à rapprocher de ce que
$\F_{iso}$ possède un foncteur de dualité, alors que $\F_{inj}$ a un
comportement \go fortement non auto-dual\gf.
\end{rem}

\begin{nota}\label{not-prz} \'Etant donné $n\in\mathbb{N}$, nous désignerons par $i^!_n
  : \,_{\kk[GL_n(\kk)]}\mathbf{Mod}\to\F_{surj}(\kk)$ le foncteur de
  prolongement par zéro obtenu en considérant la sous-catégorie
relativement connexe $\E^{n}_{surj}(\kk)\simeq\underline{GL_n(\kk)}$
de $\E^f_{surj}(\kk)$. On a ainsi un isomorphisme $$i^!_n(M)(V)\simeq\kk[{\rm
   Iso}(E_n,V)]\underset{\kk[GL_n(\kk)]}{\otimes}M$$ naturel en
 $V\in {\rm Ob}\,\E^f_\kk$ et $M\in {\rm Ob}\,_{\kk[GL_n(\kk)]}\mathbf{Mod}$.

Par abus, nous désignerons de la même façon les foncteurs analogues
dans $\F_{inj}(\kk)$.

Nous noterons enfin ${\rm Is}_n$ l'objet $i^!_n(\kk[GL_n(\kk)])$ de
  $\F_{surj}(\kk)$, et $S_n^{surj}=i^!_n(\kk)$.
\end{nota}

Ainsi, ${\rm Is}_0\simeq S_0^{surj}$ est le foncteur de $\F_{surj}$ égal à $\kk$ évalué sur l'espace
vectoriel $0$ et nul sur les espaces non nuls. 

\begin{rem}\label{prisfs} Les endofoncteurs $i^!_n\circ {\rm ev}_n$ et
  $\,\cdot\otimes S_n^{surj}$ de $\F_{surj}$ sont isomorphes.
\end{rem}

\paragraph*{Décomposition scalaire et tors de Frobenius} La catégorie
$\E^f_{surj}$ n'est pas $\kk$-linéaire, mais elle possède une action
naturelle du groupe multiplicatif $\kk^\times$. On en déduit, de façon
analogue à la proposition~\ref{pra-ds}, l'énoncé suivant.

\begin{prdef}\label{ds-fsurj} \'Etant donné un entier $i$, notons $\F_{surj, i}(\kk)$
  la sous-catégorie pleine de $\F_{surj}(\kk)$ formée des foncteurs
  $X$ tels que $X(\lambda.id_V)=\lambda^i.id_{X(V)}$ pour tous
  $\lambda\in\kk^\times$ et $V\in {\rm Ob}\,\E^f_{surj}$. Les
  inclusions induisent une équivalence de catégories 
$$\F_{surj}(\kk)\simeq\prod_{i=1}^{q-1}\F_{surj, i}(\kk).$$

On notera $X\simeq\bigoplus_{i=1}^{q-1} X_i$ la décomposition
naturelle d'un foncteur $X$ de $\F_{surj}(\kk)$ qu'on en déduit, où
$X_i\in {\rm Ob}\,\F_{surj, i}(\kk)$. On l'appelle {\em décomposition
  scalaire} de $X$.
\end{prdef}

Contrairement au cas d'une catégorie de foncteurs de source
$\kk$-linéaire, il n'y a pas de facteur correspondant à $i=0$ dans
cette décomposition. Cela est à rapprocher de ce que la sous-catégorie
pleine des foncteurs constants de $\F_{surj}$ n'est pas épaisse,
contrairement à ce qui se produit dans $\F$ (rappelons que la
sous-catégorie $\F_0$ obtenue par la décomposition scalaire est celle
des foncteurs constants de $\F$) --- cf. aussi proposition~\ref{dsf-fsf}.

\smallskip

L'automorphisme $\phi$ de la catégorie $\E^f_\kk$ déduit du morphisme
de Frobenius induit un automorphisme de la catégorie
$\E^f_{surj}(\kk)$. La précomposition par ce foncteur procure un
automorphisme de $\F_{surj}(\kk)$ que l'on appellera encore {\em tors
  de Frobenius}. Il induit une équivalence entre les catégories $\F_{surj,
  i}(\kk)$ et $\F_{surj, pi}(\kk)$.

On a des résultats analogues dans $\F_{inj}(\kk)$ et $\F_{surj}^I$,
pour $I\subset\mathbb{N}$. On obtient par exemple une variante de la
proposition~\ref{precsurj} en remplaçant $\F^{\leq i}_{surj}$ (pour
$i=n, n-1$) par $\F^{\leq i}_{surj,q-1}$ et $GL_n(\kk)$ par $PGL_n(\kk)$.

\paragraph*{Changement de corps} Si $K$ est une extension finie de
$\kk$, on dispose de foncteurs de restriction
$\F_{surj}(K)\to\F_{surj}(\kk)$ et d'induction
$\F_{surj}(\kk)\to\F_{surj}(K)$, définis comme dans $\F$. Cependant,
ces foncteurs ne sont plus adjoints : les foncteurs de
restriction $\E^f_K\to\E^f_\kk$ et d'induction  $\E^f_\kk\to\E^f_K$ se
restreignent en des foncteurs $\E^f_{surj}(K)\to\E^f_{surj}(\kk)$, mais
ces restrictions perdent les propriétés d'adjonction des foncteurs initiaux.

\paragraph*{Produit tensoriel total} Outre sa structure tensorielle usuelle, la
catégorie $\F_{surj}$ possède une structure tensorielle qui possède
l'avantage sur la structure usuelle de préserver les objets projectifs, et de
bien se comporter à l'égard du foncteur fondamental $\varpi : \F_{surj}\to\F$ que nous
introduirons au paragraphe~\ref{parfsf}.

On rappelle que la notation $\Gr$ désigne l'ensemble des sous-espaces
d'un espace vectoriel.

\begin{defi}\label{dpttfs} \'Etant donnés deux objets $F$ et
  $G$ de $\F_{surj}$, on appelle {\em produit tensoriel
    total} de $F$ et $G$ le
  foncteur noté $F\ptt G$ et défini comme suit.
\begin{itemize}\item Si $A$ est un objet de $\E^f_{surj}$, on pose
$$(F\ptt G)(A)=\underset{V+W=A}{\bigoplus_{V,W\in\Gr(A)}}F(V)\otimes G(W).$$
\item Si $u : A\twoheadrightarrow A'$ est un morphisme de $\E^f_{surj}$,
  le morphisme $$(F\ptt G)(u) : (F\ptt G)(A)\to (F\ptt G)(A')$$ est
  défini comme la somme directe sur les sous-espaces $V$ et $W$ de $A$
  tels que $V+W=A$ des morphismes
$$F(V)\otimes G(W)\xrightarrow{F(u)\otimes G(u)}F(u(V))\otimes
G(u(W))\hookrightarrow\underset{V'+W'=A'}{\bigoplus_{V',W'\in\Gr(A')}}F(V')\otimes G(W')$$
(où l'on note encore,
par abus, $u$ pour les morphismes $V\twoheadrightarrow u(V)$ et
$W\twoheadrightarrow u(W)$ induits par~$u$) --- cette définition fait
sens puisque $u(V)+u(W)=u(V+W)=u(A)=A'$.
\end{itemize}
%
\end{defi}

\begin{pr}\begin{enumerate}\item Le produit tensoriel total définit sur
    $\F_{surj}$ une structure monoïdale symétrique exacte d'unité ${\rm Is}_0$.
\item Il existe un monomorphisme naturel $F\otimes G\hookrightarrow F\ptt G$
pour $F, G\in {\rm Ob}\,\F_{surj}$.
\end{enumerate}
\end{pr}

La première assertion se vérifie par inspection. Le monomorphisme de la seconde assertion est donné par l'inclusion du facteur direct $F(V)\otimes G(V)$
de $(F\ptt G)(V)$ sur $V\in {\rm Ob}\,\E^f_{surj}$.

Par la suite, lorsque nous nous référerons à des notions dépendant
d'une structure tensorielle sur $\F_{surj}$, nous utiliserons
l'adjectif {\em total} lorsqu'il s'agira de la structure définie par
$\ptt$, l'absence de qualificatif signifiant qu'il s'agit de la
structure tensorielle usuelle définie par $\otimes$.

%

\begin{lm} Il existe un unique foncteur
$\E^f_{surj}\times\E^f_{surj}\to\E^f_{surj}$ tel que le diagramme
$$\xymatrix{\E^f_{surj}\times\E^f_{surj}\ar[r]\ar@{^{(}->}[d] & \E^f_{surj}\ar@{^{(}->}[d] \\
\E^f\times\E^f\ar[r]^-{\oplus} & \E^f
}$$
commute. Il sera encore noté $\oplus$ par abus. 
%
\end{lm}

\begin{rem} Par dualité, on en déduit un résultat analogue dans la
  catégorie $\E^f_{inj}$, où l'on peut aussi définir un produit tensoriel total.
\end{rem}

\begin{pr}\label{projptt1} Il existe un isomorphisme $P^{surj}_V\ptt
  P^{surj}_W\simeq P^{surj}_{V\oplus W}$ naturel en les objets $V$ et
  $W$ de $\E^f_{surj}$.
\end{pr}

\begin{proof} Cette proposition s'obtient en linéarisant
  l'isomorphisme ensembliste du lemme suivant.
\end{proof}

\begin{lm}\label{bijens} Il existe une bijection
$${\rm Epi}_\E (A\oplus B,E)\simeq
\underset{V+W=E}{\coprod_{V,W\in\Gr(E)}}{\rm Epi}_\E (A,V)\times
{\rm Epi}_\E (B,W)$$
naturelle en les objets $A$, $B$ et $E$ de $\E^f_{surj}$.
\end{lm}

\begin{cor}\label{cr-pttfs} Le produit tensoriel total de deux objets
  projectifs de $\F_{surj}$ est projectif.
\end{cor}

\subsection{Foncteur de décalage et objets finis}\label{parsof} Nous
étudions les objets finis de $\F_{surj}$ et $\F_{inj}$ à l'aide du
foncteur de décalage, qui joue formellement le même rôle 
que le foncteur différence dans l'étude des objets finis de~$\F$.

\begin{defi}\label{dec-surjinj} \'Etant donné un
  espace vectoriel de dimension finie $V$, on appelle {\em foncteur de
    décalage par $V$ dans $\F_{surj}$
    (resp. $\F_{inj}$)} l'endofoncteur de
  précomposition par $\cdot\oplus V$.  
\end{defi}

\begin{nota}\label{taudesi} Nous noterons $\delta_V^{surj}$
  (resp. $\delta_V^{inj}$) l'endofoncteur de décalage par un espace
  $V$ de la catégorie  $\F_{surj}$  (resp. $\F_{inj}$). Lorsque
  $V=\kk$, ce foncteur sera simplement noté $\delta^{surj}$
  (resp. $\delta^{inj}$).

De plus, l'exposant $surj$ (resp. $inj$) sera omis lorsqu'aucune
confusion ne pourra en résulter.
\end{nota}

Ainsi, l'association $V\mapsto\delta_V$ est fonctorielle ; si $V$ est
de dimension $n$, les foncteurs $\delta_V$ et $\delta^n$ sont
isomorphes.

Les foncteurs de décalage commutent aux limites, aux
colimites et au produit tensoriel (usuel), puisque ce sont des
foncteurs de précomposition.

\begin{pr}\label{agdscr} Il
  existe un isomorphisme
$${\rm hom}_{\F_{surj}} (X\ptt P^{surj}_A,Y)\simeq {\rm
  hom}_{\F_{surj}} (X,\delta^{surj}_A (Y))$$
naturel en les objets $X$, $Y$ de $\F_{surj}$ et $A$ de $\E^f_{surj}$.
\end{pr}

\begin{proof} On définit un morphisme naturel $X\to\delta^{surj}_A
  (X\ptt P^{surj}_A)$ par les inclusions
$$X(V)\xrightarrow{X(V)\otimes [id_A]}X(V)\otimes
P^{surj}_A(A)\hookrightarrow\bigoplus_{W_1+W_2=V\oplus A}X(W_1)\otimes
P^{surj}_A(W_2)\,.$$
On vérifie, grâce à la proposition~\ref{projptt1}, que l'application naturelle
$${\rm hom}_{\F_{surj}} (X\ptt P^{surj}_A,Y)\to {\rm
  hom}_{\F_{surj}} (X,\delta^{surj}_A (Y))$$
que l'on en déduit est un isomorphisme lorque $X$ est un projectif
standard. Le cas général s'en déduit par passage à la colimite.
\end{proof}

\begin{pr}\label{dec-pris} Il existe dans $\F_{surj}$ un isomorphisme
$$\delta_A(I^{surj}_B)\simeq\bigoplus_{W_1+W_2=B} I^{surj}_{W_1}\otimes I^{surj}_{W_2}(A)$$
naturel en les objets $A$ et $B$ de $\E^f_{surj}$.
\end{pr}

\begin{proof} Cette proposition s'obtient par linéarisation de la
  bijection du lemme~\ref{bijens} (ou à partir de la proposition
  précédente et du lemme de Yoneda).
\end{proof}

\begin{cor}\label{presdec-cr}\begin{enumerate}\item Les foncteurs de
    décalage de $\F_{surj}$ préservent les objets injectifs et les objets de co-type fini.
\item Les foncteurs de
    décalage de $\F_{inj}$ préservent les objets projectifs de type
    fini et les objets de type fini.
\end{enumerate}
\end{cor}

\begin{lm}\label{lm-sfpf}\begin{enumerate}\item Soit $X$ un objet de
    $\F_{surj}$ tel que $X(0)$ est de dimension finie et $\delta(X)$
    fini dans $\F_{surj}$. Alors $X$ est un objet fini de $\F_{surj}$.
\item Soient $i\in\mathbb{N}\cup\{\infty\}$ et $X$ un objet de
    $\F_{inj}$ tel que $X(0)$ est de dimension finie et $\delta(X)$
    pf$_i$ dans $\F_{inj}$. Alors $X$ est un objet pf$_i$ de $\F_{inj}$.
\end{enumerate}
\end{lm}

\begin{proof} Le foncteur $(\delta^{surj},ev_0) : \F_{surj}\to\F_{surj}\times\E$
est exact et fidèle. Le premier point résulte donc de la proposition
\ref{pr-ff}.

Pour le second, on utilise le foncteur $(\delta^{inj},ev_0) :
\F_{inj}\to\F_{inj}\times\E$, qui est exact et fidèle, commute aux
colimites et conserve les objets projectifs de type fini (cf. corollaire \ref{presdec-cr}). La conclusion est donc donnée par
la proposition \ref{crpf2}.
\end{proof}

\begin{pr}[Objets finis, de type fini de $\F_{surj}$]\label{ftffs} Soit $X$
  un objet de $\F_{surj}$. Les assertions suivantes sont équivalentes
  :
\begin{enumerate}\item l'objet $X$ est fini ;
\item l'objet $X$ est de type fini ;
\item l'objet $X$ est à valeurs de dimension finie et nilpotent pour le
  foncteur décalage~$\delta$ (cette dernière condition signifiant encore
  que $X(V)=0$ si $\dim V$ est assez grande).
\end{enumerate}
\end{pr}
%

\begin{proof} Il est clair que 1 implique 2.

Le foncteur projectif standard $P^{surj}_V$ associé à un
  espace $V$ de dimension $n$ est à valeurs de dimension finie et
  annihilé par le foncteur $\delta^{n+1}$, puisque ${\rm Epi}\,(V,W)=0$
  si $\dim W\geq n+1$. Cela montre que la deuxième assertion implique
  la troisième.

Enfin, la troisième implique la première par le lemme \ref{lm-sfpf}.
\end{proof}

Dualement, on a le résultat suivant.

\begin{pr}[Objets finis, de co-type fini de $\F_{inj}$]\label{fctfi} Soit $X$
  un objet de $\F_{inj}$. Les assertions suivantes sont équivalentes
  :
\begin{enumerate}\item l'objet $X$ est fini ;
\item l'objet $X$ est de co-type fini ;
\item l'objet $X$ est à valeurs de dimension finie et nilpotent pour le
  foncteur décalage~$\delta$.
\end{enumerate}
\end{pr}

\begin{cor}\label{pr-fsln}\begin{enumerate}\item Tout objet de $\F_{surj}$ est localement fini ; en particulier,
$\F_{surj}$ est une catégorie localement
noethérienne.
\item Tout objet de $\F_{inj}$ est co-localement fini ; en particulier,
$\F_{inj}$ est une catégorie co-localement artinienne.
\end{enumerate}
\end{cor}

\begin{cor}\label{cr-regfs}\begin{enumerate}\item Les objets finis de
    $\F_{inj}$ sont pf$_\infty$.
\item Les objets finis de $\F_{surj}$ sont co-pf$_\infty$.
\end{enumerate}
\end{cor}

\begin{proof} La première assertion s'obtient en combinant le lemme
  \ref{lm-sfpf} et la proposition \ref{fctfi}. La seconde s'en déduit
  par dualité.
\end{proof}

\begin{nota} On note $\F_{inj}^{lf}$ pour $(\F_{inj})^{lf}$ la
sous-catégorie des objets localement finis de $\F_{inj}$. Celle-ci est
{\em épaisse} parce que les objets finis de $\F_{inj}$ sont de
présentation finie, de sorte que l'on
peut appliquer la proposition~\ref{prevff}.
\end{nota}

\begin{cor}[Filtration canonique dans $\F_{surj}$]\label{fcan-fs}
  \'Etant donné $n\in\mathbb{N}$, notons $T_n$ l'endofoncteur de
  $\F_{surj}$ composé du foncteur de restriction
  $\mathcal{R} : \F_{surj}\to\F_{surj}^{\leq n}$ et du prolongement par zéro
  $\mathcal{P} : \F_{surj}^{\leq n}\to\F_{surj}$ ; la coünité de l'adjonction entre
  $\mathcal{R}$ et $\mathcal{P}$ fournit un monomorphisme $T_n\hookrightarrow id$.

Pour tout objet $F$ de $\F_{surj}$, $(T_n(F))_{n\in\mathbb{N}}$ est
une suite croissante de sous-objets de $F$ de réunion~$F$. Si $F$ est
à valeurs de dimension finie, alors les $T_n(F)$ sont finis.
\end{cor}

Cette filtration est similaire à la filtration
polynomiale dans
$\F$ (cf. par exemple \cite{GP4},~§\,$2$ à ce sujet). Notons
$T_n^{hom}=T_n/T_{n-1}$ ; il existe un isomorphisme canonique
$T_n^{hom}(X)\simeq X\otimes S^{surj}_n$, où $S^{surj}_n=i^!_n(\kk)$ est défini dans
la notation~\ref{not-prz}. Contrairement à ce qui
advient pour la filtration polynomiale dans $\F$, les foncteurs $T_n$ et
$T_n^{hom}$ sont exacts.

\begin{pr}[Objets simples de $\F_{surj}$]\label{simfsur} Un objet $S$
  de $\F_{surj}$ est simple si et seulement s'il existe
  $n\in\mathbb{N}$ tel que ${\rm ev}_n(S)$ est un
  $\kk[GL_n(\kk)]$-module simple et que ${\rm ev}_k(S)=0$ pour
  $k\neq n$. Il revient au même de dire que $S$ est isomorphe à
  $i^!_n(R)$ pour un certain  $n\in\mathbb{N}$ et un certain  $\kk[GL_n(\kk)]$-module simple $R$.
\end{pr}

\begin{proof} Un
  objet simple dans une sous-catégorie épaisse $\F^{\leq n}_{surj}$
  reste simple dans $\F_{surj}$ ; réciproquement, la proposition \ref{ftffs}
  montre qu'un objet simple de $\F_{surj}$ appartient à
  une sous-catégorie $\F^{\leq n}_{surj}$. La conclusion s'obtient alors
  en combinant les propositions \ref{precsurj} et~\ref{crgrrec}.
\end{proof}

\begin{ex} Les objets $S_n^{surj}$ sont simples dans $\F_{surj}$.
\end{ex}

\begin{rem} On peut formaliser ce raisonnement en notant que
  $\E^{f}_{surj}$ est la colimite filtrante (en un sens à préciser) des sous-catégories
  $\E^{\leq n}_{surj}$ et que, par conséquent, $\F_{surj}$ est la
  limite filtrante (id.) des $\F_{surj}^{\leq n}$. Kuhn détaille (dans le
  cadre analogue de la catégorie $\F$) cet argument dans \cite{K2}, §\,2.
\end{rem}

\begin{cor}\label{agrofs} Dans $\F_{surj}$, le produit tensoriel
  total, donc a fortiori le produit tensoriel usuel, de
  deux objets finis, est fini. Cela munit le groupe de Grothendieck
  $G_0^f(\F_{surj})$ d'une structure d'anneau commutatif sans unité
  (le foncteur constant $\kk$ n'étant pas fini) via
  $\otimes$, et d'une structure d'anneau commutatif (unitaire, le
  foncteur ${\rm Is}_0$ étant fini) via $\ptt$. Il est
  isomorphe, pour la première structure, à l'idéal
$$\bigoplus_{n\in\mathbb{N}} G_0^f(_{\kk[GL_n(\kk)]}\mathbf{Mod})$$
de l'anneau produit
$$\prod_{n\in\mathbb{N}} G_0^f(_{\kk[GL_n(\kk)]}\mathbf{Mod})$$
(chaque facteur étant muni de la structure d'anneau induite par le
produit tensoriel sur $\kk$).
\end{cor}

\subsection{Le foncteur $\varpi : \F_{surj}\to\F$}\label{parfsf} Nous
établissons dans ce paragraphe les liens élémentaires entre les catégories $\F$ et $\F_{surj}$.

\paragraph*{Le foncteur d'oubli $o : \F\to\F_{surj}$}

\begin{nota} Le foncteur de précomposition par le foncteur d'inclusion
  $\E^f_{surj}\to\E^f$ sera noté $o : \F\to\F_{surj}$.
\end{nota}

Ce foncteur est exact et fidèle.

\begin{pr}\label{cr-strffse} L'image par le foncteur $o : \F\to\F_{surj}$ d'un
  objet simple de $\F$ est un objet unisériel\,\footnote{i.e. l'inclusion est un ordre {\em total} sur
  l'ensemble de ses
  sous-objets.} de $\F_{surj}$, dont la
  filtration canonique est l'unique suite de composition.
\end{pr}

\begin{proof} Soient $S$ un objet simple de $\F$ et $X$ un quotient de
  $o(S)$. Pour tout $n\in\mathbb{N}$, la projection $\pi_n :
  S(E_n)\twoheadrightarrow X(E_n)$ est soit nulle, soit un
  isomorphisme, par la proposition~\ref{comgr-ffs2}.

Comme $o(S)$ transforme les morphismes de $\E_{surj}$ en des
épimorphismes, la commutation des diagrammes
$$\xymatrix{S(E_m)\ar@{>>}[r]\ar@{>>}[d]_{\pi_m} & S(E_n)\ar@{>>}[d]^{\pi_n} \\
X(E_m)\ar[r] & X(E_n)
}$$
pour $m\geq n$ (induits par un épimorphisme arbitraire
$E_m\twoheadrightarrow E_n$) montre que $X(E_n)=0$ si $X(E_m)=0$. Cela
démontre la proposition. 
\end{proof}

Nous spécifions maintenant le cas particulier, immédiat mais important, du foncteur constant de $\F$.

\begin{nota}\label{notstr-cs} Pour tout $i\in\mathbb{N}$, on désigne par $\kk^{\geq i}$  l'image du
foncteur constant $\kk$ de $\F^{\geq i}_{surj}$ par le foncteur de
prolongement par zéro $\mathcal{P} : \F^{\geq
  i}_{surj}\to\F_{surj}$.
\end{nota}

\begin{cor}\label{str-cs}  Pour tout $n\in\mathbb{N}$, il existe une
  suite exacte $0\to S^{surj}_n\to\kk^{\geq n}\to \kk^{\geq n+1}\to
  0$. De plus, $S^{surj}_n$ est le socle de $\kk^{\geq n}$.

L'objet $\kk$ de $\F_{surj}$ est
unisériel, de suite de
composition $S^{surj}_0, S^{surj}_1,\dots,S^{surj}_n,\dots$
\end{cor}

Ce résultat est à comparer avec le résultat fondamental sur la
structure  du foncteur $\bar{I}_\FF$ de $\F(\FF)$ : la filtration polynomiale
de ce foncteur est son unique suite de composition, et son $n$-ième
quotient est isomorphe au foncteur $n$-ième puissance extérieure $\Lambda^n$, qui correspond à la
représentation triviale de $GL_n(\FF)$. De plus, $\bar{I}_\FF$ est le \go plus
petit\gf injectif indécomposable non constant de $\F(\FF)$. Le foncteur
$\kk\simeq I^{surj}_0$ est le \go plus petit\gf
injectif indécomposable de $\F_{surj}(\kk)$ ; sa filtration canonique est
son unique suite de composition, et son $n$-ième quotient est
$T^{hom}_n(\kk)\simeq S^{surj}_n$ (cf. corollaire \ref{fcan-fs}), qui correspond à la représentation
triviale de $GL_n(\kk)$ (noter qu'ici on a également un facteur pour $n=0$).

Le corollaire \ref{str-cs} interviendra de façon cruciale dans la section~\ref{ssecp}.

\paragraph*{L'adjonction entre les foncteurs $o$ et $\varpi$}

\begin{prdef}\label{prdfffs} Il existe un foncteur exact et fidèle,
 noté $\varpi :
  \F_{surj}\to\F$, donné de la façon suivante.
\begin{itemize}\item Si $F$ est un objet de $\F_{surj}$ et $V$ un
  objet de $\E^f$, on pose
$$\big(\varpi(F)\big)(V)=\bigoplus_{W\in\Gr(V)} F(W).$$
\item Si $F$ est un objet de $\F_{surj}$ et $V\xrightarrow{f}V'$ une
  flèche de $\E^f$, et si $W$ et $W'$ sont des sous-espaces respectifs
  de $V$ et $V'$, la composante $F(W)\to F(W')$ de
  $\big(\varpi(F)\big)(f)$ est égale à $F(g)$ si $f(W)=W'$, où $g$ est
  l'épimorphisme $W\twoheadrightarrow W'$ induit par $f$, et nulle sinon.
\item Si $F\xrightarrow{u}G$ est une flèche de $\F_{surj}$, le
  morphisme $\varpi(u) : \varpi(F)\to\varpi(G)$ de $\F$ est défini,
  sur l'espace vectoriel de dimension finie $V$, comme la somme
  directe sur les sous-espaces $W$ de $V$ des morphismes $u_W : F(W)\to G(W)$.
\end{itemize} 
\end{prdef}

\begin{pr}\label{pre-vpo}\begin{enumerate}\item Le foncteur $o$ est
    adjoint à droite au foncteur $\varpi$.
\item\label{itvp1} Pour tout espace vectoriel $V$ de dimension finie, on a un isomorphisme $\varpi(P^{surj}_V)\simeq P_V$.
\item\label{itvp2} Pour tout espace vectoriel $V$ de dimension finie,
  on a un
  isomorphisme $o(I_V)\simeq\underset{W\in\Gr(V)}{\bigoplus}I^{surj}_W$.
\item Le foncteur $o$ préserve les objets de co-type fini et, plus
  généralement, les objets co-pf$_i$, où $i\in\mathbb{N}\cup\{\infty\}$.
\item Il existe un isomorphisme $\varpi(F\ptt
  G)\simeq\varpi(F)\otimes\varpi(G)$ naturel en les objets $F$ et $G$
  de $\F_{surj}$.
\end{enumerate}
\end{pr}

\begin{proof} Soient $X$ un objet de $\F_{surj}$ et $F$ un objet de
  $\F$. Un morphisme $\varpi(X)\to F$ est la donnée, pour tout $V\in
  {\rm Ob}\,\E^f$ et tout sous-espace $W$ de $V$, d'une application
  linéaire $u_{V,W} : X(W)\to F(V)$ de sorte que, pour tout $f\in {\rm
    hom}_{\E^f} (V,V')$, le diagramme
\begin{equation}\label{eqdt}\xymatrix{X(W)\ar[r]^-{u_{V,W}}\ar[d]_-{X(f)} & F(V)\ar[d]^-{F(f)} \\
X(W')\ar[r]_-{u_{V',W'}} & F(V')
}\end{equation}
commute, où l'on a posé $W'=f(W)$. Les applications linéaires
$u_{V,V}$ fournissent en particulier un morphisme  $X\to o(F)$ dans
$\F_{surj}$.

Réciproquement, si $a : X\to o(F)$ est une flèche de $\F_{surj}$,
définissons, pour $V\in {\rm Ob}\,\E^f$ et $W\in\Gr(V)$, une
application linéaire $u_{V,W} :
X(W)\xrightarrow{a_W}F(W)\hookrightarrow F(V)$, la dernière flèche
étant induite par l'inclusion. Il est clair que les diagrammes
(\ref{eqdt}) commutent, de sorte que $u$ fournit un morphisme
$\varpi(X)\to F$ dans $\F$. On vérifie aussitôt que les deux
applications naturelles entre ${\rm hom}_{\F_{surj}} (X,o(F))$ et
${\rm hom}_{\F} (\varpi(X),F)$ définies précédemment sont réciproques
l'une de l'autre, ce qui établit le premier point.

Les assertions \ref{itvp1} et \ref{itvp2}  s'obtiennent
par adjonction à partir du lemme de Yoneda : on a des isomorphismes naturels
$${\rm hom}_\F(\varpi(P^{surj}_V),F)\simeq {\rm
  hom}_{\F_{surj}}(P^{surj}_V,o(F))\simeq o(F)(V)=F(V)\simeq {\rm
  hom}_\F(P_V,F)$$
et
$${\rm hom}_{\F_{surj}}(X,o(I_V))\simeq {\rm
  hom}_\F(\varpi(X),I_V)\simeq\varpi(X)(V)=\bigoplus_{W\in\Gr(V)}X(W)$$
$$\simeq {\rm hom}_{\F_{surj}}(X,\bigoplus_{W\in\Gr(V)}I^{surj}_W).$$

La quatrième assertion provient de \ref{itvp2} et de
l'exactitude de $o$. 

La dernière s'obtient à partir de la décomposition ensembliste naturelle des
couples de sous-espaces d'un espace vectoriel $V$ suivante : 
$$\{(A,B)\,|\,A,B\in\Gr(V)\}=\coprod_{W\in\Gr(V)}\{(A,B)\,|\,A,B\in\Gr(W),\,A+B=W\}.$$
\end{proof}

\begin{rem} L'isomorphisme de la dernière assertion est adjoint au
  monomorphisme canonique $F\ptt G\hookrightarrow o(\varpi(F)\otimes\varpi(G))$
  donné sur l'espace $A$ par l'injection du facteur direct
$$\bigoplus_{V+W=A} F(V)\otimes
G(W)\hookrightarrow\bigoplus_{V,W\in\Gr(A)} F(V)\otimes G(W).$$
\end{rem}

\begin{pr}\label{offsurj} Soient $F$ un objet de $\F$ et $i\in\mathbb{N}\cup\{\infty\}$.
\begin{enumerate}\item Si $o(F)$ est artinien, alors $F$ est artinien.
\item Si $o(F)$ est co-pf$_i$, alors $F$ est co-pf$_i$.
\end{enumerate}
\end{pr}

\begin{proof} Le foncteur $o$ est exact et fidèle, ce qui fournit la
  première assertion, en appliquant la proposition~\ref{preltf}. De plus, il
  commute aux limites et transforme les cogénérateurs injectifs de
  $\F_{surj}$ en des objets injectifs de co-type fini. La variante
  duale de la proposition~\ref{crpf2} donne donc la deuxième partie de la proposition.
\end{proof}

La proposition suivante, que l'on vérifie par inspection, précise la
compatibilité entre les foncteurs $o$ et $\varpi$ d'une part, et la
décomposition scalaire et le tors de Frobenius d'autre part.

\begin{pr}\label{dsf-fsf}\begin{enumerate}\item Soit $F\in {\rm Ob}\,\F$. Il existe des isomorphismes
  naturels $o(F)_i\simeq o(F_i)$
  pour $1\leq i<q-1$ et $o(F)_{q-1}\simeq o(F_{q-1})\oplus o(F(0))$.
\item Soit $X\in {\rm Ob}\,\F_{surj}$. Il existe des isomorphismes
  naturels $\varpi(X)_0\simeq X(0)$, $\varpi(X)_i\simeq\varpi(X_i)$
  pour $1\leq i<q-1$ et $\varpi(X)_{q-1}\simeq\varpi(X_{q-1})/X(0)$.
\item Les foncteurs $o$ et $\varpi$ commutent au tors de Frobenius, à
  isomorphisme naturel près.
\end{enumerate}
\end{pr}

De façon similaire, le foncteur $o$ commute aux foncteurs de
restriction et d'induction associés à une extension finie du corps
$\kk$. En revanche, le comportement du foncteur $\varpi$ relativement
aux changements de corps est délicat.

Nous donnons maintenant  quelques traductions dans le cadre
de la catégorie $\F_{inj}$ de la proposition~\ref{pre-vpo}.

\begin{prdef}\label{pdfd-inj}\begin{enumerate}\item Nous noterons
    $o_{inj} : \F\to\F_{inj}$ le foncteur
    d'oubli (précomposition par l'inclusion $\E^f_{inj}\to\E^f$) ; il
    est exact et fidèle. De plus, il préserve les objets de type fini.
\item Il existe un foncteur exact et fidèle $\varpi_{inj} : \F_{inj}\to\F$ donné sur
  les objets par
$$\varpi_{inj}(X)(V)=\bigoplus_{W\in\Gr(V)} X(V/W)\,,$$
le morphisme $\varpi_{inj}(X)(V)\to\varpi_{inj}(X)(V')$ induit par une
application linéaire $f : V\to V'$ ayant pour composante $X(V/W)\to
X(V'/W')$ l'application $X(\bar{f})$, où $\bar{f}$ est le monomorphisme $V/W\to V'/W'$
induit par $f$, si $W=f^{-1}(W')$, $0$ sinon.
\item Le foncteur $\varpi_{inj}$ est adjoint à droite à $o_{inj}$.
\end{enumerate}
\end{prdef}

Nous avons noté $o_{inj}$ le foncteur d'oubli de $\F$ vers $\F_{inj}$,
et simplement $o$ plutôt que $o_{surj}$ le foncteur d'oubli de $\F$
vers $\F_{surj}$ (qui est dual du précédent), car ce
dernier interviendra beaucoup plus souvent. 

\paragraph*{Une propriété homologique du foncteur $\varpi$} Le
comportement homologique général de $\varpi$ semble délicat ; dans le
cas d'une action des groupes linéaires sur les foncteurs considérés, on
dispose cependant de l'utile résultat suivant.  

\begin{pr}\label{pr-extfs} Soient $n$ un entier strictement positif et
  $Z$ un objet de $\F_{surj}$ vérifiant les propriétés suivantes : 
\begin{enumerate}\item l'action du groupe linéaire sur le $GL_{n+1}(\kk)$-module $M={\rm ev}_{n+1}(Z)$
  (resp. le $GL_n(\kk)$-module $N={\rm ev}_n(Z)$) est triviale ;
\item le morphisme $u : M\to N$ induit par la projection
  $E_{n+1}\twoheadrightarrow E_n$ sur les $n$ premières coordonnées
  est injectif ;
\item on a $Z(E_k)=0$ si $k\notin\{n,n+1\}$.
\end{enumerate}

On note $X=i^!_{n+1}(M)$ et $Y=i^!_n(N)$.

Alors l'extension
$$0\to\varpi(Y)\to\varpi(Z)\to\varpi(X)\to 0$$
de $\F$ obtenue par application du foncteur exact $\varpi$ à la suite exacte
$0\to Y\to Z\to X\to 0$ de $\F_{surj}$ est essentielle.
\end{pr}

\begin{proof} Considérons un
  morphisme $f : \varpi(Z)\to F$ de $\F$
  dont la restriction à $\varpi(Y)$ est injective. On établit que
  l'application linéaire $f_V : \varpi(Z)(V)\to F(V)$ est injective pour
  tout $V\in {\rm Ob}\,\E^f$ par récurrence sur $\dim V$. C'est clair
  si $\dim V\leq n$, nous supposerons donc $\dim V>n$.

Soit $a$ un élément de $ker\,f_V$ ; si $W$ (resp. $B$) est un sous-espace de
dimension $n+1$ (resp. $n$) de $V$, notons $x_W$ (resp. $y_B$)
l'élément de $M$ (resp. $N$) correspondant à la composante de $a$ dans
$X(W)\simeq M$ (resp. $Y(B)\simeq N$). Pour tout
élément $v$ non nul de $V$, l'hypothèse de récurrence montre que
l'image de $a$ dans $\varpi(Z)(V/v)$ est nulle. Soient $H$ un sous-espace
de dimension $n$ de $V/v$ et $W$ son image réciproque dans $V$ : la
composante dans $Z(H)$ de cette image est 
$$u(x_W)+\underset{W=B\oplus v}{\sum_{B\in\Gr_n(V)}}y_B=0.$$

\'Etant donné un sous-espace $W$ de dimension $n+1$ de $V$, faisons la somme des relations ainsi
obtenues pour $v\in W\setminus\{0\}$ : on obtient
$$u(x_W)+\underset{W=B\oplus v}{\sum_{v\in
    W\setminus\{0\},B\in\Gr_n(W)}}y_B=0\,,$$
soit
$$u(x_W)=\sum_{B\in\Gr_n(W)}{\rm Card}\,(W\setminus B)\, y_B.$$

Mais comme $n>0$, les cardinaux qui apparaissent dans cette somme,
égaux à $q^{n+1}-q^n=q^n(q-1)$, sont nuls dans $\kk$, d'où $u(x_W)=0$.

Comme $u$ est supposé injectif, cela donne $x_W=0$ ; autrement dit,
$a\in\varpi(Y)(V)$. Puisque la restriction à $\varpi(Y)(V)$ de $f_V$
est par hypothèse injective, il vient finalement $a=0$.
Ainsi, l'extension de l'énoncé est bien essentielle. 
\end{proof}

\paragraph*{Foncteurs de Powell} Les duaux des foncteurs définis
ci-après ont été
introduits par Powell dans~\cite{GP2} (où la catégorie $\F_{surj}$
n'apparaît pas explicitement) sous le nom de {\em foncteurs
  co-Weyl} ; il en a montré l'intérêt dans l'étude de la conjecture artinienne.

\begin{defi}\label{fctpow} On appelle {\em foncteur de Powell} l'image
  par le foncteur $\varpi : \F_{surj}\to\F$ d'un objet simple de $\F_{surj}$.

Nous nommerons {\em filtration de Powell} d'un foncteur $F$ de $\F$ toute filtration finie
$$0=A_0\subset A_1\subset\dots\subset A_n=F$$
de $F$ dont les sous-quotients $A_i/A_{i-1}$ sont des foncteurs de
Powell.
\end{defi}

Dans~\cite{GP2}, les foncteurs possédant une filtration de Powell sont
appelés foncteurs {\em DJ-bons}. De nombreuses propriétés en sont
établies. Nous nous contenterons de démontrer la stabilité par produit tensoriel, particulièrement commode avec le formalisme de la catégorie $\F_{surj}$.

\begin{pr}\label{pregp} Si $X$ est un objet fini de $\F_{surj}$, alors
  $\varpi(X)$ admet une filtration de Powell. 
\end{pr}

Cette propriété résulte de l'exactitude du foncteur $\varpi$.

\begin{ex} Les projectifs standard $P_V\simeq\varpi(P^{surj}_V)$ de
  $\F$ admettent une filtration de Powell.
\end{ex}

\begin{cor}\label{crptop} Le produit tensoriel de deux foncteurs de $\F$ possédant une
  filtration de Powell admet  une filtration de Powell.
\end{cor}

\begin{proof} Cela provient de la dernière assertion de la proposition~\ref{pre-vpo}, du corollaire~\ref{agrofs} et de l'exactitude du
  produit tensoriel de $\F$.
\end{proof}

Ce corollaire, qui illustre l'intérêt du produit tensoriel total, est
établi dans \cite{GP2} à partir d'un critère homologique pour
l'existence d'une filtration de Powell.

Signalons également le lien entre les objets simples de $\F$ et
$\F_{surj}$ réalisé par les foncteurs de Powell, qui se déduit
de~\cite{GP2}, §\,$2.3$ (cf. aussi~\cite{K2}).

\begin{pr}[Kuhn-Powell]\label{comgr-ffs} Le cosocle d'un foncteur de Powell est
  simple. De plus, tout foncteur simple de $\F$ est le cosocle d'un
  unique foncteur de Powell, à isomorphisme près.
\end{pr}

\begin{rem} Cette propriété procure un isomorphisme de groupes
  $G^f_0(\F)\simeq G_0^f(\F_{surj})$. En revanche, les structures multiplicatives sur
  $G_0^f(\F_{surj})$ et $G_0^f(\F)$ semblent fort délicates à comparer.
\end{rem}

\paragraph*{Liens avec la conjecture artinienne} \`A l'aune du
corollaire~\ref{pr-fsln}, les problèmes de finitude dans $\F_{surj}$
se posent en les termes suivants.

\begin{conj}\label{casf} La catégorie $\F_{surj}$ est co-localement artinienne.
\end{conj}

Par dualité, cette conjecture équivaut au caractère localement noethérien de
$\F_{inj}$.

\begin{pr}\label{ca20} La conjecture \ref{casf} implique la
  conjecture artinienne.
\end{pr}

\begin{proof} Si la conjecture \ref{casf} est vérifiée, $o(I_V)$ est
  artinien pour tout $V\in {\rm Ob}\,\E^f$ (cf.~proposition~\ref{pre-vpo}.\,\ref{itvp2}) ; la proposition \ref{offsurj} montre que les
  injectifs standard $I_V$ de $\F$ sont alors également artiniens.
\end{proof}

 Il ne semble en revanche exister aucun argument formel
  pour obtenir la réciproque de la proposition \ref{ca20}.

\begin{nota}\label{not-tilff} \'Etant donné $n\in\mathbb{N}$, nous noterons $\widetilde{P}(n)$ l'objet $\varpi({\rm
    Is}_n)$ de $\F$.
\end{nota}
\nopagebreak[3]

On a ainsi un isomorphisme canonique $\widetilde{P}(n)(V)\simeq\kk[{\rm Pl}_\E (E_n,V)]$ pour
$V\in {\rm Ob}\,\E^f$.

\pagebreak[3]

\begin{lm}\label{lmfsc}\begin{enumerate}\item \'Etant donné un objet $X$
    de $\F_{surj}$, les assertions suivantes sont équivalentes.
\begin{enumerate}\item L'objet $X$ est de co-type fini.
\item L'ensemble $\{n\in\mathbb{N}\,|\,{\rm hom}\,({\rm Is}_n,X)\neq 0\}$
est fini et $X$ est à valeurs de dimension finie.
\item Le socle de $X$ est fini.
\end{enumerate}
\item Soit $X$ un objet co-tf 
    de $\F_{surj}$. Les assertions suivantes sont équivalentes.
\begin{enumerate}\item L'objet $X$ est co-pf.
\item L'ensemble $\{n\in\mathbb{N}\,|\,{\rm Ext}^1({\rm Is}_n,X)\neq
  0\}$ est fini.
\end{enumerate}
\end{enumerate}
\end{lm}

\begin{pr}\label{cafsurj} La conjecture artinienne équivaut
  à l'assertion suivante : pour tout objet de co-type fini
  $F$ de $\F$, l'ensemble
$$\{n\in\mathbb{N}\,|\,{\rm Ext}^1_\F (\widetilde{P}(n),F)\neq 0\}$$
est fini.
\end{pr}

Le lemme et la proposition précédents sont laissés au lecteur, qui
pourra en trouver la démonstration dans~\cite{these}.

Powell a émis, à partir des résultats de \cite{GP2}, la conjecture suivante, dont nous proposerons une
version renforcée ultérieurement.

\begin{conj}[Conjecture artinienne très forte]\label{ca4} Pour tout
  objet simple $S$ de $\F_{surj}$ tel que $ev_n(S)\neq 0$, le foncteur de Powell $\varpi(S)$ est
  simple noethérien de type~$n$.
\end{conj}

\subsection{Liens avec les systèmes de coefficients}\label{par-dwyer}

 Nous rappelons la définition de Dwyer de la catégorie des systèmes de
coefficients (\cite{Dw}).

\begin{conv} Dans ce paragraphe, on note ${\rm R}_n :
\,  _{\kk[GL_{n+1}(\kk)]}\mathbf{Mod}\to\, _{\kk[GL_n(\kk)]}\mathbf{Mod}$, pour
      $n\in\mathbb{N}$, le foncteur de restriction. Le
          groupe $GL_n(\kk)$ est plongé dans $GL_m(\kk)$, pour $n\leq m$, par 
$$\alpha_{n,m} : A\mapsto\left(\begin{array}{cc} A & 0 \\
0 & I_{m-n}\end{array}\right).$$

On note, si $n\leq m$ sont des entiers, $l_{n,m} : E_n\hookrightarrow E_m$
l'inclusion $x\mapsto (x,0)$.

On note enfin, pour $n\leq m$, $GL(m,n)$ le sous-groupe de
$GL_m(\kk)$ des automorphismes $g$ de $E_m$ tels que $g\circ l_{n,m}=l_{n,m}$. 
\end{conv}

\begin{rem}\label{rqagl} Comme le sous-groupe $GL_n(\kk)$ de $GL_m(\kk)$
  normalise $GL(m,n)$, on peut voir le foncteur $M\mapsto M^{GL(m,n)}$
  (invariants sous l'action de $GL(m,n)$) défini sur les $GL_m(\kk)$-modules comme un sous-foncteur de la
  restriction ${\rm R}_n\dots{\rm R}_{m-1} :
\,  _{\kk[GL_m(\kk)]}\mathbf{Mod}\to\, _{\kk[GL_{n}(\kk)]}\mathbf{Mod}$.
\end{rem}

\begin{defi} On note $\E^f_{coef}$ la catégorie dont les objets sont
  les espaces $E_n$ pour $n\in\mathbb{N}$, les morphismes étant donnés
  par ${\rm hom}_{\E^f_{coef}}(E_n,E_m)=GL_m(\kk)$ si $n\leq m$, $\varnothing$ si
  $n>m$, et où la composition est définie comme suit. Si $0\leq n\leq m\leq
  k$ sont des entiers, la composée $g\circ h$, où $g\in {\rm
    hom}_{\E^f_{coef}}(E_m,E_k)=GL_k(\kk)$ et $h\in {\rm
    hom}_{\E^f_{coef}}(E_n,E_m)=GL_m(\kk)$, est l'élément
  $g.\alpha_{m,k}(h)$ de $GL_k(\kk)= {\rm
    hom}_{\E^f_{coef}}(E_n,E_k)$.

On note $\C oef$ la catégorie $\mathbf{Fct}(\E^f_{coef},\E)$ ; ses
objets sont appelés {\em systèmes de coefficients}.
\end{defi}

\begin{rem}\label{lddsc} Dans \cite{Dw}, un système de coefficients
  est une
  suite $(M_n)_{n\in\mathbb{N}}$, où $M_n$ est un $GL_n(\kk)$-module à
  gauche, munie de morphismes de $GL_n(\kk)$-modules $r_n : M_n\to {\rm
    R}_n(M_{n+1})$.

L'équivalence entre cette définition et la nôtre  résulte de ce que la catégorie
$\E^f_{coef}$ est engendrée  (cf. \cite{Mac}, chapitre~II, §\,8) par les éléments de $GL_n(\kk)={\rm
  hom}_{\E^f_{coef}}(E_n,E_n)$ et les morphismes $i_{n,m} : E_n\to E_m$
correspondant à $I_m\in GL_m(\kk)$, pour $n\leq m$, soumis aux
relations $i_{n,n+1}g=\alpha_{n,n+1}(g)\in {\rm
  hom}_{\E^f_{coef}}(E_n,E_{n+1})$ pour $n\in\mathbb{N}$ et $g\in GL_n(\kk)$.
\end{rem}

\begin{lm}\label{prcfi-dw} On définit un foncteur plein et
  essentiellement surjectif $\E^f_{coef}\to\E^f_{inj}$ de la façon
  suivante :
\begin{itemize}\item on associe à l'objet $E_n$ de $\E^f_{coef}$
  l'objet $E_n$ de $\E^f_{inj}$ ;
\item si $0\leq n\leq m$ sont des entiers, on associe à $g\in {\rm
    hom}_{\E^f_{coef}}(E_n,E_m)=GL_m(\kk)$ le monomorphisme
  $E_n\hookrightarrow E_m$ composé de l'inclusion 
  $l_{n,m} : E_n\hookrightarrow E_m$ et de l'automorphisme $g : E_m\to E_m$.
\end{itemize} 
\end{lm}

Ce lemme se vérifie par inspection.

\begin{nota}\label{dw-finj}  On note ${\rm C} : \F_{inj}\to\C oef$ le
  foncteur de précomposition par le foncteur du lemme~\ref{prcfi-dw}. 
\end{nota}

On rappelle que la notion de sous-catégorie de Serre est introduite
dans la définition~\ref{dfserre}.

\begin{pr}\label{encdw} Le foncteur ${\rm C}$ est pleinement fidèle ;
  son image est une sous-catégorie de Serre de $\C oef$.
\end{pr}

\begin{proof} Ce résultat s'obtient en combinant le
  lemme~\ref{prcfi-dw} et la proposition~\ref{lm-form}. 
\end{proof}


On peut préciser cette propriété comme suit.

\begin{pr}\label{fd-dw}  L'image essentielle du foncteur ${\rm C}$ est constituée
des systèmes de coefficients $((M_n),(r_n))$
tels que $im\,r_m r_{m-1}\dots r_n\subset {M_{m+1}}^{GL(m+1,n)}$ pour
tous entiers $n$ et $m$ tels que $0\leq n\leq m$, où l'on utilise les
notations de la remarque~\ref{lddsc} pour les systèmes de coefficients.
\end{pr}

\begin{proof} C'est une conséquence formelle de l'observation suivante
  : le squelette de la catégorie $\E^f_{inj}$ constitué des espaces
  $E_n$ est engendré par les inclusions $l_{n,n+1}$ et les éléments des
  différents $GL_n(\kk)$, soumis aux relations $\alpha_n(g)\circ l_n=l_n\circ g$ pour
  $g\in GL_n(\kk)$, $g\circ h=gh$ pour $g,h\in GL_n(\kk)$ et $g\circ l_{n,m}=l_{n,m}$ pour $g\in GL(m+1,n)$.
\end{proof}

\begin{rem} En particulier, tous les systèmes de coefficients
  appartenant à l'image du foncteur ${\rm C}$
  sont {\em centraux} au sens de \cite{Dw}, §\,$2$.
\end{rem}

\paragraph*{Les catégories $\C oef^{lf}$ et $\C oef/\C oef^{lf}$} Le
foncteur de décalage $\delta^{inj} : \F_{inj}\to\F_{inj}$ admet un
relèvement à $\C oef$. Précisément, il existe un endofoncteur exact
$\delta^{coef}$ de $\C oef$ tel que le diagramme 
\begin{equation}\label{ldc-coinj}\xymatrix{\F_{inj}\ar[r]^-{\delta^{inj}}\ar[d]_{{\rm C}} & \F_{inj}\ar[d]^{{\rm C}} \\
\C oef\ar[r]^-{\delta^{coef}} &  \C oef
}\end{equation}
commute (à isomorphisme canonique près).

Le foncteur $\delta^{coef}$ s'obtient par précomposition par
l'endofoncteur de $\E^f_{coef}$ donné sur les objets par $E_n\mapsto
E_{n+1}$ et sur les morphismes par l'inclusion
$GL_n(\kk)\hookrightarrow GL_{n+1}(\kk)$ donnée par 
$$A\mapsto\left(\begin{array}{cc} 1 & 0 \\
0 & A\end{array}\right)
$$
($\delta^{coef}$ est le foncteur
$\Sigma$ de \cite{Dw}, §\,2).

On établit, par la même méthode que celle employée dans $\F_{inj}$, le
résultat suivant.

\begin{pr}\label{fcoef}\begin{enumerate}\item Un objet de $\C oef$ est
    fini si et seulement s'il est nilpotent pour le foncteur
    $\delta^{coef}$ et à valeurs de dimension finie.
\item Tout objet fini de $\C oef$ est pf$_\infty$. Par conséquent, la
  sous-catégorie pleine $\C oef^{lf}$ de $\C oef$ des systèmes de
  coefficients localement finis est épaisse. 
\end{enumerate}
\end{pr}

Soit $GL(\kk)$ la
colimite des groupes linéaires $GL_n(\kk)$ (relativement aux
inclusions $\alpha_{n,m}$). On peut voir $GL(\kk)$ comme le groupe
des automorphismes linéaires $g$ de l'espace vectoriel
$E_\infty=\underset{n\in\mathbb{N}}{\col} E_n$ tels que $ker\,(g-id)$
est de codimension finie.

On définit un foncteur ${\rm e}_{coef} :
  \C oef\to\,_{\kk[GL(\kk)]}\mathbf{Mod}$  par 
${\rm e}_{coef}(X)=\underset{n\in\mathbb{N}}{\col} X(E_n)$, l'action du
groupe linéaire $GL(\kk)$ provenant de l'action de $GL_n(\kk)$ sur $X(E_n)$.

\begin{pr}\label{coefslf}\begin{enumerate}\item Le foncteur ${\rm
      e}_{coef}$ est adjoint à gauche au foncteur
    ${\rm c} :\,_{\kk[GL(\kk)]}\mathbf{Mod}\simeq\mathbf{Fct}(\underline{GL(\kk)},\E)\to\C oef=\mathbf{Fct}(\E^f_{coef},\E)$ donné par la précomposition par le foncteur $\E^f_{coef}\to\underline{GL(\kk)}$ obtenu en plongeant les ensembles de morphismes $GL_n(\kk)$ dans $GL(\kk)$.
\item  Le foncteur ${\rm
      e}_{coef}$ induit une équivalence entre les catégories $\C oef/\C oef^{lf}$ et $_{\kk[GL(\kk)]}\mathbf{Mod}$.
\end{enumerate}
\end{pr}

\begin{proof} Si $P$ est un générateur projectif standard de $\C oef$,
  il existe un isomorphisme canonique ${\rm
    e}_{coef}(P)\simeq\kk[GL(\kk)]$. On en déduit, pour tout
  $GL(\kk)$-module $M$, un isomorphisme canonique
$${\rm hom}_{\kk[GL(\kk)]}({\rm
    e}_{coef}(P),M)\simeq M\simeq {\rm hom}_{\C oef}(P,{\rm c}(M)).$$

Cela donne la première assertion, en écrivant un objet de $\C oef$
comme colimite de projectifs standard, puisque le foncteur ${\rm
  e}_{coef}$ commute aux colimites.

La première partie de la proposition~\ref{fcoef} implique que le noyau
du foncteur exact $e_{coef}$ est égal à $\C oef^{lf}$ ; il induit donc
un foncteur exact et fidèle $\C oef/\C
oef^{lf}\to\,_{\kk[GL(\kk)]}\mathbf{Mod}$.

La conclusion s'obtient alors formellement à partir des deux observations suivantes
: 
\begin{enumerate}\item la coünité ${\rm e}_{coef}{\rm c}\to id$ de l'adjonction est un
isomorphisme ; 
\item l'unité $id\to {\rm c}{\rm e}_{coef}$ induit un isomorphisme
  après application du foncteur ${\rm e}_{coef}$.
\end{enumerate}
\end{proof}

\begin{rem}\label{rq-cafx} L'anneau $\kk[GL(\kk)]$ n'est pas
  noethérien à gauche. En effet, si l'on note
  $A_n=ker\,(\kk[GL(\kk)]\twoheadrightarrow\kk[GL(\kk)/GL_n(\kk)])$,
  les $A_n$ forment une suite strictement croissante d'idéaux à gauche de
  $\kk[GL(\kk)]$. Par conséquent, la proposition~\ref{coefslf} montre
  que la catégorie $\C oef$ n'est pas localement noethérienne.

Ce phénomène est à rapprocher de l'exemple~\ref{exnca} : dans cet
exemple, l'obstruction au caractère noethérien de la catégorie de
foncteurs considérée est de nature \go combinatoire\gf (la catégorie
source n'a pas d'endomorphismes non identiques). Ici, l'obstruction
vient plutôt de la théorie des groupes.
\end{rem}

\begin{nota}\label{not-e-inj} Nous noterons ${\rm e} :
  \F_{inj}(\kk)\to\,_{\kk[GL(\kk)]}\mathbf{Mod}$ le foncteur composé
$$\F_{inj}(\kk)\xrightarrow{{\rm C}}\C oef\twoheadrightarrow\C oef/\C
  oef^{lf}\simeq\,_{\kk[GL(\kk)]}\mathbf{Mod}.$$
\end{nota}

\begin{rem}
Le foncteur ${\rm e}$ est beaucoup moins
élémentaire que le foncteur analogue
$\F(\kk)\to\,_{\kk[\mathcal{M}(\kk)]}\mathbf{Mod}$, où
$\mathcal{M}(\kk)$ désigne le monoïde sans unité colimite des
$\mathcal{M}_n(\kk)$ relativement aux inclusions
$$A\mapsto\left(\begin{array}{cc} A & 0 \\
0 & 0\end{array}\right).$$

Pour une discussion à ce sujet, voir \cite{K2},~§\,$3$.
\end{rem}

\begin{pr}\label{finjmlf} Le foncteur ${\rm C} : \F_{inj}\to\C oef$
  induit une équivalence entre la catégorie $\F_{inj}/\F_{inj}^{lf}$
  et la sous-catégorie de Serre de $\C oef/\C
  oef^{lf}\simeq\,_{\kk[GL(\kk)]}\mathbf{Mod}$ des $GL(\kk)$-modules
  $M$ tels que $M=\underset{n\in\mathbb{N}}{\col} M^{GL(\infty,n)}$,
  où $GL(\infty,n)$ désigne le sous-groupe de $GL(\kk)$ des
  automorphismes $g$ tels que $g(x)=x$ pour $x\in E_n$.
\end{pr}

\begin{proof} Le diagramme commutatif (\ref{ldc-coinj}) montre qu'un
  objet $X$ de $\F_{inj}$ est fini si et seulement si son image par le
  foncteur ${\rm C}$ est finie. Comme ce foncteur commute aux
  colimites, un objet de $\F_{inj}$ est localement fini si et
  seulement s'il en est de même pour son image par ${\rm C}$. Ainsi,
  ce foncteur induit un foncteur exact et fidèle $\F_{inj}/\F_{inj}^{lf}\to\C oef/\C
  oef^{lf}\simeq\,_{\kk[GL(\kk)]}\mathbf{Mod}$.

D'autre part, le lemme de Yoneda montre que le foncteur ${\rm e}$ 
admet un adjoint à droite $r$ donné par 
$$r(M)={\rm hom}_{GL(\kk)}(\kk[{\rm Pl}_\E (V,E_\infty)],M).$$

Comme ${\rm Pl}_\E (E_n,E_\infty)$ est un $GL(\kk)$-ensemble
transitif, et que le stabilisateur de l'inclusion canonique
$E_n\hookrightarrow E_\infty$ est $GL(\infty,n)$, on a
$r(M)(E_n)\simeq M^{GL(\infty,n)}$. En particulier, la coünité de
l'adjonction ${\rm e}r\to id$ est un {\em monomorphisme}. On en déduit
formellement, comme dans la démonstration de la
proposition~\ref{coefslf}, que le foncteur
$\F_{inj}/\F_{inj}^{lf}\to\,_{\kk[GL(\kk)]}\mathbf{Mod}$ qui nous
intéresse  induit une équivalence entre  $\F_{inj}/\F_{inj}^{lf}$
et la sous-catégorie pleine de $_{\kk[GL(\kk)]}\mathbf{Mod}$ des
$GL(\kk)$-modules $M$ tels que la coünité ${\rm e}r(M)\to M$ est un
isomorphisme, d'où la proposition.
\end{proof}

\section{Catégories de comodules sur un foncteur en coalgèbres de
Boole}\label{sctccf}

Nous exposons dans cette section des résultats formels généralisant, en
termes de catégories de foncteurs, l'observation suivante. {\em Soit $X$
un ensemble, munissons l'espace vectoriel $\kk[X]$ de la structure de
coalgèbre donnée par la diagonale $[x]\mapsto [x]\otimes [x]$ et la
coünité $[x]\mapsto 1$. La catégorie des $\kk$-espaces vectoriels
$X$-gradués est naturellement équivalente à la sous-catégorie
$\mathbf{Comod}_{\kk[X]}$ des $\kk[X]$-comodules de $\E_\kk$.} Une
coalgèbre du type $\kk[X]$ est appelée {\em coalgèbre de
  Boole}, car, lorsque $\kk=\FF$, la
structure de coalgèbre sur $\FF[X]$ est duale de la structure
d'algèbre de Boole de $\FF^X$. 

L'étude de catégories de modules ou de comodules sur des objets
présentant une structure de type booléen dans un contexte proche des
catégories de foncteurs remonte aux travaux des années~$90$ sur les
modules instables sur l'algèbre de Steenrod. 

Dans la deuxième partie de~\cite{HLS}, Henn, Lannes et Schwartz
établissent un lien fondamental entre les foncteurs en algèbres de
Boole de la catégorie $\F(\mathbb{F}_p)$ et les algèbres instables sur
l'algèbre de Steenrod modulo $p$ (voir aussi \cite{Lannes}). Ils utilisent dans~\cite{HLS2} des structures de
module sur ces algèbres pour étudier la cohomologie équivariante
modulo $p$. Des considérations analogues apparaissent dans les travaux
de Lannes et Zarati (cf.~\cite{LZ2}).

\begin{conv} Dans toute cette section, on se donne une catégorie
essentiellement petite~$\I$.
\end{conv}

\begin{rem} L'hypothèse de finitude du corps $\kk$ n'interviendra pas dans cette section.
\end{rem}

\subsection{La catégorie de comodules $\mathbf{Fct}(\I_{\backslash
    X},\E_\kk)$}\label{sctcom} Ce paragraphe a pour but d'identifier
les catégories de comodules sur un foncteur obtenu par linéarisation
d'un foncteur ensembliste à des catégories de foncteurs convenables.

\begin{conv} Dans ce paragraphe, $X$ désigne un foncteur de $\I$
vers $\bf{Ens}$. 
\end{conv}

\begin{nota}\label{notfcom}  On définit une catégorie $\I_{\backslash X}$ de la manière
  suivante.
\begin{itemize}
\item Les objets de $\I_{\backslash X}$ sont les couples $(E,x)$, où $E$ est un
  objet de $\I$ et $x$ un élément de $X(E)$.
\item Les morphismes dans $\I_{\backslash X}$ sont donnés par
$${\rm hom}_{\I_{\backslash X}}((E,x),(F,y))=\{u\in {\rm hom}_\I(E,F)\,|\,X(u)(x)=y\}.$$
\item La composition des morphismes dans $\I_{\backslash X}$ est induite par celle
  de $\I$.
\end{itemize}

Le foncteur d'oubli $(E,x)\mapsto E : \I_{\backslash X}\to\I$ sera noté $\mathcal{O}_{\I,X}$.
\end{nota}

\begin{rem} Il s'agit d'un cas particulier de la construction
  catégorique classique étudiée dans \cite{Mac}, chapitre~II, §\,$6$.
\end{rem}

\begin{nota}\label{dfup} Nous désignerons par $\Upsilon_X :
  \mathbf{Fct}(\I,\E_\kk)\to\mathbf{Fct}(\I_{\backslash X},\E_\kk)$ le
  foncteur de précomposition par $\mathcal{O}_{\I,X}$.
\end{nota}

Ainsi, $\Upsilon_X$ est un foncteur
exact qui commute aux limites, aux colimites et au produit
tensoriel (cf. proposition~\ref{lm-form}). Dans le cas où $X$ prend ses valeurs dans les ensembles non
vides, $\mathcal{O}_{\I,X}$ est essentiellement surjectif, donc
$\Upsilon_X$ est fidèle.

\begin{prdef}\label{prdfom} Il existe un foncteur exact et fidèle,
  appelé {\em foncteur de $X$-intégrale} et noté $\Omega_X :
  \mathbf{Fct}(\I_{\backslash
    X},\E_\kk)\to\mathbf{Fct}(\I,\E_\kk)$,
  défini de la façon suivante.
\begin{enumerate}\item Si $F$ est un objet de $\mathbf{Fct}(\I_{\backslash X},\E_\kk)$ et $E$ un
  objet de $\I$, on pose
$$\big(\Omega_X(F)\big)(E)=\bigoplus_{x\in X(E)} F(E,x).$$
\item Si $F$ est un objet de $\mathbf{Fct}(\I_{\backslash X},\E_\kk)$ et $E\xrightarrow{f}E'$ une
  flèche de $\I$, et si $x$ et $x'$ sont des éléments respectifs
  de $X(E)$ et $X(E')$, la composante $F(E,x)\to F(E',x')$ de
  $\big(\Omega_X(F)\big)(f)$ est égale à $F(f)$ si $X(f)(x)=x'$, et à $0$ sinon.
\item Si $F\xrightarrow{u}G$ est une flèche de $\mathbf{Fct}(\I_{\backslash X},\E_\kk)$, le
  morphisme $\Omega_X(u) : \Omega_X(F)\to\Omega_X(G)$ de $\mathbf{Fct}(\I,\E_\kk)$ est défini,
  sur l'objet $E$ de $\I$, comme la somme
  directe sur les éléments $x$ de $X(E)$ des morphismes $u(E,x) : F(E,x)\to G(E,x)$.
\end{enumerate} 
\end{prdef}

\begin{ex} On a $\Omega_X(\kk)=\kk[X]$ (linéarisation du foncteur $X$)
  dans $\mathbf{Fct}(\I,\E_\kk)$, où le foncteur $\kk\in {\rm Ob}\,\mathbf{Fct}(\I_{\backslash
    X},\E_\kk)$ du membre de gauche est le foncteur constant.
\end{ex}


\begin{pr}\label{adj-intou} Le foncteur $\Omega_X$ est adjoint à
  gauche à $\Upsilon_X$.
\end{pr}

Cette proposition se vérifie  de façon similaire à
l'adjonction entre les foncteurs  $o$ et $\varpi$ dans $\F_{surj}$.

\begin{pr}\label{preomu}\begin{enumerate}\item L'endofoncteur
    $\Omega_X\Upsilon_X$ de $\mathbf{Fct}(\I,\E_\kk)$ est isomorphe à
    $\cdot\otimes \kk[X]$. Plus généralement, on a un isomorphisme
$$\Omega_X(A\otimes\Upsilon_X (F))\simeq\Omega_X(A)\otimes F$$
naturel en les objets $A$ de $\mathbf{Fct}(\I_{\backslash X},\E_\kk)$ et
$F$ de $\mathbf{Fct}(\I,\E_\kk)$.
\item La coünité  $\Omega_X\Upsilon_X\to id$ de l'adjonction de la
  proposition \ref{adj-intou} s'identifie au produit tensoriel par le
  morphisme d'augmentation $\kk[X]\to\kk$ obtenu par linéarisation de
  l'unique transformation naturelle $X\to *$.
\end{enumerate}
\end{pr}

\begin{pr}\label{ptensco} Les injections diagonales
$$\bigoplus_{x\in X(E)} F(E,x)\otimes
G(E,x)\hookrightarrow\bigoplus_{x,y\in X(E)} F(E,x)\otimes
G(E,y)$$
fournissent un monomorphisme $\Omega_X (F\otimes
G)\hookrightarrow\Omega_X(F)\otimes\Omega_X(G)$ naturel en les objets
$F$ et $G$ de $\mathbf{Fct}(\I_{\backslash X},\E_\kk)$. 
\end{pr}

Cette proposition se vérifie par inspection. Appliquée aux isomorphismes canoniques $F\xrightarrow{\simeq} F\otimes \kk$  de
$\mathbf{Fct}(\I_{\backslash X},\E_\kk)$, elle donne le corollaire suivant.

\begin{cor}\label{crf-com}\begin{enumerate}\item Le foncteur
    $\kk[X]=\Omega_X(\kk)$ est canoniquement une coalgèbre
    cocommutative dans la
    catégorie monoïdale symétrique $\mathbf{Fct}(\I,\E_\kk)$.
\item Pour tout objet $F$ de $\mathbf{Fct}(\I_{\backslash X},\E_\kk)$,
  $\Omega_X(F)$ est naturellement un $\kk[X]$-comodule. Autrement dit,
  on peut compléter le diagramme suivant.
$$\xymatrix{\mathbf{Fct}(\I_{\backslash X},\E_\kk)\ar[r]^-{\Omega_X}\ar@{.>}[dr] & \mathbf{Fct}(\I,\E_\kk) \\
 & \mathbf{Comod}_{\kk[X]}\ar[u]_-{oubli}
}$$
\end{enumerate}
\end{cor}

\begin{rem} La structure de coalgèbre sur $\kk[X]$ généralise celle de
  l'algèbre d'un groupe, puisque sa diagonale est donnée par les
  applications linéaires $[x]\mapsto [x]\otimes [x]$.
\end{rem}

\begin{pr}\label{prf-com} Le foncteur $\Omega_X$ induit une équivalence de
  catégories entre $\mathbf{Fct}(\I_{\backslash X},\E_\kk)$ et $\mathbf{Comod}_{\kk[X]}$.
\end{pr}

\begin{proof} On applique la proposition \ref{monadique}, de sorte
  qu'il suffit de constater qu'un $\kk[X]$-comodule est un comodule sur la comonade déterminée
  par l'adjonction entre $\Omega_X$ et $\Upsilon_X$. Cela découle de
  la proposition \ref{preomu} et du corollaire \ref{crf-com}.
\end{proof}

\begin{rem}\label{rqesc}\begin{enumerate}\item Dans l'équivalence de
    catégories de la proposition, le produit tensoriel de
    $\mathbf{Fct}(\I_{\backslash X},\E_\kk)$ correspond au produit
    cotensoriel de $\kk[X]$-comodules.
\item L'image du foncteur $\Upsilon_X$ correspond aux comodules libres.
\end{enumerate}
\end{rem}

\begin{rem}\label{rqesc2} Les considérations de ce paragraphe sont fonctorielles
  en $X$ en le sens suivant. Toute transformation naturelle de foncteurs ensemblistes $X\to
  X'$ induit un foncteur $\I_{\backslash X}\to\I_{\backslash X'}$,
  donc par précomposition un foncteur $i : \mathbf{Fct}(\I_{\backslash
    X'},\E_\kk)\to\mathbf{Fct}(\I_{\backslash X},\E_\kk)$ (qui généralise
  $\Upsilon_X$). Via
  l'équivalence de catégories de la proposition~\ref{prf-com}, ce foncteur
  correspond à la coïnduction : pour tous objets $F$ de $\mathbf{Fct}(\I_{\backslash
    X'},\E_\kk)$ et $E$ de $\I$, le
  $\kk[X(E)]$-comodule $\Omega_X(i(F))(E)$ est le produit cotensoriel
  du $\kk[X'(E)]$-comodule $\Omega_{X'}(F)(E)$ et de $\kk[X(E)]$ (vu
  comme un $\kk[X'(E)]$-comodule par le morphisme de coalgèbres
  $\kk[X(E)]\to\kk[X'(E)]$ déduit de la transformation naturelle
  $X\to X'$).
\end{rem}

\subsection{Recollements de catégories de comodules}\label{svnc}

Ce paragraphe traite d'une variante non coünitaire des considérations
précédentes, que l'on étudie commodément à partir d'un diagramme de
recollement (proposition~\ref{prfrc}).

\begin{conv} Dans ce paragraphe, $X$ désigne un foncteur de $\I$ vers
la catégorie $\mathbf{Ens}$ et $Y$ un sous-foncteur de $X$.

Nous désignerons par $\I_{X,Y}$ la sous-catégorie
pleine de $\I_{\backslash X}$ dont les objets sont les couples $(E,x)$
pour lesquels $x\in X(E)\setminus Y(E)$.
\end{conv}

On prendra garde que $X\setminus Y$ n'est {\em pas} en général un foncteur
ensembliste.  Nous noterons cependant, par abus, $\kk[X\setminus Y]$ le
conoyau du monomorphisme $\kk[Y]\hookrightarrow \kk[X]$ de $\mathbf{Fct}(\I,E_\kk)$
induit par l'inclusion $Y\hookrightarrow X$.

On a un scindement canonique 
\begin{equation}\label{eqscfl} \kk[X/Y]\simeq \kk[X\setminus Y]\oplus \kk
\end{equation}
obtenu en constatant que le foncteur ensembliste $X/Y$ se relève
canoniquement en un foncteur vers les ensembles pointés (on convient
que $E/\varnothing =E\amalg\{*\}$).

Comme l'épimorphisme $\kk[X/Y]\twoheadrightarrow \kk$ est la coünité du
foncteur en coalgèbres $\kk[X/Y]$, le foncteur $\kk[X\setminus Y]$ est une
coalgèbre sans coünité dans $\mathbf{Fct}(\I,\E_\kk)$ ; nous y
reviendrons à la proposition \ref{vncprf}. 

Dans le lemme suivant, nous nommons {\em sous-catégorie
  complémentaire} d'une sous-catégorie pleine $\B$ d'une catégorie
$\A$ la sous-catégorie pleine de $\A$
dont la classe d'objets est le complémentaire de celle de $\B$.

\begin{lm}\label{lmef-recc} La catégorie $\I_{\backslash Y}$
  s'identifie canoniquement à une sous-catégorie pleine complète à
  droite de
  $\I_{\backslash X}$, dont la catégorie complémentaire est
  $\I_{X,Y}$.
\end{lm}

\begin{proof} Il suffit de remarquer que, si $f : (E,x)\to (E',x')$ est
  un morphisme de $\I_{\backslash X}$, où $x\in Y(E)$, alors on a
  $x'=X(f)(x)=Y(f)(x)\in Y(E')$ puisque $Y$ est un sous-foncteur de $X$.
\end{proof}

Appliquant le corollaire \ref{crf-rec}, on en déduit le
résultat suivant.

\begin{pr}\label{prfrc} Il existe un diagramme de recollement
$$\xymatrix{\mathbf{Fct}(\I_{\backslash Y},\E_\kk)\ar[r]|-{\mathcal{P}} &
  \mathbf{Fct}(\I_{\backslash X},\E_\kk)\ar[r]|-{\mathcal{R}}\ar@/_/[l]\ar@/^/[l]^-{\mathcal{R}} &
  \mathbf{Fct}(\I_{X,Y},\E_\kk)\ar@/_/[l]\ar@/^/[l]^-{\mathcal{P}}}$$
où $\mathcal{P}$ désigne le prolongement par zéro et $\mathcal{R}$ la restriction. 
\end{pr}

Nous adaptons maintenant certains résultats de la section
précédente à la catégorie de foncteurs $\mathbf{Fct}(\I_{X,Y},\E_\kk)$.

\begin{nota}\label{notiota} Nous noterons $\Upsilon_{X,Y}$ le foncteur composé
$$\mathbf{Fct}(\I,\E_\kk)\xrightarrow{\Upsilon_X}\mathbf{Fct}(\I_{\backslash X},\E_\kk)\xrightarrow{\mathcal{R}}\mathbf{Fct}(\I_{X,Y},\E_\kk).$$
\end{nota}

Autrement dit, $\Upsilon_{X,Y}$ est le foncteur de précomposition par
le foncteur d'oubli $\I_{X,Y}\to\I$.

\begin{defi}\label{prdfom2} Le {\em foncteur de
    $X\setminus Y$-intégrale}, noté $\Omega_{X,Y}$, est défini comme la composée
$$\Omega_{X,Y} : \mathbf{Fct}(\I_{X,Y},\E_\kk)\xrightarrow{\mathcal{P}}\mathbf{Fct}(\I_{\backslash X},\E_\kk)\xrightarrow{\Omega_X}\mathbf{Fct}(\I,\E_\kk).$$
\end{defi}

On a donc $\Omega_{X,Y}(F)(E)=\underset{x\in X(E)\setminus
  Y(E)}{\bigoplus}F(E,x)$ pour $F\in {\rm
  Ob}\,\mathbf{Fct}(\I_{X,Y},\E_\kk)$ et $E\in {\rm
  Ob}\,\I$. On prendra garde que les foncteurs $\Upsilon_{X,Y}$ et
$\Omega_{X,Y}$ ne sont généralement pas adjoints.

\begin{rem}\label{rqevnc}\begin{enumerate}\item On a un isomorphisme canonique
    $\Omega_{X,Y}(\kk)\simeq \kk[X\setminus Y]$.
\item Le foncteur $\Omega_{X,Y}$ est exact et fidèle.
\end{enumerate}
\end{rem}

La proposition \ref{preomu} entraîne le résultat suivant, puisque les
foncteurs de restriction et de prolongement par zéro commutent
canoniquement au produit tensoriel.

\begin{pr}\label{preomu2} L'endofoncteur
    $\Omega_{X,Y}\Upsilon_{X,Y}$ de $\mathbf{Fct}(\I,\E_\kk)$ est isomorphe à
    $\cdot\otimes \kk[X\setminus Y]$. Plus généralement, on a un isomorphisme
$$\Omega_{X,Y}(A\otimes\Upsilon_{X,Y} (F))\simeq\Omega_{X,Y}(A)\otimes F$$
naturel en les objets $A$ de $\mathbf{Fct}(\I_{X,Y},\E_\kk)$ et
$F$ de $\mathbf{Fct}(\I,\E_\kk)$.
\end{pr}

De même, la proposition \ref{ptensco} et le corollaire \ref{crf-com} procurent le résultat
suivant.

\begin{pr}\label{vncprf}\begin{enumerate}\item Les injections diagonales fournissent un monomorphisme $\Omega_{X,Y} (F\otimes
G)\hookrightarrow\Omega_{X,Y}(F)\otimes\Omega_{X,Y}(G)$ naturel en les objets
$F$ et $G$ de $\mathbf{Fct}(\I_{X,Y},\E_\kk)$. Celui-ci est
compatible aux isomorphismes d'associativité et de commutativité du
produit tensoriel.
\item Le foncteur $\kk[X\setminus Y]$ est canoniquement une coalgèbre
    cocommutative sans coünité dans la
    catégorie monoïdale symétrique $\mathbf{Fct}(\I,\E_\kk)$ ; le morphisme
    canonique $\kk[X]\twoheadrightarrow \kk[X\setminus Y]$ est un
    morphisme de coalgèbres sans coünité.
\item Pour tout objet $F$ de $\mathbf{Fct}(\I_{X,Y},\E_\kk)$,
  $\Omega_{X,Y}(F)$ est naturellement un $\kk[X\setminus Y]$-comodule.
\end{enumerate}
\end{pr}

\begin{rem} Le scindement (\ref{eqscfl}) (page~\pageref{eqscfl})
  permet d'identifier les catégories de comodules sur $\kk[X\setminus
  Y]$ et $\kk[X/Y]$. On rappelle que les comodules sur une coalgèbre coünitaire, comme
    $\kk[X/Y]$, sont supposés compatibles à la coünité.
\end{rem}

\begin{pr}\label{prf-com2} Le foncteur $\Omega_{X,Y}$ induit une équivalence de
  catégories entre $\mathbf{Fct}(\I_{X,Y},\E_\kk)$
  et la sous-catégorie pleine de $\mathbf{Comod}_{\kk[X\setminus Y]}$
  formée des $\kk[X\setminus Y]$-comodules à droite $(C,\psi_C)$ (que
  nous nommerons {\em fidèles}) tels que la comultiplication $\psi_C$
  est injective, notée $\mathbf{Comod}^{fid}_{\kk[X\setminus Y]}$.
\end{pr}

\begin{proof} On utilise l'identification de la remarque précédente.  Le foncteur de prolongement par zéro et la proposition
  \ref{prf-com} identifient alors $\mathbf{Fct}(\I_{X,Y},\E_\kk)$ à la
  sous-catégorie pleine de $\mathbf{Comod}_{\kk[X\setminus Y]}$ des comodules
  $(C,\psi_C)$ tels que le produit fibré de l'inclusion $C\hookrightarrow
  C\oplus (C\otimes \kk[X\setminus Y])$ et du morphisme $(id_C,\psi_C) : C\to C\oplus (C\otimes \kk[X\setminus Y])$ est
  nul. Ce produit fibré étant égal à $ker\,\psi_C$, cette condition équivaut à la fidélité de
  $\psi_C$, ce qui établit la proposition.
\end{proof}

\subsection{La catégorie de modules $\mathbf{Fct}(\I_{/X},\E_\kk)$}

Nous présentons maintenant la situation duale de celle du
paragraphe~\ref{sctcom}. Toutes les propriétés énoncées se démontrent de façon analogue à
celles du paragraphe~\ref{sctcom}, c'est pourquoi nous les laissons au
lecteur.

\begin{conv} Dans ce paragraphe, $X$ désigne un foncteur de la
catégorie $\I^{op}$ vers la catégorie $\mathbf{Ens}^f$ des ensembles
{\em finis}.
\end{conv}

Nous expliquerons en fin de paragraphe l'utilité de cette restriction
aux ensembles finis.

\begin{nota}  Nous noterons $\I_{/X}$ la
  catégorie $((\I^{op})_{\backslash X})^{op}$. Le foncteur d'oubli
  $\I_{/X}\to\I$ sera noté $\mathcal{O}^{\I,X}$.
\end{nota}

Explicitement, les objets de $\I_{/X}$ sont les couples $(E,x)$ formés
d'un objet $E$ de $\I$ et d'un élément $x$ de $X(E)$. On a $\mathcal{O}^{\I,X}(E,x)=E$ sur les objets.

\begin{nota}\label{dfup3} Nous désignerons par $\Upsilon^X :
  \mathbf{Fct}(\I,\E_\kk)\to\mathbf{Fct}(\I_{/ X},\E_\kk)$ le
  foncteur de précomposition par $\mathcal{O}^{\I,X}$.
\end{nota}


\begin{prdef}\label{prdfom3} Il existe un foncteur exact et fidèle $\Omega^X :
  \mathbf{Fct}(\I_{/ X},\E_\kk)\to\mathbf{Fct}(\I,\E_\kk)$,
  appelé {\em foncteur de
    $X$-intégrale},
  défini de la façon suivante.
\begin{enumerate}\item Si $F$ est un objet de $\mathbf{Fct}(\I_{/ X},\E_\kk)$ et $E$ un
  objet de $\I$, on pose
$$\big(\Omega^X(F)\big)(E)=\bigoplus_{x\in X(E)} F(E,x).$$
\item Si $F$ est un objet de $\mathbf{Fct}(\I_{/ X},\E_\kk)$ et $E\xrightarrow{f}E'$ une
  flèche de $\I$, et si $x$ et $x'$ sont des éléments respectifs
  de $X(E)$ et $X(E')$, la composante $F(E,x)\to F(E',x')$ de
  $\big(\Omega^X(F)\big)(f)$ est égale à $F(f)$ si $X(f)(x')=x$, et à $0$ sinon.
\item Si $F\xrightarrow{u}G$ est une flèche de $\mathbf{Fct}(\I_{/ X},\E_\kk)$, le
  morphisme $\Omega^X(u) : \Omega^X(F)\to\Omega^X(G)$ de $\mathbf{Fct}(\I,\E_\kk)$ est défini,
  sur l'objet $E$ de $\I$, comme la somme
  directe sur les éléments $x$ de $X(E)$ des morphismes $u(E,x) : F(E,x)\to G(E,x)$.
\end{enumerate} 
\end{prdef}

\begin{ex} On a $\Omega^X(\kk)=\kk^X$.
\end{ex}

\begin{pr}\label{adj-intou3} Le foncteur $\Omega^X$ est adjoint à
  droite à $\Upsilon^X$.
\end{pr}

\begin{pr}\label{preomu3}\begin{enumerate}\item L'endofoncteur
    $\Omega^X\Upsilon^X$ de $\mathbf{Fct}(\I,\E_\kk)$ est isomorphe à
    $\cdot\otimes \kk^X$. Plus généralement, on a un isomorphisme
$$\Omega^X(A\otimes\Upsilon^X (F))\simeq\Omega^X(A)\otimes F$$
naturel en les objets $A$ de $\mathbf{Fct}(\I_{/ X},\E_\kk)$ et
$F$ de $\mathbf{Fct}(\I,\E_\kk)$.
\item L'unité  $id\to\Omega^X\Upsilon^X$ de l'adjonction de la
  proposition \ref{adj-intou3} s'identifie au produit tensoriel par le
  morphisme $\kk\to \kk^X$ obtenu en appliquant cette unité au foncteur
  constant $\kk$.
\end{enumerate}
\end{pr}

\begin{pr}\label{ptensco3} Les projections canoniques
$$\bigoplus_{x,y\in X(E)} F(E,x)\otimes
G(E,y) \twoheadrightarrow \bigoplus_{x\in X(E)} F(E,x)\otimes
G(E,x)$$
fournissent un épimorphisme $\Omega^X(F)\otimes\Omega^X(G)\twoheadrightarrow\Omega^X (F\otimes
G)$ naturel en les objets
$F$ et $G$ de $\mathbf{Fct}(\I_{/ X},\E_\kk)$.
\end{pr}

Cette proposition, appliquée aux isomorphismes canoniques $F\otimes \kk\xrightarrow{\simeq} F$  de
$\mathbf{Fct}(\I_{/ X},\E_\kk)$, procure le corollaire suivant.

\begin{cor}\label{crf-com3}\begin{enumerate}\item Le foncteur
    $\kk^X=\Omega^X(\kk)$ est canoniquement une algèbre
    commutative dans la
    catégorie monoïdale symétrique $\mathbf{Fct}(\I,\E_\kk)$.
\item Pour tout objet $F$ de $\mathbf{Fct}(\I_{/ X},\E_\kk)$,
  $\Omega^X(F)$ est naturellement un $\kk^X$-module. Autrement dit,
  on peut compléter le diagramme suivant.
$$\xymatrix{\mathbf{Fct}(\I_{/ X},\E_\kk)\ar[r]^-{\Omega^X}\ar@{.>}[dr] & \mathbf{Fct}(\I,\E_\kk) \\
 & \mathbf{Mod}_{\kk^X}\ar[u]_-{oubli}
}$$
\end{enumerate}
\end{cor}

\begin{pr}\label{prf-com3} Le foncteur $\Omega^X$ induit une équivalence de
  catégories entre $\mathbf{Fct}(\I_{/ X},\E_\kk)$ et $\mathbf{Mod}_{\kk^X}$.
\end{pr}

\begin{rem} L'hypothèse de finitude des
ensembles $X(E)$ permet d'assurer que, dans la définition de $\Omega^X$,
on peut remplacer la somme directe par un produit, ce qui est
nécessaire pour établir la proposition \ref{adj-intou3}. Il est
obligatoire de considérer une somme directe pour disposer des
propriétés de commutation avec le produit tensoriel.
\end{rem}

\section{Les catégories $\E^f_{\Gr,I}$, $\widetilde{\E}^f_{\Gr}$  et
  $\E^f_{\Pl,n}$}\label{sctcatb}

%

Dans cette section, nous introduisons les catégories sources des
{\em catégories de foncteurs en grassmanniennes}. Ces catégories
sources, notées $\E^f_{\Gr,I}(\kk)$, $\widetilde{\E}^f_{\Gr}(\kk)$ et
$\E^f_{\Pl,n}(\kk)$, possèdent une riche structure : de nombreux
foncteurs, avec des propriétés d'adjonction et de composition,
apparaissent très naturellement. Toutes les propriétés s'établissent
de façon directe, c'est pourquoi nous ne démontrons que certaines
d'entre elles.

Une partie des définitions et des résultats de cette section
constituent des illustrations des constructions catégoriques utilisées
dans la section~\ref{sctccf}.

\begin{conv} Dans toute cette section, $I$ désigne une partie de
  $\mathbb{N}$ et $n$ un entier naturel.
\end{conv}

\subsection{Définition des catégories et foncteurs
utilisés}\label{ssctdcfu} On rappelle que $\Gr(V)$ désigne la
grassmannienne des sous-espaces d'un sous-espace vectoriel $V$. Dans
la suite, nous regarderons $\Gr$ comme un foncteur
$\E^f_\kk\to\bf{Ens}$, en définissant $\Gr(f) : \Gr(V)\to\Gr(V')$, où $f : V\to V'$ est un
morphisme de $\E^f$, comme étant la fonction $W\mapsto f(W)$.

On rappelle également que, si $I$ est une partie de $\mathbb{N}$, $\Gr_I(V)$ désigne le sous-ensemble de
$\Gr(V)$ formé des sous-espaces de $V$ dont la dimension appartient à
$I$.

\begin{rem}\label{rqcate1} On notera que $\Gr_I$ est un
  sous-foncteur de $\Gr$ si et seulement si $I$ est de la forme $\leq
  n$ ou $\mathbb{N}$.
\end{rem}

\begin{defi}[Catégories $\E^f_{\Gr,I}(\kk)$ et $\widetilde{\E}^f_{\Gr,
    I}(\kk)$]\label{dfcate}
\begin{enumerate}\item Nous désignerons par  $\widetilde{\E}^f_{\Gr,
    I}(\kk)$ la catégorie
  donnée comme suit.
\begin{itemize}\item Les objets de $\widetilde{\E}^f_{\Gr, I}(\kk)$ sont les couples
  $(V,W)$, où $V$ est un $\kk$-espace vectoriel de dimension finie et $W$ un
  élément de $\Gr_I(V)$.
\item Les morphismes de $(V,W)$ vers $(V',W')$ dans
  $\widetilde{\E}^f_{\Gr, I}(\kk)$ sont les applications
  linéaires $f : V\to V'$ telles que $f(W)\subset W'$.
\item La composition des morphismes s'obtient par
composition des applications linéaires sous-jacentes.
\end{itemize}
\item Nous désignerons par $\E^f_{\Gr, I}(\kk)$ la sous-catégorie de $\widetilde{\E}^f_{\Gr, I}(\kk)$ ayant les
  mêmes objets et dont les morphismes sont donnés par
$${\rm hom}_{\E^f_{\Gr, I}(\kk)}((V,W),(V',W'))=\{f\in {\rm
  hom}_{\widetilde{\E}^f_{\Gr, I}(\kk)}((V,W),(V',W'))\,|\,f(W)=W'\}.$$
\item Nous noterons
  $\widetilde{incl}_I$ le foncteur (fidèle et essentiellement
  surjectif) d'inclusion de $\E^f_{\Gr, I}(\kk)$ dans $\widetilde{\E}^f_{\Gr, I}(\kk)$.
\end{enumerate}

La mention du corps $\kk$ sera souvent omise.

De plus, lorsque la partie $I$ est égale à $\mathbb{N}$, nous omettrons
l'indice $I$ dans la notation de ces catégories et des foncteurs où
elles interviennent. 

Par exemple, nous noterons $\wt{incl} : \E^f_\Gr\to\wt{\E}^f_\Gr$ pour $\wt{incl}_\mathbb{N} : \E^f_{\Gr,\mathbb{N}}(\kk)\to\wt{\E}^f_{\Gr,\mathbb{N}}(\kk)$.
\end{defi}

Dans la suite, nous ne considérerons la catégorie
$\widetilde{\E}^f_{\Gr,I}(\kk)$ que lorsque $I=\mathbb{N}$. Nous ne
donnerons donc la plupart des définitions et des propriétés des foncteurs où cette catégorie
intervient que dans ce cadre, mais beaucoup d'entre elles se
généralisent sans difficulté au cas de 
$\widetilde{\E}^f_{\Gr,I}(\kk)$ pour $I\subset\mathbb{N}$ quelconque.

\begin{nota}[Catégorie $\E^f_{\Pl,n}(\kk)$]\label{dfcatpl} 
\begin{enumerate}\item \'Etant donné un objet $V$ de $\E^f_\kk$, nous
  noterons $\Pl_n(V)$ l'ensemble ${\rm Pl}_{\E_\kk^f} (E_n,V)$ des monomorphismes
  $E_n\hookrightarrow V$.
\item Nous désignerons par
  $\E^f_{\Pl,n}(\kk)$ la catégorie définie
  ainsi.
\begin{itemize}\item Les objets de $\E^f_{\Pl,n}(\kk)$ sont les couples
  $(V,u)$ formés d'un objet $V$ de $\E^f_\kk$ et d'un élément $u$ de $\Pl_n(V)$.
\item Les morphismes dans $\E^f_{\Pl,n}(\kk)$ de $(V,u)$ vers $(V',u')$
  sont les applications linéaires $f : V\to V'$ faisant commuter le
  diagramme
$$\xymatrix{E_n\ar[r]^-u\ar[rd]_-{u'} & V\ar[d]^-f \\
 & V'
}$$
\item La composition des morphismes est induite par la
composition des applications linéaires.
\end{itemize}
\end{enumerate}
\end{nota}

\begin{rem}\label{rqfece}\begin{enumerate}\item Si $I$ est une partie du type $\leq n$ ou
  $\mathbb{N}$, alors $\E^f_{\Gr,I}(\kk)$ est la catégorie
  $(\E^f_\kk)_{\backslash\Gr_I}$ (cf. notation~\ref{notfcom}).
\item Avec les notations du paragraphe \ref{svnc}, $\E^f_{\Gr,n}$ est la catégorie
  $\E^f_{\Gr_{\leq n}, \Gr_{\leq n-1}}$.
\item L'association $V\mapsto\Pl_n(V)$ n'est pas fonctorielle, mais si
  l'on note ${\rm hom}_{\leq i} (E,V)$ le sous-ensemble de ${\rm
    hom}_\E (E,V)$ formé des applications linéaires de rang au plus $i$, alors ${\rm hom}_{\leq i} (E_n,.)$ est un
  sous-foncteur de ${\rm hom}_{\E^f} (E_n,.)$. On peut voir
  $\E^f_{\Pl,n}$ comme la catégorie $\E^f_{{\rm
      hom}_{\E^f} (E_n,.), {\rm hom}_{\leq n-1} (E_n,.)}$.
\item L'ensemble $\Gr_n(V)$ s'identifie canoniquement au quotient de
  $\Pl_n(V)$ par l'action à droite libre du groupe $GL_n(\kk)$.
\item Les catégories  $\E^f_{\Gr,0}$ et $\E^f_{\Pl,0}$ s'identifient canoniquement à $\E^f$.
\item Les catégories $\E^f_{\Gr,1}(\FF)$ et $\E^f_{\Pl,1}(\FF)$ sont isomorphes.
\item Toutes les catégories introduites vérifient l'hypothèse
  \ref{hypf3} (finitude des ensembles de morphismes), car $\kk$ est fini.
\end{enumerate}
\end{rem}


\begin{nota}[Foncteurs d'oubli]\label{notfctou}\begin{enumerate}\item
    Soit $J$ une partie de $\mathbb{N}$ telle que $I\subset
    J$. On note $incl_{I,J} : \E^f_{\Gr, I}\to\E^f_{\Gr, J}$ le foncteur (pleinement fidèle) d'oubli.
\item Nous désignerons par
  $inc^{\Pl}_n$ le foncteur (fidèle et essentiellement
  surjectif) d'oubli du plongement $\E^f_{\Pl,n}\to\E^f_{\Gr,n}$ associant à un objet
  $(V,\xymatrix{E_n\ar@{^{(}->}[r]^-u & V})$ de $\E^f_{\Pl,n}$ l'objet
    $(V,im\,u)$ de $\E^f_{\Gr,n}$ et égal à l'inclusion évidente sur les morphismes.
\end{enumerate}
\end{nota}

Comme nous le verrons en fin de section, les catégories introduites
précédemment ne sont en général pas abéliennes. Une partie de la
structure abélienne de la catégorie $\E^f$ s'y reflète cependant,
grâce à la notion suivante.

\begin{defi}[Foncteurs de translation]\label{dfftr}  On définit des  foncteurs  $\E^f\times\E^f_{\Gr,I}\to\E^f_{\Gr,I}$,
 $\E^f\times\widetilde{\E}^f_{\Gr}\to\widetilde{\E}^f_{\Gr}$ et
$\E^f\times\E^f_{\Pl,n}\to\E^f_{\Pl,n}$, appelés {\em foncteurs de
  translation}, et notés
$\boxplus$, par $V\boxplus (A,B)=(V\oplus A,B)$
(dans les deux premiers cas) et $V\boxplus (A, u : E_n\hookrightarrow
A)=(V\oplus A, E_n\xrightarrow{u}A\hookrightarrow V\oplus A)$ (dans le
troisième) sur les objets, l'action sur les morphismes se déduisant de
la fonctorialité de $\oplus : \E^f\times\E^f\to\E^f$. 
\end{defi}

\begin{rem}\label{remfe-cate} On a des isomorphismes $0\boxplus X\simeq X$ et $(V'\oplus V)\boxplus
  X\simeq V'\boxplus (V\boxplus X)$ naturels en les objets $V$, $V'$ de
  $\E^f$ et $X$ de $\E^f_{\Gr,I}$ (resp. $\widetilde{\E}^f_{\Gr}$,
  $\E^f_{\Pl,n}$), vérifiant des propriétés de cohérence qu'on
  laisse au lecteur le soin d'expliciter. Ainsi, on peut voir $\boxplus$ comme un
  foncteur d'action de la catégorie additive $\E^f$ sur $\E^f_{\Gr,I}$ (resp. $\widetilde{\E}^f_{\Gr}$,
  $\E^f_{\Pl,n}$).
\end{rem}

\begin{defi}[Foncteurs fondamentaux de source $\E^f_{\Gr,I}$,
  $\widetilde{\E}^f_{\Gr,I}$ ou $\E^f_{\Pl,n}$]\label{df-ff}
\begin{enumerate}\item Le {\em foncteur d'oubli
    principal} $\widetilde{\mathfrak{O}}_I :
  \widetilde{\E}^f_{\Gr,I}\to\E^f$ est le foncteur  associant à un objet $(V,W)$ de $\widetilde{\E}^f_{\Gr,I}$ l'espace vectoriel $V$, et à un morphisme
  l'application linéaire sous-jacente.
\item On appelle également {\em foncteurs d'oubli principaux} les
    foncteurs composés
$$\mathfrak{O}_I :
\E^f_{\Gr,I}\xrightarrow{\widetilde{incl}_I}\widetilde{\E}^f_{\Gr,I}\xrightarrow{\widetilde{\mathfrak{O}}_I}\E^f$$
et
$$\bar{\mathfrak{O}}_n :
  \E^f_{\Pl,n}\xrightarrow{inc_n^\Pl}\E^f_{\Gr,n}\xrightarrow{\mathfrak{O}_n}\E^f.$$
\item Le foncteur {\em base} $\mathfrak{B}_I
  : \E^f_{\Gr,I}\to\E^I_{surj}$ est
  défini sur les objets par $\mathfrak{B}_I(V,W)=W$, et associe à un
  morphisme $f: (V,W)\to (V',W')$ l'application linéaire (surjective par
  définition de $\E^f_{\Gr,I}$) $W\to W'$ induite par $f$.
\item Le foncteur {\em d'oubli secondaire} $\widetilde{\mathfrak{B}} :
  \widetilde{\E}^f_\Gr\to\E^f$ associe à un objet $(V,W)$ l'espace vectoriel $W$ et à un morphisme $(V,W)\to (V',W')$ l'application linéaire induite $W\to W'$.
\item Le foncteur de {\em
    réduction}
$\widetilde{\mathfrak{K}}_I : \widetilde{\E}^f_{\Gr,I}\to\E^f$ est donné par $(V,W)\mapsto V/W$ sur les objets et
  associe à un morphisme l'application linéaire induite.
\item On appelle également {\em foncteurs de réduction} les
    foncteurs composés
$$\mathfrak{K}_I :
\E^f_{\Gr,I}\xrightarrow{\widetilde{incl}_I}\widetilde{\E}^f_\Gr\xrightarrow{\widetilde{\mathfrak{K}}_I}\E^f$$
et
$$\bar{\mathfrak{K}}_n :
  \E^f_{\Pl,n}\xrightarrow{inc_n^\Pl}\E^f_{\Gr,n}\xrightarrow{\mathfrak{K}_n}\E^f.$$
\end{enumerate}

On rappelle que l'indice $I$ sera omis dans toutes ces notations
lorsque $I=\mathbb{N}$.
\end{defi}

\begin{rem}\label{reme1}\begin{enumerate}\item Soit $J$ une partie de $\mathbb{N}$ telle que $I\subset J$. Le foncteur
    composé
    $\E^f_{\Gr,I}\xrightarrow{incl_{I,J}}\E^f_{\Gr,J}\xrightarrow{\mathfrak{O}_J}\E^f$ est égal à $\mathfrak{O}_I$. De même, le diagramme
$$\xymatrix{\E^f_{\Gr,I}\ar[r]^-{\mathfrak{B}_I}\ar@{^{(}->}[d]_{incl_{I,J}} & \E^I_{surj}\ar@{^{(}->}[d] \\
\E^f_{\Gr,J}\ar[r]^-{\mathfrak{B}_J} & \E^J_{surj}
}$$
commute. On a d'autres propriétés analogues de compatibilité à l'extension de la
partie $I$ avec les différents foncteurs introduits.
\item Le diagramme
$$\xymatrix{\E^f_{\Gr}\ar[d]_{\mathfrak{B}}\ar[r]^-{\widetilde{incl}} & \widetilde{\E}^f_\Gr\ar[d]^{\widetilde{\mathfrak{B}}}
  \\
\E^f_{surj}\ar[r]_-{oubli} & \E^f
}$$
commute.
\item Les foncteurs base ou d'oubli secondaire n'ont pas d'analogue
  non trivial en termes de la catégorie $\E^f_{\Pl,n}$ ; le rôle de
  cette catégorie est justement de simplifier certaines des
  considérations relatives à $\E^f_{\Gr,n}$ en \go rendant la base
  canoniquement isomorphe à $E_n$\gf.
\end{enumerate}
\end{rem}

\begin{rem} Le foncteur
    $\mathfrak{K}_I : \E^f_{\Gr,I}\to\E^f$ est le conoyau de la transformation naturelle
    injective tautologique $\mathfrak{B}_I\to\mathfrak{O}_I$ (où l'on
    note par abus $\mathfrak{B}_I$ pour la composée de ce foncteur
    avec l'inclusion $\E^I_{surj}\to\E^f$) ; un
    constat analogue vaut pour $\widetilde{\mathfrak{K}}_I$ et $\bar{\mathfrak{K}}_n$.
\end{rem}

La proposition suivante établit des liens entre les foncteurs de
translation et les foncteurs introduits dans la définition~\ref{df-ff}.

\begin{pr}\label{actc-prf} Il existe des isomorphismes
$${\rm hom}_{\E^f_{\Gr,I}}(V\boxplus X,Y)\simeq {\rm
  hom}_{\E^f}(V,\mathfrak{O}_I(Y))\times {\rm
  hom}_{\E^f_{\Gr,I}}(X,Y)\quad\text{et}$$
$${\rm hom}_{\E^f_{\Gr,I}}(X,V\boxplus Y)\simeq {\rm
  hom}_{\E^f}(\mathfrak{K}_I(X),V)\times {\rm
  hom}_{\E^f_{\Gr,I}}(X,Y)$$ 
naturels en les objets $V$ de $\E^f$ et $X$, $Y$ de $\E^f_{\Gr,I}$.

On a des énoncés similaires dans $\widetilde{\E}^f_\Gr$ et $\E^f_{\Pl,n}$.
\end{pr}

%
\begin{proof} Les deux isomorphismes étant très analogues, nous nous
  bornerons à montrer le second, qui s'obtient par la suite
  d'isomorphismes naturels
$${\rm hom}_{\E^f_{\Gr,I}}(X,V\boxplus Y)\simeq \{f\in {\rm
  hom}_{\E^f}(\mathfrak{O}_I(X),V\oplus\mathfrak{O}_I(Y))\,|\,f(\mathfrak{B}_I(X))=\mathfrak{B}_I(Y)\}\simeq$$
$$\{(a,b)\in {\rm
  hom}_{\E^f}(\mathfrak{O}_I(X),V)\times {\rm
  hom}_{\E^f}(\mathfrak{O}_I(X),\mathfrak{O}_I(Y))\,|\,a(\mathfrak{B}_I(X))=0\;\text{et}\;b(\mathfrak{B}_I(X))=\mathfrak{B}_I(Y)\}$$
$$\simeq {\rm hom}_{\E^f} (\mathfrak{K}_I(X),V)\times {\rm
  hom}_{\E^f_{\Gr,I}}(X,Y).$$
\end{proof}

\begin{cor}\label{adjkaio} Supposons que $I$  contient $0$. Le foncteur d'inclusion
  $incl_{0,I} : \E^f\simeq\E^f_{\Gr,0}\hookrightarrow\E^f_{\Gr,I}$ est adjoint à droite à $\mathfrak{K}_I$.
\end{cor}

\begin{defi}[Foncteurs fondamentaux de but $\E^f_{\Gr,I}$,
  $\widetilde{\E}^f_\Gr$ ou $\E^f_{\Pl,n}$]\label{df-ff2}
\begin{enumerate}\item Les foncteurs de {\em plongement
    diagonal} sont les foncteurs $\mathfrak{D}_I :
  \E^I_{surj}\to\E^f_{\Gr,I}$ et $\widetilde{\mathfrak{D}} :
  \E^f\to\widetilde{\E}^f_\Gr$ donnés sur les objets par $V\mapsto
  (V,V)$ et par le plongement évident sur les morphismes.
\item Les foncteurs de {\em plongement relatif} sont les foncteurs composés
$$\mathfrak{L}_I :
\E^f\times\E^I_{surj}\xrightarrow{\E^f\times\mathfrak{D}_I}\E^f\times\E^f_{\Gr,I}\xrightarrow{\boxplus}\E^f_{\Gr,I}$$
et
$$\widetilde{\mathfrak{L}} :
  \E^f\times\E^f\xrightarrow{\E^f\times\widetilde{\mathfrak{D}}}\E^f\times\widetilde{\E}^f_\Gr\xrightarrow{\boxplus}\widetilde{\E}^f_\Gr.$$

Ces foncteurs sont donc donnés sur les objets par $(A,B)\mapsto (A\oplus B,B)$.
\item Le foncteur de {\em décalage
    pointé} $\mathfrak{S}_n : \E^f\to\E^f_{\Pl,n}$ associe à un
  objet $V$ de $\E^f$ l'objet $(V\oplus E_n,E_n\hookrightarrow V\oplus
  E_n)$ et à une application linéaire $u$ le morphisme $u\oplus E_n$.
\end{enumerate}
\end{defi}

\begin{rem}\label{rqidit}\begin{enumerate}\item Le diagramme
$$\xymatrix{\E^f_{\Gr}\ar[r]^-{\widetilde{incl}} & \widetilde{\E}^f_\Gr
  \\
\E^f_{surj}\ar[u]^{\mathfrak{D}}\ar[r]^-{oubli} & \E^f\ar[u]_{\widetilde{\mathfrak{D}}}
}$$
commute.
\item Le diagramme
$$\xymatrix{\E^f\ar@{^{(}->}[r]\ar[d]_{\mathfrak{S}_n} & \E^f\times\E^n_{surj}\ar[d]^{\mathfrak{L}_n} \\
\E^f_{\Pl,n}\ar[r]^-{inc^\Pl_n} & \E^f_{\Gr,n}
}$$
commute, où l'inclusion supérieure est donnée par $V\mapsto (V,E_n)$.
\end{enumerate}
\end{rem}

\begin{rem}\label{remdffp} Le plongement évident $\underline{GL_n(\kk)}\to\E^n_{surj}(\kk)$ étant une
  équivalence de catégories, nous commettrons parfois des abus
  de notation consistant à l'assimiler à une égalité pour les foncteurs mettant en jeu $\E^n_{surj}(\kk)$. 
\end{rem}

\subsection{Propriétés des foncteurs fondamentaux}\label{ssctpff} Nous
commençons par donner des propriétés d'adjonction entre les foncteurs
introduits précédemment qui seront utilisées de façon intensive
dans tout l'article.

\begin{pr}\label{adjfpr}
\begin{enumerate}\item Le foncteur de plongement diagonal
  $\mathfrak{D}_I : \E^I_{surj}\to\E^f_{\Gr,I}$ est adjoint
  à gauche à $\mathfrak{B}_I : \E^f_{\Gr,I}\to\E^I_{surj}$.
\item Le foncteur de plongement diagonal $\widetilde{\mathfrak{D}} : \E^f\to\wt{\E}^f_\Gr$ est adjoint
  à gauche à $\widetilde{\mathfrak{B}} :
  \widetilde{\E}^f_\Gr\to\E^f$.
\item Le foncteur de
    plongement relatif $\mathfrak{L}_{I} :
    \E^f\times\E^I_{surj}\to\E^f_{\Gr,I}$ est adjoint à gauche au
    foncteur $\mathfrak{O}_I\times\mathfrak{B}_I : \E^f_{\Gr,I}\to\E^f\times\E^I_{surj}$.
\item Le foncteur de
    plongement relatif $\widetilde{\mathfrak{L}} :
    \E^f\times\E^f\to\widetilde{\E}^f_\Gr$ est adjoint à gauche au
    foncteur $\widetilde{\mathfrak{O}}\times\widetilde{\mathfrak{B}} : \widetilde{\E}^f_\Gr\to\E^f\times\E^f$.
\item Le décalage pointé $\mathfrak{S}_n : \E^f\to\E^f_{\Pl,n}$ est
  adjoint à gauche à $\bar{\mathfrak{O}}_n : \E^f_{\Pl,n}\to\E^f$.
\end{enumerate}
\end{pr}

\begin{proof} Si $A$ est un objet de $\E^I_{surj}$ et
  $(V,W)$ un objet de $\E^f_{\Gr,I}$, on a un isomorphisme
  canonique
$${\rm hom}_{\E^f_{\Gr,I}}(\mathfrak{D}_{I}(A),(V,W))=\{f\in {\rm
  hom}_{\E^f}(A,V)\,|\,f(A)=W\}\simeq {\rm Epi}_{\E^f} (A,W),$$
ce qui démontre la première assertion. La troisième assertion s'en
déduit en utilisant la proposition \ref{actc-prf}.
 
Les deuxième et quatrième points se traitent pareillement. 

\'Etablissons le dernier : si $A$ est un objet de $\E^f$ et $(V, u
: E_n\hookrightarrow V)$ un objet de $\E^f_{\Pl,n}$, on a un
isomorphisme canonique
$${\rm hom}_{\E^f_{\Pl,n}}(\mathfrak{S}_n(A),(V,u))=\{f\in {\rm
  hom}_{\E^f} (A\oplus E_n,V)\,|\,f_{|E_n}=u\}\simeq {\rm hom}_{\E^f} (A,V).$$
Cela achève la démonstration. 
\end{proof}

La proposition suivante, laissée au lecteur, jouera un rôle fondamental par la suite.

\begin{pr}[Compositions fondamentales]\label{compf-prc}  Les foncteurs composés
$$\E^f\times\E^I_{surj}\xrightarrow{\mathfrak{L}_I}\E^f_{\Gr,I}\xrightarrow{\mathfrak{K}_I\times\mathfrak{B}_I}\E^f\times\E^I_{surj}$$
et
$$\E^f\times\E^f\xrightarrow{\widetilde{\mathfrak{L}}}\widetilde{\E}^f_\Gr\xrightarrow{\widetilde{\mathfrak{K}}\times\widetilde{\mathfrak{B}}}\E^f\times\E^f$$
sont canoniquement isomorphes aux foncteurs identités.
\end{pr}

\begin{pr}\label{crf-catec} Supposons $I$ non vide.
\begin{enumerate}\item Les foncteurs d'oubli principaux sont fidèles.
\item Le foncteur $\bar{\mathfrak{O}}_n :
  \E^f_{\Pl,n}\to\E^f$ induit un foncteur essentiellement surjectif
  $\E^f_{\Pl,n}\to\E^{\geq n}$, où le but désigne la sous-catégorie
  pleine de $\E^f$ des espaces de dimension au moins $n$. Si $I$ a pour plus
  petit élément $n$, $\mathfrak{O}_I : \E^f_{\Gr,I}\to\E^f$ induit un foncteur essentiellement surjectif
  $\E^f_{\Gr,I}\to\E^{\geq n}$. Le foncteur $\widetilde{\mathfrak{O}} :
  \widetilde{\E}^f_\Gr\to\E^f$ est essentiellement surjectif.
\item Les foncteurs de plongement relatif
  $\mathfrak{L}_I$ et $\widetilde{\mathfrak{L}}$ sont
  fidèles et essentiellement surjectifs. Il en est de même pour le
  décalage pointé $\mathfrak{S}_n$. 
\item Les foncteurs
  $\mathfrak{K}_I\times\mathfrak{B}_I$ (donc en particulier
  $\mathfrak{K}_I$ et $\mathfrak{B}_I$) et
  $\widetilde{\mathfrak{K}}\times\widetilde{\mathfrak{B}}$ (donc en particulier
  $\widetilde{\mathfrak{K}}$ et $\widetilde{\mathfrak{B}}$) sont
  pleins et essentiellement surjectifs.
\item Les foncteurs de plongement diagonal $\mathfrak{D}_I$ et
  $\widetilde{\mathfrak{D}}$ sont pleinement fidèles.
\end{enumerate}
\end{pr}

\begin{proof} Les deux premières assertions s'établissent par
  inspection. L'essentielle surjectivité des foncteurs de plongement
  relatif et de décalage pointé découle de ce que tout sous-espace d'un
  espace vectoriel est facteur direct ; leur fidélité est claire. La
  quatrième assertion s'obtient en combinant la proposition
  \ref{compf-prc} et l'essentielle surjectivité des plongements
  relatifs. La dernière s'obtient à partir des
  deux premières adjonctions de la proposition \ref{adjfpr} et du
  constat que leurs unités $id\to\mathfrak{B}_I\mathfrak{D}_I$ et
  $id\to\widetilde{\mathfrak{B}}\widetilde{\mathfrak{D}}$ sont
  des isomorphismes. 
\end{proof}

\subsection{Propriétés de structure des catégories $\E^f_\Gr$,
  $\widetilde{\E}^f_{\Gr}$ et $\E^f_{\Pl,n}$}\label{ssctpsci} Les catégories qui nous intéresseront le plus par la suite sont les
$\E^f_{\Gr,I}$. Les propriétés qui suivent montrent que, par
certains côtés, les catégories $\widetilde{\E}^f_\Gr$ de
$\E^f_{\Pl,n}$ ont une
structure plus \go régulière\gf, c'est pourquoi nous
serons parfois amenés à travailler dans ces catégories auxiliaires.

\begin{pr}\label{pr-tilad} La catégorie
  $\widetilde{\E}^f_{\Gr}$ est additive et $\kk$-linéaire. Les foncteurs
  d'oubli principal et secondaire sont additifs --- autrement dit, on
  a un isomorphisme $(V,W)\oplus (V',W')\simeq (V\oplus V',W\oplus
  W')$ naturel en les objets $(V,W)$ et $(V',W')$ de
  $\widetilde{\E}^f_{\Gr}$. En particulier, on a un
  isomorphisme $E\boxplus (V,W)\simeq (E,0)\oplus (V,W)$ naturel en
  les objets $E$ de $\E^f$ et $(V,W)$ de $\widetilde{\E}^f_{\Gr}$.
\end{pr}

\begin{rem} La catégorie $\wt{\E}^f_\Gr$ n'est pas abélienne. En
  effet, on vérifie que pour tout $V\in {\rm Ob}\,\E^f$, le morphisme
  $(V,0)\to (V,V)$ de $\wt{\E}^f_\Gr$ dont l'application linéaire
  sous-jacente est l'identité est à la fois un monomorphisme et un
  épimorphisme. En revanche, ce n'est pas un isomorphisme, si $V$ est
  non nul.
\end{rem}

\begin{prdef}[Dualité dans $\widetilde{\E}^f_\Gr$]\label{dfgrt} Le
  foncteur de dualité $(\cdot)^* : (\E^f)^{op}\to\E^f$ induit une
  équivalence de catégories $(\cdot)^\vee :
  (\widetilde{\E}_\Gr^f)^{op}\to\widetilde{\E}_\Gr^f$ donnée sur les
  objets par $(V,W)^\vee=(V^*,W^\perp)$\index{Nota}{$(V,W)^\vee$}.
\end{prdef}


\begin{pr}\label{pr-cfpl} Soit $n\in\mathbb{N}$. 
\begin{enumerate}\item\begin{enumerate}\item La catégorie  $\E^f_{\Pl,n}$ possède des sommes
  finies.
\item Son objet initial est $(E_n,id_{E_n})$.
\item La somme $(A,a)\amalg (B,b)$ de deux objets $(A,a)$ et $(B,b)$ de
  $\E^f_{\Pl,n}$ s'obtient en formant le carré cocartésien
  d'inclusions
$$\xymatrix{E_n\ar[r]^a\ar[d]_b & A\ar[d] \\
B\ar[r] & A\underset{E_n}{\oplus}B}$$
et en munissant l'espace vectoriel $A\underset{E_n}{\oplus}B$ du
plongement donné par la diagonale du carré.
\end{enumerate}
\item Deux objets $(A,a)$ et $(B,b)$ de
  $\E^f_{\Pl,n}$ possèdent toujours un  produit. Il est donné par
  $$(A,a)\times (B,b)=(A\oplus B,(a,b) : E_n\hookrightarrow A\oplus B).$$
\item On a des isomorphismes 
$$\mathfrak{S}_n(V)\amalg X\simeq V\boxplus X\quad\text{et}$$
$$\mathfrak{S}_n(V)\times X\simeq\mathfrak{S}_n(V\oplus\bar{\mathfrak{O}}_n(X))$$
naturels en les objets $V$ de $\E^f$ et $X$ de $\E^f_{\Pl,n}$.
\end{enumerate}
\end{pr}

La démonstration des propositions \ref{pr-tilad}, \ref{dfgrt}
et~\ref{pr-cfpl} est laissée au lecteur. La proposition~\ref{pr-cfpl}
ne sera d'ailleurs pas utilisée explicitement.

\begin{rem}\label{rqpathfct}\begin{enumerate}\item Si $I$  contient $0$, alors $(0,0)$ est objet final de
    $\E^f_{\Gr,I}$, et tout objet de $\E^f_{\Gr,I}$ admet un produit avec $(V,0)$ (où $V\in {\rm Ob}\,\E^f$),
    qui est donné par le foncteur $V\boxplus\cdot$. 
\item En revanche, on vérifie facilement que deux objets $(V,W)$ et
  $(V',W')$ de $\E^f_{\Gr,I}$ tels que $W$ et $W'$ sont non
  nuls ne possèdent jamais de somme ni de
  produit.
\end{enumerate}
\end{rem}

Les propriétés suivantes fournissent un substitut à l'absence de
sommes et de produits dans $\E^f_\Gr$ ; elles reposent sur le lemme
 simple et très utile suivant.

\begin{lm}\label{lmclei} Il existe une bijection
$${\rm hom}_{\E^f}(A,A')\simeq\coprod_{B'\in\Gr(A')} {\rm
  hom}_{\E^f_\Gr}((A,B),(A',B'))$$
naturelle en l'objet $(A,B)$ de $\E^f_\Gr$ et l'objet $A'$ de $\E^f$. 
\end{lm}

La fonctorialité doit être comprise dans le sens suivant :
\begin{itemize}\item pour le terme de gauche, on considère le foncteur
$$(\E^f_\Gr)^{op}\times\E^f\xrightarrow{(\mathfrak{O})^{op}\times\E^f}(\E^f)^{op}\times\E^f\xrightarrow{{\rm hom}}\bf{Ens}\,;$$
\item pour le terme de droite, la fonctorialité en $(A,B)$ provient de
  manière usuelle du foncteur hom ; pour la fonctorialité en $A'$, on
  fait correspondre à une application linéaire $u : A'\to A''$ et à un
  élément $f$ de ${\rm hom}_{\E^f_\Gr}((A,B),(A',B'))$ (où
  $B'\in\Gr(A')$) la flèche $(A,B)\to (A'',B'')$, où $B''=u(B')$,
  donnée par $u\circ f : A\to A''$.
\end{itemize}

\begin{proof} Cette bijection s'obtient en faisant correspondre à une
  application linéaire $f : A\to A'$ le sous-espace $B'=f(B)$ de $A'$
  et le morphisme $(A,B)\to (A',B')$ de $\E^f_\Gr$ induit par $f$.
\end{proof}

\begin{pr}\label{prsef1} Il existe une bijection
$${\rm hom}_{\E^f_\Gr}((A\oplus A',B\oplus
B'),(V,W))\simeq$$
$$\underset{W_1+W_2=W}{\coprod_{W_1,W_2\in\Gr(W)}}{\rm
  hom}_{\E^f_\Gr}((A,B),(V,W_1))\times {\rm
  hom}_{\E^f_\Gr}((A',B'),(V,W_2))$$
 naturelle en les objets $(A,B)$, $(A',B')$ et $(V,W)$ de $\E^f_\Gr$.
\end{pr}

Cette proposition, comme la proposition~\ref{prsef2} ci-après, est
laissée au lecteur. La fonctorialité doit être comprise dans le sens suivant :
\begin{itemize}\item sur les ensembles hom, on utilise la
  fonctorialité usuelle ;
\item pour $\E^f_\Gr\times\E^f_\Gr\to\E^f_\Gr\quad ((A,B),(A',B'))\mapsto (A\oplus A',B\oplus
B')$, on associe à un morphisme $(u,v)$ de $\E^f_\Gr\times\E^f_\Gr$ le
morphisme $u\oplus v$ de $\E^f_\Gr$ ;
\item dans le terme de droite, la fonctorialité en $(V,W)$ s'obtient
  comme suit. Si $u: (V,W)\to (V',W')$ est un morphisme de $\E^f_\Gr$
  et $W_1, W_2$ deux sous-espaces de $W$ tels que $W_1+W_2=W$, on pose
  $W'_i=f(W_i)$ ($i\in\{1,2\}$), de sorte que $W'_1+W'_2=f(W)=W'$. Le
  morphisme induit par $u$ s'obtient par somme sur les $(W_1,W_2)$ des
  morphismes ${\rm
  hom}_{\E^f_\Gr}((A,B),(V,W_1))\times {\rm
  hom}_{\E^f_\Gr}((A',B'),(V,W_2))\to {\rm
  hom}_{\E^f_\Gr}((A,B),(V',W'_1))\times {\rm
  hom}_{\E^f_\Gr}((A',B'),(V',W'_2))$ induits par $u$.
\end{itemize}

%

\begin{nota}\label{notgrd} Soient $V$ et $W$ deux espaces vectoriels
  de dimension finie. Nous noterons $Gr(V,W)$\index{Nota}{Gr@$Gr$} le sous-ensemble de $\Gr(V\oplus W)$
  formé des sous-espaces $E$ de $V\oplus W$ tels que
  les morphismes \mbox{$E\hookrightarrow V\oplus W\twoheadrightarrow V$} et
\mbox{$E\hookrightarrow V\oplus W\twoheadrightarrow W$} soient
surjectifs.

On définit ainsi un foncteur $Gr :
\E^f_{surj}\times\E^f_{surj}\to\bf{Ens}$, l'action sur les morphismes
étant obtenue par un biais analogue à celui détaillé précédemment.
\end{nota}

\begin{pr}\label{prsef2} Il existe une bijection
$${\rm hom}_{\E^f_\Gr}((V,W),(A,B))\times {\rm
  hom}_{\E^f_\Gr}((V,W),(A',B'))$$
$$\simeq\coprod_{C\in Gr(B,B')} {\rm hom}_{\E^f_\Gr}((V,W),(A\oplus A',C))$$
naturelle en les objets $(A,B)$, $(A',B')$ et $(V,W)$ de $\E^f_\Gr$.
\end{pr}

La fonctorialité repose ici sur celle de $Gr$ (d'une manière similaire
à celle explicitée pour la proposition \ref{prsef1}).


\begin{rem} Les propositions \ref{prsef1} et \ref{prsef2} sont
  spécifiques à la catégorie \go globale\gf $\E^f_{\Gr}$ ; elles n'ont
  pas d'analogue dans $\E^f_{\Gr,n}$, par exemple.
\end{rem}

\part{Les catégories de foncteurs en grassmanniennes}\label{p-dfg}

Le principal sujet de cet article réside dans l'étude de la catégorie
$\F_\Gr(\kk)$, appelée parfois {\em catégorie de foncteurs en
  grassmanniennes globale}, et de ses sous-catégories $\F_{\Gr,n}(\kk)$, où
$n\in\mathbb{N}$. Ces catégories possèdent un intérêt intrinsèque en
raison de leur riche structure algébrique ; elles s'interprètent
notamment en termes de modules ou de comodules. De plus, nous verrons dans la partie~\ref{p-omeg} que le {\em foncteur
  d'intégrale en grassmanniennes} $\omega : \F_\Gr(\kk)\to\F(\kk)$ constitue un outil
puissant d'étude de la catégorie $\F(\kk)$.

La section~\ref{cfgri} décrit la structure de base de la catégorie
$\F_\Gr$. Après avoir introduit les foncteurs fondamentaux reliant
cette catégorie aux autres catégories de foncteurs que nous avons
introduites, notamment $\F$ et $\F_{surj}$, elle décrit les groupes de Grothendieck de ses objets finis et
projectifs de type fini. Dans la section~\ref{fgm}, nous présentons
une approche monadique de la catégorie $\F_\Gr$. Enfin, la
section~\ref{fhifd} présente quelques propriétés générales des adjoints au
produit tensoriel dans les catégories $\F_\Gr$ et~$\F_{\Gr,n}$.

\medskip

Les autres catégories de foncteurs en grassmanniennes que nous
introduisons  serviront surtout d'auxiliaires dans l'étude de la catégorie
$\F_\Gr$. Dans la section~\ref{cfgrt}, nous présentons quelques
propriétés d'une catégorie notée $\wt{\F}_\Gr(\kk)$ et étudions le
foncteur~$\kk[\Gr]$, dont une propriété fondamentale, l'auto-dualité,
constitue un cas particulier d'un résultat sur la catégorie
$\wt{\F}_\Gr(\kk)$. Le rôle fondamental du foncteur~$\kk[\Gr]$
provient de ce que la catégorie de foncteurs en grassmanniennes
globale $\F_\Gr$ est équivalente à la catégorie des
$\kk[\Gr]$-comodules. 

La section~\ref{cfpl} étudie une dernière famille de catégories de
foncteurs en grassmanniennes, notées $\F_{\Pl,n}(\kk)$, qui
fournit un analogue des catégories $\F_{\Gr,n}$ obtenu en 
trivialisant l'action du groupe linéaire $GL_n(\kk)$ (assertion qui
sera précisée par la proposition~\ref{liengrpl}). En outre, la
structure monadique de la catégorie $\F_{\Pl,n}$ donne lieu à de
nouvelles structures dans les catégories de foncteurs en grassmanniennes.

%

\section{Les catégories $\F_{\Gr,I}$}\label{cfgri}

Cette section est consacrée à une première étude des catégories
suivantes. Son principal objectif consiste à en comprendre les objets finis.

\begin{defi}\label{defcatf} \'Etant donnée une partie $I$ de
$\mathbb{N}$, on introduit la catégorie de
foncteurs 
$$\F_{\Gr,I}(\kk)=\mathbf{Fct}(\E^f_{\Gr,I}(\kk),\E_\kk).$$
Nous noterons simplement $\F_\Gr(\kk)$, ou $\F_\Gr$, pour $\F_{\Gr,\mathbb{N}}(\kk)$.
\end{defi}

\begin{rem}
Les  cas les plus intéressants sont ceux où $I=\mathbb{N}$, $\leq n$
ou $n$.
\end{rem}

\subsection{Généralités}\label{subsfgr} Nous introduisons les
foncteurs obtenus par précomposition à partir de ceux du
paragraphe~\ref{ssctdcfu}, et utilisons les propriétés établies dans
le paragraphe~\ref{ssctpff} et la section~\ref{sctccf}.

\begin{nota}\label{notglfgr}\begin{enumerate}\item Soit $J$ une partie de
$\mathbb{N}$. Nous abrégerons la notation des espaces vectoriels ${\rm
  hom}_{\F_{\Gr,J}}$ en ${\rm hom}_{\Gr,J}$, et utiliserons une
convention analogue pour les groupes d'extensions ou les hom internes. De même, nous noterons
les projectifs standard $P^{\Gr,J}_A$ (où $A\in {\rm
  Ob}\,\E^f_{\Gr,J}$) plutôt que $P^{\E^f_{\Gr,J}}_A$, et les injectifs
standard $I^{\Gr,J}_A$ plutôt que $I^{\E^f_{\Gr,J}}_A$ (cf. §\,\ref{sgs-p}). L'exposant $J$
sera omis pour $J=\mathbb{N}$.
\item Soient  $I$ et $J$ deux parties de $\mathbb{N}$ telles que
  $J\subset I$. Nous noterons $\mathcal{R}_{I,J} :
  \F_{\Gr,I}\to\F_{\Gr,J}$ le foncteur de restriction
  $(incl_{J,I})^*$. Le foncteur de prolongement par zéro, lorsqu'il
  est défini (cf. paragraphe \ref{p-prec}), sera noté $\mathcal{P}_{J,I} :
  \F_{\Gr,J}\to\F_{\Gr,I}$. Les
  indices seront omis lorsqu'il n'y a pas d'ambiguïté.
\end{enumerate}
\end{nota}

La proposition \ref{prfrc} fournit le résultat suivant, dans lequel
nous omettons les indices des foncteurs de restriction et de
prolongement par zéro (de sorte que la notation $\mathcal{R}$ et
$\mathcal{P}$ désigne à chaque fois deux foncteurs différents).

\begin{pr}\label{recf10} Pour tout entier $n\geq 0$, il existe un diagramme de recollement
$$\xymatrix{\F_{\Gr,\leq n-1}\ar[r]|-{\mathcal{P}} &
  \F_{\Gr,\leq n}\ar[r]|-{\mathcal{R}}\ar@/_/[l]\ar@/^/[l]^-{\mathcal{R}} &
  \F_{\Gr,n}\ar@/_/[l]\ar@/^/[l]^-{\mathcal{P}}}.$$
\end{pr}

Nous introduisons maintenant les deux foncteurs fondamentaux de source
$\F$ et de but $\F_{\Gr,I}$.

\begin{defi}\label{dfffond1} Soit $I$ une partie de $\mathbb{N}$.
\begin{enumerate}\item On définit le foncteur de {\em plongement
    standard} $\iota_I : \F\to\F_{\Gr,I}$ comme le foncteur de
  précomposition par le foncteur d'oubli principal $\mathfrak{O}_I :
  \E^f_{\Gr,I}\to\E^f$.
  
\item Le foncteur de {\em plongement réduit} $\kappa_I : \F\to\F_{\Gr,I}$ est le foncteur de
  précomposition par le foncteur de réduction $\mathfrak{K}_I :
  \E^f_{\Gr,I}\to\E^f$.
\end{enumerate}

Autrement dit, on a 
$$\iota_I(F)(V,W)=F(V)$$
et
$$\kappa_I(F)(V,W)=F(V/W)$$
pour $F\in {\rm Ob}\,\F$ et $(V,W)\in {\rm Ob}\,\E^f_{\Gr,I}$.

L'indice $I$ sera omis quand $I=\mathbb{N}$ ; des conventions
analogues vaudront dans les notations suivantes.
\end{defi}

\begin{rem} Les foncteurs $\iota_{\leq n}$ et $\iota$ sont, avec la
  notation du paragraphe \ref{sctcom}, les foncteurs
  $\Upsilon_{\Gr_{\leq n}}$ et $\Upsilon_\Gr$ respectivement. De même,
  le foncteur $\iota_n$ est, selon la convention du paragraphe
  \ref{svnc}, le foncteur $\Upsilon_{\Gr_{\leq n},\Gr_{\leq n-1}}$.
\end{rem}

Grâce à la remarque précédente, nous pouvons considérer les foncteurs
d'intégrale introduits dans les définitions~\ref{prdfom}
et~\ref{prdfom2}. Il s'agit des seuls foncteurs de cette section qui
ne soient pas des foncteurs de précomposition. Des foncteurs d'une
catégorie de foncteurs en grassmanniennes vers la catégorie $\F$, ce
sont les plus fondamentaux (cf. partie~\ref{p-omeg}).

\begin{defi}[Foncteurs d'intégrale en grassmanniennes]\label{figras} Soit $n\in\mathbb{N}$.
\begin{enumerate}\item Nous noterons $\omega_{\leq n} : \F_{\Gr,\leq n}\to\F$ le foncteur $\Omega_{\Gr_{\leq n}}$.
\item Nous noterons $\omega : \F_{\Gr}\to\F$ le foncteur $\Omega_{\Gr}$.
\item Nous désignerons par $\omega_n : \F_{\Gr,n}\to\F$ le foncteur $\Omega_{\Gr_{\leq n},\Gr_{\leq n-1}}$.
\end{enumerate}

Ces foncteurs seront appelés {\em foncteurs d'intégrale en
  grassmanniennes}.
\end{defi}

\begin{rem} On a des isomorphismes canoniques $\omega_{\leq
    n}\simeq\omega\circ\mathcal{P}_{\leq n,\mathbb{N}}$ et $\omega_n\simeq\omega_{\leq n}\circ\mathcal{P}_{n,\leq n}$.
\end{rem}

Les résultats des paragraphes \ref{sctcom} et \ref{svnc} se
traduisent par le résultat formel mais fondamental suivant.

\begin{pr}\label{prfig} Soit $n\in\mathbb{N}$.
\begin{enumerate}\item Le foncteur $\omega_{\leq n}$ est adjoint à
  gauche à $\iota_{\leq n}$. Il induit une équivalence de catégories
  entre $\F_{\Gr,\leq n}$ et la sous-catégorie
  $\mathbf{Comod}_{\kk[\Gr_{\leq n}]}$ de $\F$.
\item Le foncteur $\omega$ est adjoint à
  gauche à $\iota$. Il induit une équivalence de catégories
  entre $\F_{\Gr}$ et la sous-catégorie
  $\mathbf{Comod}_{\kk[\Gr]}$ de $\F$.
\item Le foncteur $\omega_n$ induit une équivalence de catégories
  entre $\F_{\Gr,n}$ et la sous-catégorie
  $\mathbf{Comod}_{\bar{G}(n)}^{fid}$ de $\F$ des
  $\bar{G}(n)$-comodules fidèles.
\item On a un isomorphisme
$$\omega_I(X\otimes\iota_I(F))\simeq\omega_I(X)\otimes F$$
naturel en les objets $X$ de $\F_{\Gr,I}$ et $F$ de $\F$, où $I=\mathbb{N}$, $\leq
n$ ou $n$.
\end{enumerate}
\end{pr}

\begin{proof}  Les deux premières assertions constituent des cas
  particuliers des propositions~\ref{adj-intou} et~\ref{prf-com}, la
  troisième de la proposition~\ref{prf-com2}. Le dernier point résulte
  pour sa part des propositions~\ref{preomu} et~\ref{preomu2}. 
\end{proof}

Nous définissons maintenant les deux foncteurs fondamentaux entre les
catégories $\F_{\Gr,I}$ et $\F^I_{surj}$.

\begin{defi} Soit $I$ une partie de $\mathbb{N}$.
\begin{enumerate}\item Le foncteur de {\em plongement secondaire} $\rho_I : \F^I_{surj}\to\F_{\Gr,I}$ est le foncteur de précomposition par le foncteur base $\mathfrak{B}_I : \E^f_{\Gr,I}\to\E^I_{surj}$.
\item Le {\em foncteur d'évaluation
    généralisée} $\varepsilon_I :
  \F_{\Gr,I}\to\F^I_{surj}$ est la précomposition par le foncteur de plongement diagonal $\mathfrak{D}_I : \E^I_{surj}\to\E^f_{\Gr,I}$.
\end{enumerate}

Autrement dit, on a
$$\rho_I(F)(V,W)=F(W)$$
pour $F\in {\rm Ob}\,\F^I_{surj}$ et $(V,W)\in {\rm Ob}\,\E^f_{\Gr,I}$
, et
$$\varepsilon_I(X)(V)=X(V,V)$$
pour $X\in {\rm Ob}\,\F_{\Gr,I}$ et $V\in {\rm Ob}\,\E^I_{surj}$.
\end{defi}

\begin{rem}\label{rqnoti} Avec l'abus de la remarque \ref{remdffp}, le foncteur $\varepsilon_n :
  \F_{\Gr,n}\to\,_{\kk[GL_n(\kk)]}\mathbf{Mod}$ s'identifie au foncteur
  d'évaluation ${\rm ev}_{(E_n,E_n)}$, ce qui justifie la terminologie
  employée.

On peut de même voir $\rho_n$ comme un foncteur
$_{\kk[GL_n(\kk)]}\mathbf{Mod}\to\F_{\Gr,n}$ ; il est donné par $\rho_n(X)(V,W)={\rm
   Iso}(E_n,W)\underset{\kk[GL_n(\kk)]}{\otimes}M$. 
\end{rem}

On rappelle que l'on désigne par $\F\otimes\F^I_{surj}$ la catégorie de
foncteurs $\mathbf{Fct}(\E^f\times\E^I_{surj},\E)$, selon la notation
\ref{ptcat}. 

\begin{defi}[Foncteurs fondamentaux entre $\F_{\Gr,I}$ et
  $\F\otimes\F^I_{surj}$] Soit $I$ une partie de $\mathbb{N}$.
\begin{enumerate}\item Le foncteur de {\em plongement
    complet} $\xi_I : \F\otimes\F^I_{surj}\to\F_{\Gr,I}$
  est la précomposition par le foncteur
  $\mathfrak{O}_I\times\mathfrak{B}_I :
  \E^f_{\Gr,I}\to\E^f\times\E^I_{surj}$.
\item Le foncteur de {\em plongement
    total} $\theta_I :
  \F\otimes\F^I_{surj}\to\F_{\Gr,I}$ est
  la précomposition par le foncteur
  $\mathfrak{K}_I\times\mathfrak{B}_I :
  \E^f_{\Gr,I}\to\E^f\times\E^I_{surj}$.
\item On définit le foncteur de {\em décalage en
    grassmanniennes} $\sigma_I : \F_{\Gr,I}\to\F\otimes\F^I_{surj}$ comme le foncteur de précomposition par le foncteur de plongement relatif $\mathfrak{L}_I : \E^f\times\E^I_{surj}\to\E^f_{\Gr,I}$.
\end{enumerate}

Ainsi, on a
$$\xi_I(F)(V,W)=F(V,W)\,,$$
$$\theta_I(F)(V,W)=F(V/W,W)$$
pour $F\in {\rm Ob}\,\F\otimes\F^I_{surj}$ et $(V,W)\in {\rm Ob}\,\E^f_{\Gr,I}$
, et
$$\sigma_I(X)(A,B)=X(A\oplus B,B)$$
pour $X\in {\rm Ob}\,\F_{\Gr,I}$, $A\in {\rm Ob}\,\E^f$ et $B\in {\rm Ob}\,\E^I_{surj}$.
\end{defi}

Intuitivement, la catégorie $\F\otimes\F^I_{surj}$ ne doit pas être vue comme
très différente de la catégorie $\F$. En effet, dans le cas où $I$ est
réduite à un élément $n$, on a des équivalences de catégories canoniques
$\F\otimes\F^n_{surj}\simeq\mathbf{Fct}(\E^f\times\underline{GL_n(\kk)},\E)\simeq\mathbf{Fct}(\underline{GL_n(\kk)},\mathbf{Fct}(\E^f,\E))=\F_{GL_n(\kk)}$
(cf. notation~\ref{notmo-i} donnée à la fin de l'introduction). Si $I$
est une partie quelconque de $\mathbb{N}$, la catégorie
$\F\otimes\F^I_{surj}$ s'obtient par recollement de telles catégories
(cf. proposition~\ref{precsurj}). 

\begin{rem} Compte-tenu des identifications que nous venons de
  mentionner, nous considérerons $\xi_n$ et $\theta_n$ comme des foncteurs
  $\F_{GL_n(\kk)}\to\F_{\Gr,n}$ (et de même $\sigma_{n}
    : \F_{\Gr, n}\to\F_{GL_n(\kk)}$).
\end{rem}

\paragraph*{Premières propriétés} 

Nous commençons par relier les foncteurs projectifs ou injectifs
standard des catégories $\F_{\Gr,I}$, $\F$ et $\F^I_{surj}$ à l'aide
des foncteurs introduits plus haut.

\begin{pr}\label{proj-cfg} Soit $J$ une partie de $\mathbb{N}$. On a des isomorphismes
\begin{equation}\label{eq-cfg1} P^{\Gr,J}_{\mathfrak{L}_J(V,W)}\simeq\iota_J(P_V)\otimes\rho_J(P_W^{\E^J_{surj}})
\end{equation}
et
\begin{equation}\label{eq-cfg2}
  I^{\Gr,J}_{\mathfrak{L}_J(V,W)}\simeq\kappa_J(I_V)\otimes I^{\Gr,J}_{\mathfrak{D}_J(W)}
\end{equation}
naturels en les objets $V$ de $\E^f$ et $W$ de $\E^J_{surj}$.
\end{pr}

\begin{proof}   Par la troisième adjonction de la proposition \ref{adjfpr}, on
 dispose d'un isomorphisme canonique
$${\rm hom}_{\E^f_{\Gr,J}}(\mathfrak{L}_J(V,W), \cdot)\simeq ({\rm hom}_{\E^f}(V,\cdot)\circ\mathfrak{O}_J)\times ({\rm hom}_{\E^J_{surj}}(W,\cdot)\circ\mathfrak{B}_J).$$
On en déduit le premier isomorphisme, en linéarisant.

De même, on obtient formellement l'isomorphisme (\ref{eq-cfg2}) à
partir du second isomorphisme de la proposition \ref{actc-prf}.
\end{proof}

\begin{rem}\label{efrq1}\begin{enumerate}\item Comme le foncteur de
    plongement relatif $\mathfrak{L}_J : \E^f\times\E^J_{surj}\to\E^f_{\Gr,J}$ est essentiellement surjectif,
cette proposition décrit tous les projectifs et tous les injectifs
standard de $\F_{\Gr,J}$.
\item Les injectifs du type $I^{\Gr,J}_{\mathfrak{D}_J(W)}=I^{\Gr,J}_{(W,W)}$ ne
se ramènent pas facilement à des injectifs de catégories plus
simples. Illustrons-le pour $J=\mathbb{N}$, afin de simplifier les
notations.

En effet, le comportement du foncteur $I^{\Gr}_{(W,W)}$ se rapproche de
l'injectif $I^{surj}_W$ de $\F_{surj}$ dans la mesure où il existe un isomorphisme canonique
$\varepsilon_J(I^{\Gr}_{(W,W})\simeq I^{surj}_W$ (car ${\rm
  hom}_{\E^f_\Gr}((V,V),(W,W))\simeq {\rm Epi}_\E (V,W)$), et aussi
parce que $I^{\Gr}_{(W,W)}$ comme $I^{surj}_W$ ont un anneau
d'endomorphismes isomorphe à $\kk[GL(W)]$. Cependant,
$I^{\Gr}_{(W,W)}$ n'est pas isomorphe à l'image de $I^{surj}_W$ par le foncteur
$\rho : \F_{surj}\to\F_\Gr$.
\end{enumerate}
\end{rem}

La propriété suivante des foncteurs injectifs et projectifs standard de la
catégorie $\F_\Gr$ (qui n'a pas d'équivalent dans $\F_{\Gr,I}$ pour
$I\neq\mathbb{N}$) est analogue aux assertions \ref{itvp1} et
\ref{itvp2} de la proposition \ref{pre-vpo}.

\begin{pr}\label{prfif}\begin{enumerate}\item Il existe un isomorphisme
    $\omega (P^\Gr_{(V,W)})\simeq P_V$ naturel en l'objet $(V,W)$ de $\E^f_\Gr$.
\item Il existe un isomorphisme
    $\iota(I_V)\simeq\underset{W\in\Gr(V)}{\bigoplus}I^\Gr_{(V,W)}$
    naturel en l'objet $V$ de $\E^f$.
\end{enumerate}
\end{pr}

\begin{proof} Ces deux assertions proviennent, par linéarisation, du
  lemme \ref{lmclei}.
\end{proof}

La proposition qui suit donne les propriétés d'usage courant des foncteurs
de précomposition introduits en début de section.

\begin{pr}\label{prfac-ff} Soit $I$ une partie non vide de $\mathbb{N}$.
\begin{enumerate}\item Les foncteurs $\iota_I : \F\to\F_{\Gr,I}$,
  $\sigma_I : \F_{\Gr,I}\to\F\otimes\F^I_{surj}$, $\xi_I : \F\otimes\F^I_{surj}\to\F_{\Gr,I}$ et
  $\varepsilon_I : \F_{\Gr,I}\to\F^I_{surj}$ sont exacts ; ils commutent au produit tensoriel, aux limites et aux
  colimites. De plus, les foncteurs $\iota_I$, $\sigma_I$ et $\xi_I$ sont fidèles.
\item Les foncteurs $\kappa_I : \F\to\F_{\Gr,I}$, $\rho_I :
  \F^I_{surj}\to\F_{\Gr,I}$ et $\theta_I : \F\otimes\F^I_{surj}\to\F_{\Gr,I}$ sont exacts et
  pleinement fidèles ; ils commutent au produit tensoriel, aux limites
  et aux colimites. De surcroît,
  leurs images sont des sous-catégories de Serre de $\F_{\Gr,I}$.
\end{enumerate}
\end{pr}

\begin{proof} Tous les foncteurs de l'énoncé sont des foncteurs de
  précomposition, ils commutent donc aux limites, aux colimites et au
  produit tensoriel, par la proposition~\ref{lm-form}. 

Comme le foncteur $\mathfrak{L}_I$ est essentiellement
  surjectif, par la proposition~\ref{crf-catec}, le foncteur de
  précomposition associé, $\sigma_I$, est fidèle (par la proposition~\ref{lm-form}). La
  fidélité du foncteur $\iota_I$ s'établit de même, car l'inclusion
  $\E^{\leq n}\hookrightarrow\E^f$ (cf. la deuxième assertion de la proposition~\ref{crf-catec}), où $n$ désigne le plus petit
  élément de $I$, induit par précomposition une équivalence de
  catégories
  $\F=\mathbf{Fct}(\E^f,\E)\xrightarrow{\simeq}\mathbf{Fct}(\E^{\leq
    n},\E)$. La fidélité de $\xi_I$ est analogue, en remarquant que
  tout objet de $\E^f\times\E^I_{surj}$ est rétracte d'un objet de
  l'image du foncteur
  $\mathfrak{O}_I\times\mathfrak{B}_I$.

Comme les foncteurs $\mathfrak{K}_I$, $\mathfrak{B}_I$ et
  $\mathfrak{K}_I\times\mathfrak{B}_I$ sont pleins et essentiellement
  surjectifs, par la proposition~\ref{crf-catec}, la
  proposition~\ref{lm-form} montre qu'ils sont pleinement fidèles et
  que leurs images sont des sous-catégories de Serre de $\F_{\Gr,I}$.
\end{proof}

La partie de la proposition~\ref{prfac-ff} relative au foncteur
$\theta_I$, ainsi que la composition fondamentale donnée par la
proposition suivante, apparaîtront sous un jour nouveau dans la section~\ref{fgm}.

\begin{pr}\label{prcompf3} Soit $I$ une partie de $\mathbb{N}$. Le
  foncteur composé
$$\F\otimes\F^I_{surj}\xrightarrow{\theta_I}\F_{\Gr,I}\xrightarrow{\sigma_I}\F\otimes\F^I_{surj}$$
est canoniquement isomorphe au foncteur identique. 
\end{pr}

Cette propriété découle de la proposition \ref{compf-prc}.

\begin{pr}\label{pradjf3} Soit $I$ une partie de $\mathbb{N}$. Le
  foncteur $\xi_I : \F\otimes\F^I_{surj}\to\F_{\Gr,I}$ est adjoint à gauche au foncteur $\sigma_I$.
\end{pr}

Ce résultat et le suivant s'obtiennent en combinant les propositions \ref{adjfpr} et
\ref{lm-form}.\,5.

\begin{pr}\label{prcompa4} Soit $I$ une partie de $\mathbb{N}$. Le
  foncteur $\rho_I$ est adjoint à gauche au foncteur
  $\varepsilon_I$. De surcroît, l'unité $id\to\varepsilon_I\rho_I$ de l'adjonction est un isomorphisme.
\end{pr}

Le corollaire \ref{adjkaio} et la proposition \ref{prfrec} fournissent les adjonctions suivantes.

\begin{pr}\label{adjf5} Soit $I$ une partie de $\mathbb{N}$
  contenant $0$. Le foncteur $\mathcal{R}_{I,0} : \F_{\Gr,I}\to\F$  est adjoint  à gauche à $\kappa_I$ et à droite à $\mathcal{P}_{0,I}$.
\end{pr}

Nous terminons ce sous-paragraphe  par deux propriétés
de compatibilité  relatives aux foncteurs entre $\F$, $\F_{surj}$ et $\F_\Gr$.

\begin{pr}\label{cofod} Le foncteur composé
  $\F\xrightarrow{\iota}\F_{\Gr}\xrightarrow{\varepsilon}\F_{surj}$ s'identifie canoniquement au foncteur d'oubli $o$.
\end{pr}

\begin{proof} En effet, le foncteur composé
$\E^f_{surj}\xrightarrow{\mathfrak{D}}\E^f_{\Gr}\xrightarrow{\mathfrak{O}}\E^f$
s'identifie au foncteur d'inclusion. \end{proof}

Par adjonction, on en déduit
(cf. propositions \ref{prfig}, \ref{prcompa4} et \ref{pre-vpo}) le corollaire suivant, que l'on peut évidemment établir par
une vérification directe.

\begin{cor}\label{cofo} Le foncteur $\omega\circ\rho :
  \F_{surj}\to\F$ est canoniquement isomorphe à~$\varpi$.
\end{cor}

\paragraph*{Décomposition scalaire et tors de
  Frobenius} La catégorie $\E^f_{\Gr,I}(\kk)$ est munie d'une action
naturelle du groupe $\kk^\times$. De façon analogue à la
proposition/définition~\ref{ds-fsurj}, on en déduit le résultat
suivant.

\begin{prdef}\label{ds-fgr} \'Etant donné un entier $i$, notons $(\F_{\Gr, I})_i(\kk)$
  la sous-catégorie pleine de $\F_{\Gr,I}(\kk)$ formée des foncteurs
  $X$ tels que $X(\lambda.id_V)=\lambda^i.id_{X(V)}$ pour tous
  $\lambda\in\kk^\times$ et $V\in {\rm Ob}\,\E^f_{\Gr,I}$. Les
  inclusions induisent une équivalence de catégories 
$$\F_{\Gr,I}(\kk)\simeq\prod_{i=1}^{q-1}(\F_{\Gr, I})_i(\kk).$$

On notera $X\simeq\bigoplus_{i=1}^{q-1} X_i$ la décomposition
naturelle d'un foncteur $X$ de $\F_{\Gr,I}(\kk)$ qu'on en déduit, où
$X_i\in {\rm Ob}\,(\F_{\Gr, I})_i(\kk)$. On l'appelle {\em décomposition
  scalaire} de $X$.
\end{prdef}

De même que $\F_{surj}(\kk)$, la catégorie $\F_{\Gr,I}(\kk)$ possède
un {\em tors de Frobenius}, parce que l'automorphisme $\phi$ de
la catégorie $\E^f_\kk$ déduit du morphisme de Frobenius se prolonge
naturellement en un automorphisme de la catégorie $\E^f_{\Gr,I}(\kk)$.

Les foncteurs introduits précédemment possèdent des propriétés de
commutation vis-à-vis de la décomposition scalaire et du tors de
Frobenius qu'on laisse au lecteur le soin d'expliciter (pour les
foncteurs de précomposition, elles se lisent sur le foncteur d'origine
entre les catégories sources).

\paragraph*{Changement de corps} Si $K$ est une extension finie de
$\kk$, les foncteurs d'induction et de restriction induisent des
foncteurs $\E^f_{\Gr,I}(\kk)\to\E^f_{\Gr,I}(K)$ et
$\E^f_{\Gr,I}(K)\to\E^f_{\Gr,I}(\kk)$, d'où l'on déduit, comme dans
$\F$ et $\F_{surj}$, des foncteurs (encore dits d'induction et de
restriction) $\F_{\Gr,I}(\kk)\to\F_{\Gr,I}(K)$ et
$\F_{\Gr,I}(K)\to\F_{\Gr,I}(\kk)$. Néanmoins, comme dans le cas de
$\F_{surj}$, on perd la propriété d'adjonction mutuelle entre
induction et restriction.

Tous les foncteurs de précomposition introduits dans ce paragraphe
commutent aux changements de corps. En revanche, ce n'est pas le cas du foncteur
d'intégrale en grassmanniennes.

\subsection{Structures tensorielles} Il est naturel, comme le montrera
la proposition~\ref{prpttgr2}, d'introduire sur la catégorie de foncteurs en
grassmanniennes globale $\F_\Gr$, à côté de la structure tensorielle
usuelle donnée par $\otimes$ (à laquelle nous nous référerons
lorsque nous parlerons de structure tensorielle sans plus de précision) une seconde structure tensorielle,
très analogue à celle donnée par le produit tensoriel total de
$\F_{surj}$ (cf. paragraphe~\ref{parsgs}).

\begin{defi}[Produit tensoriel total] \'Etant donnés deux objets $X$ et
$Y$ de $\F_\Gr$, on appelle {\em produit tensoriel
    total} de $X$ et $Y$ le
  foncteur noté $X\ptt Y$ et défini comme suit.
\begin{enumerate}\item Si $(V,W)$ est un objet de $\E^f_\Gr$, on pose
$$(X\ptt Y)(V,W)=\underset{A+B=W}{\bigoplus_{A,B\in\Gr(W)}}X(V,A)\otimes Y(V,B).$$
\item Si $u : (V,W)\twoheadrightarrow (V',W')$ est un morphisme de $\E^f_\Gr$,
  le morphisme $(X\ptt Y)(u) : (X\ptt Y)(V,W)\to (X\ptt Y)(V',W')$ est
  défini comme la somme directe sur les sous-espaces $A$ et $B$ de $W$
  tels que $A+B=W$ des morphismes
$$X(V,A)\otimes Y(V,B)\xrightarrow{X(u)\otimes Y(u)}X(V',u(A))\otimes
Y(V',u(B))$$
$$\hookrightarrow\underset{A'+B'=W'}{\bigoplus_{A',B'\in\Gr(W')}}X(V',A')\otimes Y(V',B')$$
où l'on note encore,
par abus, $u$ pour les morphismes $(V,A)\to (V',u(A))$ et
$(V,B)\to (V',u(B))$ induits par $u$.

Cette définition fait
sens puisque $u(A)+u(B)=u(A+B)=u(W)=W'$.
\end{enumerate}
\end{defi}

Avant de préciser les liens entre ce produit tensoriel total et celui
défini dans $\F_{surj}$, définissons un produit
tensoriel total dans la catégorie auxiliaire $\F\otimes\F_{surj}$.

\begin{defi}\label{pttptf} Le {\em produit tensoriel total} sur
  $\F\otimes\F_{surj}$ est le foncteur $\ptt : (\F\otimes\F_{surj})\times
  (\F\otimes\F_{surj})\to\F\otimes\F_{surj}$ donné, via l'isomorphisme
  $\F\otimes\F_{surj}\simeq\mathbf{Fct}(\E^f,\F_{surj})$, par 
$$\mathbf{Fct}(\E^f,\F_{surj})\times\mathbf{Fct}(\E^f,\F_{surj})\simeq\mathbf{Fct}(\E^f,\F_{surj}\times\F_{surj})\xrightarrow{\ptt_*}\mathbf{Fct}(\E^f,\F_{surj}).$$

Autrement dit,
$$(X\ptt Y)(A,B)=\underset{V+W=B}{\bigoplus_{V,W\in\Gr(B)}}
X(A,V)\otimes Y(A,W)\qquad (A\in{\rm Ob}\,\E^f, B\in{\rm Ob}\,\E^f_{surj}).$$
\end{defi}

\begin{pr}\label{prpttgr}\begin{enumerate}\item Le produit tensoriel total définit sur
    $\F_\Gr$ une structure monoïdale symétrique exacte d'unité $\rho({\rm Is}_0)$.
\item Il existe un monomorphisme $X\otimes Y\hookrightarrow X\ptt Y$
  naturel en les objets $X$ et $Y$ de $\F_\Gr$.
\item Les produits tensoriels totaux sur les catégories $\F_\Gr$,
  $\F_{surj}$ et $\F\otimes\F_{surj}$ vérifient les propriétés de
  compatibilité suivantes : on a des isomorphismes naturels
 \begin{eqnarray*}\rho(A\ptt  B) & \simeq & \rho(A)\ptt\rho(B)  \qquad (A, B\in {\rm
  Ob}\,\F_{surj}), \\
\xi(F\ptt
  G) & \simeq & \xi(F)\ptt\xi(G)  \qquad (F,
  G\in {\rm Ob}\,\F\otimes\F_{surj}).
\end{eqnarray*}
\end{enumerate}
\end{pr}

Cette proposition se vérifie par inspection.

\begin{rem} En revanche, les foncteurs $\varepsilon$, $\sigma$ et
  $\theta$ ne commutent pas au produit tensoriel total. Par exemple,
  on a
$$\theta(F\ptt G)(V,W)\simeq (F\ptt
G)(V/W,W)=\underset{A+B=W}{\bigoplus_{A,B\in\Gr(W)}}F(V/W,A)\otimes
G(V/W,B)$$
tandis que
$$(\theta(F)\ptt\theta(G))(V,W)=\underset{A+B=W}{\bigoplus_{A,B\in\Gr(W)}}\theta(F)(V,A)\otimes\theta(G)(V,B)$$
$$\simeq\underset{A+B=W}{\bigoplus_{A,B\in\Gr(W)}}F(V/A,A)\otimes
G(V/B,B).$$
On voit cependant ainsi qu'il existe un épimorphisme canonique $\theta(F)\ptt\theta(G)\twoheadrightarrow\theta(F\ptt G)$.
\end{rem}

On dispose d'un foncteur $\E^f_\Gr\times\E^f_\Gr\to\E^f_\Gr$ donné sur les objets
par $((A,B),(A',B'))\mapsto (A\oplus A',B\oplus B')$, dont l'action
sur les morphismes se déduit de la fonctorialité de $\oplus :
\E^f\times\E^f\to\E^f$. La proposition suivante montre l'utilité du
produit tensoriel total à l'aide de ce foncteur.

\begin{pr}\label{prpttgr2} Il existe dans $\F_\Gr$ un isomorphisme $P^\Gr_{(A,B)}\ptt
  P^\Gr_{(A',B')}\simeq P^\Gr_{(A\oplus A',B\oplus B')}$ naturel en
  les objets $(A,B)$ et $(A',B')$ de $\E^f_\Gr$.
\end{pr}

\begin{proof} Cet énoncé découle de la proposition \ref{prsef1}.
\end{proof}

\begin{cor}\label{presptt} Le produit tensoriel total de $\F_\Gr$
  préserve les objets projectifs, les objets de type fini et les
  objets pf$_n$.
\end{cor}

\begin{rem}\begin{enumerate}\item Le produit tensoriel usuel de $\F_\Gr$ préserve également les objets de type fini et
  pf$_n$.
\item Cependant, le produit tensoriel usuel de deux objets projectifs de
  $\F_\Gr$ n'est généralement pas projectif, car il en est ainsi dans
  $\F_{surj}$ --- par exemple, $P^{surj}_\kk\otimes P^{surj}_{E_2}$
  n'est pas projectif. 
\item Le produit tensoriel de deux objets projectifs d'une catégorie
  $\F_{\Gr,n}$ est en revanche projectif. Cela provient de ce que le produit tensoriel
  d'un $\kk[GL_n(\kk)]$-module projectif et d'un $\kk[GL_n(\kk)]$-module fini
  est projectif, de ce que le produit tensoriel dans $\F$ de deux
  projectifs est projectif, et de l'isomorphisme (\ref{eq-cfg1}) de la
  proposition~\ref{proj-cfg}.
\end{enumerate}
\end{rem}

Contrairement à la situation que nous venons d'observer pour les
projectifs, le produit tensoriel total possède un comportement
déplaisant sur les injectifs. En revanche, le produit tensoriel ordinaire
est adapté à leur étude. 

On rappelle que le symbole $Gr$
qui intervient dans la proposition suivante a été introduit dans la
notation \ref{notgrd}.

\begin{pr}\label{prtifgr} Il existe un isomorphisme
$$I^\Gr_{(A,B)}\otimes I^\Gr_{(A',B')}\simeq\bigoplus_{C\in Gr(B,B')}
I^\Gr_{(A\oplus A',C)}$$
naturel en les objets $(A,B)$ et $(A',B')$ de $\E^f_\Gr$.
\end{pr}

\begin{proof} C'est la version linéarisée de la proposition~\ref{prsef2}.
\end{proof}

\begin{cor}\label{pgri} Le produit tensoriel de $\F_\Gr$ préserve :
\begin{enumerate}\item les objets injectifs de co-type fini ;
\item les objets de co-type fini ;
\item les objets co-pf$_n$.
\end{enumerate}
\end{cor}

La proposition suivante (qui généralise la dernière assertion de la
proposition \ref{pre-vpo}) constitue la principale motivation de
l'introduction du produit tensoriel total dans $\F_\Gr$.

\begin{pr}\label{comom} Il existe dans $\F$ un isomorphisme 
$$\omega(X\ptt Y)\simeq\omega(X)\otimes\omega(Y)$$
naturel en les objets $X$ et $Y$ de $\F_\Gr$.
\end{pr}

\begin{proof} Cela provient, par linéarisation, de la décomposition
  ensembliste
$$\Gr(V)\times\Gr(V)\simeq\coprod_{W\in\Gr(V)}\{(A,B)\in\Gr(W)\times\Gr(W)\,|\,A+B=W\}$$
naturelle en l'objet $V$ de $\E^f$.
\end{proof}

\subsection{Le foncteur différence}\label{subs-fd} Comme dans la catégorie $\F$, il
existe dans les catégories $\F_{\Gr,I}$ un foncteur différence,
fondamental dans l'étude de ces catégories, notamment de leurs objets
finis. L'analogie avec le cas de $\F$ est particulièrement étroite du
fait que le foncteur différence des $\F_{\Gr,I}$~\guillemotleft~n'agit pas sur la base~\guillemotright.

\begin{conv} Dans tout ce paragraphe, $J$ désigne une partie de
$\mathbb{N}$.
\end{conv}

\begin{defi}[Foncteurs de décalage et foncteur différence dans
  $\F_{\Gr,J}$]\label{dfgrp} 
\begin{enumerate}\item Soit $V$ un objet de
    $\E^f$. Le foncteur $\F_{\Gr,J}\to\F_{\Gr,J}$ de précomposition
    par le foncteur de translation $V\boxplus\cdot :
    \E^f_{\Gr,J}\to\E^f_{\Gr,J}$ (cf. définition \ref{dfftr}) est
    appelé {\em foncteur de décalage} par $V$ et se note
    $\Delta_V^{\Gr,J}$.

On a ainsi $\Delta_V^{\Gr,J} X(A,B)=X(V\oplus A,B)$
pour $X\in {\rm Ob}\,\F_{\Gr,J}$ et $(A,B)\in {\rm Ob}\,\E^f_{\Gr,J}$.

 La bifonctorialité
    de $\boxplus$ rend l'association $V\mapsto\Delta^{\Gr,J}_V$
    fonctorielle. 
\item Le {\em foncteur différence} $\Delta^{\Gr,J}$ de $\F_{\Gr,J}$ est le noyau de la transformation naturelle
  $\Delta_\kk^{\Gr,J}\to id$ induite par le morphisme $\kk\to
  0$ de $\E^f$.
\end{enumerate}
\end{defi}


La proposition suivante, laissée au lecteur, montre que le
comportement de ces foncteurs ne dépend guère de la partie $J$ de $\mathbb{N}$.

\begin{pr}\label{rqedgh}Si $I$ est un sous-ensemble de $J$, le diagramme 
$$\xymatrix{\F_{\Gr,I}\ar[d]_{\mathcal{P}_{I,J}}\ar[r]^-{\Delta^{\Gr,I}}
    & \F_{\Gr,I}\ar[d]^{\mathcal{P}_{I,J}} \\
\F_{\Gr,J}\ar[r]^-{\Delta^{\Gr,J}}
    & \F_{\Gr,J}
}$$
commute à isomorphisme canonique près, lorsque le prolongement par
zéro $\mathcal{P}_{I,J}$ est défini. Il existe une propriété
analogue relative au foncteur de restriction $\mathcal{R}_{J,I}$.
\end{pr}

Une grande part des considérations relatives aux foncteurs de décalage
et différence de $\F$ se transcrivent  dans
$\F_{\Gr,J}$. Notons tout d'abord que la transformation
naturelle $id\to\Delta^{\Gr,J}_\kk$ procure un scindement
canonique
$$\Delta^{\Gr,J}_\kk\simeq\Delta^{\Gr,J}\oplus id\,,$$
de sorte que $\Delta^{\Gr,J}$ commute aux limites et colimites, comme les foncteurs de décalage.

Une autre propriété fondamentale, similaire à la
proposition~\ref{adj-fctd} relative à $\F$, réside dans l'existence d'adjoints à
gauche et à droite exacts et explicites aux foncteurs de décalage et différence.

\begin{pr}\label{adelgr} Soit $V$ un objet de $\E^f$. Le foncteur
  $\Delta_V^{\Gr,J}$ est adjoint :
\begin{enumerate}\item à droite à $\cdot\otimes\iota_J(P_V)$ ;
\item à gauche à $\cdot\otimes\kappa_J(I_V)$.
\end{enumerate}
Ces adjonctions sont naturelles en $V$.
\end{pr}

\begin{proof} On combine les propositions \ref{prdffdec} et
  \ref{prdffdec2} avec les isomorphismes de la proposition~\ref{actc-prf}.
\end{proof}

\begin{cor}\label{cretgr} Le foncteur différence 
  $\Delta^{\Gr,J}$ est adjoint :
\begin{enumerate}\item à droite à $\cdot\otimes\iota_J(\bar{P}_\kk)$ ;
\item à gauche à $\cdot\otimes\kappa_J(\bar{I}_\kk)$.
\end{enumerate}
\end{cor}

\begin{cor}\label{crdifgr} Les foncteurs de décalage et le foncteur
  différence de $\F_{\Gr,J}$ conservent les objets projectifs et les
  objets injectifs.
\end{cor}

\begin{proof} C'est une conséquence formelle de la proposition et du
  corollaire précédents et de l'exactitude des adjoints.
\end{proof}

Nous évaluons maintenant les foncteurs de décalage sur les objets
projectifs et injectifs standard de la catégorie $\F_{\Gr,J}$.

\begin{pr}\label{decip} Il existe des isomorphismes
$$\Delta^{\Gr,J}_V (P^{\Gr,J}_{(A,B)})\simeq P_{A/B}(V)\otimes
P^{\Gr,J}_{(A,B)}$$
et
$$\Delta^{\Gr,J}_V (I^{\Gr,J}_{(A,B)})\simeq I_A(V)\otimes
I^{\Gr,J}_{(A,B)}$$
naturels en les objets $V$ de $\E^f$ et $(A,B)$ de $\E^f_{\Gr,J}$.
\end{pr}

\begin{proof} C'est une conséquence formelle de la proposition~\ref{actc-prf}.
\end{proof}

Comme dans le cas de $\F$, cette propriété a l'utile conséquence
suivante. 

\begin{cor}\label{cr-reggr} Les foncteurs de décalage et le foncteur
  différence de $\F_{\Gr,J}$ conservent les objets de type fini, de
  co-type fini, pf$_n$ et co-pf$_n$.
\end{cor}

\begin{proof} Ce corollaire se déduit de la proposition~\ref{decip} et
  de la commutation des
  foncteurs de décalage et différence aux limites et colimites.
\end{proof}

La proposition suivante donne les propriétés de compatibilité entre les
deux produits tensoriels de $\F_\Gr$ et son foncteur différence.

\begin{pr}\label{prtnfgr} Soient $X$ et $Y$ deux objets de
  $\F_{\Gr,J}$ et $V$ un objet de $\E^f$. Il existe des isomorphismes
  naturels
 $$\Delta^\Gr_V(X\otimes Y)\simeq\Delta^\Gr_V
    X\otimes\Delta^\Gr_V Y,$$ 
 $$\Delta^\Gr(X\otimes Y)\simeq (\Delta^\Gr
    X\otimes Y)\oplus (X\otimes\Delta^\Gr Y)\oplus (\Delta^\Gr
    X\otimes\Delta^\Gr Y)\,;$$
et, dans le cas où $J=\mathbb{N}$,
 $$\Delta^\Gr_V(X\ptt Y)\simeq\Delta^\Gr_V
    X\ptt\Delta^\Gr_V Y,$$ 
 $$\Delta^\Gr(X\ptt Y)\simeq (\Delta^\Gr
    X\ptt Y)\oplus (X\ptt\Delta^\Gr Y)\oplus (\Delta^\Gr
    X\ptt\Delta^\Gr Y).$$
\end{pr}

\begin{proof} \'Etablissons la première assertion relative au produit
  tensoriel total : si $(A,B)$ est un objet de $\E^f_\Gr$, on a
$$\Delta_V(X\ptt Y)(A,B)=(X\ptt Y)(V\oplus
A,B)=\underset{W+W'=B}{\bigoplus_{W,W'\in\Gr(B)}}X(V\oplus A,W)\otimes
Y(V\oplus A,W')$$
$$=\underset{W+W'=B}{\bigoplus_{W,W'\in\Gr(B)}}\Delta_ V(X)(A,W)\otimes
\Delta_ V(Y)(A,W')=\Delta_ V(X)\ptt\Delta_ V(Y)(A,B).$$

L'inclusion canonique $$X\ptt Y\hookrightarrow\Delta^\Gr_\kk(X\ptt
Y)=\Delta^\Gr_\kk(X)\ptt\Delta^\Gr_\kk(Y)\simeq
(X\oplus\Delta^\Gr(X))\ptt (Y\oplus\Delta^\Gr(Y))$$
$$\simeq (X\ptt Y)\oplus
(X\ptt\Delta^\Gr(Y))\oplus (\Delta^\Gr(X)\ptt Y)\oplus
(\Delta^\Gr(X)\ptt\Delta^\Gr(Y))$$
s'identifie à l'inclusion du facteur direct, d'où le dernier point.

La première assertion relative au produit tensoriel ordinaire résulte
de ce que $\Delta^\Gr_V$ est un foncteur de précomposition ; la
seconde s'en déduit comme dans le cas du produit tensoriel total.
\end{proof}

\paragraph*{Commutation des foncteurs fondamentaux aux foncteurs
  différences} Con\-for\-mé\-ment à nos conventions générales (cf. notation~\ref{postpre}),
l'endofoncteur $\Delta_*$ de
$\F\otimes\F^J_{surj}\simeq\mathbf{Fct}(\E^J_{surj},\F)$ qui
apparaît dans la proposition suivante est donné par la
  postcomposition par le foncteur différence $\Delta$ de $\F$.

\begin{pr}\label{pr-comd2}
\begin{enumerate}\item On a des isomorphismes
    canoniques $\Delta^{\Gr,J}\circ\iota_J\simeq\iota_J\circ\Delta$ et
    $\Delta^{\Gr,J}\circ\kappa_J\simeq\kappa_J\circ\Delta$ de
    foncteurs $\F\to\F_{\Gr,J}$.
\item On a des isomorphismes
    canoniques $\Delta^{\Gr,J}\circ\xi_J\simeq\xi_J\circ\Delta_*$ et
    $\Delta^{\Gr,J}\circ\theta_J\simeq\theta_J\circ\Delta_*$ de
    foncteurs $\F\otimes\F^J_{surj}\to\F_{\Gr,J}$.
\item On a un isomorphisme canonique
    $\Delta_*\circ\sigma_J\simeq\sigma_J\circ\Delta^{\Gr,J}$ de
    foncteurs $\F_{\Gr,J}\to\F\otimes\F^J_{surj}$.
\end{enumerate}
\end{pr}

\begin{proof}\'Etablissons par exemple le dernier point, les autres se
  montrant de manière analogue. Le diagramme
$$\xymatrix{\E^f\times\E^J_{surj}\ar[d]_{\mathfrak{L}_J}\ar[rr]^{(\FF\oplus\cdot)\times\E^J_{surj}}
  & & 
  \E^f\times\E^J_{surj}\ar[d]^{\mathfrak{L}_J} \\
\E^f_{\Gr,J}\ar[rr]^-{\kk\boxplus\cdot} & & \E^f_{\Gr,J}
}$$
commute à isomorphisme canonique près, d'où un isomorphisme canonique
$\wt{\Delta}_*\circ\sigma_J\simeq\sigma_J\circ\Delta^{\Gr,J}_\kk$, puis $\Delta_*\circ\sigma_J\simeq\sigma_J\circ\Delta^{\Gr,J}$.
\end{proof}

\subsection{Foncteurs polynomiaux}\label{s-fp} Le foncteur différence
de $\F_{\Gr,J}$ permet d'introduire les foncteurs polynomiaux de cette
catégorie. On commence, à l'aide du lemme et de la proposition suivants, par identifier le
noyau du foncteur différence, qui, contrairement au cas de la catégorie
$\F$, ne se réduit pas aux foncteurs constants.

\begin{lm}\label{lmev999} Le foncteur composé $\Delta^{\Gr,J}\rho_J : \F^J_{surj}\to\F_{\Gr,J}$ est nul.
\end{lm}

\begin{proof} Pour tout $V\in {\rm Ob}\,\E^f$, l'injection naturelle
$\rho_J\hookrightarrow\Delta^{\Gr,J}_V\circ\rho_J$ est un isomorphisme, puisque le
  diagramme
$$\xymatrix{\E^f_{\Gr,J}\ar[r]^-{\mathfrak{B}_J}\ar[d]_{V\boxplus\,\cdot} & \E^J_{surj} \\
\E^f_{\Gr,J}\ar[ru]_-{\mathfrak{B}_J} &
}$$
commute (à isomorphisme naturel près).
\end{proof}

\begin{prdef}[Foncteurs pseudo-constants]\label{prdfpcst} Soit $X$ un objet de $\F_{\Gr,J}$. Les
  conditions suivantes sont équivalentes :
\begin{enumerate}\item le foncteur $\Delta^{\Gr,J} X$ est nul ;
\item la coünité $\rho_J\varepsilon_J (X)\to X$ de l'adjonction de la
  proposition \ref{prcompa4} est un isomorphisme ;
\item il existe un objet $R$ de $\F^J_{surj}$ tel que $X$ est
  isomorphe à $\rho_J(R)$.
\end{enumerate}

Lorsqu'elles sont vérifiées, nous dirons que $X$ est un foncteur {\em pseudo-constant}.
\end{prdef}

\begin{proof} La coünité $\rho_J\varepsilon_J (X)\to X$ est donnée sur
  l'objet $(V,W)$ de $\E^f_{\Gr,J}$ par l'injection
  $X(W,W)\hookrightarrow X(V,W)$ induite par le morphisme canonique
  $(W,W)\to (V,W)$ de $\E^f_{\Gr,J}$. La première assertion signifie que le
  monomorphisme canonique $X(A,B)\hookrightarrow X(A\oplus\kk,B)$ est
  un isomorphisme pour tout objet $(A,B)$ de $\E^f_{\Gr,J}$, on en
  déduit donc que le morphisme précédent est un isomorphisme par
  récurrence sur la codimension de $W$ dans $V$. Ainsi, la première
  assertion implique la deuxième.

Il est clair que la deuxième implique la troisième. Enfin, la troisième entraîne la première par le lemme~\ref{lmev999}.
\end{proof}

Nous introduisons maintenant la notion générale de foncteur polynomial
et de foncteur analytique de $\F_{\Gr,J}$, de la même manière que dans
$\F$. Nous avons précédemment identifié les foncteurs polynomiaux de
degré $0$. 

\begin{defi}[Foncteurs polynomiaux, analytiques dans $\F_{\Gr,J}$] Un foncteur $X\in
  {\rm Ob}\,\F_{\Gr,J}$ est dit :
\begin{enumerate}\item {\em polynomial} s'il existe un entier naturel
  $n$ tel que $(\Delta^{\Gr,J})^n X=0$ ;
\item {\em analytique} s'il est
  colimite de ses sous-foncteurs polynomiaux.
\end{enumerate}

Le {\em degré} d'un foncteur polynomial $X$ est le plus grand entier
positif $n$ tel que $(\Delta^{\Gr,J})^n X\neq 0$ si $X$ est non nul, on le note
$\deg X$ ; on pose également $\deg 0=-\infty$.
\end{defi}

\begin{nota} Nous désignerons par $\F_{\Gr,J}^k$ la sous-catégorie
  pleine de $\F_{\Gr,J}$ formée des foncteurs polynomiaux de degré au
  plus $k$. 
\end{nota}

La proposition~\ref{prtnfgr} a la conséquence suivante.

\begin{pr}\label{prtenfgr}\begin{enumerate}
\item Si $X$ et $Y$ sont deux foncteurs polynomiaux de $\F_{\Gr,J}$,
  alors $X\otimes Y$ et $X\ptt Y$ sont
  polynomiaux et $\deg (X\otimes Y)=\deg (X\ptt Y)=\deg X+\deg Y$.
\item Le produit tensoriel total (donc aussi usuel) de deux foncteurs analytiques de
  $\F_{\Gr,J}$ est analytique.
\item Un foncteur $X$ de $\F_{\Gr,J}$ est polynomial si et seulement si $\Delta^{\Gr,J} X$ est
  polynomial ; dans ce cas $\deg (\Delta^{\Gr,J} X)=\deg X-1$ si $X$ n'est pas pseudo-constant.
\end{enumerate}
\end{pr}


Afin de donner les propriétés de conservation des foncteurs
polynomiaux ou analytiques par les foncteurs fondamentaux de source ou
de but $\F_{\Gr,J}$, nous avons besoin d'introduire la notion de
foncteur polynomial ou analytique dans la catégorie $\F\otimes\F^J_{surj}$.

\begin{prdef}\label{prdfdelot} Soient $X$ un objet de
  $\F\otimes\F^J_{surj}$ et $i\in\mathbb{N}$. 
\begin{enumerate}\item
Les assertions suivantes sont équivalentes.
\begin{enumerate}\item Le foncteur $\Delta_*^i X$ est nul. 
\item Le foncteur $X$ appartient à la sous-catégorie épaisse
  $\mathbf{Fct}(\E^J_{surj},\F^{i-1})$ de $\mathbf{Fct}(\E^J_{surj},\F)$.
\end{enumerate}

Lorsque ces conditions sont vérifiées, nous dirons que $X$ est un
foncteur polynomial de degré strictement inférieur à $i$.
\item On dit que $X$ est analytique s'il est colimite de ses
  sous-foncteurs polynomiaux.
\end{enumerate}
\end{prdef}

\begin{proof} C'est une application directe de la proposition~\ref{lmfora}.
\end{proof}

L'énoncé qui suit découle quant à lui de la proposition~\ref{pr-comd2}.

\begin{pr}\label{pr-comd}
\begin{enumerate}\item\begin{enumerate}\item Si $F$ est un objet de $\F$, $F$ est polynomial si et seulement
  si $\iota_J (F)$ est un objet polynomial de $\F_{\Gr,J}$. Ils ont
  alors même degré.
\item Si $F$ est un objet de $\F$, $F$ est polynomial si et seulement
  si $\kappa_J (F)$ est un objet polynomial de $\F_{\Gr,J}$. Ils ont
  alors même degré.
\end{enumerate}
\item\begin{enumerate}\item Un objet $F$ de $\F\otimes\F^J_{surj}$ est polynomial si et seulement
  si l'objet $\xi_J(F)$ de $\F_{\Gr,J}$ est polynomial. Ils ont
  alors même degré.
\item Si $F$ est un objet de $\F\otimes\F^J_{surj}$, $F$ est polynomial si et seulement
  si $\theta_J (F)$ est un objet polynomial de $\F_{\Gr,J}$. Ils ont
  alors même degré.
\end{enumerate}
\item Un objet $X$ de $\F_{\Gr,J}$ est polynomial si et
  seulement si l'objet $\sigma_J(X)$ de
  $\F\otimes\F^J_{surj}$ est polynomial. Ils ont alors même degré.
\end{enumerate}
\end{pr}

Comme les foncteurs considérés commutent aux colimites, on en déduit
le corollaire suivant.

\begin{cor}\label{rqana} Les foncteurs $\iota_J$, $\kappa_J$, $\xi_J$,
  $\theta_J$ et $\sigma_J$ préservent les foncteurs analytiques.
\end{cor}

\begin{rem} En revanche, le foncteur d'intégrale en grassmanniennes
  $\omega : \F_\Gr\to\F$ ne préserve {\em pas} les foncteurs
  analytiques. Nous reviendrons en détail sur cette observation essentielle dans
  la partie~\ref{p-omeg}.
\end{rem}

\paragraph*{Quotients de la filtration polynomiale de $\F_{\Gr,J}$}
Grâce à la proposition~\ref{prdfpcst}, la catégorie $\F^0_{\Gr,J}$ est
équivalente à la catégorie $\F^J_{surj}$. \`A l'aide des
lemmes~\ref{lmcpd1} et~\ref{lmcpd12}, nous allons généraliser ce
résultat en identifiant les catégories quotients
$\F^k_{\Gr,J}/\F^{k-1}_{\Gr,J}$ pour tout entier $k$. Dans le
paragraphe suivant, nous
appliquerons ces résultats à la description des objets simples de~$\F_{\Gr,J}$.

\begin{lm}\label{lmcpd1} Soit $X$ un objet polynomial de degré $i\geq 0$ de
  $\F_{\Gr,J}$. Le noyau de l'épimorphisme canonique
  $p_X : \xi_J\,\sigma_J(X)\twoheadrightarrow X$ (coünité
  de l'adjonction de la proposition \ref{pradjf3}) est de degré
  strictement inférieur à $i$.
\end{lm}

\begin{proof} La proposition~\ref{pr-comd2} montre que les
  endofoncteurs $\xi_J\,\sigma_J$ et $\Delta^{\Gr,J}$ de $\F_{\Gr,J}$
  commutent à isomorphisme naturel près ; de plus,
  $\Delta^{\Gr,J}(p_X)\simeq p_{\Delta^{\Gr,J} X}$. Comme $p_X$ est un isomorphisme lorsque
  le foncteur $X$ est pseudo-constant, cela donne la conclusion.
\end{proof}

\begin{lm}\label{exfk} Si $F$ est un objet polynomial de degré $i\geq 0$ de $\F\otimes\F^J_{surj}$, le noyau de l'épimorphisme canonique $\xi_J(F)\twoheadrightarrow\theta_J(F)$ est de degré strictement inférieur à~$i$.
\end{lm}

\begin{proof} Ce lemme s'obtient par application du précédent à
  l'objet $\theta_J(F)$ de $\F_{\Gr,J}$, qui est de degré $i$ par la
  proposition~\ref{pr-comd}, puisque $\sigma_J\circ\theta_J\simeq id$ (proposition~\ref{prcompf3}). 
\end{proof}

Le lemme suivant se démontre de la même façon que le lemme~\ref{lmcpd1}.

\begin{lm}\label{lmcpd12} Soit $X$ un objet polynomial de degré $i\geq 0$ de
  $\F_{\Gr,J}$. Le conoyau du monomorphisme canonique
  $X\hookrightarrow\sigma_J\,\xi_J (X)$ (unité
  de l'adjonction de la proposition \ref{pradjf3}) est de degré
  strictement inférieur à~$i$.
\end{lm}

Les lemmes \ref{lmcpd1} et~\ref{lmcpd12} impliquent l'important
résultat suivant.

\begin{pr}\label{prcqd}
\begin{enumerate}\item Le foncteur
  $\xi_J : \mathbf{Fct}(\E^f\times\E^J_{surj},\E)\simeq\mathbf{Fct}(\E^J_{surj},\F)\to\F_{\Gr,J}$  induit un foncteur
  $\mathbf{Fct}(\E^J_{surj},\F^n/\F^{n-1})\to\F_{\Gr,J}^n/\F_{\Gr,J}^{n-1}$ pour tout $n\in\mathbb{N}$.
\item Ce foncteur est une équivalence de catégories dont un inverse est
donné par le foncteur $\F_{\Gr,J}^n/\F_{\Gr,J}^{n-1}\to\mathbf{Fct}(\E^J_{surj},\F^n/\F^{n-1})$
induit par $\sigma_J :
\F_{\Gr,J}\to\mathbf{Fct}(\E^f\times\E^J_{surj},\E)\simeq\mathbf{Fct}(\E^J_{surj},\F)$.
\end{enumerate}
\end{pr}

Nous utiliserons également la description suivante de l'équivalence de
catégories donnée par la proposition~\ref{prcqd}.

\begin{pr}\label{rqkio} Le foncteur
  $\mathbf{Fct}(\E^J_{surj},\F^n/\F^{n-1})\to\F_{\Gr,J}^n/\F_{\Gr,J}^{n-1}$ induit par $\theta_J$  est le même que celui qu'induit $\xi_J$.
\end{pr}

\begin{proof} Il s'agit d'une conséquence directe du lemme~\ref{exfk}.
\end{proof}

Les propositions \ref{prcqd} et \ref{prec-f} fournissent le corollaire
suivant.

\begin{cor}\label{crggr} On
  suppose $d=1$, i.e. le corps $\kk$ premier. 

Pour tout  entier naturel $n$, les catégories
  $\F_{\Gr,J}^n/\F_{\Gr,J}^{n-1}$ et
  $\mathbf{Fct}(\E^J_{surj},\mathbf{Mod}_{\kk[\Sigma_n]})$ sont équivalentes.
\end{cor}

\subsection{Foncteurs finis}\label{s-ffg} Les résultats du paragraphe précédent permettent de décrire les
objets simples et les objets finis des catégories $\F_{\Gr,I}$ à partir de ceux de $\F$ et $\F^I_{surj}$. Comme pour la
catégorie $\F$, l'un des résultats les plus importants réside dans le
caractère polynomial des objets finis.

\begin{conv} Dans ce paragraphe, $I$ désigne une partie de $\mathbb{N}$.
\end{conv}

\begin{lm}\label{lmprce} Soit $X$ un objet de $\F_{\Gr,I}$ tel que
$\Delta^{\Gr,I} X$ est un objet fini de $\F_{\Gr,I}$ et $\varepsilon_I(X)$ un
objet fini de $\F^I_{surj}$. Alors $X$ est fini.
\end{lm}

\begin{proof} Le foncteur exact
$$(\Delta^{\Gr,I},\varepsilon_I) : \F_{\Gr,I}\to\F_{\Gr,I}\times\F^I_{surj}\,.$$
est fidèle : son noyau est constitué des foncteurs pseudo-constants
$X$ tels que $\varepsilon_I(X)=0$. Mais un tel foncteur $X$ est
isomorphe à $\rho_I\varepsilon_I(X)=0$, par la
proposition~\ref{prdfpcst}, d'où la fidélité annoncée. La
proposition~\ref{pr-ff} donne donc la conclusion.
\end{proof}

\begin{rem}\label{remut3} On a des résultats semblables en remplaçant
  dans cet énoncé {\em fini} par {\em de type fini}, {\em pf$_n$} ou
  {\em co-pf$_n$}.
\end{rem}

\begin{defi}[Niveau]\label{nivcon} On appelle {\em niveau}
  d'un objet $X$ de $\F_{\Gr,I}$ l'élément
$${\rm niv} (X)=sup\,\{\dim W\,|\,(V,W)\in {\rm
  Ob}\,\E^f_{\Gr,I}\;\,X(V,W)\neq 0\}$$
de $I\cup\{-\infty,+\infty\}$. On dit que $X$ est de {\em niveau fini}
si ${\rm niv} (X)<+\infty$.

On définit de même le {\em coniveau} de $X$ comme l'élément
$${\rm coniv} (X)=inf\,\{\dim W\,|\,(V,W)\in {\rm Ob}\,\E^f_{\Gr,I}\;\,X(V,W)\neq 0\}$$
de $I\cup\{+\infty\}$.
\end{defi}

Ainsi, lorsque l'ensemble $I$ est fini, tous les foncteurs de
$\F_{\Gr,I}$ sont de niveau fini.

\begin{pr}\label{prf-grof}  Un objet de $\F_{\Gr,I}$ est fini si et seulement s'il est
polynomial, à valeurs de dimension finie et de niveau fini.
\end{pr}

\begin{proof} Soient $S$ un objet simple de $\F_{\Gr,I}$ et $(A,B)$ un
  objet de $\E^f\times\E^I_{surj}$ tel que $S(\mathfrak{L}_I(A,B))\neq
  0$ (le foncteur $\mathfrak{L}_I$ est essentiellement surjectif ---
  cf. proposition \ref{crf-catec}), de sorte qu'il
  existe un épimorphisme
  $P^{\Gr,I}_{\mathfrak{L}_I(A,B)}\simeq\iota_I(P_A)\otimes\rho_I(P^{\E^I_{surj}}_B)\twoheadrightarrow S$ (cf. proposition \ref{proj-cfg}). Cela montre d'une part que $S$ est à valeurs de dimension finie, et que $S(V,W)=0$ pour $\dim W>\dim B$. D'autre part, comme l'objet $P_A$ de $\F$ est co\-analytique (proposition~\ref{prf-fpol}), il existe un objet polynomial $F$ de $\F$ et un épimorphisme $\iota_I(F)\otimes\rho_I(P^{surj}_B)\twoheadrightarrow S$. La proposition \ref{pr-comd} montre que $\iota_I(F)\otimes\rho_I(P^{surj}_B)$, donc $S$, est polynomial. On en déduit les mêmes propriétés pour les foncteurs finis par un argument d'épaisseur.

La réciproque résulte du lemme \ref{lmprce} et de la proposition \ref{ftffs} par récurrence sur le
degré polynomial.
\end{proof}

Cette proposition essentielle permet d'établir la préservation des
foncteurs finis par les foncteurs $\xi_I$, $\theta_I$ et $\sigma_I$,
donnée par le corollaire suivant.

\begin{cor}\label{cr-ffe2} Soient $F$ un objet de $\F$ et $A$ un objet
  de $\F^I_{surj}$.
\begin{enumerate}\item Si $F$ et $A$ sont finis, il en est de même pour
  $\iota_I(F)\otimes\rho_I(A)$. La réciproque est vraie si $A$ et $F$
  sont non nuls.  
\item Si $F$ et $A$ sont finis, il en est de même pour
  $\kappa_I(F)\otimes\rho_I(A)$. La réciproque est vraie si $A$ et $F$
  sont non nuls.
\item Un objet $X$ de $\F_{\Gr,I}$ est fini si et
  seulement si l'objet $\sigma_I(X)$ de $\F\otimes\F^I_{surj}$ est
  fini.
\end{enumerate}
\end{cor}

\begin{proof} Traitons le premier cas, en écartant le cas trivial où
  $F$ ou $A$ est nul. Le foncteur
  $\iota_I(F)\otimes\rho_I(A)=\xi_I(F\boxtimes A)$ (on rappelle que le
  produit tensoriel extérieur $\boxtimes$ est introduit au
  §\,\ref{pap-pte}) est polynomial si
  et seulement s'il en est de même pour $F\boxtimes A$, par la
  proposition~\ref{pr-comd}, ou encore de $F$, puisque
  $\Delta_*(F\boxtimes A)=\Delta(F)\boxtimes A$. Le foncteur
  $\iota_I(F)\otimes\rho_I(A)$ est à valeurs de dimension finie si et
  seulement s'il en est de même pour $F$ et $A$. Enfin, le foncteur
  $F\boxtimes A$ est de niveau fini si et seulement si $A(E_n)=0$ pour
  $n$ assez grand. Les propositions~\ref{prf-grof} et~\ref{ftffs}
  fournissent donc la première assertion.

Les autres s'établissent de façon analogue.\end{proof}

\begin{nota} Nous désignerons par la suite par $\F_{\Gr,I}^{lf}$ la sous-catégorie
pleine de $\F_{\Gr,I}$ formée des objets localement finis. 
\end{nota}

Le corollaire suivant donne les principales propriétés de régularité
des objets finis et localement finis de
$\F_{\Gr,I}$. 

\begin{cor}\label{crfff3}\begin{enumerate}\item Si $X$ est un objet fini
    de $\F_{\Gr,I}$, il en est de même pour $\Delta^{\Gr,I} X$.
\item Le produit tensoriel total (donc, a fortiori, usuel)  de deux objets finis de $\F_{\Gr,I}$ est fini.
\item Les objets finis de $\F_{\Gr,I}$ sont pf$_\infty$ et
  co-pf$_\infty$. 
\item La sous-catégorie $\F_{\Gr,I}^{lf}$ de $\F_{\Gr,I}$ est
  épaisse. Ses objets sont les foncteurs analytiques.
\item Un objet de $\F_{\Gr,I}$ est de co-type fini si et seulement
  s'il est analytique et de socle fini.
\end{enumerate} 
\end{cor}

\begin{proof} Elle est similaire à celle du corollaire
  \ref{crf-fpolf}. Pour le dernier point, on établit que les injectifs
  standard de $\F_{\Gr,I}$ sont analytiques de la façon suivante : le
  corollaire~\ref{rqana} et la proposition~\ref{prf-fpol} montrent que
 l'image par le foncteur $\iota_I$ d'un injectif standard de $\F$ est
 analytique. Comme tout injectif standard de $\F_{\Gr,I}$ est facteur
 direct d'un tel foncteur, par la proposition~\ref{prfif}, cela donne
 la conclusion.
\end{proof}

\paragraph*{Description des objets simples} Les objets simples de la
catégorie $\F_{\Gr,I}$ seront ramenés à ceux de la catégorie
$\F\otimes\F^I_{surj}$. Le lemme suivant les décrit à partir des
objets simples de $\F$ et $\F^I_{surj}$, qui sont eux-mêmes assez bien
compris (cf. sections~\ref{s-rf} et~\ref{sct-surj}).

\begin{lm}\label{prelmspl} Les objets simples de $\F\otimes\F^I_{surj}$ sont, à
  isomorphisme près, les $S\boxtimes R$, où $S$ est un objet simple de
  $\F$ et $R$ un objet simple de $\F^I_{surj}$. Le produit tensoriel
  extérieur induit de plus un isomorphisme d'anneaux (sans
  unité si $I$ est infini) $G_0^f(\F\otimes\F^I_{surj})\simeq
G_0^f(\F)\otimes G_0^f(\F^I_{surj})$.
\end{lm}

\begin{proof} Cela découle de la proposition
  \ref{simple-pte} et de
  son corollaire \ref{gro-pte}, dont les hypothèses sont vérifiées
  grâce à la proposition~\ref{crf-sf}.
\end{proof}

Avec la proposition~\ref{prf-grof}, le résultat suivant est le plus
important de cette section.

\begin{pr}\label{prfffgr}\begin{enumerate}\item \'Etant donné un objet $X$ de $\F_{\Gr,I}$, les
assertions suivantes sont équivalentes :
\begin{enumerate}\item l'objet $X$ de $\F_{\Gr,I}$ est simple ;
\item l'objet $\sigma_I(X)$ de $\F\otimes\F^I_{surj}$ est simple ;
\item il existe un objet simple $F$ de $\F$ et un objet simple $R$ de
  $\F^I_{surj}$ tel que $X$ est isomorphe à $\kappa_I(F)\otimes\rho_I(R)$.
\end{enumerate}
\item Les foncteurs exacts $\sigma_I$ et
  $\theta_I$ induisent des isomorphismes d'anneaux (sans
  unité si $I$ est infini)
entre $G_0^f(\F_{\Gr,I})$ et $G_0^f(\F\otimes\F^I_{surj})\simeq
G_0^f(\F)\otimes G_0^f(\F^I_{surj})$ réciproques l'un de l'autre. 
\end{enumerate}
\end{pr}

\begin{proof} La proposition \ref{prfac-ff}.\,2 montre que le
  foncteur $\theta_I$ transforme un objet simple de
  $\F\otimes\F^I_{surj}$ en un objet simple de $\F_{\Gr,I}$ et induit
  un {\em monomorphisme} de groupes abéliens, compatible au produit,
  $G_0^f(\F\otimes\F^I_{surj})\hookrightarrow G_0^f(\F_{\Gr,I})$. Ces
  deux groupes sont naturellement gradués par le degré polynomial, qui
  est respecté par $\theta_I$ (cf. proposition \ref{pr-comd}), de
  sorte que ce monomorphisme s'identifie à la somme directe sur
  $k\in\mathbb{N}$ des
  morphismes
$$G_0^f(\mathbf{Fct}(\E^I_{surj},\F^k/\F^{k-1}))\to\F_{\Gr,J}^k/\F_{\Gr,J}^{k-1}$$
induits par $\theta_I$, lesquels sont des isomorphismes par la
proposition \ref{prcqd}. La compatibilité au produit tensoriel
extérieur pour $I=\mathbb{N}$ provient de la proposition~\ref{prpttgr}.

Cela démontre la seconde assertion de l'énoncé, et la première
compte-tenu du lemme \ref{prelmspl} et de ce que l'isomorphisme
$G_0^f(\F_{\Gr,I})\xrightarrow{\simeq}G_0^f(\F\otimes\F^I_{surj})$
inverse du précédent est induit par $\sigma_I$ (par les propositions~\ref{pr-comd} et~\ref{rqkio}).
\end{proof}

\begin{cor}\label{ncrptt} Le foncteur exact $\xi_I$ induit un isomorphisme d'anneaux (sans
  unité si $I$ est infini)
entre $G_0^f(\F)\otimes G_0^f(\F^I_{surj})\simeq
G_0^f(\F\otimes\F^I_{surj})$ et $G_0^f(\F_{\Gr,I})$.

Dans le cas où $I=\mathbb{N}$, ces foncteurs induisent également un
isomorphisme d'anneaux pour la structure induite par le produit tensoriel total sur $\F_\Gr$
et $\F_{surj}$ (et le produit tensoriel sur $\F$).
\end{cor}

\begin{proof} Pour obtenir l'isomorphisme de groupes, il suffit de
  reprendre la démonstration de la proposition~\ref{prfffgr}, compte-tenu de la proposition~\ref{rqkio}. La proposition~\ref{prpttgr} procure la compatibilité au
  produit tensoriel total.
\end{proof}

\paragraph*{Quelques conséquences} Avant d'aborder la description des
objets projectifs indécomposables de $\F_{\Gr,I}$ (proposition~\ref{crkz}), qui sera cependant
moins explicite que celle de ses objets simples, nous tirons de la
proposition~\ref{prfffgr} le corollaire suivant relatif aux cosocles
--- rappelons que cette notion est introduite dans la définition~\ref{defsocs}.

\begin{cor}\label{cr-cosoc} Pour tout objet $F$ de $\F\otimes\F^I_{surj}$, on a
  des isomorphismes naturels
$${\rm cosoc}\,\xi_I(F)\simeq {\rm
  cosoc}\,\theta_I(F)\simeq \theta_I({\rm cosoc}\,F).$$
Le premier isomorphisme est induit par
  la projection canonique
  $\xi_I(F)\twoheadrightarrow\theta_I(F)$ et le
  second par $F\twoheadrightarrow {\rm cosoc}\,F$.
\end{cor}

\begin{proof} Grâce à la proposition \ref{prfffgr}, l'isomorphisme
  ${\rm cosoc}\,\xi_I(F)\simeq\theta_I({\rm cosoc}\,F)$ provient des
  isomorphismes naturels
$${\rm hom}_\Gr (\xi_I(F),\theta_I(S))\simeq {\rm
  hom}_{\F\otimes\F^I_{surj}} (F,\sigma_I\theta_I(S))\simeq {\rm
  hom}_{\F\otimes\F^I_{surj}} (F,S)$$
déduits des propositions \ref{pradjf3} et
\ref{prcompf3}.

L'isomorphisme ${\rm cosoc}\,\theta_I(F)\simeq \theta_I({\rm cosoc}\,F)$ se déduit de la seconde assertion de la proposition~\ref{prfac-ff}.
\end{proof}

\begin{rem} Le premier isomorphisme du corollaire~\ref{cr-cosoc} implique, que pour tout objet de type fini
  $F$ de $\F$, l'épimorphisme canonique
  $\iota_I(F)\twoheadrightarrow\kappa_I(F)$ est essentiel. En
  effet, un épimorphisme entre objets co-localement finis et à valeurs
  de dimension finie de $\F_{\Gr,I}$ qui induit
  un isomorphisme entre les cosocles est essentiel --- cet énoncé
  s'obtient, par dualité, à partir de l'observation que dans une
  catégorie de Grothendieck, un monomorphisme entre objets localement
  finis qui induit un isomorphisme entre les socles est essentiel.

Le comportement de l'épimorphisme canonique
  $\iota_I\twoheadrightarrow\kappa_I$ est totalement différent sur les objets de co-type fini
  ; ainsi, la proposition~\ref{prfif} montre que, dans le cas où $I=\mathbb{N}$, cet épimorphisme se
  scinde toujours lorsque $F$ est un injectif standard $I_V$ de $\F$
  (alors que $\iota(I_V)\twoheadrightarrow\kappa(I_V)$ n'est un
  isomorphisme que pour $V=0$).
\end{rem}

Rappelons que $K_0$ désigne le groupe de Grothendieck des classes d'objets
projectifs indécomposables de type fini (cf. notation~\ref{not-grot}).

\begin{lm}\label{llmqt}\begin{enumerate}\item Le produit tensoriel extérieur induit un isomorphisme $K_0(\F)\otimes
K_0(\F^I_{surj})\xrightarrow{\simeq} K_0(\F\otimes\F^I_{surj})$.
\item Si $\I$ est une catégorie essentiellement petite vérifiant
  l'hypothèse~\ref{hypf3} (page~\pageref{hypf3}), on obtient des
  isomorphismes réciproques l'un de l'autre entre
  $G^f_0(\mathbf{Fct}(\I,\E_\kk))$ et $K_0(\mathbf{Fct}(\I,\E_\kk))$ en
  associant à un objet simple sa couverture projective et à un objet
  projectif de type fini son cosocle.
\end{enumerate}
\end{lm}

\begin{proof} La première assertion se démontre de manière analogue
  au lemme~\ref{prelmspl}.

La seconde est un résultat classique en théorie des représentations
(cf. \cite{CR}) ; dans le cas général, on le déduit aisément en
dualisant (à l'aide de l'hypothèse~\ref{hypf3}) des propriétés des
enveloppes injectives dans les catégories de Grothendieck
(cf. \cite{Gab}) --- on pourra se reporter à \cite{these} pour une
démonstration détaillée.
\end{proof}

\begin{pr}\label{crkz} Le foncteur $\xi_I :
  \F\otimes\F^I_{surj}\to\F_{\Gr,I}$ induit un isomorphisme de groupes abéliens
  $K_0(\F\otimes\F^I_{surj})\xrightarrow{\simeq}K_0(\F_{\Gr,I})$.

Ainsi, le produit tensoriel des morphismes induits par $\iota_I$ et $\rho_I$
fournit un isomorphisme $K_0(\F)\otimes
K_0(\F^I_{surj})\xrightarrow{\simeq}K_0(\F_{\Gr,I})$.

Dans le cas où $I=\mathbb{N}$, le produit tensoriel total sur
$\F_\Gr$ et $\F_{surj}$ et le produit tensoriel sur $\F$ munissent les
deux membres de structures d'anneau commutatif ; l'isomorphisme
précédent est un isomorphisme d'anneaux pour ces structures. 
\end{pr}

\begin{proof} La proposition \ref{proj-cfg} montre que le foncteur
  exact $\xi_I$ conserve les objets projectifs de
  type fini, il induit donc un morphisme de groupes
  $K_0(\F\otimes\F^I_{surj})\to K_0(\F_{\Gr,I})$. Par le
  corollaire~\ref{cr-cosoc}, le diagramme
$$\xymatrix{K_0(\F\otimes\F^I_{surj})\ar[r]^-{(\xi_I)_*}\ar[d]_\simeq &
  K_0(\F_{\Gr,I})\ar[d]^\simeq \\
G^f_0(\F\otimes\F^I_{surj})\ar[r]^-{(\theta_I)_*} &
  G^f_0(\F_{\Gr,I})
}$$
commute, où les flèches verticales, induites par le cosocle, sont des
isomorphismes par la seconde assertion du lemme~\ref{llmqt}. Comme le
morphisme $(\theta_I)_*$ est un isomorphisme par la
proposition~\ref{prfffgr}, il en est de même pour le morphisme $(\xi_I)_*$.

La première assertion du lemme~\ref{llmqt} permet d'en déduire l'isomorphisme $K_0(\F)\otimes
K_0(\F^I_{surj})\xrightarrow{\simeq}K_0(\F_{\Gr,I})$.

Le corollaire~\ref{presptt} établit enfin la  compatibilité au produit tensoriel extérieur.
\end{proof}

\begin{rem} L'inverse de l'isomorphisme de la proposition~\ref{crkz} ne se décrit pas simplement
  en termes des foncteurs fondamentaux depuis $\F_{\Gr,I}$.
\end{rem}

\section{La catégorie $\widetilde{\F}_\Gr$}\label{cfgrt}

Nous introduisons dans cette section une catégorie de foncteurs $\wt{\F}_\Gr(\kk)$, qui est intuitivement à
$\F_\Gr(\kk)$ ce que $\F(\kk)$ est à $\F_{surj}(\kk)$. Nous en
étudions quelques propriétés au paragraphe~\ref{substil} en vue
d'applications à d'autres catégories. C'est la propriété de dualité
démontrée au paragraphe~\ref{pduf}, qui n'a pas d'analogue dans
$\F_\Gr(\kk)$, qui nous permettra de compléter les renseignements sur
la structure du foncteur~$\kk[\Gr]$ déduits de la
section~\ref{sct-surj} (§\,\ref{pskg}).

Cette propriété interviendra également dans la partie~\ref{p-omeg}, afin d'appliquer des
propriétés du foncteur $\omega$ à la catégorie $\F_{inj}(\kk)$.

\begin{nota} Nous noterons $\wt{\F}_\Gr(\kk)$ la catégorie $\mathbf{Fct}(\wt{\E}^f_\Gr(\kk),\E_\kk)$.
\end{nota}

\subsection{Généralités}\label{substil} La catégorie $\wt{\F}_\Gr$ possède un comportement assez
différent des catégories $\F_{\Gr,I}$ : elle n'entre
pas dans le cadre étudié à la section~\ref{sctccf}, de sorte qu'elle
ne s'interprète pas en termes de (co)modules. 

Plutôt que de donner une description complète des foncteurs qui
apparaissent naturellement, par précomposition, à partir des foncteurs
étudiés dans la section \ref{sctcatb}, nous nous focaliserons sur les
liens entre  $\wt{\F}_\Gr$ et $\F_{\Gr}$.

\begin{nota}\label{notfraktil} Le foncteur de précomposition
  ${\widetilde{incl}}^* :
  \widetilde{\mathcal{F}}_\Gr\to\mathcal{F}_\Gr$ sera noté $\mathfrak{R}$.
\end{nota}

Nous introduisons maintenant un foncteur très analogue au foncteur
$\varpi$.

\begin{prdef}\label{defj} Il existe un foncteur exact et fidèle
  $\mathfrak{J} :
  \mathcal{F}_\Gr\to\widetilde{\mathcal{F}}_\Gr$ défini de la manière suivante :
\begin{itemize}
\item{\em  action sur les objets : }
$$\mathfrak{J}(X)(V,B)=\bigoplus_{W\in\Gr(B)} X(V,W)$$
\item {\em action sur les morphismes : } si $f : (V,B)\to (V',B')$ est un morphisme de $\widetilde{\mathcal{E}}^f_\Gr$, $\mathfrak{J}(X)(f)$ a pour composante $X(V,W)\to X(V',W')$ (où $W\in\Gr(B)$ et $W'\in\Gr(B')$) le morphisme induit par $f$ si $W'=f(W)$ et $0$ sinon,
\item {\em fonctorialité : } si $t : X\to Y$ est un morphisme de $\F_{\Gr}$, $\mathfrak{J}(t) : \mathfrak{J}(X)\to \mathfrak{J}(Y)$ s'obtient sur l'objet $(V,B)$ par somme directe des $t(V,W)$ pour $W\in\Gr(B)$.
\end{itemize}
\end{prdef}

\begin{pr}\label{pratil}  Le foncteur $\mathfrak{J} :
  \mathcal{F}_\Gr\to\widetilde{\mathcal{F}}_\Gr$ est adjoint à gauche à  $\mathfrak{R}$.
\end{pr}

La vérification de ces propriétés, analogue à celle des
propositions~\ref{prdfffs} et~\ref{pre-vpo} (ou~\ref{prdfom}
et~\ref{adj-intou}), est laissée au lecteur.


\begin{nota}\label{notil3}
L'endofoncteur  $\mathfrak{R}\circ\mathfrak{J}$ de
$\mathcal{F}_\Gr$ sera noté $\mathcal{I}$.

Explicitement, on a
\begin{equation}\label{ident-ieq}\I(X)(V,W)=\bigoplus_{B\in\Gr(W)} X(V,B)\end{equation}
pour $X\in {\rm Ob}\,\F_{\Gr}$ et $(V,W)\in {\rm Ob}\,\E^f_{\Gr}$.
\end{nota}

L'unité de l'adjonction de la proposition \ref{pratil} procure une
transformation naturelle injective $id\to\I$. Elle est
donnée sur un objet $X$ par l'inclusion
$X(V,W)\hookrightarrow\I(X)(V,W)$ correspondant au facteur direct de (\ref{ident-ieq})
obtenu pour $B=W$.

%
%
%
%
%
%

\begin{nota}\label{prdftnf} Nous noterons $j^\omega :
  \mathcal{I}\to\iota\circ\omega$ la
  transformation naturelle donnée par les inclusions
$$\bigoplus_{W\in\Gr(B)}
X(V,W)\hookrightarrow\bigoplus_{W\in\Gr(V)} X(V,W)\qquad (X\in {\rm
  Ob}\,\F_{\Gr},\, (V,B)\in {\rm
  Ob}\,\E^f_{\Gr}).$$ 
\end{nota}

Cette transformation naturelle jouera un rôle essentiel dans le paragraphe~\ref{sctisof}.

\begin{rem} Le diagramme suivant d'endofoncteurs de $\F_\Gr$, dans lequel les flèches non
  spécifiées sont les unités des adjonctions des propositions~\ref{pratil}
  et~\ref{prfig}.\,$2$, commute.
$$\xymatrix{id\ar[r]\ar[dr] & \iota\omega \\
 & \I\ar[u]_{j^\omega}
}$$
\end{rem}

Nous introduisons à présent l'auto-dualité de la
catégorie $\widetilde{\F}_\Gr$ que l'on déduit de la
proposition/définition~\ref{dfgrt}, dont on conserve la notation $(\cdot)^\vee$.

\begin{defi} Le {\em foncteur de dualité} $(\widetilde{\F}_\Gr)^{op}\to\widetilde{\F}_\Gr$, noté
  $D_\Gr$, est la composée de
  la précomposition par $(\cdot)^\vee :
  (\wt{\E}^f_\Gr)^{op}\to\wt{\E}^f_\Gr$ et de la postcomposition par
  $(\cdot)^* : \E^{op}\to\E$. 

On a ainsi  $(D_\Gr X)(V,W)=X((V,W)^\vee)^*$ pour tous $X\in {\rm
  Ob}\,\wt{\F}_\Gr$ et $(V,W)\in {\rm Ob}\,\wt{\E}^f_\Gr$.\end{defi}

La proposition suivante, que l'on déduit aussitôt de la proposition~\ref{prdcff}, justifie l'appellation de dualité.

\begin{pr}\begin{enumerate}\item Le foncteur $D_\Gr :
  (\widetilde{\F}_\Gr)^{op}\to\widetilde{\F}_\Gr$ est adjoint à
  droite à~$D_\Gr^{op}$.
\item Le foncteur $D_\Gr$ induit une équivalence entre la
  sous-catégorie pleine de $\wt{\F}_\Gr$ des foncteurs à valeurs de
  dimension finie et sa catégorie opposée.
\end{enumerate}
\end{pr}

Cela permet de définir la notion d'objet {\em auto-dual} de
$\wt{\F}_\Gr$, de façon similaire à la définition~\ref{dfad}.

\paragraph*{Décomposition scalaire, tors de Frobenius et changement de
corps} Comme la catégorie $\wt{\E}^f_\Gr(\kk)$ est $\kk$-linéaire, il
existe dans $\wt{\F}_\Gr$ une décomposition scalaire au sens de la
proposition/définition~\ref{pra-ds}. On dispose par ailleurs, comme
dans la catégorie $\F_\Gr(\kk)$, d'un tors de Frobenius et de foncteurs
de restriction et d'induction relativement à une extension finie
de~$\kk$. On laisse au lecteur le soin d'écrire les détails et les
propriétés de compatibilité de ces différents foncteurs.

\paragraph*{Objets finis} La catégorie $\wt{\E}^f$ étant additive, on
dispose d'un endofoncteur de décalage $\Delta^{\wt{\Gr}}_A$ dans
$\wt{\F}_\Gr$ associé à chaque objet $A$ de $\wt{\E}^f$, et l'on a des
inclusions canoniques scindées
$id\hookrightarrow\Delta^{\wt{\Gr}}_{(\kk,0)}$ et
$id\hookrightarrow\Delta^{\wt{\Gr}}_{(\kk,\kk)}$, de conoyaux notés
$\bar{\Delta}^{\wt{\Gr}}_{(\kk,0)}$ et
$\bar{\Delta}^{\wt{\Gr}}_{(\kk,\kk)}$  respectivement. On peut alors
montrer le
résultat suivant.

\begin{pr} Un objet de $\wt{\F}_\Gr$ est fini si et seulement s'il est
  nilpotent pour les deux foncteurs $\bar{\Delta}^{\wt{\Gr}}_{(\kk,0)}$ et
$\bar{\Delta}^{\wt{\Gr}}_{(\kk,\kk)}$ et à valeurs de dimension finie.
\end{pr}

Esquissons la démonstration de cette proposition, que nous
n'utiliserons pas :
\begin{enumerate}\item le foncteur
  $(\bar{\Delta}^{\wt{\Gr}}_{(\kk,0)},\bar{\Delta}^{\wt{\Gr}}_{(\kk,\kk)}) : \wt{\F}_\Gr\to\wt{\F}_\Gr\times\wt{\F}_\Gr$ est exact et son noyau se réduit aux foncteurs constants. On en déduit facilement qu'un foncteur de $\wt{\F}_\Gr$ nilpotent pour $\bar{\Delta}^{\wt{\Gr}}_{(\kk,0)}$ et
$\bar{\Delta}^{\wt{\Gr}}_{(\kk,\kk)}$ et à valeurs de dimension finie
est fini ;
\item à partir du caractère analytique des injectifs standard de $\F$,
  on obtient que les injectifs standard de $\wt{\F}_\Gr$ sont colimite
  de foncteurs nilpotents pour $\bar{\Delta}^{\wt{\Gr}}_{(\kk,0)}$ et
$\bar{\Delta}^{\wt{\Gr}}_{(\kk,\kk)}$ et à valeurs de dimension finie,
ce qui montre la réciproque.
\end{enumerate}

\begin{rem} On peut déduire de la proposition précédente un
  isomorphisme d'anneaux $G_0^f(\wt{\F}_\Gr)\simeq G_0^f(\F)^{\otimes
    2}$, d'une manière analogue à la proposition~\ref{prfffgr}.
\end{rem}

\subsection{Propriété de dualité du foncteur
  $\widetilde{\omega}$}\label{pduf} Ce paragraphe est consacré à
l'étude d'une propriété de dualité liée au foncteur $\omega :
\F_\Gr\to\F$. Son cadre naturel est la catégorie $\wt{\F}_\Gr$, par
l'intermédiaire du foncteur introduit dans la notation suivante. En
effet, la catégorie $\wt{\F}_\Gr$ est
étroitement liée à $\F_\Gr$ et possède un foncteur dualité,
contrairement à $\F_\Gr$.

\begin{nota}\label{ndw} Nous désignerons par $\widetilde{\omega} :
\widetilde{\F}_\Gr\to\F$ le foncteur composé
$$\widetilde{\F}_\Gr\xrightarrow{\mathfrak{R}}\F_\Gr\xrightarrow{\omega}\F.$$
\end{nota}

Le foncteur $\wt{\omega}$ est donné explicitement par
$\widetilde{\omega}(X)(V)=\bigoplus_{B\in\Gr(V)}X(V,B)$ ;  un
morphisme $f : V\to V'$ de $\E^f$ induit l'application linéaire
$\wt{\omega}(X)(V)\to\wt{\omega}(V')$ dont  la composante
$X(V,B)\to X(V',B')$ est le morphisme induit par $f$ si
$B'=f(B)$, $0$ sinon.

L'introduction de la variante suivante du foncteur $\wt{\omega}$ est
motivée par les propositions~\ref{pr-oop} et~\ref{comdom2}.

\begin{nota}\label{nwtp} On note $\widetilde{\omega}' :
  \widetilde{\F}_\Gr\to\F$ le foncteur défini par 
$$\widetilde{\omega}'(X)(V)=\widetilde{\omega}(X)(V)=\bigoplus_{B\in\Gr(V)}X(V,B)$$
et tel que pour tout morphisme $f : V\to V'$ de $\E^f$, la composante
$X(V,B)\to X(V',B')$ de $\widetilde{\omega}'(X)(f)$ est l'application
induite par le morphisme $(V,B)\to (V',B')$ induit par $f$ si
$f^{-1}(B')=B$, $0$ sinon.
\end{nota}

La proposition suivante exprime que les foncteurs $\wt{\omega}$ et
$\wt{\omega}'$ sont {\em duaux}.

\begin{pr}\label{pr-oop} Il existe des isomorphismes 
  $\wt{\omega}\circ D_\Gr\simeq D\circ\wt{\omega}'$ et
  $\wt{\omega}'\circ D_\Gr\simeq D\circ\wt{\omega}$ de foncteurs
  contravariants de $\wt{\F}_\Gr$ vers $\F$.
\end{pr}

\begin{proof} Cela découle de l'observation suivante : si $f : V\to
  V'$ est un morphisme de $\E^f$, $W$ un élément de $\Gr(V)$ et $W'$
  un élément de $\Gr(V')$, les conditions $W^\perp=\,^t f(W'^\perp)$ et
  $W=f^{-1}(W')$ sont équivalentes, où $^t f : V'^*\to V^*$ désigne la
  transposée de~$f$.
\end{proof}

La proposition suivante constitue le résultat principal de ce paragraphe.

\begin{pr}\label{comdom2} Les foncteurs $\wt{\omega}$ et
  $\wt{\omega}'$ sont isomorphes.
\end{pr}

\begin{proof} Soit $u_{X,V} : \widetilde{\omega}(X)(V)\to\widetilde{\omega}'(X)(V)$, pour $X\in {\rm Ob}\,\widetilde{\F}_\Gr$ et $V\in {\rm
  Ob}\,\E^f$, l'application linéaire linéaire dont la composante
$X(V,W)\to X(V,B)$ (où $W,B\in\Gr(V)$) est induite par l'inclusion $(V,W)\hookrightarrow
(V,B)$ si $W\subset B$, $0$ sinon. Alors $u_{X,V}$ est un isomorphisme d'espaces vectoriels, car
si l'on munit $\Gr(V)$ d'un ordre total $\leq$ tel que $W\subset B$
implique $W\leq B$, on obtient pour $u_{X,V}$ une matrice triangulaire
par blocs, avec des blocs diagonaux identiques.

De plus, pour toute application linéaire $f : V\to V'$,
le diagramme
$$\xymatrix{\widetilde{\omega}(X)(V)\ar[rr]^-{\widetilde{\omega}(X)(f)}\ar[d]_-{u_{X,V}}
  & & \widetilde{\omega}(X)(V')\ar[d]^-{u_{X,V'}} \\
\widetilde{\omega}'(X)(V)\ar[rr]^-{\widetilde{\omega}'(X)(f)} & & \widetilde{\omega}'(X)(V')
}$$
commute. En effet, la composante $X(V,W)\to X(V',B')$ de
l'application $\widetilde{\omega}(X)(V)\to \widetilde{\omega}'(X)(V')$ obtenue en suivant la composée
supérieure est la somme des applications induites par le morphisme
$(V,W)\to (V',B')$ induit par $f$ indexée sur les
$W'\in\Gr(V')$ tels que $f(W)=W'$ et $W'\subset B'$. Autrement dit, cette composante est l'application induite par le morphisme
$(V,W)\to (V',B')$ induit par $f$ si $f(W)\subset B'$ et $0$ sinon.

De même, la composante $X(V,W)\to X(V',B')$ de
l'application $\widetilde{\omega}(X)(V)\to \widetilde{\omega}'(X)(V')$ obtenue en suivant la composée
inférieure est  l'application induite par le morphisme
$(V,W)\to (V',B')$ induit par $f$ si $W\subset
f^{-1}(B')$ et $0$ sinon, d'où la commutativité recherchée.

Par conséquent, les applications linéaires $u_{X,V}$ définissent un
isomorphisme $u_X : \widetilde{\omega}(X)\xrightarrow{\simeq}\widetilde{\omega}'(X)$ de
$\F$, qui est naturel en $X$, d'où un isomorphisme de foncteurs
$u : \widetilde{\omega}\xrightarrow{\simeq}\widetilde{\omega}'$.
\end{proof}

Les propositions~\ref{comdom2} et~\ref{pr-oop} fournissent le
corollaire suivant.

\begin{cor}\label{comdom} Il existe un isomorphisme $\alpha : \widetilde{\omega}\circ
  D_\Gr\xrightarrow{\simeq} D\circ\widetilde{\omega}$
  tel que  pour tout objet $X$ de
  $\widetilde{\F}_\Gr$, le diagramme 
$$\xymatrix{\widetilde{\omega}(D^2_\Gr X)\ar[r]^-{\alpha_{D_\Gr X}} &
  D\widetilde{\omega}(D_\Gr X) \\
\widetilde{\omega}(X)\ar@{^{(}->}[r]\ar@{^{(}->}[u] & D^2\widetilde{\omega}(X)\ar[u]_-{D\alpha_X}
}$$
dont les monomorphismes non spécifiés sont les unités de l'adjonction commute.
\end{cor}

\begin{proof} Soient $a : \widetilde{\omega}\circ D_\Gr\xrightarrow{\simeq} D\circ\widetilde{\omega}'$ et
$b : \widetilde{\omega}'\circ D_\Gr\xrightarrow{\simeq}
D\circ\widetilde{\omega}$ les deux isomorphismes canoniques de la
proposition~\ref{pr-oop}, et $u :
\widetilde{\omega}\xrightarrow{\simeq}\widetilde{\omega}'$
l'isomorphisme de la proposition~\ref{comdom2}. Les deux composées $\widetilde{\omega}\circ D_\Gr\xrightarrow{a}D\circ \widetilde{\omega}'\xrightarrow{Du}D\circ \widetilde{\omega}$
et
$\widetilde{\omega}\circ D_\Gr\xrightarrow{u_{D_\Gr}}\widetilde{\omega}'\circ D_\Gr\xrightarrow{b}D\circ \widetilde{\omega}$
coïncident, car si $W$ et $B$ sont deux sous-espaces
de $V\in {\rm Ob}\,\E^f$, les conditions $W\subset B$ et
$B^\perp\subset W^\perp$ sont équivalentes. Notons $\alpha$
l'isomorphisme donné par ces
composées. La commutation du diagramme de
l'énoncé se ramène à celle de
$$\xymatrix{\widetilde{\omega}(D_\Gr^2 X)\ar[r]^-{a_{D_\Gr X}} &
  D\widetilde{\omega}'(D_\Gr X) \\
\widetilde{\omega}(X)\ar@{^{(}->}[r]\ar@{^{(}->}[u] & D^2\widetilde{\omega}(X)\ar[u]_-{Db_X}
}$$
qui se vérifie par inspection.
\end{proof}

Dans le paragraphe suivant, nous utiliserons le
corollaire~\ref{comdom} par le biais de sa conséquence directe suivante.

\begin{cor}\label{autd-drtf} Si $X$ est un objet auto-dual de
  $\widetilde{\F}_\Gr$, alors $\widetilde{\omega}(X)$ est un objet
  auto-dual de $\F$.
\end{cor}

\subsection{Structure de $\kk[\Gr]$}\label{pskg} La
proposition~\ref{prfig} identifiant la catégorie $\F_\Gr$ à celle des
$\kk[\Gr]$-comodules de $\F$ illustre l'importance de ce foncteur dans
l'étude des catégories de foncteurs en grassmanniennes.

\paragraph*{Structure fondamentale} La proposition~\ref{dsf-fsf}
montre que la décomposition scalaire de $\kk[\Gr]=\varpi(\kk)$ se réduit à
$\kk[\Gr]\simeq\kk\oplus\overline{\kk[\Gr]}$, où
$\overline{\kk[\Gr]}\in {\rm Ob}\,\F_{q-1}$. La proposition~\ref{pr-extfs} fournit quant à elle le résultat
  suivant.

\begin{pr}\label{str-fgrb} Le foncteur
$\kk[\Gr]\simeq\kk\oplus\overline{\kk[\Gr]}$ est colimite filtrante des sous-foncteurs
$\varpi(T_n(\kk))\simeq\kk[\Gr_{\leq n}]$ ; on a des extensions essentielles
$$0\to\overline{\kk[\Gr_{\leq n-1}]}\to\overline{\kk[\Gr_{\leq n}]}\to\kk[\Gr_n]\to 0$$
pour tout $n\in\mathbb{N}^*$. De plus, $\overline{\kk[\Gr_{\leq 1}]}\simeq\kk[\Gr_1]\simeq P_{\kk,q-1}$.
\end{pr}

Comme le foncteur $P_{\kk,q-1}$ est indécomposable, on en déduit le
corollaire suivant.

\begin{cor} Le foncteur $\overline{\kk[\Gr]}$ est indécomposable.
\end{cor}

Appliqué au foncteur constant
$\kk$, le corollaire~\ref{autd-drtf} fournit l'important résultat suivant.

\begin{pr}\label{autd-ins} Le foncteur $\kk[\Gr]$ est auto-dual.
\end{pr}

\`A partir de cette auto-dualité et de la structure de comodule de
$\kk[\Gr]$, on obtient le corollaire suivant.

\begin{cor}\label{cor-hopfgr} Munissons le foncteur
  $\kk[\Gr]$ du produit et de l'unité obtenus en dualisant sa
  structure de coalgèbre de Boole.

Le foncteur $\kk[\Gr]$ devient ainsi un {\em foncteur en algèbres de Hopf}.
\end{cor}

\begin{proof} En reprenant la
  démonstration de la proposition~\ref{comdom2} et en notant que $B\subset W_1^\perp$ et
  $B\subset W_2^\perp$ équivaut à $B\subset (W_1+W_2)^\perp$, où $B$,
  $W_1$ et $W_2$ sont des sous-espaces d'un  espace vectoriel de
  dimension finie $V$, on voit que la structure d'algèbre sur $\kk[\Gr]$ est donnée sur  $V$ par 
$$\kk[\Gr(V)]\otimes\kk[\Gr(V)]\to\kk[\Gr(V)]\qquad [W_1]\otimes
[W_2]\mapsto [W_1+W_2].$$

Ainsi, le produit de deux générateurs canoniques de $\kk[\Gr(V)]$ est encore un
générateur canonique de $\kk[\Gr(V)]$, ce qui montre que les
structures d'algèbre et de coalgèbre de Boole sur cet espace vectoriel
sont compatibles : c'est une algèbre de Hopf.
\end{proof}

\`A l'aide de la proposition~\ref{autd-ins} et du lemme suivant,
laissé au lecteur, nous allons décrire les
foncteurs ${\rm hom}_\F(\kk[\Gr],\cdot)$ et ${\rm
  Ext}^i_\F(\cdot,\kk[\Gr])$. On rappelle que ${\rm ev}_n :
\F_{surj}\to\,_{GL_n(\kk)}\mathbf{Mod}$ désigne le foncteur
d'évaluation sur $E_n$.

\begin{lm}\label{ide-fhcfs} Le foncteur ${\rm hom}_{\F_{surj}}
  (\kk,\cdot)$ est isomorphe à $\underset{n\in\mathbb{N}}{\lim}\,
  {\rm ev}_n^{GL_n(\kk)}$, où :
\begin{itemize}\item l'on désigne par ${\rm ev}_n^{GL_n(\kk)}$ le
  foncteur obtenu en prenant les invariants de ${\rm ev}_n$ sous
  l'action naturelle de $GL_n(\kk)$ ;
\item la limite est relative aux transformations naturelles
  ${\rm ev}_n^{GL_n(\kk)}\to {\rm ev}_m^{GL_m(\kk)}$ (pour $n\geq m$)
  induites par la projection $E_n\twoheadrightarrow E_m$ sur les $m$
  premières coordonnées.
\end{itemize}
\end{lm}

\begin{pr} Soient $F\in {\rm Ob}\,\F$ et $i\in\mathbb{N}$. Il existe
  des isomorphismes naturels
$${\rm hom}_\F(\kk[\Gr],F)\simeq\underset{n\in\mathbb{N}}{\lim}\,
F(E_n)^{GL_n(\kk)}$$
et
$${\rm Ext}^i_\F(F,\kk[\Gr])\simeq \big(\underset{n\in\mathbb{N}}{\col}
H_i(GL_n(\kk),F(E_n))\big)^*.$$
\end{pr}

\begin{proof} Le lemme \ref{ide-fhcfs} et l'adjonction entre $\varpi$ et
  $o$ (cf. proposition \ref{pre-vpo}) donnent le premier
  isomorphisme. Par auto-dualité de $\kk[\Gr]$ (proposition~\ref{autd-ins}), on en
déduit des isomorphismes naturels
$${\rm hom}_\F(F,\kk[\Gr])\simeq {\rm hom}_\F(\kk[\Gr],DF)\simeq\underset{n\in\mathbb{N}}{\lim}\,
DF(E_n)^{GL_n(\kk)}\simeq\big(\underset{n\in\mathbb{N}}{\col}
F(E_n)_{GL_n(\kk)}\big)^*.$$
L'intérêt de la dernière écriture est de remplacer la limite par une
colimite filtrante, exacte. Cela permet d'en déduire le second
isomorphisme, en dérivant les foncteurs considérés.
\end{proof}

\paragraph*{Anneau d'endomorphismes} Nous poursuivons l'étude du
foncteur $\kk[\Gr]$ par la détermination de son anneau
d'endomorphisme, donnée par le corollaire~\ref{agr-sf}.

\begin{defi}[Algèbre en grassmanniennes]\label{dfagr} On appelle {\em
    algèbre en grassmanniennes}, et l'on note $\agr(\kk)$, la $\kk$-algèbre
  ${\rm End}_\F (\kk[\Gr])$. L'{\em algèbre en grassmanniennes réduite}, notée $\agrr(\kk)$, est
  définie comme ${\rm End}_\F (\overline{\kk[\Gr]})$.
\end{defi}

Le scindement $\varpi(\kk)\simeq\kk\oplus\overline{\varpi(\kk)}$
fournit donc un isomorphisme d'algèbres $\agr(\kk)\simeq\kk\oplus\agrr(\kk)$.

\smallskip

Pour déterminer explicitement ces algèbres, nous emploierons le lemme
combinatoire élémentaire suivant. On rappelle que $p$ désigne la
caractéristique de $\kk$ et $q$ son cardinal.

\begin{lm}\label{lmcomb1} Soient $l$ et $i$ deux entiers naturels.
\begin{enumerate}\item  Soit $a$ un élément
    non nul de l'espace vectoriel $E_l$. Les classes d'équivalences de la relation définie par $W\sim
    W'$ si $W+a=W'+a$ sur l'ensemble des éléments $W$ de $\Gr_i(E_l)$ ne contenant pas $a$ sont
    de cardinal $q^i$, donc multiple de $p$ si $i>0$.
\item Si $i\leq l$, le cardinal de l'ensemble  $\Gr_i(E_l)$ est congru à $1$ modulo $p$.
\end{enumerate}
\end{lm}

\begin{proof} Pour le premier point, on note que le groupe $GL(W\oplus
  a,a)$ des automorphismes $u$ de $W\oplus a$ tels que $u(a)=a$ opère transitivement sur la classe d'équivalence de $W$, le
  stabilisateur de $W$ étant $GL(W)$. Le cardinal de cette classe est
  donc celui de  $GL(W\oplus a,a)/GL(W)$, soit $q^i$.

Pour la seconde assertion, on remarque que le cardinal de l'ensemble ${\rm
  Pl}_\E(E_i,E_l)$ est
$\prod_{j=0}^{i-1}(q^l-q^j)$. Comme l'ensemble $\Gr_i(E_l)$
s'identifie au quotient de ${\rm
  Pl}_\E(E_i,E_l)$ par l'action libre de $GL_i(\kk)$, on a pour $i\leq l$
$${\rm Card}\,\Gr_i(l)=\frac{{\rm Card}\,{\rm
  Pl}_\E(E_i,E_l)}{{\rm
    Card}\,GL_i(\kk)}=\prod_{j=0}^{i-1}\frac{q^{l-j}-1}{q^{i-j}-1}\equiv
1\quad (mod\; p).$$
\end{proof}

La proposition suivante fournit une première description de l'algèbre
en grassmanniennes.

\begin{pr}\label{idenagr} La $\kk$-algèbre $\agr(\kk)$ est isomorphe
  au $\kk$-espace vectoriel $\kk^\mathbb{N}$ muni de loi multiplicative $*$  définie par
$$(f*g)(n)=\underset{i,j\geq 0}{\sum_{i+j=n}}f(i)g(j)+\underset{i,j>0}{\sum_{i+j=n+1}}f(i)g(j).$$
\end{pr}

\begin{proof}  Par la proposition~\ref{pre-vpo}, on a un isomorphisme (linéaire) d'adjonction $\agr(\kk)\simeq {\rm hom}_{\F_{surj}}(\kk,o(\kk[\Gr]))$. Le $GL_n(\kk)$-module
  ${\rm ev}_n(o(\kk[\Gr]))$ est librement
  engendré comme $\kk$-espace vectoriel par $\Gr(E_n)$, sur lequel
  $GL_n(\kk)$ agit tautologiquement. Le sous-espace vectoriel ${\rm
    ev}_n(o(\kk[\Gr]))^{GL_n(\kk)}$ de
  $\kk[\Gr(E_n)]$ a pour base $s_0^n,s_1^n,\dots,s^n_n$, où $s_i^n$
  désigne la somme des générateurs canoniques $[B]$ associés à un
  sous-espace $B$ de dimension $i$ de $E_n$. L'application linéaire
  induite par la projection $E_{n+1}\twoheadrightarrow E_n$ envoie
  $s^{n+1}_i$ sur $s^n_{i-1}$ si $i>0$ et $s^{n+1}_0$ sur $s^n_0$,
  grâce à la première assertion du lemme~\ref{lmcomb1}. 

Par le lemme~\ref{ide-fhcfs}, on en déduit une identification
entre $\agr(\kk)$ et la limite $L$ des espaces vectoriels $\kk^{n+1}$ (dont
nous continuerons à noter $s_0^n,s_1^n,\dots,s^n_n$ une base
privilégiée)
relativement aux applications linéaires $f_n : \kk^{n+1}\to\kk^n$ données par $s^{n}_i\mapsto
s^{n-1}_{i-1}$ pour $i>0$ et $s^{n}_0\mapsto s^{n-1}_0$.

Soit $l_0^n,l_1^n,\dots,l^n_n$ la base duale de
$s_0^n,s_1^n,\dots,s^n_n$. L'application linéaire $a : L\to\kk^\mathbb{N}\qquad (v_n)_{n\in\mathbb{N}}\mapsto
(l_0^n(v_n))_{n\in\mathbb{N}}$ est bijective ; sa réciproque est donnée par
$$b : \kk^\mathbb{N}\to L\qquad (t_n)_{n\in\mathbb{N}}\mapsto \Big(t_n
s_0^n+\sum_{i=1}^n (t_{n-i}+t_{n-i+1})s_i^n\Big)_{n\in\mathbb{N}}\,.$$
En effet, $b$ prend bien ses valeurs dans $L$ parce que
$$f_n\Big(t_n
s_0^n+\sum_{i=1}^n (t_{n-i}+t_{n-i+1})s_i^n\Big)=t_n
s_0^{n-1}+\sum_{i=1}^n (t_{n-i}+t_{n-i+1})s_{i-1}^{n-1}$$
$$=t_{n-1} s_0^{n-1}+\sum_{j=1}^{n-1} (t_{n-1-j}+t_{n-j})s_j^{n-1}\,,$$
l'égalité $a\circ b=id$ est immédiate, et $b\circ a=id$ se déduit
facilement des remarques précécentes.

On a ainsi obtenu une identification de $\agr(\kk)$ et $\kk^\mathbb{N}$
comme {\em espaces vectoriels} ; il reste à lire la structure
multiplicative de $\agr(\kk)$ dans les isomorphismes précédents.

Soient $u$ un endomorphisme de $\kk[\Gr]$ et
$(t_n)_{n\in\mathbb{N}}$ l'élément de $\kk^\mathbb{N}$ correspondant. Par ce qui
précède, pour tout sous-espace $W$ de $V\in {\rm
  Ob}\,\E^f$, le générateur $[W]$ de $\kk[\Gr(V)]$ est envoyé par
$u(V)$ sur $t_m s_0(W)+\sum_{i=1}^m (t_{m-i}+t_{m-i+1})s_i(W)$, où $m=\dim
W$ et $s_i(W)$ est la somme des générateurs  de $\kk[\Gr(V)]$
associés aux sous-espaces de dimension $i$ de $V$. On en déduit que
$u_{E_n}$ envoie $s^n_m$ sur $t_m s_0^n+\sum_{i=1}^m
(t_{m-i}+t_{m-i+1})s_i^n$, puisque le coefficient d'un générateur $[B]$ (où $\dim
B=i$) dans cette image est égal au cardinal de  l'ensemble des sous-espaces $W$ de
$E_n$ de dimension $m$ contenant $B$, multiplié par $t_0$ si $i=0$,
$t_{m-i}+t_{m-i+1}$ si $1\leq i\leq m$, et $0$ sinon. Le cardinal en question n'est autre que celui des sous-espaces de dimension $m-i$
de $E_n/B$, égal à $1$ dans $\kk$ par la seconde assertion du lemme
\ref{lmcomb1}, d'où notre assertion.

Soient $u'$ un autre endomorphisme de $\kk[\Gr]$ et $u''=u'u$ ;
notons $(t'_n)_{n\in\mathbb{N}}$ et $(t''_n)_{n\in\mathbb{N}}$ les
éléments de $\kk^\mathbb{N}$ correspondant à $u'$ et $u''$
respectivement. On a
$$u''(E_n)([E_n])=u'\Big(t_n s_0^n+\sum_{i=1}^n
(t_{n-i}+t_{n-i+1})s_i^n\Big)$$
$$=t_n t'_0 s^n_0+\sum_{i=1}^n
(t_{n-i}+t_{n-i+1})\Big(t'_i
s_0^n+\sum_{j=1}^i(t'_{i-j}+t'_{i-j+1})s^n_j\Big)\,,$$
d'où en identifiant le coefficient de $s^n_0$
$$t''_n=t_n t'_0+\sum_{i=1}^n
(t_{n-i}+t_{n-i+1})t'_i=\underset{i,j\geq 0}{\sum_{i+j=n}}t_i
t'_j+\underset{i,j>0}{\sum_{i+j=n+1}}t_i t'_j\,,$$
ce qui achève la démonstration.
\end{proof}

\begin{cor}\label{cragrr} La loi $*$ sur l'espace vectoriel
  $\kk^{\mathbb{N}^*}$ définie par
$$(f*g)(n)=\underset{i,j\geq
  1}{\sum_{i+j=n}}f(i)g(j)+\underset{i,j\geq 1}{\sum_{i+j=n+1}}f(i)g(j)$$
fait de $\kk^{\mathbb{N}^*}$ une $\kk$-algèbre isomorphe à $\agrr(\kk)$.
\end{cor}

\begin{proof} Dans l'isomorphisme précédent, $\agrr(\kk)$ correspond à l'idéal de
  $\kk^\mathbb{N}$ des fonctions nulles en $0$.
\end{proof}

\begin{cor}\label{agr-sf} L'algèbre $\agrr(\kk)$ est une algèbre de séries
  formelles sur l'élément $\tau$ donné par la fonction
  $\mathbb{N}^*\to\kk$ associant $0$ à $1$ et $1$ à $n\geq 2$.
\end{cor}

\begin{proof} Pour tout $f\in\agrr(\kk)$, on a $\tau f(n)=f(n-1)$ si $n\geq
  2$ et $\tau f(1)=0$. On en déduit que $\tau^k$ est la fonction
  $n\mapsto 1$ si $n\geq k+1$, $0$ sinon. Cela fournit aussitôt le résultat.
\end{proof}

Explicitement, l'isomorphisme est donné par $\sum_{i\in\mathbb{N}}a_i\tau^i\mapsto \Big(\sum_{i=0}^{n-1}a_i\Big)_{n\in\mathbb{N}^*}$.

\begin{rem}\label{identt}\begin{enumerate}\item En reprenant la démonstration de la
  proposition~\ref{idenagr}, on voit que l'endomorphisme $\tau$ envoie
  un générateur canonique $[W]$ sur 
$$\underset{{\rm codim}_W B=1}{\sum_{B\in\Gr(W)}} [B].$$
\item Malgré le corollaire~\ref{agr-sf}, il pourra être utile de
  conserver la description de $\agr(\kk)$ donnée à la proposition~\ref{idenagr}. Cela apparaîtra clairement à la remarque~\ref{hopf-boole}.\,\ref{agr-ia}.
\end{enumerate}
\end{rem}


Nous terminons cette section par quelques considérations combinant
l'auto-dualité de $\kk[\Gr]$ (proposition~\ref{autd-ins}) et le corollaire~\ref{agr-sf}.

\begin{cor} L'involution de l'algèbre $\agr(\kk)$ induite par
  l'auto-dualité du foncteur $\kk[\Gr]$ est triviale.
\end{cor}

\begin{proof} Par le corollaire~\ref{agr-sf}, il suffit de le vérifier
  sur l'endomorphisme $\tau$ de $\kk[\Gr]$.

Pour cela, on note que la structure auto-duale de $\kk[\Gr]$ est
donnée par les formes bilinéaires
$$b_V : \kk[\Gr(V)]\times\kk[\Gr(V^*)]\to\kk\qquad ([W],[H])\mapsto 1
\text{ si } H\subset W^\perp, 0 \text{ sinon.}$$

On calcule alors, pour $W\in\Gr_n(V)$,
$$b_V(\tau([W]),[H])=\sum_{B\in\Gr_{n-1}(W)}b_V([B],[H])=$$
$${\rm
  Card}\,\{B\in\Gr_{n-1}(W)\,|\,H\subset B^\perp\}={\rm
  Card}\,\Gr_{n-1}(W\cap H^\perp).$$

En utilisant le lemme~\ref{lmcomb1}, on voit que cet élément de $\kk$ vaut
$1$ si $\dim (W\cap H^\perp)\geq n-1$, i.e. $\dim W-\dim (W\cap
H^\perp)\leq 1$, et $0$ sinon. Comme cette condition est symétrique en
$W$ et $H$ (modulo l'identification de $V$ à $V^{**}$), cela donne la conclusion.
\end{proof}

\begin{rem}\label{hopf-boole}\begin{enumerate}\item\label{agr-ia} On déduit la structure d'algèbre de Hopf sur $\kk[\Gr]$
donnée par le corollaire~\ref{cor-hopfgr} une nouvelle structure d'algèbre sur $\agr(\kk)$, donnée comme suit. Si $u$ et $v$
  sont deux éléments de $\agr(\kk)$, leur produit $u.v$ est le morphisme
  composé
$$\kk[\Gr]\to\kk[\Gr]\otimes\kk[\Gr]\xrightarrow{u\otimes
  v}\kk[\Gr]\otimes\kk[\Gr]\to\kk[\Gr]\,,$$
où la première flèche est le coproduit et la dernière le produit
(cf. \cite{K3}, §\,$5$).

Dans l'isomorphisme $\agr(\kk)\simeq\kk^\mathbb{N}$ de la proposition
\ref{idenagr}, le produit $.$ est le produit usuel d'algèbre de Boole
de l'anneau produit $\kk^\mathbb{N}$. Pour le voir, il suffit de
reprendre la démonstration de la proposition \ref{idenagr} et de
constater que, en conservant ses notations, le produit de $s^n_i$
et de $s^n_j$ (pour la structure d'algèbre de Boole de
$\kk[\Gr(E_n)]$) est $s^n_0$ si $i=j=0$, et que sinon le coefficient
de $s^n_0$ dans sa décomposition dans la base $s^n_0,\dots,s^n_n$ de
$\kk[\Gr(E_n)]^{GL_n(\kk)}$ est nul (en fait, on a
$s^n_i.s^n_j=s^n_{max(i,j)}$). Avec la description du
corollaire~\ref{agr-sf}, on a $\tau^i.\tau^j=\tau^{max(i,j)}$.
\item Si $X$ est un objet de $\F_\Gr$, on peut munir naturellement le
  groupe abélien ${\rm hom}_\Gr(\omega(X),\kk[\Gr])$ d'une structure de
  module sur l'algèbre de Boole $(\agr(\kk),.)$, en utilisant la structure
  de $\kk[\Gr]$-comodule de $\omega(X)$ et en procédant comme au point
  précédent. Cela suggère une structure algèbrique très riche sur
  ces objets, qui sont également, comme tous les groupes ${\rm hom}_\Gr(F,\kk[\Gr])$, des modules sur l'algèbre de séries
  formelles $(\agr(\kk),*)\simeq\kk\oplus\kk[[\tau]]$.
\end{enumerate}
\end{rem}

\section{La catégorie $\F_{\Gr,I}$ comme catégorie de modules}\label{fgm}

L'équivalence entre la catégorie $\F_\Gr$ et la catégorie des
$\kk[\Gr]$-comodules repose sur l'adjonction entre les foncteurs
$\iota$ et $\omega$. Nous donnons maintenant une équivalence
entre la catégorie $\F_\Gr$ (et plus généralement, toutes les
catégories $\F_{\Gr,I}$) et une catégorie de modules sur une monade
explicite à partir de l'adjonction entre les foncteurs $\xi$ et~$\sigma$.

\begin{conv} Dans cette section, on se donne une partie $I$ de $\mathbb{N}$.
\end{conv}

On rappelle que le foncteur $\xi_I : \F\otimes\F^I_{surj}\to\F_{\Gr,I}$ est adjoint
à gauche à $\sigma_I$ (proposition~\ref{pradjf3}). Les transformations
naturelles introduites dans la notation suivante seront identifiées
dans la proposition~\ref{adjidis}.

\begin{nota}\label{nds} Nous désignerons par $\T_{\Gr,I}=(\widetilde{\Delta}^I_{surj},u_{\Gr,I},\mu_{\Gr,I})$ la
  monade sur $\F\otimes\F^I_{surj}$ associée à l'adjonction entre les foncteurs
  $\xi_I$ et $\sigma_I$
  conformément à la proposition \ref{mon-adj}. Ainsi :
\begin{enumerate}\item on a $\widetilde{\Delta}^I_{surj}=\sigma_I\circ
\xi_I$, soit $\widetilde{\Delta}^I_{surj}
(F)(A,B)=F(A\oplus B,B)$ sur les objets ($F\in {\rm
  Ob}\,\F\otimes\F^I_{surj}$, $A\in {\rm Ob}\,\E^f$, $B\in {\rm Ob}\,\E^I_{surj}$).
\item La transformation naturelle $u_{\Gr,I} : id\to \widetilde{\Delta}^I_{surj}$ est l'unité de
l'adjonction. 
\item La
transformation naturelle $\mu_{\Gr,I} :
(\widetilde{\Delta}^I_{surj})^2\to\widetilde{\Delta}^I_{surj}$ est
donnée par $\sigma_I (v_{\xi_I})$, où $v$ désigne la coünité de l'adjonction.
\end{enumerate}
\end{nota}

\begin{pr}\label{adjidis} Soit $F$ un objet de $\F\otimes\F^I_{surj}$.
\begin{enumerate}\item L'unité $u_{\Gr,I} :
    id\hookrightarrow\widetilde{\Delta}^I_{surj}$
    est la transformation naturelle injective telle que $((u_{\Gr,I})_F)_{(A,B)} : F(A,B)\to
    F(A\oplus B,B)$ est induit par le monomorphisme canonique
    $(A,B)\hookrightarrow (A\oplus B,B)$ pour tous $F\in {\rm
      Ob}\,\F\otimes\F^I_{surj}$, $A\in {\rm
      Ob}\,\E^f$ et $B\in {\rm
      Ob}\,\E^I_{surj}$.
\item La multiplication
  $\mu_{\Gr,I} : (\widetilde{\Delta}^I_{surj})^2\to\widetilde{\Delta}^I_{surj}$  est fournie par le morphisme
$F(A\oplus B\oplus B,B)\to F(A\oplus B,B)$ induit le morphisme
$A\oplus B\oplus B\to A\oplus B$ somme directe de $id_A$ et de la
somme $B\oplus B\to B$, et par le morphisme identique $B\to B$.
\item Il existe un scindement naturel
$$F\xrightarrow{(u_{\Gr,I})_F}\widetilde{\Delta}^I_{surj} (F)\xrightarrow{(p_{\Gr,I})_F} F$$
où $(p_{\Gr,I})_F$ est donné sur l'objet $(A,B)$ par le morphisme induit par
l'épimorphisme canonique $(A\oplus
B,B)\twoheadrightarrow (A,B)$, pour tous $F\in {\rm Ob}\,\F\otimes\F^I_{surj}$, $A\in {\rm
      Ob}\,\E^f$ et $B\in {\rm Ob}\,\E^I_{surj}$. 

De plus, $(F,(p_{\Gr,I})_F)$ est un module sur la monade $\T_{\Gr,I}$.
\end{enumerate}
\end{pr}

\begin{proof} Analysons la
  monade associée à l'adjonction entre les foncteurs
  $\mathfrak{O}_I\times\mathfrak{B}_I : \E^f_{\Gr,I}\to\E^f\times\E^I_{surj}$ et
  $\mathfrak{L}_I : \E^f\times\E^I_{surj}\to\E^f_{\Gr,I}$
  (cf. proposition~\ref{adjfpr}). Son unité est la transformation
  naturelle $id\to (\mathfrak{O}_I\times\mathfrak{B}_I)\mathfrak{L}_I$
  donnée par l'inclusion $(A,B)\hookrightarrow (A\oplus B,B)$ (de
  composantes l'inclusion du facteur direct $A$ et $id_B$). Sa
  multiplication est donnée par le morphisme
$(A\oplus B\oplus B,B)\to (A\oplus B,B)$ induit le morphisme
$A\oplus B\oplus B\to A\oplus B$ somme directe de $id_A$ et de la
somme $B\oplus B\to B$, et par le morphisme identique $B\to B$. En
effet, la coünité de l'adjonction est donnée sur l'objet $(V,W)$ de
$\E^f_{\Gr,I}$ par le morphisme $(V\oplus W,W)\to (V,W)$ de
composantes $id_V$ et $W\hookrightarrow V$.

En utilisant la proposition~\ref{lm-form}, on en déduit les deux
premières assertions. La dernière est immédiate.
\end{proof}

Le dernier point de cette proposition conduit à donner la définition
suivante, qui introduit une sorte de foncteur de différence dans
$\F\otimes\F^I_{surj}$, distinct de l'endofoncteur $\Delta_*$
considéré à la fin du paragraphe~\ref{subs-fd}.

\begin{defi}\label{notdsurj} Le noyau de $(p_{\Gr,I})_F$, qui s'identifie donc au conoyau de $(u_{\Gr,I})_F$, sera noté
$\Delta^I_{surj}(F)$. On définit ainsi un endofoncteur exact
$\Delta^I_{surj}$ de $\F\otimes\F^I_{surj}$.
\end{defi}

Le reste de cette section s'emploie à tirer les conséquences du
résultat suivant.

\begin{pr}\label{fgr-mon} La catégorie $\F_{\Gr,I}$ est équivalente à
  la catégorie des modules sur la monade $\T_{\Gr,I}$ de
  $\F\otimes\F^I_{surj}$.
\end{pr}

\begin{proof} Il s'agit d'un cas particulier de la proposition~\ref{dualmonadique}.
\end{proof}

\begin{conv}
Dans la suite de cette section, nous {\em identifierons} la catégorie
$\F_{\Gr,I}$ avec la sous-catégorie des modules sur $\T_{\Gr,I}$ de
$\F\otimes\F^I_{surj}$. Autrement dit, un objet de $\F_{\Gr,I}$ sera
désormais un couple $(X,\widetilde{\Delta}^I_{surj}
X\xrightarrow{\tilde{m}_X}X)$, où $X$ est un objet de $\F\otimes\F^I_{surj}$ et
$\tilde{m}_X$ un morphisme tel que :
\begin{enumerate}\item la composée $X\xrightarrow{(u_{\Gr,I})_X}\widetilde{\Delta}^I_{surj}
X\xrightarrow{\tilde{m}_X}X$ est le morphisme identique ;
\item les composées $(\widetilde{\Delta}^I_{surj})^2 X\xrightarrow{(\mu_{\Gr,I})_X}\widetilde{\Delta}^I_{surj}
X\xrightarrow{\tilde{m}_X}X$ et  $(\widetilde{\Delta}^I_{surj})^2 X\xrightarrow{\widetilde{\Delta}^I_{surj}\tilde{m}_X}\widetilde{\Delta}^I_{surj}
X\xrightarrow{\tilde{m}_X}X$ coïncident.
\end{enumerate}

Par abus, nous noterons souvent simplement $X$ pour
$(X,\tilde{m}_X)$. Nous
désignerons aussi par $m_X$ le morphisme $\Delta^I_{surj}X\hookrightarrow\widetilde{\Delta}^I_{surj}
X\xrightarrow{\tilde{m}_X}X$.

Avec ces conventions, les morphismes $(X,\tilde{m}_X)\to (Y,\tilde{m}_Y)$ de $\F_{\Gr,I}$ sont
les morphismes $f : X\to Y$ de $\F\otimes\F^I_{surj}$ tels que le
diagramme
$$\xymatrix{\widetilde{\Delta}^I_{surj}X\ar[r]^-{\widetilde{\Delta}^I_{surj}f}\ar[d]_{\tilde{m}_X}
  & \widetilde{\Delta}^I_{surj}Y\ar[d]^{\tilde{m}_Y} \\
X\ar[r]^f & Y
}$$
commute. Cette condition est équivalente à la commutation du diagramme analogue
sans tilde.
\end{conv}

Le lien avec la définition originelle de la catégorie $\F_{\Gr,I}$
s'obtient à partir des remarques suivantes :
\begin{enumerate}\item le $\T_{\Gr,I}$-module associé à un foncteur $X
  : \E^f_{\Gr,I}\to\E$ est $\sigma_I(X)$ (muni de la multiplication
  dérivant de l'adjonction de la proposition~\ref{pradjf3}) ;
\item le foncteur $\E^f_{\Gr,I}\to\E$ associé à un $\T_{\Gr,I}$-module
  $X$ est le coégalisateur de $\xi_I(\tilde{m}_X)$ et de la flèche canonique
 $\xi_I(\widetilde{\Delta}^I_{surj}X)\to\xi_I(X)$ (adjointe à $id_{\widetilde{\Delta}^I_{surj}X}$). 
\end{enumerate}

%

\smallskip

Nous identifions maintenant le foncteur $\theta_I$ en termes de
modules sur la monade $\T_{\Gr,I}$.

\begin{lm}\label{lmthe} Soient $F$ un objet de $\F\otimes\F^I_{surj}$,
et  $(A,B)$ un objet de $\E^f_{\Gr,I}$. Notons
  $p_{A,B} : A\oplus B\twoheadrightarrow A$ et $\pi_{A,B} :
  A\twoheadrightarrow A/B$ les projections canoniques et
  $q_{A,B} : A\oplus B\to A$ le morphisme dont les
  composantes sont $id_A$ et l'inclusion $B\hookrightarrow A$. La suite
\begin{equation}\label{sekp} F(A\oplus B,B)\xrightarrow{F(p_{A,B},id_B)+F(q_{A,B},id_B)}
  F(A,B)\xrightarrow{F(\pi_{A,B},id_B)} F(A/B,B)\to 0
\end{equation}
de $\E$ est exacte.
\end{lm}

\begin{proof} Notons $i : B\hookrightarrow A$ l'inclusion : on a
  $\pi_{A,B}\circ i=0$, donc $\pi_{A,B}\circ q_{A,B}=\pi_{A,B}\circ
  p_{A,B}$, ce qui montre que la suite (\ref{sekp}) est un
  complexe. La surjectivité de $F(\pi_{A,B},id_B)$ provient de ce que
  $(\pi_{A,B},id_B) : (A,B)\to (A/B,B)$ admet une section.

Pour établir l'exactitude en $F(A,B)$, considérons une rétraction $r :
A\twoheadrightarrow B$ de $i$ et notons $u : (A,B)\to (A\oplus B, B)$
le morphisme $(id_A\oplus r,id_B)$. Alors $(p_{A,B},id_B)\circ
u=id_{(A,B)}$, tandis que $(q_{A,B},id_B)\circ u$ est nul sur $(B,0)$,
donc se factorise par $(\pi_{A,B},id_B)$. Par conséquent, la
restriction à $N=ker\,F(\pi_{A,B},id_B)$ de $F(u)$ est une section du
morphisme $F(A\oplus B,B)\to N$ induit par
$F(p_{A,B},id_B)+F(q_{A,B},id_B)$, ce qui achève la démonstration.
\end{proof}

\begin{rem}\label{rqkap-ch6} Cette suite exacte est une partie d'une
  suite exacte longue dépendant d'une structure simpliciale (cf. proposition~\ref{blpr}).
\end{rem}

On  déduit du lemme~\ref{lmthe}, compte-tenu des remarques précédentes sur le lien entre
les deux descriptions de $\F_{\Gr,I}$, la proposition suivante.

\begin{pr}\label{identk} Le foncteur $\theta_I :
  \F\otimes\F^I_{surj}\to\F_{\Gr,I}$ est donné par
  $\theta_I(F)=(F,(p_{\Gr,I})_F : \widetilde{\Delta}^I_{surj} F\twoheadrightarrow F)$
  sur les objets ---
  cf. proposition~\ref{adjidis}.\,3 --- et par l'égalité
  ${\rm hom}_{\Gr,I}(\theta_I(F),\theta_I(G))={\rm
    hom}_{\F\otimes\F^I_{surj}}(F,G)$ sur les morphismes.

Autrement dit, $\theta_I$ identifie $\F\otimes\F^I_{surj}$ à la
sous-catégorie pleine de $\F_{\Gr,I}$ formée des objets $X$ tels que $m_X=0$.
\end{pr}

\begin{rem} Ce résultat fournit une seconde démonstration de la pleine
  fidélité du foncteur $\theta_I$ et de ce que son image est une
  sous-catégorie de Serre de $\F_{\Gr,I}$ (cf. proposition~\ref{prfac-ff}).
\end{rem}

\begin{ex}\label{sntun} Considérons le cas où $I$ est réduit à
  l'entier $1$. On a $\Delta^1_{surj}=\Delta$ ; le foncteur $\omega_1
  : \F_{\Gr,1}\to\F$ associe à un objet $(X,\tilde{m}_X)$ le
  coégalisateur de la flèche canonique $\bar{P}\otimes\Delta
  X\to\bar{P}\otimes X$ (composée de $\bar{P}\otimes\Delta X\xrightarrow{j\otimes\Delta
    X}\bar{P}^{\otimes 2}\otimes\Delta X$, où $j$ est le coproduit de
  $\bar{P}$, et de $\bar{P}^{\otimes 2}\otimes\Delta
  X\xrightarrow{\bar{P}\otimes f}\bar{P}\otimes X$, où $f$ est la
  coünité de l'adjonction) et de $\bar{P}\otimes m_X$. Cette
  description de $\omega_1$ se déduit aisément de l'exactitude de ce
  foncteur de l'isomorphisme canonique $\omega_1\circ\iota_1\simeq\bar{P}\otimes\cdot$.

Par conséquent, la proposition~\ref{identk} montre que, pour tout objet $F$ de $\F$, le foncteur
  $\omega_1\kappa_1(F)$ est canoniquement isomorphe au conoyau de
  l'application naturelle $\bar{P}\otimes\Delta F\to\bar{P}\otimes F$.
\end{ex}

\subsection{Le foncteur $\eta_I$}\label{sfe-a} Les considérations précédentes
permettent de décrire très naturellement l'adjoint à gauche au
foncteur $\theta_I$.

\begin{defi} On définit un foncteur $\eta_I :
  \F_{\Gr,I}\to\F\otimes\F^I_{surj}$ par
  $\eta_I(X)=coker\,m_X$ sur les objets ; si $f : X\to Y$ est un
  morphisme de $\T_{\Gr,I}$-modules, $\eta_I(f) :
  \eta_I(X)\to\eta_I(Y)$ est le morphisme induit par $f$. 
\end{defi}

Ainsi, le foncteur $\eta_I : \F_{\Gr,I}\to\F\otimes\F^I_{surj}$ est un
quotient du foncteur $\sigma_I$. 

\begin{lm}\label{lmeta} Pour tout foncteur $X$ de $\F_{\Gr,I}$, il
  existe un épimorphisme naturel $X\twoheadrightarrow\theta_I\eta_I
  (X)$. De plus, $\theta_I\eta_I
  (X)$ est le plus grand quotient de $X$ appartenant à l'image du
  foncteur $\theta_I : \F\otimes\F^I_{surj}\to\F_{\Gr,I}$.
\end{lm}

\begin{proof} Le diagramme
$$\xymatrix{\Delta^I_{surj} \sigma_I(X)\ar[r]^-{m_X}\ar@{>>}[d] & \sigma_I(X)\ar@{>>}[d] \\
\Delta^I_{surj} \eta_I(X)\ar[r]_-{0} & \eta_I(X)
}$$
de $\F\otimes\F^I_{surj}$ commute, de sorte que la projection $\sigma_I(X)\twoheadrightarrow\eta_I(X)$
induit un épimorphisme $X\twoheadrightarrow\theta_I\eta_I(X)$ dans $\F_{\Gr,I}$. La proposition
\ref{identk} et  la propriété universelle du conoyau montrent ensuite que tout
épimorphisme de $X$ sur un foncteur dans l'image de $\theta_I$ se
factorise par $X\twoheadrightarrow\theta_I\eta_I(X)$.
\end{proof}

Le principal résultat de ce paragraphe est le suivant.

\begin{pr}\label{adjeta} Le foncteur $\eta_I$ est adjoint à gauche au
  foncteur $\theta_I : \F\otimes\F^I_{surj}\to\F_{\Gr,I}$.
\end{pr}

\begin{proof} C'est une conséquence directe du lemme précédent,
  puisque le foncteur $\theta_I$ est pleinement fidèle (cf. proposition~\ref{prfac-ff}).
\end{proof}

Les deux énoncés suivants sont à rapprocher du
corollaire~\ref{cr-cosoc} relatif aux cosocles.

\begin{cor}\label{cosoceta} Il existe un isomorphisme
$${\rm cosoc}\,X\simeq\theta_I({\rm cosoc}\,\eta_I(X))$$
naturel en l'objet $X$ de $\F_{\Gr,I}$.
\end{cor}

\begin{proof} Ce corollaire se déduit des propositions~\ref{prfffgr}
  et~\ref{adjeta}.
\end{proof}

\begin{pr}\label{propeta} Il existe des isomorphismes d'endofoncteurs
  de $\F\otimes\F^I_{surj}$
$$\eta_I\circ\xi_I\simeq\eta_I\circ\theta_I\simeq id.$$
\end{pr}

\begin{proof} L'isomorphisme $\eta_I\xi_I\simeq id$ vient de ce que le foncteur
  $\eta_I\xi_I$ est adjoint à gauche à
  $\sigma_I\theta_I\simeq id$ (proposition~\ref{prcompf3}).

L'isomorphisme $\eta_I\theta_I\simeq id$ découle de la proposition
\ref{identk} et de la définition de~$\eta_I$.
\end{proof}

Nous signalons dans la remarque qui suit un autre lien entre les
foncteurs $\eta_I$, $\xi_I$, $\sigma_I$ et $\theta_I$.

\begin{rem}\label{diageta} On vérifie qu'il existe un diagramme commutatif
 {\em cocartésien}
$$\xymatrix{\xi_I\,\sigma_I (X)\ar@{>>}[r]\ar@{>>}[d] & X\ar@{>>}[d] \\
\theta_I \sigma_I (X)\ar@{>>}[r] & \theta_I \eta_I(X)
}$$
naturel en l'objet $X$ de $\F_{\Gr,I}$.
\end{rem}

La proposition suivante, qui ne découle pas formellement de
l'adjonction entre les foncteurs $\eta_I$ et $\theta_I$, illustre
l'utilité de la description monadique explicite de la catégorie $\F_{\Gr,I}$.

\begin{pr}\label{pretapt} Il existe un isomorphisme
$$\eta_I(X\otimes Y)\simeq \eta_I(X)\otimes\eta_I(Y)$$
naturel en les objets $X$ et $Y$ de $\F_{\Gr,I}$.
\end{pr}

\begin{proof} La commutation du foncteur $\widetilde{\Delta}^I_{surj}$ au
  produit tensoriel procure un isomorphisme canonique $\Delta^I_{surj}
  (X\otimes Y)\simeq (\Delta^I_{surj} X\otimes Y)\oplus
  (X\otimes\Delta^I_{surj} Y)\oplus(\Delta^I_{surj}
  X\otimes\Delta^I_{surj} Y)$, dans lequel $m_{X\otimes Y}$ se lit
  comme le morphisme de composantes $\Delta^I_{surj} X\otimes
  Y\xrightarrow{m_X\otimes Y}X\otimes Y$, $X\otimes\Delta^I_{surj}
  Y\xrightarrow{X\otimes m_Y}X\otimes Y$ et $\Delta^I_{surj} X\otimes\Delta^I_{surj}
  Y\xrightarrow{m_X\otimes m_Y}X\otimes Y$. Par conséquent, l'image de
  $m_{X\otimes Y}$ est la somme des sous-objets $im\,(m_X)\otimes Y$,
  $X\otimes im\,(m_Y)$ et $im\,(m_X)\otimes im\,(m_Y)$ de $X\otimes Y$,
  i.e. $im\,(m_X)\otimes Y+X\otimes im\,(m_Y)=ker\,(X\otimes
  Y\twoheadrightarrow\eta_I(X)\otimes\eta_I(Y))$. Ainsi,
  $coker\,m_{X\otimes Y}=\eta_I(X\otimes Y)$ s'identifie à  $\eta_I(X)\otimes\eta_I(Y)$.
\end{proof}

\subsection{Résolution canonique et algèbre homologique monadique}
L'intérêt majeur de l'adjonction entre les foncteurs $\xi_I : \F\otimes\F^I_{surj}\to\F_{\Gr,I}$ et
$\sigma_I$, qui sous-tend toute cette section, réside dans la
possibilité de construire une résolution naturelle d'un objet de
$\F_{\Gr,I}$ par des foncteurs dans l'image du foncteur $\xi_I$. Cette
résolution, donnée par la proposition suivante, permet de ramener le comportement homologique de la
catégorie $\F_{\Gr,I}$ à celui de la catégorie $\F\otimes\F^I_{surj}$.

\begin{pr}\label{blpr} Il existe une suite exacte
$$\cdots\to
   \xi_I(\Delta^I_{surj})^n\sigma_I(X)\to
  \xi_I(\Delta^I_{surj})^{n-1}\sigma_I(X)\to\cdots$$
\begin{equation}\label{beqrc}\cdots\to\xi_I\Delta^I_{surj}\sigma_I(X)\to\xi_I\,\sigma_I(X)\to X\to 0
\end{equation}
naturelle en l'objet $X$ de $\F_{\Gr,I}$.
\end{pr}

\begin{proof}  Nous revenons à la définition originelle de
  $\F_{\Gr,I}$, et notons $\bf{E}$ la sous-catégorie pleine de
  $\mathbf{Ens}^f$ formée des ensembles finis non vides.

On définit un foncteur $a : \E^f_{\Gr,I}\times
\mathbf{E}^{op}\to\E^f_{\Gr,I}$ en associant à un objet $(V,W)$ de
$\E^f_{\Gr,I}$ et à un ensemble fini non vide $E$ la somme amalgamée de $V$ et de $W^E$ relativement aux monomorphismes
$W\hookrightarrow V$ (inclusion) et $W\hookrightarrow W^E$ (plongement
diagonal --- on utilise ici la non-vacuité de $E$) muni de la base
$W$. L'action sur les morphismes se déduit de la fonctorialité de
l'association $\E^f\times (\mathbf{Ens}^f)^{op}\to\E^f\quad (V,E)\mapsto
V^E\simeq V\otimes\kk^E$. On remarque que si $E$ est de cardinal $n$,
l'endofoncteur $a(\cdot,E)$ de $\E^f_{\Gr,I}$ est isomorphe à la
$n-1$-ième itérée du foncteur $\mathfrak{L}_I\circ
  (\mathfrak{O}_I\times\mathfrak{B}_I)$.

Par restriction à la sous-catégorie simpliciale $\mathbf{\Delta}$ de $\bf{E}$ (squelette
de la sous-catégorie des ensembles totalement ordonnés, les morphismes
étant les applications croissantes), on en déduit un foncteur
$\E^f_{\Gr,I}\times\mathbf{\Delta}^{op} \to\E^f_{\Gr,I}$, puis
$c^{\mathbf{\Delta}}_{\Gr,I} : \E^f_{\Gr,I}\to\mathbf{Fct}(\mathbf{\Delta}^{op},\E^f_{\Gr,I})$.
Ce foncteur vérifie les propriétés suivantes :
\begin{enumerate}\item en degré zéro, on a
  $(c^{\mathbf{\Delta}}_{\Gr,I})_0\simeq id$ ;
\item en degré un, on a
  $(c^{\mathbf{\Delta}}_{\Gr,I})_1\simeq\mathfrak{L}_I\circ
  (\mathfrak{O}_I\times\mathfrak{B}_I)$ ;
\item plus généralement, en degré $n\in\mathbb{N}^*$, on a  $(c^{\mathbf{\Delta}}_{\Gr,I})_n\simeq\mathfrak{L}_I\circ\big((\mathfrak{O}_I\times\mathfrak{B}_I)\mathfrak{L}_I\big)^{n-1}\circ
  (\mathfrak{O}_I\times\mathfrak{B}_I)$.
\end{enumerate}

Par précomposition, on en déduit un foncteur
$$C^{\mathbf{\Delta}}_{\Gr,I} : \F_{\Gr,I}=\mathbf{Fct}(\E^f_{\Gr,I},\E)\to\mathbf{Fct}(\E^f_{\Gr,I}\times
\mathbf{\Delta}^{op},\E)\simeq\mathbf{Fct}(\mathbf{\Delta}^{op},\F_{\Gr,I}).$$
Les remarques précédentes montrent que
$(C^{\mathbf{\Delta}}_{\Gr,I})_0\simeq id$ et
$(C^{\mathbf{\Delta}}_{\Gr,I})_n\simeq\xi_I(\widetilde{\Delta}^I_{surj})^{n-1}\sigma_I$
pour $n>0$. 

Le complexe de Moore associé à cet objet simplicial fournit dans
$\F_{\Gr,I}$ un complexe
$$\cdots\to
   \xi_I(\widetilde{\Delta}^I_{surj})^n\sigma_I(X)\to
  \xi_I(\widetilde{\Delta}^I_{surj})^{n-1}\sigma_I(X)\to\cdots$$
$$\cdots\to\xi_I\widetilde{\Delta}^I_{surj}\sigma_I(X)\to\xi_I\,\sigma_I(X)\to X\to 0$$
naturel en $X$. Nous allons montrer qu'il est acyclique. \'Evaluée sur
un objet $(V,W)$ de $\E^f_{\Gr,I}$, sa
différentielle $\partial_n : (C^{\mathbf{\Delta}}_{\Gr,I})_{n}\to
(C^{\mathbf{\Delta}}_{\Gr,I})_{n-1}$ est la somme alternée des $n+1$ morphismes
induits par $[v,w_0,\dots,w_n]\mapsto [v,w_0,\dots,\hat{w}_i,\dots,w_n]$ (le chapeau
signifiant que le terme considéré doit être omis), où l'on désigne par
$[v,w_0,\dots,w_n]$ la classe dans $V\underset{W}{\oplus}W^{n+1}$
(base de $(c^{\mathbf{\Delta}}_{\Gr,I})_n(V,W)$) de
$(v,w_0,\dots,w_n)\in V\oplus W^{n+1}$. On obtient donc une homotopie
entre les endomorphismes nul et identique de ce complexe en considérant les morphismes
induits par $[v,w_0,\dots,w_n]\mapsto [v,-\pi(v),w_0,\dots,w_n]$, où
$\pi : V\to W$ est un projecteur.

La suite exacte de l'énoncé s'obtient en considérant le complexe
normalisé associé à $C^{\mathbf{\Delta}}_{\Gr,I}$, qui est
homotopiquement équivalent au complexe de Moore, donc également acyclique.
\end{proof}

\begin{rem}\label{rqos1}\begin{enumerate}\item La proposition~\ref{blpr} repose uniquement sur la
  proposition~\ref{adjidis}.\,3 (relative au scindement naturel
  de la monade $\T_{\Gr,I}$).  Nous avons préféré en donner
  une démonstration directe car les objets simpliciaux de
  $\E^f_{\Gr,I}$ qui apparaissent sont particulièrement naturels, et
  plus parlants que la construction générale d'un objet simplicial
  canonique à partir d'une monade scindée. 
\item Les monades et les comonades (en particulier, celles
  qui proviennent d'adjonctions) fournissent un cadre général efficace
  pour faire de l'algèbre homologique, y compris dans un contexte non
  abélien ; la proposition~\ref{blpr} et les quelques conséquences que
  nous développons en sont un cas particulier. On trouvera dans
  \cite{BB} une exposition systématique de ces
  notions. 

On prendra garde au fait que la notion de résolution canonique dans ce cadre général est légèrement différente (on part d'une
{\bf co}monade pour obtenir une résolution homologique). En appliquant
le foncteur $\sigma_I$ à la résolution de la proposition~\ref{blpr},
on obtient la résolution canonique pour la {\bf co}monade associée à
l'adjonction entre $\xi_I$ et $\sigma_I$ qui sert de point de départ à
\cite{BB}. La possibilité de \go relever\gf cette résolution dans
$\F\otimes\F^I_{surj}$ en une résolution  dans $\F_{\Gr,I}$ provient
du scindement de la monade $\T_{\Gr,I}$ ; alors que la résolution
initiale dans $\F\otimes\F^I_{surj}$ n'apporte essentiellement rien,
la proposition~\ref{blpr} constitue un résultat important sur la
structure de $\F_{\Gr,I}$.
\end{enumerate}
\end{rem}

\begin{defi}[Résolution canonique]\label{rescan} Le complexe concentré
  en degrés positifs
$$\cdots\to
   \xi_I(\Delta^I_{surj})^n\sigma_I(X)\to
  \xi_I(\Delta^I_{surj})^{n-1}\sigma_I(X)\to\cdots\to\xi_I\Delta^I_{surj}\sigma_I(X)\to\xi_I\,\sigma_I(X)$$
de la proposition~\ref{blpr}  est appelé {\em résolution
  canonique}\index{Gal}{r\'esolution canonique} de l'objet $X$ de $\F_{\Gr,I}$.
Nous la noterons $\Re^{\Gr,I}_\bullet (X)$. Ainsi, $\Re^{\Gr,I}_n
  (X)=\xi_I(\Delta^I_{surj})^n\sigma_I(X)$ si
  $n\geq 0$, $0$ sinon. 
\end{defi}

Le résultat suivant constitue l'une des conséquences les plus
notables de la proposition~\ref{blpr}.

\begin{cor}\label{crss} Il existe une suite spectrale du premier
  quadrant (donc convergente) naturelle en les objets $X$ et $Y$ de
  $\F_{\Gr,I}$ dont le terme $E^1$ est donné par
$$E^1_{p,q}={\rm Ext}^p_{\F\otimes\F^I_{surj}}
  \big((\Delta^I_{surj})^q \sigma_I(X),\sigma_I(Y)\big)$$
et dont l'aboutissement est ${\rm Ext}^*_{\Gr,I}(X,Y)$.
\end{cor}

\begin{proof} On considère les deux suites spectrales associées au
  bicomplexe ${\rm hom}_{\Gr,I} (\Re^{\Gr,I}_\bullet (X),J^* (Y))$ (cf. \cite{ML}, chapitre~XI, §\,6),
  où $J^* (Y)$ désigne une résolution injective de $Y$, que l'on
  peut choisir fonctorielle en $Y$, car la catégorie $\F_{\Gr,I}$
  possède un cogénérateur injectif.  La suite spectrale obtenue en
  prenant d'abord la différentielle de $\Re^{\Gr,I}_\bullet (X)$
  dégénère au terme $E^2$, donné par ${\rm Ext}^*_{\Gr,I}(X,Y)$,
  tandis que celle obtenue en considérant d'abord la différentielle de
  $J^* (Y)$ a le terme $E^1$ indiqué dans l'énoncé, par adjonction entre les foncteurs exacts $\xi_I$ et $\sigma_I$.
\end{proof}

La proposition suivante montre que l'on peut ramener théoriquement le
calcul des groupes d'extension dans $\F_{\Gr,I}$ entre des foncteurs {\em
  finis} à un nombre fini de calculs de groupes d'extension entre
foncteurs finis de~$\F\otimes\F^I_{surj}$.

\begin{pr}\label{rqimp} La résolution canonique  d'un foncteur polynomial $X$
  de $\F_{\Gr,I}$ est finie ; sa longueur est majorée par $\deg X+1$ si $X\neq 0$. 
\end{pr}

\begin{proof} Grâce à la proposition~\ref{pr-comd}, il suffit de voir
  que l'endo\-fonc\-teur $\Delta^I_{surj}$ de $\F\otimes\F^I_{surj}$
  diminue strictement le degré des foncteurs polynomiaux non
  nuls. Ce fait est une traduction du lemme~\ref{lmcpd12}.
\end{proof}

Nous appliquons maintenant la proposition~\ref{blpr} à des
considérations relatives aux foncteurs dérivés gauches du foncteur $\eta_I$, qui
  est exact à droite par la proposition~\ref{adjeta}.

\begin{nota}\label{fcth} \'Etant donné $n\in\mathbb{Z}$, nous noterons
  $h_n^{\Gr,I} : \F_{\Gr,I}\to\F\otimes\F^I_{surj}$ le $n$-ième foncteur dérivé gauche de $\eta_I$.
\end{nota}

\begin{rem} Les foncteurs dérivés de $\eta_I$ mesurent, intuitivement,
le défaut d'essentielle surjectivité du foncteur $\theta_I$, donc la
différence homologique entre les catégories $\F\otimes\F^I_{surj}$ et
$\F_{\Gr,I}$. Cela motive la notation employée pour ces foncteurs.

En effet, si l'on se restreint aux foncteurs finis, nous avons vu (proposition~\ref{prfffgr}) que
tout objet de $\F_{\Gr,I}$ s'obtient par extensions successives de
foncteurs appartenant à l'image du foncteur $\theta_I$. La description
de ces objets équivaut donc essentiellement au calcul de groupes
d'extensions entre deux objets de l'image de $\theta_I$. Ces calculs
font naturellement intervenir les foncteurs $h_n^{\Gr,I}$ : il existe
une suite spectrale du premier quadrant
$$E^2_{p,q}={\rm
    Ext}^p_{\F\otimes\F^I_{surj}}(h^{\Gr,I}_q(X),F)\Rightarrow {\rm Ext}^*_{\Gr,I} (X,\theta_I(F))$$
naturelle en les objets $F$ de $\F\otimes\F^I_{surj}$ et $X$ de $\F_{\Gr,I}$.
\end{rem}

La proposition et le corollaire suivants établissent le lien entre le
foncteur gradué $h_*^{\Gr,I}$ et la résolution canonique.

\begin{pr}\label{iotacycl} Pour  tout entier $n>0$ et tout objet $F$ de $\F\otimes\F^I_{surj}$, on a $h_n^{\Gr,I}(\xi_I(F))=0$.
\end{pr}

\begin{proof} Soit $P_\bullet$ une résolution projective de $F$. Comme
  le foncteur $\xi_I$ est exact et préserve les projectifs (son
  adjoint à droite $\sigma_I$ est exact), $\xi(P_\bullet)$ est une
  résolution projective de $\xi_I(F)$. La proposition~\ref{propeta}
  montre que $\eta_I\xi_I(P_\bullet)$ s'identifie à $P_\bullet$. Ce
  complexe, dont l'homologie est (isomorphe à)
  $h_*^{\Gr,I}(\xi_I(F))$, est donc acyclique en degrés strictement positifs.
\end{proof}

\begin{cor}\label{rcfd} Le foncteur gradué $h_*^{\Gr,I}$ est
  canoniquement isomorphe à l'homologie du complexe $\eta_I\big(\Re^{\Gr,I}_\bullet\big)$ de foncteurs $\F_{\Gr,I}\to\F\otimes\F^I_{surj}$.
\end{cor}

Nous illustrons à présent la proposition~\ref{iotacycl} par un calcul élémentaire.

\begin{ex}[Calcul sur les injectifs standard --- cas $I=\mathbb{N}$]
  La proposition~\ref{prfif} permet un calcul rapide des objets
  $h^\Gr_*(I^\Gr_{(V,W)})$. En effet, par les propositions~\ref{propeta} et~\ref{iotacycl}, on a $h_0^\Gr(\iota(I_V))\simeq I_V$, où l'on plonge
  $\F$ dans $\F\otimes\F_{surj}$ par le foncteur $\cdot\boxtimes\kk$,
  et  $h_n^\Gr(\iota(I_V))=0$ si $n\neq 0$. Comme
  $\eta(I_{(V,0)})\simeq\eta(\kappa(I_V))\simeq I_V$ (par les
  propositions~\ref{proj-cfg} et~\ref{propeta}), on a finalement : 
\begin{enumerate}\item l'objet gradué $h_*^\Gr(I^\Gr_{(V,W)})$ est nul
  si $\dim W>0$ ;
\item l'objet gradué $h_*^\Gr(I^\Gr_{(V,0)})$ est concentré en degré
  $0$, où il est naturellement isomorphe à $I_V$.
\end{enumerate}
\end{ex}

La proposition~\ref{pretapt} fournit la formule de Künneth suivante.

\begin{pr}\label{kunneth} Il existe un isomorphisme d'objets gradués
  de $\F\otimes\F^I_{surj}$
$$h^{\Gr,I}_*(X\otimes Y)\simeq h^{\Gr,I}_*(X)\otimes h^{\Gr,I}_*(Y)$$
naturel en les objets $X$ et $Y$ de
$\F_{\Gr,I}$.
\end{pr}

\begin{proof} Comme le foncteur $\xi_I$ commute au produit tensoriel,
  la proposition~\ref{iotacycl} montre que le complexe total du produit tensoriel des
  résolutions canoniques de $X$ et $Y$ est une résolution
  $\eta_I$-acyclique de $X\otimes Y$.  L'homologie de ce complexe
  étant naturellement isomorphe au produit tensoriel de
  $h^{\Gr,I}_*(X)$ et $h^{\Gr,I}_*(Y)$, 
on en déduit la proposition.
\end{proof}

\section{Les catégories $\F_{\Pl,n}$}\label{cfpl}

Le diagramme de recollement de la proposition~\ref{recf10} permet de
dévisser la catégorie de foncteurs en grassmanniennes globale $\F_\Gr(\kk)$
à l'aide des catégories $\F_{\Gr,n}(\kk)$. Il est patent que l'on peut
pousser plus loin ce dévissage, en raison de l'intervention du groupe
linéaire $GL_n(\kk)$ dans un grand nombre de considérations relatives
à cette catégorie. Nous étudions dans ce qui suit des catégories
réalisant un tel dévissage, dans un sens qui sera précisé et illustré à la fin du paragraphe~\ref{fpli}.

\begin{conv} Dans toute cette section, $n$ désigne un entier naturel. 

Nous noterons encore $E_n$, pour abréger, l'objet initial $(E_n,id)$
  de $\E^f_{\Pl,n}(\kk)$.
\end{conv}

\begin{defi} La catégorie de foncteurs en grassmanniennes $\F_{\Pl,n}(\kk)$
  est la catégorie définie par
$$\F_{\Pl,n}(\kk)=\mathbf{Fct}(\E^f_{\Pl,n}(\kk),\E^f_\kk).$$
\end{defi}

\subsection{Généralités}\label{subsfpl}

Nous introduisons maintenant des foncteurs analogues à ceux du §\,\ref{subsfgr}
dans le cas $\F_{\Pl,n}$ et en donnons les propriétés
élémentaires. Les démonstrations, semblables à celles dudit
paragraphe, sont laissées au lecteur.

\begin{nota} Nous abrégerons
  respectivement en ${\rm
    hom}_{\Pl,n}$, $P^{\Pl,n}_X$, $I^{\Pl,n}_X$ les expressions ${\rm hom}_{\F_{\Pl,n}}$, $P^{\E^f_{\Pl,n}}_X$ et
  $I^{\E^f_{\Pl,n}}_X$.
\end{nota}

\begin{defi}\label{dfffpl} 
\begin{enumerate}\item Le foncteur de {\em restriction sans plongement} $\gamma_n :
  \F_{\Gr,n}\to\F_{\Pl,n}$ est défini comme étant le foncteur de
  précomposition par le foncteur d'oubli du plongement $inc^\Pl_n : \E^f_{\Pl,n}\to\E^f_{\Gr,n}$.
\item Le foncteur de {\em plongement
    standard} $\iota^\Pl_n : \F\to\F_{\Pl,n}$ est le foncteur de
  précomposition par le foncteur d'oubli principal
  $\bar{\mathfrak{O}}_n : \E^f_{\Pl,n}\to\E^f$.
\item Le foncteur de {\em plongement
    réduit} $\kappa_n^\Pl : \F\to\F_{\Pl,n}$ est le foncteur de
  précomposition par le foncteur de réduction $\bar{\mathfrak{K}}_n : \E^f_{\Pl,n}\to\E^f$. 
\item  On définit le foncteur de {\em décalage en
    grassmanniennes} $\sigma_n^\Pl :
  \F_{\Pl,n}\to\F$ comme le foncteur
  de précomposition par le foncteur de décalage pointé $\mathfrak{S}_n : \E^f\to\E^f_{\Pl,n}$.
\end{enumerate}

Explicitement, on a
$$\gamma_n(X)(V,u)=X(V,
im\,u),$$
$$\iota^\Pl_n(F)(V,u)=F(V),$$
$$\kappa_n^\Pl(F)(V,u)=F(coker\,u),$$
$$\sigma_n^\Pl(A)(W)=A(W\oplus
E_n,E_n\hookrightarrow W\oplus E_n)$$
pour $X\in {\rm Ob}\,\F_{\Gr,n}$, $F\in {\rm Ob}\,\F$, $A\in {\rm Ob}\,\F_{\Pl,n}$, $(V,u : E_n\hookrightarrow
V)\in {\rm Ob}\,\E^f_{\Pl,n}$ et $W\in {\rm Ob}\,\E^f$.
\end{defi}

\begin{rem}\label{rqcofpl}\begin{enumerate}\item On a
    $\iota^\Pl_n=\gamma_n\iota_n$ et $\kappa_n^\Pl=\gamma_n\kappa_n$.
\item Nous ne donnons pas de notation pour les foncteurs d'intégrale
$\F_{\Pl,n}\to\F$, qui ne revêtent pas la même importance que les foncteurs
d'intégrale en grassmanniennes des catégories de type $\F_{\Gr,I}$.
\end{enumerate}
\end{rem}

La proposition suivante donne les principales variantes des propriétés
du paragraphe~\ref{subsfgr} en termes des foncteurs introduits précédemment.

\begin{pr}\label{prffpl} 
\begin{enumerate}\item Les foncteurs $\gamma_n$, $\iota^\Pl_n$, $\sigma^\Pl_n$ et
  $\kappa^\Pl_n$ sont exacts et fidèles ; ils commutent au produit
  tensoriel, aux limites et aux colimites. 
\item Le foncteur $\kappa^\Pl_n$ est de plus plein, et son image est
  une sous-catégorie de Serre de $\F_{\Pl,n}$.
\item La composition
  $\F\xrightarrow{\kappa^\Pl_n}\F_{\Pl,n}\xrightarrow{\sigma^\Pl_n}\F$
  est canoniquement isomorphe au foncteur identique.
\item\label{lndlad} Le foncteur $\iota_n^\Pl$ est adjoint à gauche à $\sigma^\Pl_n$.
\item On a des isomorphismes
  $P^{\Pl,n}_{\mathfrak{S}_n(V)}\simeq\iota^\Pl_n(P_V)$ et
  $I^{\Pl,n}_{\mathfrak{S}_n(V)}\simeq\kappa^\Pl_n(I_V)\otimes I^{\Pl,n}_{E_n}$ naturels en
  l'objet $V$ de $\E^f$.
\item Le foncteur composé $_{\kk[GL_n(\kk)]}\mathbf{Mod}\xrightarrow{\rho_n}\F_{\Gr,n}\xrightarrow{\gamma_n}\F_{\Pl,n}$ est canoniquement isomorphe à $_{\kk[GL_n(\kk)]}\mathbf{Mod}\to\E\to\F_{\Pl,n}$, où la première flèche est le foncteur d'oubli de l'action de $GL_n(\kk)$ et le second le plongement canonique (donné par les foncteurs constants).
\end{enumerate}
\end{pr}

\paragraph*{Décomposition scalaire, tors de Frobenius et changement de
  corps} Ces notions se définissent comme dans $\F_{\Gr,n}$ et
possèdent des propriétés tout à fait analogues.

\paragraph*{Les objets polynomiaux et finis de la catégorie
  $\F_{\Pl,n}$} L'endofoncteur de $\F_{\Pl,n}$ donné par la
précomposition par le foncteur de translation $V\boxplus\cdot :
\E^f_{\Pl,n}\to\E^f_{\Pl,n}$ (où $V\in {\rm Ob}\,\E^f$) s'appelle {\em foncteur de décalage} et
se note $\Delta^{\Pl,n}_V$ ; il existe un scindement canonique
$\Delta^{\Pl,n}_\kk\simeq id\oplus\Delta^{\Pl,n}$, où $\Delta^{\Pl,n}$
est le {\em foncteur différence} de $\F_{\Pl,n}$.

Les {\em foncteurs polynomiaux} de $\F_{\Pl,n}$ sont ses objets
$\Delta^{\Pl,n}$-nilpotents ; on a une notion de degré. On note
$\F^k_{\Pl,n}$ la sous-catégorie épaisse de $\F_{\Pl,n}$ des foncteurs
polynomiaux de degré au plus $k$. La catégorie $\F^0_{\Pl,n}$ est réduite aux foncteurs constants,
comme dans la catégorie~$\F$ (contrairement à ce qui advient
dans~$\F_{\Gr,n}$). Nous nous contentons d'énoncer les analogues les
plus importants des propriétés des §\,\ref{s-fp} et~\ref{s-ffg}, qui
s'adaptent sans difficulté à~$\F_{\Pl,n}$.

Le résultat suivant, similaire à la proposition~\ref{prcqd} identifie les quotients de la filtration
polynomiale de $\F_{\Pl,n}$.

\begin{pr} Le foncteur $\sigma^\Pl_n : \F_{\Pl,n}\to\F$ induit pour
  tout $k\in\mathbb{N}$ une équivalence de catégories
  $\F^k_{\Pl,n}/\F^{k-1}_{\Pl,n}\to\F^k/\F^{k-1}$.

Les foncteurs $\kappa^\Pl_n$ et $\iota^\Pl_n : \F\to\F_{\Pl,n}$
induisent chacun une équivalence de catégories
$\F^k/\F^{k-1}\to\F^k_{\Pl,n}/\F^{k-1}_{\Pl,n}$ inverse de la précédente.
\end{pr}

L'analogue suivant de la proposition~\ref{prf-grof} sous-tend, avec la
proposition précédente, la proposition~\ref{jlaou} ci-dessous.

\begin{pr} Un foncteur de $\F_{\Pl,n}$ est fini si et seulement s'il
  est polynomial et à valeurs de dimension finie.
\end{pr}

Nous en venons maintenant à la description explicite des objets
simples de la catégorie $\F_{\Pl,n}$ (cf. proposition~\ref{prfffgr}).

\begin{pr}\label{jlaou}\begin{enumerate}\item \'Etant donné un objet $X$ de
    $\F_{\Pl,n}$, les assertions suivantes sont équivalentes.
\begin{enumerate}\item L'objet $X$ de $\F_{\Pl,n}$ est simple.
\item L'objet $\sigma_n^\Pl(X)$ de $\F$ est simple.
\item Il existe un objet simple $S$ de $\F$ tel que $X\simeq\kappa^\Pl_n(S)$.
\end{enumerate}
\item Les foncteurs $\kappa^\Pl_n$ et $\sigma^\Pl_n$ induisent des
isomorphismes d'anneaux inverses l'un de l'autre entre $G^f_0(\F)$ et $G^f_0(\F_{\Pl,n})$.
\end{enumerate}
\end{pr}

\paragraph*{Premiers liens entre les catégories $\F$, $\F_{\Gr,n}$ et
  $\F_{\Pl,n}$ en termes de (co)modules} L'un des thèmes de la
section~\ref{cfpl} consiste à compléter l'équivalence de catégories
fondamentale entre $\F_\Gr$ et les $\kk[\Gr]$-comodules
(cf. proposition~\ref{prfig}) par l'identification d'autres catégories
de modules ou de comodules à des catégories de foncteurs. Dans ce
sous-paragraphe, nous donnons deux telles propriétés, qui reposent sur
la section~\ref{sctccf}.

\begin{pr}\label{fgrfpl}  La catégorie
  $\F_{\Pl,n}$ est équivalente à la sous-catégorie
  $\mathbf{Comod}_{P^{\Gr,n}_{(E_n,E_n)}}$ de~$\F_{\Gr,n}$. 
\end{pr}

\begin{proof} Le foncteur
$$(\E^f_{\Gr,n})_{\backslash {\rm hom}\,((E_n,E_n),.)}\to\E^f_{\Pl,n}\quad
\big((V,W),u\big)\mapsto (V,u)$$
est une équivalence de catégories, où la catégorie
source dérive de la notation \ref{notfcom}. La proposition est
donc un cas
particulier de la proposition \ref{prf-com}.
\end{proof}

\begin{rem}\label{rqida} Le foncteur ${\rm hom}_{\Gr,n}((E_n,E_n),.)$
  est canoniquement isomorphe à la composée
$\F_{\Gr,n}\xrightarrow{\mathfrak{B}_n}\E^n_{surj}\xrightarrow{{\rm hom}\,(E_n,.)}\bf{Ens}$. Comme $\E^n_{surj}\simeq\underline{GL_n(\kk)}$ est
  équivalente à sa catégorie opposée, on peut remplacer dans la
  démonstration précédente ce foncteur par un foncteur {\em
    contravariant}, et obtenir ainsi sur $P^{\Gr,n}_{(E_n,E_n)}$ une
  structure d'algèbre telle que  $\F_{\Pl,n}$ est équivalente à la sous-catégorie
  $\mathbf{Mod}_{P^{\Gr,n}_{(E_n,E_n)}}$ de $\F_{\Gr,n}$. Ce phénomène
  est à rapprocher de la remarque suivante : $P^{\Gr,n}_{(E_n,E_n)}$
  s'identifie à $\rho_n(\kk[GL_n(\kk)])$, et $\kk[GL_n(\kk)]$ est un objet
  auto-dual de $_{\kk[GL_n(\kk)]}\bf{Mod}$ (la dualité consistant à
  associer à une représentation linéaire la représentation contragrédiente).
\end{rem}

Grâce au lemme suivant, nous établirons, à la proposition~\ref{prffpl1}, un lien direct entre les
catégories $\F$ et $\F_{\Pl,n}$.

\begin{lm} La catégorie $(\E^f_{\Pl,n})_{/{\rm hom}\,(.,E_n)}$ est
  équivalente à $\E^f$.
\end{lm}

\begin{proof} On vérifie aussitôt que les deux foncteurs introduits
  ci-après sont des équivalences de catégories réciproques l'une de
  l'autre.
\begin{itemize}\item On définit un foncteur
  $\E^f\to(\E^f_{\Pl,n})_{/{\rm hom}\,(.,E_n)}$ en associant à $V\in
  {\rm Ob}\,\E^f$ l'objet $\mathfrak{S}_n(V)$ de $\E^f_{\Pl,n}$ muni
  du morphisme vers $E_n$ donné par la projection $V\oplus
  E_n\twoheadrightarrow E_n$, et à une application linéaire $f : V\to
  W$ le morphisme $f\oplus E_n$.
\item On définit un foncteur $(\E^f_{\Pl,n})_{/{\rm
      hom}\,(.,E_n)}\to\E^f$ en associant à un objet $\big((V,u :
  E_n\hookrightarrow V), r : V\to E_n\big)$ de $(\E^f_{\Pl,n})_{/{\rm
      hom}\,(.,E_n)}$ l'espace vectoriel $coker\,u\simeq ker\,r$, et à
  un morphisme l'application linéaire induite.
\end{itemize}
\end{proof}

\begin{pr}\label{prffpl1} La catégorie $\F$ est équivalente à la
  sous-catégorie $\mathbf{Mod}_{I^{\Pl,n}_{E_n}}$ de $\F_{\Pl,n}$.
\end{pr}

\begin{proof} On combine le lemme précédent avec la proposition~\ref{prf-com3}.
\end{proof}

\subsection{L'équivalence de catégories $\F_{\Pl,n}\simeq\mathbf{Comod}_{I_{E_n}}$}\label{fpli}
L'analogue des considérations de la section~\ref{fgm} dans les
catégories $\F_{\Pl,n}$ permet d'obtenir l'équi\-va\-lence de
catégories éponyme de ce paragraphe. Les structures que l'on en déduit constituent
une justification essentielle à l'étude des catégories
$\F_{\Pl,n}$, dont nous verrons comment elle peut
intervenir dans les catégories de type $\F_{\Gr,I}$.

Le point de départ de ce paragraphe réside dans le fait que le foncteur
$\iota^\Pl_n$ est adjoint à gauche au foncteur $\sigma_n^\Pl$ (proposition~\ref{prffpl}.\,\ref{lndlad}).

\begin{prdef}\label{mon-pl}  La monade de $\F$ associée à l'adjonction entre les foncteurs
  $\iota^\Pl_n$ et $\sigma_n^\Pl$ conformément à la proposition \ref{mon-adj}, que nous désignerons par
  $\T_{\Pl,n}=(\Delta_{E_n},u_{\Pl,n},\mu_{\Pl,n})$, est donnée comme suit.
\begin{itemize}\item Le foncteur $\Delta_{E_n}$ est le foncteur de
  décalage de $\F$ (cf. section~\ref{s-rf}).
\item La transformation naturelle $u_{\Pl,n} : id\to \Delta_{E_n}$ (unité de
l'adjonction)  est induite par l'application linéaire $0\to E_n$,
compte-tenu de l'identification entre $id$ et $\Delta_0$.
\item La
multiplication $\mu_{\Pl,n} :
(\Delta_{E_n})^2\simeq\Delta_{E_n\oplus E_n}\to\Delta_{E_n}$ est la transformation naturelle
induite par la somme $E_n\oplus E_n\to E_n$.
\end{itemize}

En outre, il existe un scindement naturel
$id\xrightarrow{u_{\Pl,n}}\Delta_{E_n}\xrightarrow{p_{\Pl,n}} id$, où
$p_{\Pl,n}$ est induit par l'application linéaire $E_n\to 0$.
\end{prdef}

\begin{pr}\label{fpl-mon} La catégorie $\F_{\Pl,n}$ est équivalente à
  la catégorie des modules sur la monade $\T_{\Pl,n}$ de
  $\F$.
\end{pr}

Les deux propositions précédentes se démontrent comme les propositions \ref{adjidis} et~\ref{fgr-mon}.

\begin{conv}
Dans la suite de ce paragraphe, nous {\em identifierons} la catégorie
$\F_{\Pl,n}$ avec la sous-catégorie des modules sur $\T_{\Pl,n}$ de
$\F$. Autrement dit, un objet de $\F_{\Pl,n}$ sera
désormais un couple $(X,\Delta_{E_n}X\xrightarrow{m_X}X)$, où $X$ est un objet de $\F$ et
$m_X$ un morphisme tel que :
\begin{itemize}\item la composée $X\xrightarrow{(u_{\Pl,n})_X}
\Delta_{E_n}X\xrightarrow{m_X}X$ est le morphisme identique ;
\item les composées $(\Delta_{E_n})^2 X\xrightarrow{(\mu_{\Pl,n})_X}
\Delta_{E_n}X\xrightarrow{m_X}X$ et  $(\Delta_{E_n})^2 X\xrightarrow{\Delta_{E_n}m_X}\Delta_{E_n}
X\xrightarrow{m_X}X$ coïncident.
\end{itemize}

Par abus, nous noterons parfois simplement $X$ pour
$(X,m_X)$.
\end{conv}

\begin{nota}
Le morphisme $X\to X\otimes I_{E_n}$ adjoint à $m_X$ (cf. proposition~\ref{adj-fctd}) sera noté $\psi_X$ dans la suite
de ce paragraphe.
\end{nota}

Nous indiquons, dans la proposition suivante, l'analogue du foncteur
$\eta_I$. Nous omettons la démonstration, similaire à celle des propositions \ref{adjeta}, \ref{propeta} et~\ref{pretapt}.

\begin{prdef}[Foncteur
  $\eta_n^\Pl$]\label{etapl}\begin{enumerate}\item On définit le foncteur
    $\eta_n^\Pl : \F_{\Pl,n}\to\F$ comme le coégalisateur des deux
    transformations naturelles $p_{\Pl,n}$ (cf. notation de la proposition~\ref{mon-pl}) et $m :
    \Delta_{E_n}\sigma^\Pl_n\twoheadrightarrow\sigma^\Pl_n$. 
\item Le foncteur $\eta_n^\Pl$  est adjoint à gauche au foncteur $\kappa^\Pl_n$.
\item Il existe des isomorphismes
$$\eta_n^\Pl\circ\kappa_n^\Pl\simeq\eta_n^\Pl\circ\iota_n^\Pl\simeq id.$$
\item Le foncteur $\eta^\Pl_n$ commute au produit tensoriel à
  isomorphisme naturel près.
\end{enumerate}  
\end{prdef}

On peut introduire, comme dans la section~\ref{fgm}, une notion de
résolution canonique dans $\F_{\Pl,n}$, grâce à laquelle on peut
calculer les foncteurs dérivés du foncteur $\eta^\Pl_n$, et relier les
groupes d'extension dans $\F_{\Pl,n}$ à ceux de $\F$ par une suite spectrale. Plutôt que de détailler
ces considérations, nous abordons des constructions plus spécifiques à la
catégorie $\F_{\Pl,n}$.

\begin{rem}\label{comod-i} Pour tout espace vectoriel $V$ de dimension
  finie, le foncteur $P_V$ est muni d'une structure naturelle
  d'algèbre commutative, dont la multiplication est le morphisme
  $P_V\otimes P_V\simeq P_{V\oplus V}\to P_V$ induit par la diagonale
  $V\to V\oplus V$ et l'unité $P_V\to\kk\simeq P_0$ par $0\to
  V$. Explicitement, la structure d'algèbre sur les espaces vectoriels
  $P_V(W)=\kk[{\rm hom}\,(V,W)]$ qui s'en déduit est celle de
  l'algèbre du groupe abélien ${\rm hom}\,(V,W)$.

En utilisant le foncteur de dualité $D$, on en déduit une
structure naturelle de coalgèbre cocommutative sur les injectifs
standard $I_V$ de $\F$.
\end{rem}

La proposition fondamentale suivante repose sur l'idenfication de
l'adjoint à droite au foncteur $\Delta_{E_n}$ de la monade
$\T_{\Pl,n}$. Elle n'a pas d'analogue dans les catégories
$\F_{\Gr,I}$, car l'adjoint à droite au foncteur $\Delta^I_{surj}$ (qui existe
par le corollaire~\ref{cr-adjp}) n'est
généralement pas un foncteur qui admet une description simple.

\begin{pr}\label{fplcom2} Il existe une équivalence de
    catégories
    $\F_{\Pl,n}\xrightarrow{\simeq}\mathbf{Comod}_{I_{E_n}}$ qui factorise le foncteur $\sigma^\Pl_n :
    \F_{\Pl,n}\to\F$ à travers le foncteur d'oubli
    $\mathbf{Comod}_{I_{E_n}}\to\F$.
$$\xymatrix{\F_{\Pl,n}\ar[r]^-{\sigma^\Pl_n}\ar@{.>}[dr]_\simeq & \F \\
& \mathbf{Comod}_{I_{E_n}}\ar[u]
}$$
\end{pr}

\begin{proof} La proposition découle de ce qu'une flèche
  $\Delta_{E_n}X\to X$ de $\F$ fait de $X$ un module sur $\T_{\Pl,n}$
  si et seulement si la flèche adjointe $X\to X\otimes I_{E_n}$
  définit une structure de $I_{E_n}$-comodule sur $X$. En effet, la
  multiplication de $\T_{\Pl,n}$ comme la comultiplication de
  $I_{E_n}$ sont induites par la somme $E_n\oplus E_n\to E_n$, et
  l'unité de $\T_{\Pl,n}$ comme la coünité de $I_{E_n}$ proviennent du
  morphisme $0\to E_n$.
\end{proof}

\begin{rem}\label{compcom} La proposition \ref{fplcom2} est à comparer
  avec l'équivalence de catégories
  $\F_{\Pl,n}\xrightarrow{\simeq}\mathbf{Comod}^{fid}_{\widetilde{P}(n)}$ (fournie par un foncteur d'intégrale) déduite de la proposition~\ref{prf-com2}, où $\wt{P}(n)$ est le foncteur introduit dans la notation~\ref{not-tilff}.
%
\end{rem}

\paragraph*{Produit cotensoriel et foncteur $\chi^\Pl_n$} Dans la suite de ce paragraphe, nous identifierons les catégories $\F_{\Pl,n}$ et $\mathbf{Comod}_{I_{E_n}}$. En
  particulier, la coalgèbre $I_{E_n}$ étant cocommutative, on dispose
  dans $\F_{\Pl,n}$ d'un produit cotensoriel
  $\underset{I_{E_n}}{\square}$, qui sera noté simplement $\square$ par la suite ; c'est un bifoncteur exact à gauche.

Avant d'utiliser ce produit cotensoriel, nous avons besoin
de décrire quelques foncteurs usuels de source ou de but
$\F_{\Pl,n}$ dans l'identification de $\F_{\Pl,n}$ à $\mathbf{Comod}_{I_{E_n}}$.

Le foncteur $\sigma^\Pl_n : \F_{\Pl,n}\to\F$ correspond au foncteur
d'oubli (cf. proposition~\ref{fplcom2}).

Nous identifions, dans la proposition qui suit, le produit tensoriel
de $\F_{\Pl,n}$ en termes de $I_{E_n}$-comodules. 

\begin{pr}\label{rqtenspl} Le coproduit $\psi_{X\otimes Y}$ d'un
  produit tensoriel de deux objets $X$ et $Y$ de $\F_{\Pl,n}$ est égal
  à la composée
$$X\otimes Y\xrightarrow{\psi_X\otimes\psi_Y}(X\otimes I_{E_n})\otimes
(Y\otimes I_{E_n})\simeq (X\otimes Y)\otimes (I_{E_n}\otimes
I_{E_n})\xrightarrow{(X\otimes Y)\otimes a_n} (X\otimes Y)\otimes
I_{E_n}$$
où $a_n : I_{E_n}\otimes I_{E_n}\to I_{E_n}$ est le produit du
foncteur en algèbres de Boole $I_{E_n}$.
\end{pr}

 Ainsi, le produit tensoriel
de la catégorie de comodules $\F_{\Pl,n}$ provient de la structure
d'algèbre sur $I_{E_n}$. Le fait que le produit tensoriel de deux
$I_{E_n}$-comodules est naturellement un $I_{E_n}$-comodule provient
de la compatibilité des structures d'algèbre et de coalgèbre sur
$I_{E_n}$, qui est une {\em algèbre de Hopf} de la catégorie
symétrique monoïdale~$\F$.

\begin{proof} La propriété provient, par adjonction, de ce que le morphisme $m_{X\otimes Y} :
  \Delta_{E_n}\sigma^\Pl_n(X\otimes Y)\to\sigma^\Pl_n(X\otimes Y)$
  s'identifie, modulo les isomorphismes canoniques
  $\sigma^\Pl_n(X\otimes
  Y)\simeq\sigma^\Pl_n(X)\otimes\sigma^\Pl_n(Y)$ et $\Delta_{E_n}\sigma^\Pl_n(X\otimes
  Y)\simeq\Delta_{E_n}\sigma^\Pl_n(X)\otimes\Delta_{E_n}\sigma^\Pl_n(Y)$, à $m_X\otimes m_Y : \Delta_{E_n}\sigma^\Pl_n(X)\otimes\Delta_{E_n}\sigma^\Pl_n(Y)\to\sigma^\Pl_n(X)\otimes\sigma^\Pl_n(Y)$.
\end{proof}

\begin{nota} Dans ce paragraphe, nous noterons $j_F : F\hookrightarrow
  F\otimes I_{E_n}$, pour $F\in {\rm Ob}\,\F$, l'inclusion canonique
  déduite de $\kk\hookrightarrow I_{E_n}$.
\end{nota}

Nous identifions maintenant le foncteur $\kappa^\Pl_n$ en termes de $I_{E_n}$-comodules.

\begin{pr}\label{idtk-cm} Le foncteur $\kappa^\Pl_n$ associe à un
  objet $F$ de $\F$ le $I_{E_n}$-comodule $(F,j_F)$.

Ainsi, le foncteur $\kappa^\Pl_n$ identifie $\F$ à la
sous-catégorie de Serre des $I_{E_n}$-comodules $X$ tels que
  $\psi_X=j_X$.
\end{pr}

\begin{proof} De manière semblable à la proposition~\ref{identk}, on
  voit que $\kappa^\Pl_n(F)$ est le $\T_{\Pl,n}$-module
  $(F,(p_{\Pl,n})_F)$. La proposition se déduit alors de ce que le
  morphisme $(p_{\Pl,n})_F : \Delta_{E_n}F\to F$ est adjoint à $j_F$.
\end{proof}

\begin{cor}\label{cr-kcp} Soient $F$ un objet de $\F$ et $X$ un objet
  de $\F_{\Pl,n}$. Les morphismes $\psi_{\kappa^\Pl_n(F)\otimes X}$ et $F\otimes\psi_X :
  F\otimes\sigma^\Pl_n(X)\to F\otimes\sigma^\Pl_n(X)\otimes I_{E_n}$
  de $\F$ sont égaux.
\end{cor}

\begin{proof} C'est une conséquence directe des
  propositions~\ref{idtk-cm} et~\ref{rqtenspl}.
\end{proof}

\begin{rem}\label{rqptct} Il existe un monomorphisme canonique
  $\sigma^\Pl_n(X\square Y)\hookrightarrow\sigma^\Pl_n(X\otimes
  Y)\simeq\sigma^\Pl_n(X)\otimes \sigma^\Pl_n(Y)$ dans $\F$ ; il ne provient
  {\em pas} d'un morphisme naturel $X\square Y\to X\otimes Y$ de $\F_{\Pl,n}$. En
  revanche, le corollaire~\ref{cr-kcp} montre qu'il est induit par un monomorphisme naturel
  $X\square Y\hookrightarrow\kappa^\Pl_n\sigma^\Pl_n(X)\otimes Y$.
\end{rem}

Nous introduisons à présent un nouveau foncteur déduit de la proposition~\ref{fplcom2}.

\begin{defi}[Foncteur $\chi^\Pl_n$] On définit le foncteur $\chi^\Pl_n
  : \F_{\Pl,n}\to\F$ par la composition suivante.
$$\chi^\Pl_n : \F_{\Pl,n}\xrightarrow{\cdot\,\square\kk}\F_{\Pl,n}\xrightarrow{\sigma^\Pl_n}\F$$
\end{defi}

Ainsi, $\chi^\Pl_n$ est l'égalisateur des transformations naturelles
$\psi$ et $j_{\sigma^\Pl_n} : \sigma^\Pl_n\to\sigma^\Pl_n\otimes I_{E_n}$.

La proposition suivante constitue le  résultat principal de ce sous-paragraphe.

\begin{pr}\label{chi-pl} Le foncteur $\chi^\Pl_n$ est adjoint à droite
  à $\kappa^\Pl_n : \F\to\F_{\Pl,n}$.
\end{pr}

\begin{proof} Soient $F$ un objet de $\F$ et $X$ un
  $I_{E_n}$-comodule. Par la proposition~\ref{idtk-cm}, ${\rm
    hom}_{\Pl,n}(\kappa^\Pl_n(F),X)$ s'identifie à l'ensemble des
  morphismes $f : F\to\sigma^\Pl_n(X)$ de $\F$ tels que le diagramme
  suivant commute.
$$\xymatrix{F\ar[d]_{j_F}\ar[r]^-f & \sigma^\Pl_n(X)\ar[d]^{\psi_X}\\
F\otimes I_{E_n}\ar[r]^-{f\otimes I_{E_n}} & \sigma^\Pl_n(X)\otimes I_{E_n}
}$$

Comme le diagramme
$$\xymatrix{F\ar[d]_{j_F}\ar[r]^-f & \sigma^\Pl_n(X)\ar[d]^{j_{\sigma^\Pl_n(X)}}\\
F\otimes I_{E_n}\ar[r]^-{f\otimes I_{E_n}} & \sigma^\Pl_n(X)\otimes I_{E_n}
}$$
commute, la condition précédente revient à dire que le morphisme $f$
est à valeurs dans l'égalisateur des morphismes $j_{\sigma^\Pl_n(X)}$
et $\psi_X$, qui est $\chi_n^\Pl(X)$. Cela démontre la proposition.
\end{proof}

La proposition ci-après donne une propriété de compatibilité entre les
foncteurs $\chi^\Pl_n$ et $\kappa^\Pl_n$.

\begin{pr}\label{chikap} Il existe des isomorphismes naturels
  $\chi^\Pl_n(X\otimes\kappa^\Pl_n(F))\simeq\chi^\Pl_n(X)\otimes F$ et
$X\square\kappa^\Pl_n(F)\simeq\kappa^\Pl_n(\chi^\Pl_n(X)\otimes F)$
pour $X\in {\rm Ob}\,\F_{\Pl,n}$ et $F\in {\rm Ob}\,\F$.
\end{pr}

\begin{proof} Le foncteur $\chi^\Pl_n(X\otimes\kappa^\Pl_n(F))$ est
  l'égalisateur des flèches
  $\psi_{X\otimes\kappa^\Pl_n(F)}=\psi_X\otimes F$ (par le corollaire~\ref{cr-kcp})
  et $j_{\sigma^\Pl_n(X\otimes\kappa^\Pl_n(F))}=j_{\sigma^\Pl_n(X)}\otimes F$, il s'identifie donc canoniquement au
  produit tensoriel de $\chi^\Pl_n(X)$ et $F$.

Par ailleurs, $X\square\kappa^\Pl_n(F)$ est l'égalisateur des flèches $X\otimes F\xrightarrow{\psi_X\otimes
F}X\otimes I_{E_n}\otimes F$ et $X\otimes
F\xrightarrow{X\otimes
  \psi_{\kappa^\Pl_n(F)}}X\otimes F\otimes I_{E_n}\simeq
X\otimes I_{E_n}\otimes F$ ; comme $X\otimes
  \psi_{\kappa^\Pl_n(F)}=X\otimes j_{F}\simeq j_{\sigma^\Pl_n(X)}\otimes
F$ (modulo l'isomorphisme d'échange des facteurs du produit tensoriel), ce foncteur s'identifie (dans~$\F$) à
$\chi^\Pl_n(X)\otimes F$. Il reste à voir que sa structure de
$I_{E_n}$-comodule est triviale, ce qui provient de l'inclusion
$X\square
\kappa^\Pl_n(F)\hookrightarrow\kappa^\Pl_n\sigma^\Pl_n(X)\otimes\kappa^\Pl_n(F)\simeq\kappa^\Pl_n(\sigma^\Pl_n(X)\otimes
F)$ (cf. remarque~\ref{rqptct}).
\end{proof}

\begin{rem} On peut expliciter l'équivalence de catégories de la
  proposition~\ref{prffpl1} à l'aide des foncteurs $\kappa^\Pl_n$ et
  $\chi^\Pl_n$. Ainsi, on vérifie que les foncteurs
$$\F\xrightarrow{\kappa^\Pl_n}\F_{\Pl,n}\xrightarrow{I^{\Pl,n}_{E_n}\otimes\cdot}\mathbf{Mod}_{I^{\Pl,n}_{E_n}}$$
et
$$\mathbf{Mod}_{I^{\Pl,n}_{E_n}}\xrightarrow{oubli}\F_{\Pl,n}\xrightarrow{\chi^\Pl_n}\F$$
sont des équivalences réciproques l'une de l'autre. 
\end{rem}

\paragraph*{Le foncteur $\chi_n : \F_{\Gr,n}\to\F_{GL_n(\kk)}$}

Nous revenons maintenant à la catégorie $\F_{\Gr,n}$, en montrant comment les résultats précédents peuvent être utilisés pour en étudier
certaines propriétés.
Nous commençons par préciser le lien entre $\F_{\Pl,n}$ et
$\F_{\Gr,n}$.

\begin{prdef}\label{liengrpl}\begin{enumerate}\item Le groupe linéaire $GL_n(\kk)$ agit à
    droite sur la
    catégorie $\E^f_{\Pl,n}$ : pour tout
    $g\in GL_n(\kk)$, on définit un foncteur $\tau_g :
    \E^f_{\Pl,n}\to\E^f_{\Pl,n}$ par $(V,u)\mapsto (V,u\circ
    g)$ sur les objets et par l'égalité ${\rm
      hom}_{\Pl,n}((V,u),(V',u'))={\rm hom}_{\Pl,n}((V,u\circ
    g),(V',u'\circ g))$ sur les morphismes, et l'on a $\tau_1=id$ et 
    $\tau_g\circ\tau_h=\tau_{hg}$ pour tous $g,h\in GL_n(\kk)$. 
\item Par précomposition, on en déduit une action à gauche de $GL_n(\kk)$
  sur $\F_{\Pl,n}$ : les foncteurs $\tau_g^* :
  \F_{\Pl,n}\to\F_{\Pl,n}$ vérifient $\tau^*_1=id$ et $\tau_g^*\circ\tau_h^*=\tau_{gh}^*$.
\item On appelle {\em $GL_n(\kk)$-module} dans $\F_{\Pl,n}$ tout objet $X$ de
  $\F_{\Pl,n}$ muni de flèches $\tau^*_g(X)\xrightarrow{t_g}X$ (dites
  d'action de $GL_n(\kk)$) telles
  que le diagramme
$$\xymatrix{\tau_g^*\tau_h^*(X)\ar@{=}[d]\ar[r]^-{\tau^*_g(t_h)} &
  \tau_g^* (X)\ar[d]^{t_g} \\
\tau_{gh}^* (X)\ar[r]^{t_{gh}} & X
}$$
commute pour tous $g,h\in GL_n(\kk)$. Un morphisme de $GL_n(\kk)$-modules de
$\F_{\Pl,n}$ est un morphisme de $\F_{\Pl,n}$ compatible aux
morphismes d'action de $GL_n(\kk)$. On définit ainsi la sous-catégorie  des
  $GL_n(\kk)$-modules de $\F_{\Pl,n}$.
\item Le foncteur $\gamma_n : \F_{\Gr,n}\to\F_{\Pl,n}$ induit une
  équivalence de catégories entre $\F_{\Gr,n}$ et la sous-catégorie des
  $GL_n(\kk)$-modules de $\F_{\Pl,n}$.
\item Le foncteur $\gamma_n$ admet un adjoint à gauche $\Phi_n :
  \F_{\Pl,n}\to\F_{\Gr,n}$ donné sur les objets par $$\Phi_n(
  X)(V,W)=\bigoplus_{u\in {\rm Iso}_\E (E_n,W)}
    X(V,E_n\xrightarrow{u}W\hookrightarrow V).$$
\end{enumerate}
\end{prdef}

\begin{proof} Les deux premières assertions sont claires. Le dernier point  est un cas particulier de la
  proposition \ref{adj-intou} (cf. démonstration de la proposition
  \ref{fgrfpl}). Comme le foncteur $\gamma_n$ est exact et fidèle, la
  proposition \ref{dualmonadique} montre que $\F_{\Gr,n}$ est
  équivalente à la sous-catégorie de $\F_{\Pl,n}$ des modules sur la
  monade associée à cette adjonction. Cette monade se décrit
  comme suit.
\begin{itemize}\item L'endofoncteur $\gamma_n\Phi_n$ de $\F_{\Pl,n}$
  est $\underset{g\in GL_n(\kk)}{\bigoplus}\tau^*_g$.
\item L'unité est l'inclusion $id\hookrightarrow\underset{g\in
    GL_n(\kk)}{\bigoplus}\tau^*_g$ du facteur direct correspondant à $g=1$.
\item La multiplication est la transformation naturelle
$$\bigoplus_{h,h'\in GL_n(\kk)}\tau^*_h\tau^*_{h'}\to\bigoplus_{g\in
  GL_n(\kk)}\tau^*_g$$
dont la composante $\tau^*_h\tau^*_{h'}\to\tau^*_g$ est l'identité si
$hh'=g$, $0$ sinon.
\end{itemize}

La proposition en résulte.
\end{proof}

\begin{conv} Dans la suite de ce paragraphe, on {\em identifie}
  $\F_{\Gr,n}$ avec la catégorie des $GL_n(\kk)$-modules de $\F_{\Pl,n}$.
\end{conv}

\begin{rem}\label{rqgln} Si $\A$ est une catégorie abélienne, on peut voir la catégorie
  $\A_{GL_n(\kk)}$ comme une catégorie de $GL_n(\kk)$-modules dans $\A$. La catégorie $\F_{\Gr,n}$ est quant à elle une catégorie de $GL_n(\kk)$-modules \go tordus\gf par l'action du groupe linéaire sur $\F_{\Pl,n}$.
\end{rem}

Pour mener des raisonnements sur la catégorie $\F_{\Gr,n}$ des
$GL_n(\kk)$-modules tordus de $\F_{\Pl,n}$ analogues à ceux relatifs
aux $GL_n(\kk)$-modules ordinaires, nous sommes conduits à introduire
la notion suivante.

\begin{defi} On appelle {\em trivialisation} sur $GL_n(\kk)$ d'un foncteur $\alpha$ de source $\F_{\Pl,n}$ la donnée d'isomorphismes
$T_{g} : \alpha\circ\tau^*_{g}\to\alpha$ pour tout $g\in GL_n(\kk)$ tels que le diagramme
$$\xymatrix{\alpha\tau^*_{g}\tau^*_{h}\ar@{=}[d]\ar[rr]^-{(T_{g})_{\tau^*_{h}}}
  & & \alpha\tau^*_{h}\ar[d]^{T_{h}} \\
\alpha\tau^*_{gh}\ar[rr]^-{T_{gh}} & & \alpha
}$$
commute pour tous $g,h\in GL_n(\kk)$.
\end{defi}

\begin{prdef}\label{prdf-ch6pl} Soient $\A$ une catégorie abélienne et $\alpha : \F_{\Pl,n}\to\A$ un foncteur muni d'une trivialisation sur $GL_n(\kk)$. On définit un
foncteur $\alpha^{\Gr}_{GL_{n}} : \F_{\Gr,n}\to\A_{GL_n(\kk)}$ comme suit.
\begin{itemize}\item Sur les objets : si $X$ est un $GL_n(\kk)$-module de $\F_{\Pl,n}$, alors $\alpha^{\Gr}_{GL_{n}}(X)$ est l'objet $\alpha(X)$ de $\A$
muni de l'action de $GL_n(\kk)$ donnée par les flèches $\alpha(X)\simeq\alpha(\tau^*_{g}(X))\to\alpha(X)$ (pour $g\in GL_n(\kk)$), où la première flèche est la trivisalisation
et la seconde est induite par la structure de $GL_n(\kk)$-module de $X$.
\item La flèche qu'induit un morphisme de $GL_n(\kk)$-modules de $\F_{\Pl,n}$ via $\alpha$ est un morphisme de $\A_{GL_n(\kk)}$, ce qui permet de déduire la
fonctorialité de $\alpha^{\Gr}_{GL_{n}}$ de celle de $\alpha$.
\end{itemize}

De plus, le diagramme suivant commute (à isomorphisme canonique près)
$$\xymatrix{\F_{\Gr,n}\ar[d]_{\gamma_{n}}\ar[r]^-{\alpha^{\Gr}_{GL_{n}}} & \A_{GL_n(\kk)}\ar[d]^{O^\A_{GL_{n}}} \\
\F_{\Pl,n}\ar[r]_-{\alpha} & \A
}$$
où $O^\A_{GL_{n}}$ désigne le
foncteur d'oubli, conformément à la notation~\ref{notfcom}.
\end{prdef}

Nous illustrons cette construction élémentaire sur le foncteur
$\chi^\Pl_n : \F_{\Pl,n}\to\F$, qui aboutit à une description
explicite de l'adjoint à droite au foncteur $\theta_n : \F_{GL_n(\kk)}\to\F_{\Gr,n}$. 

\begin{prdef}\label{lmid-ch6}\begin{enumerate}\item Le foncteur
    $\sigma^{\Pl}_{n} : \F_{\Pl,n}\to\F$ admet une trivialisation sur $GL_n(\kk)$ pour laquelle le foncteur $(\sigma^{\Pl}_{n})^{\Gr}_{GL_{n}} : \F_{\Gr,n}\to\F_{GL_n(\kk)}$ 
s'identifie à $\sigma_{n}$.
\item La trivialisation précédente induit, via le monomorphime canonique $\chi^{\Pl}_{n}\hookrightarrow\sigma^{\Pl}_{n}$, une trivialisation sur $\chi^{\Pl}_{n}$.
\item Nous noterons $\chi_{n} : \F_{\Gr,n}\to\F_{GL_n(\kk)}$ le foncteur
  $(\chi^{\Pl}_{n})^{\Gr}_{GL_{n}}$. 
\item Le foncteur $\chi_n$ est adjoint à droite à $\theta_{n}$. 
\end{enumerate}
\end{prdef}

\begin{proof} La trivialisation canonique de $\sigma^{\Pl}_{n}$ se lit
  sur le foncteur de décalage pointé $\mathfrak{S}_n$ : pour tout
  $g\in GL_n(\kk)$, $\tau_g\mathfrak{S}_n$ est donné sur un objet $V$ de
  $\E^f$ par $(V\oplus E_n,E_n\xrightarrow{0\oplus g}V\oplus E_n)$, et
  les transformations naturelles
  $\tau_g\mathfrak{S}_n\to\mathfrak{S}_n$ données par le diagramme
  commutatif suivant fournissent la trivialisation recherchée. 
$$\xymatrix{E_n\ar[r]^-{0\oplus g}\ar[rd]_{0\oplus id} & V\oplus E_n\ar[d]^{id\oplus g^{-1}} \\
& V\oplus E_n
}$$

Pour en déduire le second point, on remarque que le diagramme
$$\xymatrix{\sigma^\Pl_n\tau^*_g\ar[r]\ar[d]_{T_g} & \sigma^\Pl_n\tau^*_g\otimes
  I_{E_n}\ar[d]^{T_g\otimes g^*} \\
\sigma^\Pl_n\ar[r] & \sigma^\Pl_n\otimes
  I_{E_n}
}$$
commute pour tout $g\in GL_n(\kk)$, où $T_g$ est l'isomorphisme de
trivialisation et les flèches horizontales sont donnéees par la somme
de $\psi$ et de l'inclusion canonique.

La dernière assertion s'obtient à partir des trois observations
suivantes.
\begin{enumerate}\item La composée
  $\F_{\Gr,n}\xrightarrow{\chi_n}\F_{GL_n(\kk)}\xrightarrow{O^\F_{GL_{n}}}\F$ est adjointe à droite à $\F\to\F_{GL_n(\kk)}\xrightarrow{\theta_n}\F_{\Gr,n}$, où la première flèche est la postcomposition par le foncteur $\E\to\,_{\kk[GL_n(\kk)]}\mathbf{Mod}$ d'extension des scalaires. En effet, le diagramme commutatif de la proposition \ref{prdf-ch6pl}, la proposition \ref{chi-pl} et  la dernière assertion de la proposition \ref{liengrpl} montrent que la première composée est adjointe à droite à $\F\xrightarrow{\kappa^\Pl_n}\F_{\Pl,n}\xrightarrow{\Phi_n}\F_{\Gr,n}$, qui coïncide avec $\kappa_n\otimes P^{\Gr,n}_{(E_n,E_n)}$, flèche à laquelle s'identifie également $\F\to\F_{GL_n(\kk)}\xrightarrow{\theta_n}\F_{\Gr,n}$.
\item Le foncteur $\chi_n$ est un sous-foncteur de $\sigma_n$. Cela
  découle de la première partie de la démonstra\-tion.
\item Il existe un diagramme commutatif
$$\xymatrix{{\rm hom}_{\F_{GL_n(\kk)}} (F,\chi_n(X))\ar@{-->}[d]_\simeq\ar@{^{(}->}[r] & {\rm
    hom}_{\F_{GL_n(\kk)}} (F,\sigma_n(X))\ar[d]^\simeq \\
{\rm hom}_{\Gr,n} (\theta_n(F),X)\ar@{^{(}->}[r] & {\rm hom}_{\Gr,n} (\xi_n(F),X)
}$$
naturel en les objets $F$ de $\F_{GL_n(\kk)}$ et $X$ de $\F_{\Gr,n}$, où la
flèche verticale de droite provient de la proposition~\ref{pradjf3}. L'isomorphisme en pointillé résulte des deux remarques
précédentes lorsque $F$ appartient à l'image du foncteur
$\F\to\F_{GL_n(\kk)}$ induit par l'extension des scalaires ; le cas général
s'en déduit par un argument de colimite.
\end{enumerate}
\end{proof}

\section{Foncteurs hom internes et foncteurs de division dans $\F_{\Gr,I}$}\label{fhifd}

Nous présentons dans cette section des propriétés élémentaires des
foncteurs hom internes et des foncteurs de division dans les
catégories $\F_{\Gr,I}$
(cf. proposition/définition~\ref{prdf-hid}). Ces propriétés
reposent sur les différentes adjonctions établies dans les
paragraphes~\ref{subsfgr} et~\ref{sfe-a} et la description explicite des
projectifs et injectifs standard de $\F_{\Gr,I}$.

Nous nous attacherons surtout aux foncteurs hom internes, notre but
principal étant d'aboutir, dans la partie~\ref{p-omeg}, à la
proposition~\ref{hl1} et ses corollaires, relatifs au comportement
mutuel des foncteurs hom internes du foncteur $\omega$. Les foncteurs de division n'interviendront que comme
auxiliaires, dans le paragraphe~\ref{scind-do}.

\subsection{Comparaison entre les différentes catégories
$\F_{\Gr,I}$}\label{parccfg} Il est parfois commode de passer de la
catégorie globale $\F_\Gr$ aux catégories $\F_{\Gr,n}$, par exemple,
pour traiter des adjoints au produit tensoriel. En effet, si les
foncteurs $\iota$ et $\omega$ sont adjoints, il n'en est pas de même
des foncteurs $\iota_n$ et $\omega_n$, ce qui induit une difficulté
dans $\F_{\Gr,n}$ nouvelle par rapport à $\F_\Gr$. En revanche, dès
que l'on traite de foncteurs pseudo-constants, la catégorie globale se
trouve beaucoup moins maniable que $\F_{\Gr,n}$, puisque les adjoints
au produit tensoriel entre $GL_n(\kk)$-modules sont faciles à décrire,
contrairement à ce qui advient dans $\F_{surj}(\kk)$.

\begin{conv} Dans ce paragraphe, on se donne deux parties $I$ et $J$ de $\mathbb{N}$
telles que $I\subset J$.
\end{conv}

On rappelle que $\mathcal{R}_{J,I}$ est le foncteur de restriction et
$\mathcal{P}_{I,J}$ le foncteur de prolongement par zéro --- défini sous
certaines conditions sur $I$ et $J$ --- introduits dans la notation~\ref{notglfgr}.

\begin{pr}\label{cprhe} Supposons que tout élément de $J$ supérieur ou égal à
un élément de $I$ appartient à $I$. Alors il existe dans $\F_{\Gr,J}$ un isomorphisme
$$\mathbf{Hom}_{\Gr,J}(X,\mathcal{P}_{I,J}(Y))\simeq\mathcal{P}_{I,J}(\mathbf{Hom}_{\Gr,I}(\mathcal{R}_{J,I}(X),Y))$$
naturel en les objets $X$ de $\F_{\Gr,J}$ et $Y$ de $\F_{\Gr,I}$,
et dans $\F_{\Gr,I}$ un isomorphisme naturel
$$\mathcal{R}_{J,I}(\mathcal{P}_{I,J}(Y) : X)\simeq (Y :
\mathcal{R}_{J,I}(X))$$
si $X$ est à valeurs de dimension finie.

En particulier, pour tout $n\in\mathbb{N}$, on a dans $\F_{\Gr,\leq n}$ un isomorphisme naturel 
$$\mathbf{Hom}_{\Gr,\leq n}(X,\mathcal{P}_{n,\leq
  n}(Y))\simeq\mathcal{P}_{n,\leq
  n}(\mathbf{Hom}_{\Gr,n}(\mathcal{R}_{\leq n,n}(X),Y))$$
et, si $X$ est à valeurs de dimension finie, un isomorphisme naturel dans $\F_{\Gr,n}$
$$\mathcal{R}_{\leq n,n}(\mathcal{P}_{n,\leq n}(Y) : X)\simeq (Y :
\mathcal{R}_{\leq n,n}(X)).$$
\end{pr}

\begin{proof} L'hypothèse assure que $\E^f_{\Gr,I}$ est une
  sous-catégorie {\em complète à gauche} de $\E^f_{\Gr,J}$. Par
  conséquent, le foncteur $\mathcal{P}_{I,J}$ est adjoint à droite à
  $\mathcal{R}_{J,I}$ (cf. proposition \ref{prfrec}). Le premier
  isomorphisme de la proposition
  provient alors de la commutation du foncteur de
  restriction au produit tensoriel, comme le montrent les
  isomorphismes naturels
$${\rm
  hom}_{\Gr,J}(A,\mathbf{Hom}_{\Gr,J}(X,\mathcal{P}_{I,J}(Y)))\simeq{\rm hom}_{\Gr,J}(A\otimes X,\mathcal{P}_{I,J}(Y))\simeq {\rm hom}_{\Gr,I}(\mathcal{R}_{J,I}(A\otimes X),Y)$$
$$\simeq{\rm hom}_{\Gr,I}(\mathcal{R}_{J,I}(A)\otimes
\mathcal{R}_{J,I}(X),Y)\simeq {\rm
  hom}_{\Gr,I}(\mathcal{R}_{J,I}(A),\mathbf{Hom}_{\Gr,I}(\mathcal{R}_{J,I}(X),Y))$$
$$\simeq{\rm
  hom}_{\Gr,J}(A,\mathcal{P}_{I,J}(\mathbf{Hom}_{\Gr,I}(\mathcal{R}_{J,I}(X),Y)))$$
(où $A\in{\rm Ob}\,\F_{\Gr,J}$) et le lemme de Yoneda. Le second se
traite de façon analogue.
\end{proof}

Nous donnons maintenant un résultat plus précis de commutation entre
foncteurs hom internes et foncteurs de restriction ou de prolongement
par zéro, dans le cas où la source du foncteur hom interne appartient
à l'image du foncteur $\iota_I : \F\to\F_{\Gr,I}$.

\begin{pr}\label{pcorp} Il existe des isomorphismes naturels
$$\mathbf{Hom}_{\Gr,J}(\iota_J(F),\mathcal{P}_{I,J}(A))\simeq\mathcal{P}_{I,J}(\mathbf{Hom}_{\Gr,I}(\iota_I(F),A))$$
dans $\F_{\Gr,J}$, si le prolongement par zéro est défini,  et 
$$\mathbf{Hom}_{\Gr,I}(\iota_I(F),\mathcal{R}_{J,I}(X))\simeq\mathcal{R}_{J,I}(\mathbf{Hom}_{\Gr,J}(\iota_J(F),X))$$
 dans $\F_{\Gr,I}$, où $F\in {\rm Ob}\,\F$, $A\in {\rm Ob}\,\F_{\Gr,I}$ et $X\in {\rm Ob}\,\F_{\Gr,J}$.
\end{pr}

\begin{proof} Les morphismes naturels
$$\mathcal{P}_{I,J}(\mathbf{Hom}_{\Gr,I}(\iota_I(F),A))\otimes\iota_J(F)\simeq
\mathcal{P}_{I,J}(\mathbf{Hom}_{\Gr,I}(\iota_I(F),A)\otimes\iota_I(F))\to\mathcal{P}_{I,J}(A)$$
(lorsque le prolongement par zéro est défini) et
$$\mathcal{R}_{J,I}(\mathbf{Hom}_{\Gr,J}(\iota_J(F),X))\otimes\iota_I(F)\simeq
\mathcal{R}_{J,I}(\mathbf{Hom}_{\Gr,J}(\iota_J(F),X)\otimes\iota_J(F))\to\mathcal{R}_{J,I}(X)$$
dont les secondes flèches s'obtiennent par application du foncteur
$\mathcal{P}_{I,J}$ ou $\mathcal{R}_{J,I}$ à la coünité de
l'adjonction fournissent, par
adjonction, des morphismes naturels
$\mathcal{P}_{I,J}(\mathbf{Hom}_{\Gr,I}(\iota_I(F),A))\to\mathbf{Hom}_{\Gr,J}(\iota_J(F),\mathcal{P}_{I,J}(A))$ et
$\mathcal{R}_{J,I}(\mathbf{Hom}_{\Gr,J}(\iota_J(F),X))\to\mathbf{Hom}_{\Gr,I}(\iota_I(F),\mathcal{R}_{J,I}(X))$.
Ce sont des isomorphismes dans le cas où $F$ est un projectif standard
grâce aux propositions~\ref{adelgr} et~\ref{rqedgh}. Le cas général s'en déduit par passage à la colimite.
\end{proof}

La notation {\em ad hoc} suivante n'interviendra que dans les deux
lemmes techniques ci-dessous.

\begin{nota} Dans ce paragraphe, nous noterons $r : \F_{\Gr,I}\to\F_{\Gr,J}$ le
foncteur adjoint à gauche au foncteur $\mathcal{R}_{J,I} :
\F_{\Gr,J}\to\F_{\Gr,I}$, et $R : \F^I_{surj}\to\F^J_{surj}$ l'adjoint à
  gauche au foncteur de restriction $\F^J_{surj}\to\F^I_{surj}$ (de tels adjoints existent par le corollaire \ref{cr-adjp}).
\end{nota}

\begin{lm}\label{lmrr} Le diagramme
$$\xymatrix{\F_{\Gr,I}\ar[r]^-{\sigma_I}\ar[d]_{r} & \F\otimes\F^I_{surj}\simeq\mathbf{Fct}(\E^f,\F^I_{surj})\ar[d]^-{R_*} \\
\F_{\Gr,J}\ar[r]^-{\sigma_J} & \F\otimes\F^J_{surj}\simeq\mathbf{Fct}(\E^f,\F^J_{surj})
}$$
commute à isomorphisme naturel près.
\end{lm}

\begin{proof} Soit $a : \F\otimes\F^J_{surj}\to\F\otimes\F^I_{surj}$ le foncteur de
restriction. On définit
une transformation naturelle $R_*\sigma_I\to \sigma_J r$ comme
l'adjointe (cf. proposition \ref{lmfora}) à la flèche $\sigma_I\to
a\sigma_J r\simeq\sigma_I\mathcal{R}_{J,I} r$ obtenue en composant
$\sigma_I$ et l'unité $id\to\mathcal{R}_{J,I} r$ de
l'adjonction. Cette transformation naturelle est un isomorphisme, car elle induit un isomorphisme ${\rm hom}_{\Gr,J}(\sigma_J
r(X),I^{\E^f\times\E^J_{surj}}_{(A,B)})\xrightarrow{\simeq}{\rm
  hom}_{\Gr,J}(R_*\sigma_I(X),I^{\E^f\times\E^J_{surj}}_{(A,B)})$ pour
tous $X\in {\rm Ob}\,\F_{\Gr,I}$ et $(A,B)\in {\rm Ob}\,\E^f\times\E^J_{surj}$
grâce aux isomorphismes (\ref{eq-cfg2}) (page~\pageref{eq-cfg2}).
\end{proof}

\begin{lm}\label{lmkac}   Il existe dans $\F_{\Gr,J}$ un isomorphisme
$r(X\otimes\kappa_I(F))\simeq r(X)\otimes\kappa_J(F)$ naturel en les
objets $X$ de $\F_{\Gr,I}$ et $F$ de $\F$.
\end{lm}

\begin{proof} On établit, à l'aide du lemme précédent, que la transformation naturelle
$r(X\otimes\kappa_I(F))\to r(X)\otimes\kappa_J(F)$  adjointe au
morphisme
$X\otimes\kappa_I(F)\to\mathcal{R}_{J,I}r(X)\otimes\kappa_I(F)\simeq\mathcal{R}_{J,I}(r(X)\otimes\kappa_J(F))$
obtenu en tensorisant par $\kappa_I(F)$ l'unité de l'adjonction
procure un isomorphisme lorsqu'on lui applique $\sigma_J$. Cela
provient du lemme~\ref{lmrr} et des deux observations suivantes :
\begin{itemize}\item il existe dans $\F\otimes\F^I_{surj}$ un isomorphisme canonique
  $\sigma_I(X\otimes\kappa_I(F))\simeq\sigma_I(X)\otimes
  (F\boxtimes\kk)$, que l'on peut encore
  voir comme l'image de $\sigma_I(X)$ par le foncteur
  $\F\otimes\F^I_{surj}\simeq\mathbf{Fct}(\E^I_{surj},\F)\xrightarrow{(\cdot\otimes F)_*}\mathbf{Fct}(\E^I_{surj},\F)$ ;
\item le diagramme
$$\xymatrix{\F\otimes\F^I_{surj}\simeq\mathbf{Fct}(\E^I_{surj},\F)\ar[d]_-{R_*}\ar[rr]^-{(\cdot\otimes F)_*} 
  & & \mathbf{Fct}(\E^I_{surj},\F)\simeq\F\otimes\F^I_{surj}\ar[d]^-{R_*}
  \\
\F\otimes\F^J_{surj}\simeq \mathbf{Fct}(\E^J_{surj},\F)\ar[rr]^-{(\cdot\otimes F)_*} 
  & & \mathbf{Fct}(\E^J_{surj},\F)\simeq\F\otimes\F^J_{surj}
}$$
commute (à isomorphisme canonique près).
\end{itemize}
Comme le
foncteur $\sigma_J$ est exact et fidèle, cela donne la conclusion.
\end{proof}

Nous sommes désormais en mesure d'établir la commutation des foncteurs
hom internes et des foncteurs de restriction, dans le cas où la source
du foncteur hom interne appartient à l'image du foncteur $\kappa_I : \F\to\F_{\Gr,I}$.

\begin{pr}\label{pcork}  Il existe un isomorphisme
$$\mathbf{Hom}_{\Gr,I}(\kappa_I(F),\mathcal{R}_{J,I}(X))\simeq\mathcal{R}_{J,I}\big(\mathbf{Hom}_{\Gr,J}(\kappa_J(F),X)\big)$$
naturel en les objets $F$ de $\F$ et $X$ de $\F_{\Gr,J}$.
\end{pr}

\begin{proof} Par le lemme \ref{lmkac} (dont on conserve la notation),
  on a des isomorphismes naturels
$${\rm
  hom}_{\Gr,I}(Y,\mathbf{Hom}_{\Gr,I}(\kappa_I(F),\mathcal{R}_{J,I}(X)))\simeq{\rm hom}_{\Gr,I}(Y\otimes\kappa_I(F),\mathcal{R}_{J,I}(X))\simeq$$
$${\rm hom}_{\Gr,J}(r(Y\otimes\kappa_I(F)),X)\simeq {\rm
  hom}_{\Gr,J}(r(Y)\otimes \kappa_J(F),X)\simeq$$
$${\rm hom}_{\Gr,J}(r(Y),\mathbf{Hom}_{\Gr,J}(\kappa_J(F),X))\simeq {\rm
  hom}_{\Gr,I}\big(Y,\mathcal{R}_{J,I}\big(\mathbf{Hom}_{\Gr,J}(\kappa_J(F),X)\big)\big)$$
(où $Y\in {\rm Ob}\,\F_{\Gr,I}$), ce qui démontre la proposition.
\end{proof}

\subsection{Propriétés formelles}\label{parhipf} Nous étudions à
présent le comportement mutuel des foncteurs de division ou hom
internes et des foncteurs fondamentaux du §\,\ref{subsfgr}.

\paragraph*{Foncteurs de division} C'est le cas des foncteurs de
division par des objets injectifs de la catégorie $\F_\Gr$ qui nous
intéresse le plus (cf.  proposition~\ref{divif}).

Nous donnons d'abord une propriété de compatibilité entre les foncteurs $\omega$, $\iota$ et
les foncteurs de division.


\begin{pr}\label{domeg} Il existe dans $\F$ un isomorphisme $\omega (X : \iota
  (F))\simeq (\omega(X) : F)$ naturel en les objets $X$ de $\F_\Gr$ et
  $F$ de $\F^{df}$. On a un résultat analogue dans $\F_{\Gr,\leq n}$.
\end{pr}

\begin{proof} Il s'agit d'une conséquence formelle de l'adjonction
  entre $\iota$ et $\omega$ (cf. proposition~\ref{prfig}) et de la commutation de $\iota$ au produit
  tensoriel (cf. proposition~\ref{cprhe}).
\end{proof}

La proposition élémentaire suivante interviendra dans \cite{art3}.

\begin{pr}\label{divrho} Pour toute partie $I$ de $\mathbb{N}$, il
  existe dans $\F_{\Gr,I}$ un isomorphisme canonique $(\rho_I(A) : X)\simeq\rho_I(A :
  \varepsilon_I(X))$ pour $A\in {\rm Ob}\,\F^I_{surj}$, et $X\in {\rm
    Ob}\,\F_{\Gr,I}$ à valeurs de dimension finie.
\end{pr}

\begin{proof} Cette proposition s'établit de la même manière que la
  précédente, à partir de la proposition~\ref{prcompa4} et de la commutation
  de $\varepsilon_I$ au produit tensoriel.
\end{proof}

On rappelle que le symbole $Gr$ qui apparaît dans l'importante proposition
suivante a été introduit au notation~\ref{notgrd}.

\begin{pr}\label{divif} Il existe un isomorphisme
$$(X : \iota(I_V))(A,B)\simeq\underset{im\,(C\hookrightarrow V\oplus
  A\twoheadrightarrow A)=B}{\bigoplus_{C\in\Gr(V\oplus A)}} X(V\oplus
A,C)$$
naturel en les objets $V$ de $\E^f$, $(A,B)$ de $\E^f_\Gr$ et $X$ de
$\F_\Gr$.

De plus, pour tout $W\in\Gr(V)$, le monomorphisme scindé
naturel $(X:I^\Gr_{(V,W)})\hookrightarrow (X : \iota(I_V))$ induit par
l'épimorphisme scindé $\iota(I_V)\twoheadrightarrow I_{(V,W)}^\Gr$
fourni par la proposition \ref{prfif} identifie
$(X:I^\Gr_{(V,W)})(A,B)$ au sous-espace 
$$\bigoplus_{C\in Gr(W,B)} X(V\oplus
A,C)$$
de $(X : \iota(I_V))(A,B)$.
\end{pr}


\begin{proof} On a des isomorphismes naturels
$$\big((X:I^\Gr_{(V,W)})(A,B)\big)^*\simeq {\rm hom}_\Gr
((X:I^\Gr_{(V,W)}),I^\Gr_{(A,B)})\simeq {\rm hom}_\Gr
(X,I^\Gr_{(V,W)}\otimes I^\Gr_{(A,B)})$$
$$\simeq\bigoplus_{C\in Gr(W,B)}{\rm hom}_\Gr
(X,I^\Gr_{(V\oplus A,C)}) \simeq\bigoplus_{C\in Gr(W,B)} X(V\oplus
A,C)^*$$
par la proposition \ref{prtifgr}, ce qui donne le résultat, 
pour $X$ à valeurs de dimension finie, en dualisant ; le cas
quelconque s'en déduit en écrivant $X$ comme colimite de
sous-foncteurs à valeurs de dimension finie.
\end{proof}


\paragraph*{Foncteurs hom internes} Nous ramenons dans ce
sous-paragraphe le calcul de nombreux foncteurs hom internes dans la
catégorie $\F_{\Gr,I}$ à celui de foncteurs hom internes des
catégories $\F$ ou $\F\otimes\F^I_{surj}$.

\begin{conv} Dans la suite de ce paragraphe, $I$ désigne une partie de $\mathbb{N}$.
\end{conv}

Le principal objectif de la proposition suivante est d'établir la
commutation entre le foncteur $\kappa_I$ et les foncteurs hom internes (corollaire~\ref{hiki}).

\begin{pr}\label{hieta} Il existe un isomorphisme
$$\mathbf{Hom}_{\Gr,I}(X,\theta_I(A))\simeq\theta_I\big(\mathbf{Hom}_{\F\otimes\F^I_{surj}} (\eta_I(X),A)\big)$$
naturel en les objets $X$ de $F_{\Gr,I}$ et $A$ de $\F\otimes\F^I_{surj}$.
\end{pr}

\begin{proof} C'est une conséquence formelle de l'adjonction entre les
  foncteurs $\theta_I$ et $\eta_I$ (proposition~\ref{adjeta}) et de la
  commutation du foncteur $\eta_I$ au produit tensoriel (proposition~\ref{pretapt}).
\end{proof}

\begin{cor}\label{crxet} Il existe un isomorphisme
$$\mathbf{Hom}_{\Gr,I}(\theta_I(B),\theta_I(A))\simeq\mathbf{Hom}_{\Gr,I}(\xi_I(B),\theta_I(A))\simeq\theta_I\big(\mathbf{Hom}_{\F\otimes\F^I_{surj}} (B,A)\big)$$
naturel en les objets $A$ et $B$ de $\F\otimes\F^I_{surj}$.
\end{cor}

\begin{proof} On combine les propositions \ref{hieta} et~\ref{propeta}.
\end{proof}

\begin{cor}\label{hiki} Il existe des isomorphismes
$$\mathbf{Hom}_{\Gr,I}(\kappa_I(F),\kappa_I(G))\simeq\mathbf{Hom}_{\Gr,I}(\iota_I(F),\kappa_I(G))\simeq\kappa_I\big(\mathbf{Hom}_\F(F,G)\big)$$
naturels en les objets $F$ et $G$ de $\F$.
\end{cor}

\begin{proof} C'est le cas particulier $A=G\boxtimes\kk$ et
  $B=F\boxtimes\kk$ du corollaire précédent.\end{proof}

\begin{rem} Une approche alternative de ce corollaire consiste à
  traiter d'abord le cas où la partie $I$ contient $0$ à l'aide de la
  proposition \ref{adjf5}, puis d'en déduire le cas général grâce à la
  proposition \ref{pcork}.
\end{rem}


Nous traitons à présent de la commutation du foncteur $\iota_I$ aux
foncteurs hom internes.

\begin{pr}\label{hiio} Il existe
  un isomorphisme naturel
  $\mathbf{Hom}_{\Gr,I}(\iota_I(F),\iota_I(G))\simeq\iota_I\big(\mathbf{Hom}_\F(F,G)\big)$ pour $F$, $G\in {\rm Ob}\,\F$.
\end{pr}

\begin{proof} Comme
  $\iota_I(F)=\mathcal{R}_{\mathbb{N},I}\big(\iota(F)\big)$, la proposition \ref{pcorp} permet de ne traiter que le
  cas où $I=\mathbb{N}$. La conclusion provient alors des assertions
  $2$ et $4$ de la proposition~\ref{prfig}.
\end{proof}

Il existe aussi une propriété de commutation relative au foncteur
$\sigma_I$ :

\begin{pr}\label{adjhis} Il existe dans $\F\otimes\F^I_{surj}$ un isomorphisme
$$\sigma_I\big(\mathbf{Hom}_{\Gr,I} (\xi_I(F),X)\big)\simeq\mathbf{Hom}_{\F\otimes\F^I_{surj}}(F,\sigma_I(X))$$
naturel en les objets $F$ de $\F\otimes\F^I_{surj}$ et $X$ de $\F_{\Gr,I}$.
\end{pr}

\begin{proof} Cette propriété s'obtient à partir de l'adjonction entre
  les foncteurs $\xi_I$ et $\sigma_I$ (proposition~\ref{pradjf3}) et de la commutation du foncteur
  $\xi_I$ au produit tensoriel.
\end{proof}

La généralisation de la description très simple des
foncteurs hom internes de la catégorie $_{\kk[GL_n(\kk)]}\bf{Mod}$
(cf. \cite{CR}, §\,10 D) aux
foncteurs pseudo-constants de $\F_{\Gr,n}$ est donnée par la
proposition suivante, laissée au lecteur (cf. \cite{these} pour les détails).

\begin{pr}\label{aglhi} Soient $n\in\mathbb{N}$ et $M$ un
  $\kk[GL_n(\kk)]$-module fini. Les endofoncteurs
  $\cdot\otimes\rho_n(M^*)$, $\mathbf{Hom}_{\Gr,n}(\rho_n(M),\cdot)$
  et $(\,\cdot\,:\rho_n(M))$ de $\F_{\Gr,n}$ sont naturellement isomorphes.
\end{pr}

%

On en déduit aussitôt la propriété de commutation suivante :

\begin{cor}\label{crcomdh} Soient $n\in\mathbb{N}$ et $M$ un
  $\kk[GL_n(\kk)]$-module fini. L'endofoncteur
  $\cdot\otimes\rho_n(M)$ de $\F_{\Gr,n}$ commute naturellement aux
  foncteurs hom internes et aux foncteurs de division.
\end{cor}

\paragraph*{Le foncteur $\mathbf{Ext}^*_\Gr(\iota(F),\cdot)$}
Commençons par introduire une nouvelle notation. Soit $A$ un objet
de $\E^f_\Gr$ ; nous désignerons par $\tau_A : \F_\Gr\to\F$ le
foncteur de précomposition par $\E^f\to\E^f_\Gr\quad E\mapsto
E\boxplus A$. Cette construction est fonctorielle en $A$.

\begin{pr}\label{extint} Il existe un isomorphisme
  $\mathbf{Ext}^*_\Gr(\iota(F),X)(A)\simeq {\rm
    Ext}^*_\F(F,\tau_A(X))$ naturel en les objets $F$ de $\F$,
  $X$ de $\F_\Gr$ et $A$ de $\E^f_\Gr$.
\end{pr}

\begin{proof} Il existe des isomorphismes naturels
  $\tau_A\circ\iota\simeq\Delta_{\mathfrak{O}(A)}$ et
  $\tau_A(P^\Gr_A)\simeq\kk[{\rm End}_\Gr(A)]\otimes P_{\mathfrak{K}(A)}$
  (par la proposition~\ref{actc-prf}). On en déduit une transformation
  naturelle $F\to\tau_A(\iota(F)\otimes
  P^\Gr_A)\simeq\Delta_{\mathfrak{O}(A)}(F)\otimes\kk[{\rm
    End}_\Gr(A)]\otimes P_{\mathfrak{K}(A)}$ par produit tensoriel des
  injections canoniques $F\hookrightarrow\Delta_{\mathfrak{O}(A)}(F)$,
  $\kk\hookrightarrow P_{\mathfrak{K}(A)}$ et $\kk\hookrightarrow\kk[{\rm
    End}_\Gr(A)]$ (donnée par $[id_A]$).

L'application naturelle ${\rm hom}_\Gr(\iota(F)\otimes P^\Gr_A,X)\to
{\rm hom}_\F(F,\tau_A(X))$ qu'on en déduit est bijective. En effet, il
suffit de le voir pour $F$ projectif standard de $\F$, auquel cas
c'est une conséquence de la proposition~\ref{adelgr}. Cet isomorphisme
s'étend, grâce à la proposition~\ref{adj-ah}, en un isomorphisme
naturel gradué ${\rm Ext }^*_\Gr(\iota(F)\otimes P^\Gr_A,X)\xrightarrow{\simeq} {\rm
    Ext}^*_\F(F,\tau_A(X))$. La conclusion provient alors de
  l'isomorphisme naturel (\ref{eq-hi}) de l'appendice~\ref{apfct}.
\end{proof}

\part{Propriétés du foncteur $\omega$. Applications}\label{p-omeg}

Le foncteur d'intégrale en grassmanniennes $\omega : \F_\Gr\to\F$
possède un comportement qui diffère notablement de celui des autres
foncteurs fondamentaux entre $\F_\Gr$ et $\F$ ou
$\F\otimes\F_{surj}$. En particulier, il ne préserve pas les objets localement
finis ; en fait, le foncteur $\omega(X)$ de $\F$ n'est localement fini que si le
foncteur $X$ de $\F_\Gr$ appartient à la sous-catégorie
$\F_{\Gr,0}\simeq\F$ --- cela résultera du théorème~\ref{prfondo}. Une
autre façon d'illustrer ce phénomène consiste à étudier la composée
$\Delta\circ\omega$ : elle est \go nettement plus grosse\gf que
$\omega\circ\Delta$ --- cf.~§\,\ref{scind-do}.

C'est dans ces observations que réside tout l'intérêt du
foncteur~$\omega$, en vue de l'étude de la structure de la
catégorie~$\F$. Remarquons que les progrès significatifs 
obtenus à ce sujet par Powell dans \cite{GP2} constituent essentiellement des traductions de
propriétés du foncteur $\varpi : \F_{surj}\to\F$ (cf. §\,\ref{parfsf})
; il n'est pas étonnant que le foncteur $\omega$, qui généralise
$\varpi$, permette d'aller plus loin dans la compréhension de la
catégorie~$\F$.

\medskip

L'une des avancées importantes obtenues grâce aux catégories de foncteurs
en grassmanniennes réside en le théorème~\ref{prfondo}, qui constitue le
résultat principal de cet article. Il donne une propriété d'annulation
cohomologique très générale, qui trouve deux applications
essentielles.

L'une d'entre elle, traitée dans le paragraphe~\ref{par-inj-f},
concerne le lien entre cohomologie du groupe linéaire $GL(\kk)$ (ou
$K$-théorie stable de $\kk$) et cohomologie fonctorielle.  Le résultat
obtenu, qui généralise le théorème de Betley-Suslin
donné dans l'appendice de \cite{FFSS}, fournit un isomorphisme entre la cohomologie de
$GL(\kk)$ à coefficients convenables et des groupes d'extensions dans
la catégorie $\F_\Gr(\kk)$ entre des objets finis, nettement plus faciles
d'accès. 

Les calculs cohomologiques dans la catégorie $\F(\kk)$ donnés par le
théorème~\ref{prfondo} possèdent également un intérêt
intrinsèque. Ainsi, ce théorème suggère une conjecture décrivant la filtration de Krull de la catégorie $\F$, qui renforce toutes
les formes antérieurement formulées de la conjecture artinienne ; il est l'un
des outils fondamentaux de la démonstration des formes partielles que
nous en établirons dans~\cite{art3}. Nous discutons ces questions dans
la section~\ref{s-fkf}.

L'application du théorème~\ref{prfondo} exposée dans la
section~\ref{s-rhid} traite de propriétés de commutation entre le
foncteur~$\omega$ et des foncteurs hom internes. Son intérêt est
illustré par~\cite{art1}, comme il est discuté à la fin du §\,\ref{partnhi}.

\section{Théorème d'annulation cohomologique}\label{sscep}

L'objectif de cette section consiste à établir le
théorème~\ref{prfondo} et le corollaire~\ref{fondcr} qui s'en
déduit. Leur démonstration repose sur des considérations explicites
liées à la catégorie $\F_{surj}$ et des arguments d'adjonction. On
emploie également une catégorie de foncteurs auxiliaire, variante de
$\F_\Gr$ obtenue en considérant deux éléments de la grassmannienne
d'un espace vectoriel. Hormis l'utilisation de ce nouveau type de
catégories, notre démarche procède des mêmes idées conceptuelles que
celles inaugurées par Pirashvili (cf. remarque~\ref{rq-pira}\,.\ref{rqitpira}).

\subsection{Préliminaires}\label{ssecp}

On rappelle que $o : \F\to\F_{surj}$ désigne le foncteur d'oubli (cf. section~\ref{sct-surj})  et que le
foncteur $\kk^{\geq i}$ est défini dans la notation~\ref{notstr-cs} (page~\pageref{notstr-cs}).

Les résultats de cette section reposent sur la structure du foncteur
constant $\kk=I^{surj}_0$ de $\F_{surj}$ donnée par le corollaire
\ref{str-cs}. On rappelle que ce foncteur est {\em unisériel}.

\begin{lm}\label{lmcle} Il existe dans $\F_{surj}$ un monomorphisme
$\kk^{\geq 1}\hookrightarrow o(\bar{I}_\kk)$, unique à homothétie près.
\end{lm}

\begin{proof} Par l'assertion \ref{itvp2} de la proposition
  \ref{pre-vpo}, il existe dans $\F_{surj}$ un diagramme commutatif
$$\xymatrix{o(I_0)\ar@{^{(}->}[d]\ar[r]^{\simeq} & I^{surj}_0\ar@{^{(}->}[d] \\
o(I_\kk) \ar[r]^-{\simeq} & I^{surj}_0\oplus I^{surj}_\kk
}$$
On en déduit un
isomorphisme entre $o(\bar{I}_\kk)$ et $I^{surj}_\kk$.

Par conséquent, ${\rm hom}\,(\kk^{\geq 1},o(\bar{I}_\kk))\simeq\kk^{\geq
  1}(\kk)^*=\kk$. Le morphisme correspondant à $1$ est injectif, car sa
restriction au socle de $\kk^{\geq 1}$, l'objet simple $S_1^{surj}$ (cf. corollaire
\ref{str-cs}), est non nulle, l'inclusion $S_1^{surj}\hookrightarrow
\kk^{\geq 1}$ induisant un isomorphisme par évaluation sur
$\kk$. Cela démontre le lemme.
\end{proof}

%

\begin{rem}\label{rqinvi} Le foncteur $\kk^{\geq 1}$ est {\em
    idempotent} pour le produit tensoriel : $(\kk^{\geq 1})^{\otimes
    2}\simeq\kk^{\geq 1}$. Il induit donc un endofoncteur idempotent
  $\cdot\otimes\kk^{\geq 1}$ de $\F_{surj}$.
\end{rem}

Avant d'exploiter le lemme \ref{lmcle}, nous aurons besoin d'introduire des
catégories de foncteurs auxiliaires (qui n'interviendront que dans ce
paragraphe).

\begin{nota}\begin{enumerate}\item Nous désignerons par $\E^f_{bi-\Gr}$
    la catégorie $(\E^f_\Gr)_{\backslash\Gr}$, où l'on note encore  $\Gr$,
    par abus de notation, le
    foncteur composé
$$\E^f_\Gr\xrightarrow{\mathfrak{O}}\E^f\xrightarrow{\Gr}\bf{Ens},$$
qui est donc donné sur les objets par $(V,W)\mapsto\Gr(V)$.

Autrement dit, les objets de $\E^f_{bi-\Gr}$ sont les triplets
$(V,B,W)$ formés d'un espace vectoriel de dimension finie $V$ et de
deux sous-espaces $B$ et $W$ de $V$, et les flèches $(V,B,W)\to
(V',B',W')$ de $\E^f_{bi-\Gr}$ sont les applications linéaires $f :
V\to V'$ telles que $f(B)=B'$ et $f(W)=W'$.

La catégorie de foncteurs $\mathbf{Fct}(\E^f_{bi-\Gr},\E)$ sera notée $\F_{bi-\Gr}$.
\item  Nous désignerons par $\E^f_{2-\D r}$ la catégorie
  $(\E^f_\Gr)_{\backslash\Gr'}$, où l'on désigne par $\Gr'$ le
  foncteur
$$\E^f_\Gr\xrightarrow{\mathfrak{B}}\E^f\xrightarrow{\Gr}\bf{Ens},$$
qui est donc donné sur les objets par $(V,W)\mapsto\Gr(W)$.

La catégorie de foncteurs $\mathbf{Fct}(\E^f_{2-\D r},\E)$ sera notée
$\F_{2-\D r}$.
\item L'inclusion naturelle de foncteurs ensemblistes
  $\mathfrak{B}\hookrightarrow\mathfrak{O}$ identifie $\E^f_{2-\D r}$ à
  une sous-catégorie pleine de $\E^f_{bi-\Gr}$. Nous noterons $i_{\D
    r} : \F_{bi-\Gr}\to\F_{2-\D r}$ le foncteur de précomposition par l'inclusion.
\end{enumerate}
\end{nota}

\begin{rem}  L'indice $2-\D r$ utilisé est une abréviation de {\em
    $2$-drapeau}. On pourrait naturellement généraliser les
  considérations de ce paragraphe à des catégories correspondant aux
  espaces vectoriels munis de drapeaux de longueur fixée arbitraire.
\end{rem}

On dispose, conformément aux résultats de la section \ref{sctccf},
d'un foncteur de plongement $\Upsilon_\Gr : \F_\Gr\to\F_{bi-\Gr}$ et
d'un foncteur d'intégrale $\Omega_\Gr : \F_{bi-\Gr}\to\F_\Gr$. Sur les
objets, on a $\Upsilon_\Gr (X)(V,B,W)=X(V,B)$ et $\Omega_\Gr(A)(V,B)=\underset{W\in\Gr(V)}{\bigoplus}A(V,B,W)$.

\begin{defi} Nous noterons $\N$ la sous-catégorie épaisse de $\F_{bi-\Gr}$
  noyau du foncteur exact $i_{\D r}$. Explicitement, $\N$ contient les
  objets $X$ de $\F_{bi-\Gr}$ tels que $X(V,B,W)=0$ pour tout
  objet $(V,B,W)$ de $\E^f_{bi-\Gr}$ tel que $W\subset B$.
\end{defi}

\begin{rem}\begin{enumerate}\item La catégorie $\E^f_{bi-\Gr}$ s'identifie canoniquement à
  $\E^f_{\backslash (\Gr,\Gr)}$ : c'est la catégorie des espaces
  vectoriels de dimension finie munie de deux sous-espaces. Nous
  avons préféré donner la présentation ci-dessus~\guillemotleft~dessymétrisant les
  deux bases~\guillemotright~de manière à introduire plus naturellement les
  foncteurs $\Upsilon_\Gr$, $\Omega_\Gr$ et $i_{\D r}$.
\item La catégorie $\F_{2-\D r}$ n'interviendra pas en elle-même,
c'est la sous-catégorie $\N$ de $\F_{bi-\Gr}$ qui nous intéressera.
\end{enumerate}
\end{rem}

Le lemme suivant donne les faits concrets dont nous aurons besoin pour
établir nos résultats d'annulation cohomologique.

\begin{lm}\label{lme33}\begin{enumerate}\item Il existe un foncteur
    $\beta : \E^f_{bi-\Gr}\to\E^f_{surj}$ tel que
    $\beta(V,B,W)=W/(W\cap B)$ pour tout objet $(V,B,W)$ de $\E^f_{bi-\Gr}$.
\item Il existe dans $\F_{bi-\Gr}$ un épimorphisme du foncteur
  constant $\kk$ vers le foncteur $\beta^*(\kk^{\geq 1})$, où $\beta^* : \F_{surj}\to\F_{bi-\Gr}$ désigne le foncteur de
  précomposition par $\beta$ ; $\N$ est la sous-catégorie pleine des objets $X$ de  $\F_{bi-\Gr}$ tels que
  la projection canonique $X\twoheadrightarrow X\otimes\beta^*(\kk^{\geq 1})$
  qui s'en déduit est un isomorphisme. En particulier, $\N$ est un
  {\bf idéal} de $\F_{bi-\Gr}$ : le produit tensoriel d'un objet de $\N$ et
  d'un objet de $\F_{bi-\Gr}$ appartient à~$\N$.
\item Le foncteur $\beta^*(\kk^{\geq
    1})$ se plonge dans $\Upsilon_\Gr(\kappa(\bar{I}_\kk))$.
\end{enumerate}
\end{lm}

\begin{proof} 
La vérification du premier point est immédiate. Le second point vient de ce que la projection $\kk\twoheadrightarrow\kk^{\geq 1}$ de
$\F_{surj}$ (cf. lemme \ref{str-cs}) induit un épimorphisme
$\kk=\beta^*(\kk)\twoheadrightarrow\beta^*(\kk^{\geq 1})$ dans $\F_{bi-\Gr}$
dont l'image est nulle sur l'objet $(V,B,W)$ si et seulement si
$\beta(V,B,W)$ est nul, i.e. si $W\subset B$.

On remarque ensuite que le foncteur composé
$\E^f_{bi-\Gr}\xrightarrow{\beta}\E^f_{surj}\hookrightarrow\E^f$ se
plonge dans le foncteur
$\E^f_{bi-\Gr}\xrightarrow{\mathcal{O}_{\E^f_\Gr,\Gr}}\E^f_\Gr\xrightarrow{\mathfrak{K}}\E^f$
(on rappelle que $\mathcal{O}_{\E^f_\Gr,\Gr}$ désigne le foncteur $(V,B,W)\mapsto (V,B)$
--- cf. notation~\ref{notfcom}),
via l'inclusion canonique $\beta(V,B,W)=W/(W\cap B)\hookrightarrow
V/B=\mathfrak{K}(V,B)$. Comme les foncteurs de $\F$ préservent les
monomorphismes, on en déduit que le foncteur de pré\-com\-po\-si\-tion
$\beta^*\circ o : \F\to\F_{bi-\Gr}$ s'injecte dans
$\Upsilon_\Gr\circ\kappa : \F\to\F_{bi-\Gr}$. Le lemme \ref{lmcle}
fournit alors un monomorphisme
$\beta^*(\kk^{\geq 1})\hookrightarrow\beta^*(o(\bar{I}_\kk))\hookrightarrow\Upsilon_\Gr(\kappa(\bar{I}_\kk))$,
ce qui achève la démonstration.
\end{proof}

\paragraph*{Résultats d'annulation cohomologique~\guillemotleft~abstraits~\guillemotright} La
première propriété d'annulation cohomologique de cette section, dont toutes les autres
se déduiront formellement, s'obtient à l'aide du lemme précédent et du
résultat formel donné par le corollaire~\ref{aux-adah}.

\begin{lm}\label{prf-pfo} Soient $A$ un foncteur analytique de
  $\F_\Gr$ et $X$ un objet de $\N$. Le groupe d'extensions ${\rm Ext}^i_{\Gr} (A,\Omega_\Gr
  (X))$ est nul pour tout $i\in\mathbb{N}$.
\end{lm}

\begin{proof}  Un argument de colimite filtrante (cf. \cite{Jen},
  th. 4.2) permet de se ramener au cas où $A$ est
  {\em polynomial} : il existe $d\in\mathbb{N}^*$ tel que $\Delta^d
  A=0$. 

Notons $F$ l'endofoncteur $\cdot\otimes\kappa(\bar{I}_\kk)$ de $\F_\Gr$,
$G=\Delta^\Gr$, $H : \N\to\F_\Gr$ la restriction à $\N$ de
$\Omega_\Gr$, et $K$ l'endofoncteur de $\N$ induit par
$\cdot\otimes\Upsilon_\Gr(\kappa(\bar{I}_\kk))$ (on utilise
que $\N$ est un idéal de $\F_{bi-\Gr}$). Alors :
\begin{enumerate}\item tous ces foncteurs sont exacts ;
\item $G$ est adjoint à gauche à $F$ par le corollaire \ref{cretgr} ;
\item la proposition \ref{preomu} fournit un isomorphisme $F\circ
  H\simeq H\circ K$ ;
\item grâce au lemme \ref{lme33}, il existe un monomorphisme
  $id_\N\hookrightarrow K$. 
\end{enumerate}
Le corollaire \ref{aux-adah} donne alors la conclusion.
\end{proof}

Nous pouvons maintenant énoncer le résultat principal de ce
paragraphe. 

On rappelle que ${\rm Is}_0$ est une autre écriture pour $S^{surj}_0$
(cf. notation \ref{not-prz}).

\begin{pr}\label{anco1} Soient $X$ un foncteur analytique de
  $\F_\Gr$ et $A$ un objet (quelconque) de $\F_{bi-\Gr}$. Le morphisme
gradué  naturel  ${\rm Ext}^*_{\Gr} (X,\Omega_\Gr
  (A\otimes\beta^*({\rm Is}_0)))\to {\rm Ext}^*_{\Gr} (X,\Omega_\Gr
  (A))$ induit par le
  monomorphisme canonique $A\otimes\beta^*({\rm Is}_0)\hookrightarrow A$ fourni par ${\rm Is}_0\hookrightarrow\kk$
  (cf. lemme \ref{str-cs}) est un isomorphisme.
\end{pr}

\begin{proof} On écrit la suite exacte longue de cohomologie associée
  à la suite exacte courte $0\to A\otimes\beta^*({\rm Is}_0)\to A\to
  A\otimes\beta^*(\kk^{\geq 1})\to 0$ et l'on applique le lemme
  précédent à l'objet $A\otimes\beta^*(\kk^{\geq 1})$ de $\N$
  (cf. lemme~\ref{lme33}).
\end{proof}

\subsection{Résultats fondamentaux}\label{sctisof}

Rappelons que le foncteur $\mathfrak{J} : \F_\Gr\to\wt{\F}_\Gr$, l'endofoncteur $\I$ de $\F_\Gr$ et la transformation naturelle $j^\omega :
\I\to\iota\omega$ ont été définis au paragraphe~\ref{substil}.

\begin{theo}\label{prfondo}  Soient $X$ un objet analytique de
$\mathcal{F}_{\mathcal{G}r}$ et $Y$ un objet quelconque de
$\mathcal{F}_{\mathcal{G}r}$. Le morphisme  gradué naturel
$(j^\omega_Y)_* : {\rm
  Ext}^*_\Gr(X,\mathcal{I} (Y))\to {\rm
  Ext}^*_\Gr(X,\iota\omega (Y))$ induit par $j^\omega_Y :
\I(Y)\to\iota\omega (Y)$ est un isomorphisme.

En particulier, les adjonctions des propositions \ref{prfig} et
\ref{pratil} fournissent des isomorphismes gradués naturels 
$${\rm Ext}^*_\F(\omega (X),\omega (Y))\simeq {\rm Ext}^*_\Gr(X,\mathcal{I} (Y))\simeq {\rm Ext}^*_{\widetilde{\mathcal{F}}_{\mathcal{G}r}}(\mathfrak{J} (X),\mathfrak{J} (Y)).$$
\end{theo}

\begin{proof} Notons $\Upsilon'_\Gr : \F_\Gr\to\F_{bi-\Gr}$ la
  composée de $\Upsilon_\Gr$ et du foncteur de précomposition par
  l'endofoncteur $\E^f_{bi-\Gr}$~\guillemotleft~échangeant les deux bases~\guillemotright,
  i.e. donné par $(V,B,W)\to (V,W,B)$ sur les objets et sur les
  morphismes par l'égalité
$${\rm hom}_{\E^f_{bi-\Gr}}((V,B,W),(V',B',W'))={\rm
    hom}_{\E^f_{bi-\Gr}}((V,W,B),(V',W',B')).$$ 
Autrement dit, on a  $\Upsilon'_\Gr(X)(V,B,W)=X(V,W)$ (tandis que $\Upsilon_\Gr(X)(V,B,W)=X(V,B)$).

Alors les endofoncteurs $\Omega_\Gr\Upsilon'_\Gr$ et $\iota\omega$ de
$\F_\Gr$ coïncident, car il existe des isomorphismes canoniques
$$\Omega_\Gr\Upsilon'_\Gr(X)(V,B)\simeq\bigoplus_{W\in\Gr(V)}\Upsilon'_\Gr(X)(V,B,W)$$
$$\simeq\bigoplus_{W\in\Gr(V)}X(V,W)\simeq\omega(X)(V)\simeq\iota\omega(X)(V,B).$$ De plus, il existe un isomorphisme 
$\I(Y)\simeq \Omega_\Gr\Upsilon'_\Gr(Y\otimes\beta^*({\rm Is}_0))$
naturel en $Y$, par lequel $j^\omega_Y$ s'identifie à la
transformation naturelle induite par $Y\otimes\beta^*({\rm
  Is}_0)\hookrightarrow Y$. Ainsi, $(j^\omega_Y)_* : {\rm
  Ext}^*_\Gr(X,\mathcal{I} (Y))\to {\rm
  Ext}^*_\Gr(X,\iota\omega (Y))$ est un
isomorphisme, par la proposition~\ref{anco1}.

 Les autres flèches considérées sont des isomorphismes
pour des raisons formelles (corollaire~\ref{adj-ah}).
\end{proof}

En utilisant l'adjonction entre foncteurs de restriction et de
prolongement par zéro, on en déduit un résultat analogue dans les
  catégories $\F_{\Gr,\leq n}$. Ce même principe est mis en \oe uvre
  dans le corollaire suivant.

\begin{cor}\label{fondcr} Soient $k$ et $n$ deux entiers naturels, $X$
  un objet analytique de $\F_{\Gr,k}$
  et $Y$ un objet quelconque de $\F_{\Gr,n}$.
\begin{enumerate}\item Si $k<n$, alors ${\rm
    Ext}^*_{\mathcal{F}}(\omega_k (X),\omega_n(Y))=0$.
\item Si $k=n$, alors le morphisme  naturel ${\rm
    Ext}^*_{\Gr,n}(X,Y)\to {\rm
    Ext}^*_{\mathcal{F}}(\omega_n (X),\omega_n(Y))$ induit par $\omega_n$ est un isomorphisme.
\end{enumerate}
\end{cor}

On rappelle que les foncteurs de restriction $\mathcal{R}$ et de
prolongement par zéro $\mathcal{P}$ qui interviennent dans la
démonstration ci-dessous ont été définis dans la notation \ref{notglfgr}.

\begin{proof} Comme $\omega_k(X)\simeq\omega\mathcal{P}_{k,\mathbb{N}}(X)$, le théorème \ref{prfondo} fournit un isomorphisme
  naturel gradué entre ${\rm Ext}^*_{\mathcal{F}}(\omega_k
  (X),\omega_n(Y))$ et ${\rm
    Ext}^*_{\Gr}(\mathcal{P}_{k,\mathbb{N}}(X),\I\,\mathcal{P}_{n,\mathbb{N}}(Y))$. On remarque que $\I\,\mathcal{P}_{n,\mathbb{N}}(Y)$ appartient à l'image de la catégorie $\F_{\Gr,\geq n}$ dans $\F_\Gr$ par le foncteur de prolongement par zéro. Comme $\E^f_{\Gr,\geq n}$ est une sous-catégorie complète à gauche de $\E^f_\Gr$, la proposition \ref{prfrec} (et le corollaire \ref{adj-ah}) montrent que ce groupe d'extensions est canoniquement isomorphe à $${\rm Ext}^*_{\Gr,\geq n}(\mathcal{R}_{\mathbb{N},\geq n}\mathcal{P}_{k,\mathbb{N}}(X),\mathcal{R}_{\mathbb{N},\geq n}\,\I\,\mathcal{P}_{n,\mathbb{N}}(Y)).$$

Si $k<n$, $\mathcal{R}_{\mathbb{N},\geq n}\mathcal{P}_{k,\mathbb{N}}(X)$ est nul, ce qui établit la première assertion. Si $k=n$, cet objet s'identifie à  $\mathcal{P}_{n,\geq
  n}(X)$. Puisque $\E^f_{\Gr,n}$ est une sous-catégorie complète à
droite de $\E^f_{\Gr,\geq n}$, la proposition~\ref{prfrec} montre
cette fois que ce groupe d'extensions est canoniquement isomorphe à
${\rm
  Ext}^*_{\Gr,n}(X,\mathcal{R}_{\mathbb{N},n}\,\I\,\mathcal{P}_{n,\mathbb{N}}(Y))$.
On conclut en constatant que le morphisme naturel $Y\to\mathcal{R}_{\mathbb{N},n}\,\I\,\mathcal{P}_{n,\mathbb{N}}(Y)$ est un isomorphisme.
\end{proof}

%

\begin{rem}\label{rq-pira}\begin{enumerate}\item\label{rqitpira} Les premiers résultats d'annulation
    cohomologique dans des catégories de foncteurs remontent à
    Pirashvili. Ainsi, le lemme 0.4 de \cite{FLS} (dû à Pirashvili)
    s'est avéré l'un des premiers outils efficaces pour simplifier des
    calculs de groupes d'extensions dans $\F$.
\item  Le cas $k=0$ du corollaire~\ref{fondcr} constitue une
  généralisation du résultat (énoncé d'ordinaire dans sa variante
  duale) dû à Franjou selon lequel ${\rm Ext}^*_\F(F,\bar{P})=0$ si
  $F$ est un foncteur fini de $\F$ (cf. \cite{GP5}, appendice, pour
  une démonstration due à Schwartz). Les
  lemmes techniques de \cite{GP5} illustrent l'intérêt de ce genre de
  résultat cohomologique pour aborder la conjecture artinienne.

Gaudens et Schwartz ont établi dans~\cite{GS} une généralisation du
résultat de Franjou que le théorème~\ref{prfondo} recoupe
partiellement. 
\end{enumerate}
\end{rem}

\section{Foncteur $\omega$ et  foncteurs hom internes}\label{s-rhid}

Le foncteur $\omega$ étant adjoint à gauche à $\iota$, on dispose pour
des raisons formelles de propriétés de compatibilité entre les
foncteurs de division et $\omega$ --- cf. proposition~\ref{domeg}. En
revanche, l'adjoint à gauche à $\omega$ ne se décrivant pas aisément,
il est plus délicat d'étudier le comportement mutuel des foncteurs hom
internes et $\omega$. Nous utiliserons les foncteurs de
décalage, qui sont à la fois des foncteurs hom internes et des
foncteurs de division (dans $\F_\Gr$ comme dans $\F$), pour contourner
cette difficulté.

\subsection{Scindement de $\Delta_V\circ\omega$}\label{scind-do} Par
les propositions~\ref{domeg} et~\ref{prfif}, il existe des isomorphismes naturels
\begin{equation}\label{scindom-eq}\Delta_V\omega(X)\simeq (\omega(X) :
\iota(I_V))\simeq\underset{W\in\Gr(V)}{\bigoplus}\omega(X : I^\Gr_{(V,W)}),\end{equation}
où l'on a utilisé l'identification de $\Delta_V$ à $(\cdot : I_V)$
(proposition~\ref{adj-fctd}) ; la proposition~\ref{divif} décrit
explicitement ce scindement, qui sous-tend nombre d'aspects de la
structure des catégories $\F_\Gr$ et $\F$ (cf.~\cite{these} et~\cite{art3}).

Précisons l'effet des foncteurs de division $(\cdot : I^\Gr_{(V,W)})$
qui apparaissent relativement au niveau et aux coniveau (notions
introduites dans la définition~\ref{nivcon}). On remarque que ${\rm niv}(X : I^\Gr_{(V,W)})\leq {\rm niv}(X)$, mais
les foncteurs $(\cdot : I^\Gr_{(V,W)})$ abaissent en général le
coniveau. La proposition suivante, que l'on utilisera au paragraphe
suivant, décrit la restriction à $\F_{\Gr,0}\simeq\F$ de $(X :
I^\Gr_{(V,W)})$. Elle fait
usage des foncteurs $\tau_A$ introduits en fin de section~\ref{fhifd}
et résulte de la proposition~\ref{divif} .

\begin{pr}\label{descri-nz} Il existe dans $\F$ un isomorphisme
  naturel $\mathcal{R}_{\mathbb{N},0}(X : I^\Gr_A)\simeq\tau_A(X)$
  pour $X\in {\rm Ob}\,\F_\Gr$ et $A\in {\rm Ob}\,\E^f_\Gr$.
\end{pr}

Signalons enfin que, pour élémentaire qu'il soit, le scindement~(\ref{scindom-eq}) jouera un rôle essentiel dans les arguments de~\cite{art3}.

\subsection{Le morphisme $h^*_{X,F} :
\omega\big(\mathbf{Ext}^*_\Gr(\iota(F),X)\big)\to\mathbf{Ext}^*_\F(F,\omega(X))$}\label{partnhi}

Nous donnons d'abord la définition et quelques propriétés générales
d'un morphisme gradué $h^*_{X,F}$ naturel en un objet $X$ de $\F_\Gr$ et $F$ de
$\F$. Notre objectif principal est d'établir que $h^*_{X,F}$ est un
isomorphisme lorsque $F$ est localement fini.

\begin{defi} Soit $h_{X,F}^0 : \omega\big(\mathbf{Hom}_\Gr(\iota(F),X)\big)\to\mathbf{Hom}_\F(F,\omega(X))$
le morphisme naturel en les objets $F$ de $\F$ et $X$ de $\F_\Gr$ adjoint à
$$\omega\big(\mathbf{Hom}_\Gr(\iota(F),X)\big)\otimes F\xrightarrow{\simeq}\omega\big(\mathbf{Hom}_\Gr(\iota(F),X)\otimes\iota(F)\big)\xrightarrow{\omega(u_{F,X})}\omega(X)\,,$$
où la première flèche est l'isomorphisme de la dernière assertion de
la proposition \ref{prfig} et $u_{F,X}$ la coünité de l'adjonction.
\end{defi}

Comme le foncteur $\mathbf{Hom}_\Gr(\iota(F),\cdot)$ est exact lorsque
$F$ est un objet projectif de $\F$ (par la proposition~\ref{adelgr})
et que le foncteur $\omega$ est exact, ce morphisme naturel s'étend en
un morphisme naturel gradué
$$h^*_{X,F} :
\omega\big(\mathbf{Ext}^*_\Gr(\iota(F),X)\big)\to\mathbf{Ext}^*_\F(F,\omega(X)).$$

\begin{nota}\label{notadhoc} Soit $V$ un espace vectoriel de dimension
  finie. Dans ce paragraphe, nous noterons $\pi_V$ l'épimorphisme
  canonique $\iota(I_V)\twoheadrightarrow\kappa(I_V)$ de $\F_\Gr$.
\end{nota}

\begin{rem}\label{rclhd} L'épimorphisme $\pi_V$ peut se lire comme la
  projection canonique
$$\iota(I_V)\simeq\bigoplus_{W\in\Gr(V)} I^\Gr_{(V,W)}\twoheadrightarrow
I^\Gr_{(V,0)}\simeq\kappa(I_V)$$
(cf. propositions \ref{prfif} et \ref{proj-cfg}) ; il est donc {\em
  scindé}.
\end{rem}

\begin{lm}\label{lmhio} Pour tous objets $V$ de $\E^f$ et $X$ de $\F_\Gr$, le diagramme
$$\xymatrix{\omega\big(\mathbf{Hom}_\Gr(\iota(P_V),X)\big)\ar[r]^-\simeq\ar[d]_{h_{X,P_V}^0}
  & \omega(\Delta^\Gr_V X) \ar[rr]^-\simeq & & \omega(X:\kappa(I_V))\ar@{^{(}->}[d]^{\omega(X:\pi_V)} \\
\mathbf{Hom}_\F(P_V,\omega(X))\ar[r]^-\simeq &
\Delta_V\omega(X)\ar[r]^-\simeq &  (\omega(X):I_V)\ar[r]^-\simeq & \omega(X:\iota(I_V))
}$$
commute, où :
\begin{itemize}\item les flèches supérieures sont les isomorphismes
  fournis par la proposition~\ref{adelgr} ; 
\item les flèches inférieures sont les isomorphismes fournis par les
  propositions~\ref{adj-fctd} et~\ref{domeg}.
\end{itemize}
\end{lm}

\begin{proof} Le morphisme
$$\omega(\Delta^\Gr_V
  X)\simeq\omega\big(\mathbf{Hom}_\Gr(\iota(P_V),X)\big)\xrightarrow{h_{X,P_V}^0}\mathbf{Hom}_\F(P_V,\omega(X))\simeq\Delta_V\omega(X)$$
est l'adjoint du morphisme $P_V\otimes\omega(\Delta^\Gr_V X)\to\omega(X)$ donné, sur l'espace
vectoriel $A$, par $[f]\otimes a\mapsto X(i_f)(a)$ pour $f\in {\rm
  hom}_\E(V,A)$, $B\in\Gr(A)$ et $a\in X(A,B)$, où l'on désigne par
$i_f : (V\oplus A,B)\to (A,B)$ le morphisme de composante $f$ sur $V$
et $id_A$ sur~$A$. 

Par conséquent, la flèche $$\omega(\Delta^\Gr_V
  X)\simeq\omega\big(\mathbf{Hom}_\Gr(\iota(P_V),X)\big)\xrightarrow{h_{X,P_V}^0}\mathbf{Hom}_\F(P_V,\omega(X))\simeq\Delta_V\omega(X)$$ est donnée sur l'espace $A$ par l'inclusion 
$$\omega(\Delta^\Gr_V
  X)(A)=\bigoplus_{B\in\Gr(A)}X(V\oplus
  A,B)\hookrightarrow\bigoplus_{B\in\Gr(V\oplus A)}X(V\oplus A,B)=\Delta_V\omega(X)(A).$$ 

La conclusion s'obtient alors par la proposition~\ref{divif}.
\end{proof}

\begin{pr}\label{prinjdo} Pour tous objets $F$ de $\F$ et $X$ de
  $\F_\Gr$, le morphisme $h_{X,F}^0$ est injectif.
\end{pr}

\begin{proof} Le lemme~\ref{lmhio} montre l'assertion lorsque $F$ est
  un projectif standard $P_V$.  On en déduit (par commutation de $\omega$ aux
  limites et de $\iota$ aux colimites) que $h_{X,P}^0$ est également
  injectif lorsque $P$ est une somme directe de projectifs
  standard. Le cas général s'en déduit en considérant un épimorphisme $p$
  d'un tel projectif $P$ sur $F$ et en considérant le diagramme
  commutatif
$$\xymatrix{\omega\big(\mathbf{Hom}_\Gr(\iota(F),X)\big)\ar@{^{(}->}[r]^-{p^*}\ar[d]_{h_{X,F}^0} & \omega\big(\mathbf{Hom}_\Gr(\iota(P),X)\big)\ar@{^{(}->}[d]^{h_{X,P}^0}\\
\mathbf{Hom}_\F(F,\omega(X))\ar[r]^-{p^*} & \mathbf{Hom}_\F(P,\omega(X)).
}$$
\end{proof}

\begin{rem}\label{rqos} Il existe même un endofoncteur $T_F$ de $\F_\Gr$, un
  isomorphisme naturel
  $\mathbf{Hom}_\F (F,\omega(X))\simeq\omega(T_F(X))$ et une
  transformation naturelle injective $j_F : \mathbf{Hom}_\Gr (\iota(F),
  \cdot)\hookrightarrow T_F$ telle que $h^0_{\cdot,F}$ s'identifie à
  $\omega(j_F)$. On le voit par un argument analogue, en considérant une
  présentation de $F$ par des sommes directes de projectifs
  standard. Ce qui suit précise ces considérations.
\end{rem}

\paragraph*{Description explicite du morphisme $h^*_{X,F}$} Soit $V\in
{\rm Ob}\,\E^f$. On a d'une part
$$\omega\big(\mathbf{Ext}^*_\Gr(\iota(F),X)\big)(V)\simeq\bigoplus_{W\in\Gr(V)}\mathbf{Ext}^*_\Gr(\iota(F),X)(V,W)\simeq\bigoplus_{W\in\Gr(V)}{\rm
  Ext}^*_\F(F,\tau_{(V,W)}(X))$$
par la proposition~\ref{extint}.

D'autre part,
$$\mathbf{Ext}^*_\F(F,\omega(X))(V)\simeq{\rm Ext}^*_\F(F\otimes
P_V,\omega(X))\simeq{\rm Ext}^*_\F(F,\Delta_V\omega(X))$$
$$\simeq\bigoplus_{W\in\Gr(V)}{\rm
  Ext}^*_\F(F,\omega(X:I^\Gr_{(V,W)}))$$
grâce à l'isomorphisme~(\ref{scindom-eq}) et à la
proposition~\ref{prdf-hid}.

\begin{pr}\label{ident-mh} Via les isomorphismes précédents et la proposition~\ref{divif}, le morphisme $h^*_{X,F}$ est induit par
  l'inclusion naturelle
$$\tau_{(V,W)}(X)(E)=X(E\oplus
V,W)\hookrightarrow\underset{B\in Gr(E,A)}{\bigoplus_{A\in\Gr(W)}} X(E\oplus V,B)=\omega(X:I^\Gr_{(V,W)})(E) .$$
\end{pr}

\begin{proof} Par naturalité, il suffit de vérifier cette
  identification lorsque $F$ est un projectif standard de $\F$, auquel
  cas elle résulte du lemme~\ref{lmhio}.
\end{proof}

Les préliminaires précédents nous permettent, à l'aide du théorème
fondamental de la section précédente, d'établir le résultat principal
de cette section.

\begin{pr}\label{hl1} Si $F$ est un objet localement fini de $\F$, le
  morphisme naturel gradué
$$h^*_{X,F} :
\omega\big(\mathbf{Ext}^*_\Gr(\iota(F),X)\big)\to\mathbf{Ext}^*_\F(F,\omega(X))$$
est un isomorphisme pour tout objet $X$ de $\F_\Gr$.
\end{pr}

\begin{proof} Notons $i_E^X(V,W)$ l'inclusion naturelle de la
  proposition~\ref{ident-mh} : l'énoncé équivaut à l'annulation des
  groupes d'extensions ${\rm Ext}^*_\F(F,\omega(coker\,i_E^X))$ pour
  tout $E\in {\rm Ob}\,\E^f$. La proposition~\ref{descri-nz} montre
  que l'image par le foncteur de restriction
  $\mathcal{R}_{\mathbb{N},0}$ du foncteur $coker\,i_E^X$ de $\F_\Gr$
  est nulle. Le théorème~\ref{prfondo} (ou son
  corollaire~\ref{fondcr}) donne alors la conclusion.
\end{proof}

Dans les conséquences qui suivent, nous nous contentons du cas du
degré $0$ (i.e. des foncteurs hom internes), le plus
significatif.

\begin{cor}\label{vind2} Soient $F$ un objet localement fini de $\F$, $n$ et $k$ des entiers naturels et $I$ une partie de $\mathbb{N}$ du type
  $\leq n$ ou $n$. Les
  foncteurs $\omega_I\circ\mathbf{Hom}_{\Gr,I}
  (\iota_I(F),\cdot)$ et $\mathbf{Hom}_\F
  (F,\cdot)\circ\omega_I$ de $\F$ vers $\F_{\Gr,I}$ sont isomorphes.
\end{cor}

\begin{proof} Par la proposition~\ref{hl1}, il existe un morphisme 
$$\mathbf{Hom}_\F (F,\cdot)\circ\omega_I\simeq\mathbf{Hom}_\F (F,\cdot)\circ\omega\circ\mathcal{P}_{I,\mathbb{N}}\simeq\omega\circ\mathbf{Hom}_\Gr
  (\iota(F),\cdot)\circ\mathcal{P}_{I,\mathbb{N}}.$$

Le corollaire s'en déduit parce que $\mathbf{Hom}_\Gr
  (\iota(F),\cdot)\circ\mathcal{P}_{I,\mathbb{N}}\simeq\mathcal{P}_{I,\mathbb{N}}\circ\mathbf{Hom}_{\Gr,I}
  (\iota_I(F),\cdot)$, par la proposition~\ref{pcorp}.
\end{proof}

Nous terminons cette section en donnant un corollaire des résultats
précédents qui décrit notamment l'image par un foncteur
$\mathbf{Hom}_\F(F,\cdot)$ (où $F\in\F_\omega$) d'un $\bar{G}(n)$-comodule
simple (cf. propositions~\ref{prfig} et~\ref{prfffgr}).

\begin{cor}\label{crfph} Soient $n$ et $k$ des entiers naturels, $F$
  un objet localement fini de $\F$, $G$ un objet de $\F$ et $M$ un
  $\kk[GL_n(\kk)]$-module fini. Il existe dans $\F$ un isomorphisme naturel
$$\mathbf{Hom}_\F\big(F,\omega_n\big(\kappa_n(G)\otimes\rho_n(M)\big)\big)\simeq
\omega_n\big(\kappa_n(\mathbf{Hom}_\F(F,G))\otimes\rho_n(M)\big).$$
\end{cor}

\begin{proof} Ce résultat s'obtient en combinant les corollaires~\ref{hiki}, \ref{crcomdh} et~\ref{vind2}.
\end{proof}

\paragraph*{Retour sur la méthode de \cite{art1}} L'article
\cite{art1} étudie, dans le cas $\kk=\FF$, la structure de foncteurs
de type fini\,\footnote{En fait, cet article travaille sur les objets
  duaux, de co-type fini. Nous avons traduit, pour la cohérence de
  l'exposition, ses énoncés en termes de foncteurs de type fini.} de
$\F(\FF)$ (précisément, le produit tensoriel entre le projectif
$P_{\FF^2}$ et une puissance extérieure) par un argument de récurrence utilisant le foncteur $\mathbf{Hom}_{\F(\FF)}(\Lambda^1,\cdot)$, où $\Lambda^1$
désigne le foncteur d'inclusion $\E^f\hookrightarrow\E$. 

Le corollaire~\ref{crfph} montre, en particulier, que les $\bar{G}(n)$-comodules finis sont nilpotents pour
le foncteur $\mathbf{Hom}_{\F(\FF)}(\Lambda^1,\cdot)$
(cf. proposition~\ref{pnil-hi} ci-après). Cet énoncé sous-tend toute la démarche de~\cite{art1}
; cependant, sans le formalisme des catégories de foncteurs en
grassmanniennes, il est malaisé d'en donner une forme et une
démonstration générales. La proposition~$5.29$ de~\cite{art1}, qui
fournit une estimation de l'image par le foncteur
$\mathbf{Hom}_{\F(\FF)}(\Lambda^1,\cdot)$ des foncteurs étudiés dans l'article en question, sert de
succédané élémentaire au corollaire~\ref{crfph}. Dans le dernier
chapitre de~\cite{these}, nous étendons la méthode de~\cite{art1}
en employant le corollaire~\ref{crfph} de façon beaucoup plus générale.

\section{La filtration de Krull de la catégorie~$\F$}\label{s-fkf}

La conjecture que nous introduisons dans ce paragraphe, appelée {\em
conjecture artinienne extrêmement forte} (elle implique la conjecture
artinienne très forte~\ref{ca4}), affirme que
la restriction  aux objets localement finis du
foncteur d'intégrale en grassmanniennes $\omega : \F_\Gr\to\F$ envoie
la filtration de $\F_\Gr^{lf}$ donnée par les sous-catégories $\F^{lf}_{\Gr,\leq
  n}$ dans la filtration de Krull de $\F$, et induit une équivalence
entre les quotients associés. Après des préliminaires présentés dans
le paragraphe~\ref{pusfkf}, nous exposons les différentes formes de
la conjecture artinienne extrêmement forte et les cas particuliers
qu'on peut en établir (§\,\ref{sdsfkf}) puis en examinons des
conséquences importantes (§\,\ref{ptsfkf}).

\subsection{Foncteurs oméga-adaptés}\label{pusfkf}  Nous nous intéresserons à la
propriété suivante d'une sous-catégorie pleine~$\A$ de $\F$ : 
\begin{hyp}\label{semep} Pour toute suite exacte courte $0\to A\to B\to
  C\to 0$ de $\F$, si deux des objets $A, B, C$ appartiennent à $\A$, il en
  est de même du troisième.
\end{hyp}

Conjecturalement, les foncteurs oméga-adaptés de hauteur au plus $n$
introduits ci-après sont exactement les foncteurs noethériens de type
$n$ de~$\F$. L'étude de ceux-là se ramène essentiellement à celle
d'objets {\em finis} de $\F_\Gr$, ce qui les rend assez facilement maniables.

\begin{defi}\label{dffog} Soit $n\in\mathbb{Z}$. Nous noterons $\F^{\omega - ad(n)}$ la plus petite sous-catégorie pleine de $\F$ vérifiant
  \ref{semep} et contenant l'image de la restriction à la
  sous-catégorie $\F^f_{\Gr,\leq n}$ des objets finis de $\F_{\Gr,\leq n}$ du foncteur $\omega_{\leq n} : \F_{\Gr,\leq n}\to\F$.

Nous dirons que $F$ est {\em oméga-adapté de hauteur au plus $n$} s'il est objet
  de $\F^{\omega - ad(n)}$.
\end{defi}

\begin{pr}\label{pre-omd}\begin{enumerate}\item\label{apr4} Le produit
    tensoriel d'un objet de $\F^{\omega - ad(n)}$ et d'un objet de
    $\F^{\omega - ad(m)}$ est un objet de $\F^{\omega - ad(n+m)}$.
\item\label{apr6} Les sous-catégories $\F^{\omega-ad(n)}$ sont stables par le foncteur différence de $\F$.
\item\label{apr5} Tout foncteur oméga-adapté est pf$_\infty$.
\end{enumerate}
\end{pr}

\begin{proof} Si $X$ est un objet fini de $\F_{\Gr,\leq k}$, alors le
produit tensoriel de $\omega_{\leq k}(X)$ et d'un objet oméga-adapté de
hauteur au plus~$n$ est oméga-adapté de hauteur au plus~$n+k$. En
effet, la sous-catégorie pleine des objets de $\F$ dont le produit
tensoriel par $\omega_{\leq k}(X)$ est oméga-adapté de hauteur au plus $n+k$ contient
les $\omega_{\leq n}(A)$, pour $A\in {\rm Ob}\,\F^f_{\Gr,\leq n}$,  par la
proposition~\ref{comom}, et vérifie
l'hypothèse~\ref{semep}. On démontre ensuite par le même raisonnement
l'assertion~\ref{apr4}.

La proposition~\ref{prfif} montre, par exactitude de
$\omega$, que $\omega$ transforme un objet pf$_\infty$, donc en
particulier un objet fini (cf. corollaire~\ref{crfff3}), en un objet
pf$_\infty$. Comme la sous-catégorie pleine des objets pf$_\infty$ de
$\F$ vérifie l'hypothèse~\ref{semep}, par la proposition~\ref{pr:pfi}, on en déduit l'assertion~\ref{apr5}.

Pour l'assertion \ref{apr6}, on considère la sous-catégorie pleine
$\C_n$ des objets $F$ de $\F$ tels que $\Delta F$ appartient à
$\F^{\omega - ad(n)}$. Les propositions~\ref{domeg} et~\ref{divif} montrent
que $\C_n$ contient $\omega_{\leq n}(X)$ pour $X\in {\rm
  Ob}\,\F^f_{\Gr,\leq n}$. D'autre part, l'exactitude du foncteur
différence montre que $\C_n$ vérifie l'hypothèse~\ref{semep}, ce qui achève la démonstration.\end{proof}

\begin{ex}\label{exidtr} Les projectifs standard
  $P_V=\omega(P^\Gr_{(V,V)})=\omega(\rho(P^{surj}_V))$ de $\F$ sont
  oméga-adaptés, de hauteur $\dim V$ (le fait que cette hauteur n'est
  pas strictement inférieure à $\dim V$ n'est pas tout à fait immédiat
  ; il découlera des considérations
  de~\cite{art3}).
\end{ex}

Nous utiliserons, dans cette section, le théorème~\ref{prfondo} par
l'intermédiaire du résultat d'annulation cohomologique suivant.

\begin{pr}\label{prhom} Soient $n\in\mathbb{N}$, $X$ un objet de
  $\F_{\Gr,\geq n}$ et $F$ un objet de $\F^{\omega-ad(n-1)}$. On a ${\rm Ext}^*_\F(F,\omega_{\geq n}(X))=0$.
\end{pr}

\begin{proof} Cette propriété s'obtient à partir du corollaire~\ref{fondcr} et de l'observation que la sous-catégorie pleine de
  $\F$ formée des objets $F$ tels que ${\rm Ext}^*_\F(F,\omega_{\geq
    n}(X))=0$ vérifie l'hypothèse~\ref{semep}.
\end{proof}

\subsection{La conjecture artinienne extrêmement forte}\label{sdsfkf}

Nous introduisons la conjecture artinienne extrêmement forte sous une
forme globale utilisant la notion de foncteur oméga-adapté.

\begin{conj}[Conjecture artinienne extrêmement forte]\label{caef} Pour tout $n\in\mathbb{N}$, un quotient d'un
  foncteur oméga-adapté de hauteur au plus $n$ est oméga-adapté de
  hauteur au plus~$n$.
\end{conj}

\begin{rem}\label{rqcafbl} Cet énoncé  implique
déjà la conjecture artinienne sous sa forme minimale en raison de l'assertion~\ref{apr5} de la
proposition~\ref{pre-omd}.
\end{rem}

La proposition suivante lie la conjecture artinienne extrêmement forte
à la filtration de Krull de~$\F$.

\begin{pr}\label{prvcaf} Supposons la conjecture~\ref{caef} vérifiée.
\begin{enumerate}\item\label{itca1} Pour tout $n\in\mathbb{N}$, la
  sous-catégorie $\F^{\omega - ad(n)}$ de $\F$ est épaisse. Elle
  est égale à la sous-catégorie $\F_{\mathbf{NT}(n)}$ des objets
  noethériens de type $n$ de $\F$.
\item\label{itca2} De plus, le foncteur $\omega_n : \F_{\Gr,n}\to\F$
  induit une équivalence de catégories entre $\F^f_{\Gr,n}$ et $\F_{\mathbf{NT}(n)}/\F_{\mathbf{NT}(n-1)}$.
\item\label{itca3} Désignons par $\overline{\F^{\omega - ad(n)}}$ la plus petite sous-catégorie épaisse de $\F$ stable
par colimites contenant  $\F^{\omega - ad(n)}$. C'est aussi la
catégorie des foncteurs qui sont colimite de leurs sous-foncteurs
noethériens de type $n$. Alors la filtration de Krull de $\F$ est
donnée par $\K_n(\F)=\overline{\F^{\omega - ad(n)}}$, et le foncteur
$\omega_n$ induit une équivalence de catégories entre $\F^{lf}_{\Gr,n}$ et $\K_n(\F)/\K_{n-1}(\F)$.
\end{enumerate}
\end{pr}

\begin{proof} Les sous-catégories $\F^{\omega - ad(n)}$ de $\F$
  sont par hypothèse stables par quotients et elles
  vérifient~\ref{semep}, elles sont donc épaisses. La description de $\overline{\F^{\omega - ad(n)}}$ comme sous-catégorie pleine des objets colimite de leurs
sous-objets oméga-adaptés de hauteur au plus $n$ vient de ce que les
foncteurs oméga-adaptés sont de présentation finie (assertion~\ref{apr5} de la
proposition~\ref{pre-omd}) et de la proposition~\ref{pr:pf}.

 On note ensuite que $\F^{\omega - ad(n)}$ est l'image réciproque par le foncteur
canonique $\F\twoheadrightarrow\F/\F^{\omega - ad(n-1)}$ de l'image, notée $\C_n$,
du foncteur $\F^f_{\Gr,n}\to\F/\F^{\omega - ad(n-1)}$ induit par
$\omega_n$. En effet, la sous-catégorie des objets de $\F$ dont
l'image dans $\F/\F^{\omega - ad(n-1)}$ appartient à $\C_n$ vérifie
\ref{semep} et contient bien les $\omega_i(X)$ pour $i\leq n$ et $X\in
{\rm Ob}\,\F^f_{\Gr,i}$, donc elle contient  $\F^{\omega -
  ad(n)}$. Réciproquement, si $F$ est un objet de $\F$ isomorphe à
$\omega_n(X)$ dans $\F/\F^{\omega - ad(n-1)}$, avec $X\in
{\rm Ob}\,\F^f_{\Gr,n}$, la proposition~\ref{prhom}
prouve\,\footnote{Utiliser les résultats de base sur les catégories
  abéliennes quotients donnés dans \cite{Gab}.} qu'il existe dans $\F$ un morphisme $F\to\omega_n(X)$ dont le
noyau et le conoyau appartiennent à $\F^{\omega - ad(n-1)}$, ce qui
entraîne que $F$ appartient à $\F^{\omega - ad(n)}$.

On établit les autres résultats par récurrence sur
  l'entier $n$. Ils sont clairs pour $n=0$. On suppose donc $n>0$ et la
proposition vérifiée au rang $n-1$.

Le foncteur $\omega_n$ induit une équivalence entre
$\F_{\Gr,n}^f$ et la sous-catégorie {\em épaisse} $\F^{\omega -
  ad(n)}/\F^{\omega - ad(n-1)}$ de $\F/\F^{\omega - ad(n-1)}$, et
aussi de $\F_{\Gr,n}^{lf}$ vers la sous-catégorie épaisse
$\overline{\F^{\omega - ad(n)}}/\overline{\F^{\omega - ad(n-1)}}$ de
$\F/\overline{\F^{\omega - ad(n-1)}}$. Cela montre en particulier, compte-tenu de l'hypothèse
de récurrence, que ces images sont respectivement incluses dans
$\F_{\mathbf{NT}(n)}/\F_{\mathbf{NT}(n-1)}$ et
$\K_n(\F)/\K_{n-1}(\F)$.

Si $X$ est un objet fini non nul de $\F_{\Gr,k}$, avec $k>n$, alors $\omega_k(X)$
est limite de ses quotients appartenant à  $\F^{\omega - ad(n)}$ et n'y
appartient pas lui-même, par la proposition~\ref{prhom}. On
en déduit que son image dans $\F/\F^{\omega - ad(n-1)}$ est infinie,
donc que $\omega_k(X)$ n'est pas noethérien de type $n$, puisque
$\F^{\omega - ad(n-1)}=\F_{\mathbf{NT}(n-1)}$ par
l'hypothèse de récurrence. Cela entraîne qu'il n'y a pas, dans $\F$,  d'objet
noethérien de type $n$ qui ne soit pas oméga-adapté de hauteur au plus
$n$. On a donc $\F_{\mathbf{NT}(n)}=\F^{\omega - ad(n)}$.

Comme la catégorie $\F$ est localement noe\-théri\-en\-ne, par la
remarque~\ref{rqcafbl}, la proposition~\ref{prkrull} donne $\K_n(\F)=\overline{\F_{\mathbf{NT}(n)}}=\overline{\F^{\omega - ad(n)}}$, d'où la proposition.
\end{proof}

Les deux énoncés qui suivent constituent les versions~\guillemotleft~locales~\guillemotright~de la conjecture artinienne extrêmement forte.

\begin{conj}\label{caef2} Soit $n\in\mathbb{N}$. Le foncteur exact $\F_{\Gr,n}\xrightarrow{\omega_n}\F\twoheadrightarrow\F/\K_{n-1}(\F)$
induit un isomorphisme $G^f_0(\F_{\Gr,n})\xrightarrow{\simeq}G^f_0(\F/\K_{n-1}(\F))$
qui préserve les classes des objets simples. 
\end{conj}

\begin{conj}\label{caef3} Soit $n\in\mathbb{N}$. Pour tout objet simple
  $S$ de $\F_{\Gr,n}$, le foncteur $\omega_n(S)$ de $\F$ est simple noethérien
  de type $n$. De plus, un foncteur de $\F$ est
  simple noethérien de type $n$ si et seulement s'il est isomorphe à
  un quotient non nul d'un tel objet.
\end{conj}

La proposition~\ref{prvcaf} fournit l'équivalence entre les
différentes conjectures de ce paragraphe.

\begin{cor}\label{eqcaf} Les conjectures~\ref{caef},~\ref{caef2}
  et~\ref{caef3} sont
  équivalentes.
\end{cor}

\begin{proof}La proposition~\ref{prvcaf} montre que
  la conjecture~\ref{caef} implique la conjecture~\ref{caef3}. 

La proposition~\ref{prkrull} montre que la conjecture~\ref{caef3} entraîne la conjecture~\ref{caef2}.

Supposons maintenant la conjecture~\ref{caef2} satisfaite. Alors tout foncteur
oméga-adapté de hauteur au plus $n$ appartient à $\K_n(\F)$. On montre par
récurrence sur $n\in\mathbb{N}$ que la sous-catégorie $\F^{\omega -
  ad(n)}$ de $\F$ est épaisse et que les objets de $\K_n(\F)$ sont
colimite de sous-objets oméga-adaptés de hauteur au plus $n$. Si cette assertion
est vérifiée, tout quotient strict d'un foncteur du type $\omega_{n+1}(S)$,
où $S\in {\rm Ob}\,\F_{\Gr,n+1}$ est simple, est oméga-adapté de
hauteur au plus $n$ : un tel quotient appartient à $\K_n(\F)$ (parce que la
conjecture~\ref{caef2} est satisfaite), et est de type fini. En
particulier, tous les quotients de $\omega_{n+1}(S)$ sont dans $\F^{\omega -
  ad(n+1)}$. On en
déduit que la sous-catégorie $\F^{\omega -
  ad(n+1)}$ de $\F$ est épaisse. Son image dans le quotient
$\F/\K_n(\F)$ contient les objets simples de cette catégorie (car la
conjecture~\ref{caef2} est satisfaite) ; elle est épaisse. Le fait que les foncteurs oméga-adaptés sont de présentation
finie implique, par la proposition~\ref{pr:pf}, que la
sous-catégorie des foncteurs de $\F$ qui sont colimite de leurs
sous-foncteurs oméga-adaptés de hauteur au plus $n+1$ est aussi épaisse. On
en déduit que cette sous-catégorie coïncide avec $\K_{n+1}(\F)$, par
définition de la filtration de Krull. Cela établit la
conjecture~\ref{caef} et achève la démonstration. 
\end{proof}

En considérant les objets simples pseudo-constants des $\F_{\Gr,n}$
dans la conjecture~\ref{caef3}, on obtient le résultat suivant.

\begin{cor}\label{creqcama} La conjecture artinienne extrêmement forte
  implique la conjecture artinienne très forte~\ref{ca4}.
\end{cor}

\paragraph*{Résultats partiels sur la filtration de Krull de~$\F$}
Dans~\cite{art3}, nous établirons les résultats suivants. Le
premier démontre la moitié de la conjecture~\ref{caef3} pour $n=1$ et
$\kk=\FF$ ;
le second, qui généralise l'important théorème de simplicité de Powell
(cf.~\cite{GP2}), consistue une forme affaiblie de cette conjecture
valable pour tout~$n$.

\begin{theo}\label{tha31} Le foncteur $\omega_1 : \F_{\Gr,1}(\FF)\to\F(\FF)$
  induit une équivalence entre la sous-catégorie $\F^f_{\Gr,1}(\FF)$ des
  objets finis de $\F_{\Gr,1}(\FF)$ et une sous-catégorie
  épaisse de $\F(\FF)/\F_\omega(\FF)$. 
\end{theo}

Dans le théorème~\ref{tha32} ci-dessous, les endofoncteurs $\tilde{\nabla}_n$ de
$\F(\FF)$ sont les
{\em duaux} de ceux introduits par Powell
dans~\cite{GP4}. Explicitement, on a
$$\tilde{\nabla}_n(F)(V)=im\,\big(F(V\oplus\FF^{\oplus
  n})\xrightarrow{\sum_{l\in (\FF^{\oplus n})^*}F(V\oplus l)}F(V\oplus\FF)\big).$$
On désigne par
$\overline{\N il}_{\tilde{\nabla}_n}$  la plus petite sous-catégorie épaisse
stable par colimites de $\F$ contenant les foncteurs
$\tilde{\nabla}_n$-nilpotents.

\begin{theo}\label{tha32} Pour tout $n\in\mathbb{N}^*$, le foncteur $\omega_n : \F_{\Gr,n}(\FF)\to\F(\FF)$
  induit une équivalence entre la sous-catégorie $\F^f_{\Gr,n}(\FF)$ des
  objets finis de $\F_{\Gr,n}(\FF)$ et une sous-catégorie
  épaisse de $\F(\FF)/\overline{\N il}_{\tilde{\nabla}_n}$.
\end{theo}

L'intérêt de ce théorème réside dans le fait qu'un foncteur
oméga-adapté de hauteur strictement inférieure à $n$ est
$\tilde{\nabla}_n$-nilpotent.

Nous déduirons du théorème~\ref{tha32} le résultat suivant.

\begin{theo} Pour tout foncteur fini $F$ de $\F(\FF)$, le foncteur
  $P_{\FF}^{\otimes 2}\otimes F$ est noethérien de type~$2$.
\end{theo}

\subsection{Conséquences de la conjecture artinienne
  extrêmement forte}\label{ptsfkf}
Nous donnons des propriétés des foncteurs oméga-adaptés qui montrent que la conjecture artinienne
  extrêmement forte implique des conjectures profondes sur la structure
de la catégorie $\F$ que les formes plus faibles de la conjecture
artinienne ne suffisent pas à résoudre.

\begin{pr}\label{pnil-hi}  Soient $k\in\mathbb{N}$ et $F$ un foncteur de oméga-adapté de
  $\F$. Il existe $N\in\mathbb{N}$ tel que pour tout $n\geq N$  et
  tout foncteur analytique $A$  de $\F$ tel que $A(0)=0$, on ait $\mathbf{Ext}^*_\F(A^{\otimes n},F)=0$.
\end{pr}

\begin{proof} Comme la sous-catégorie pleine des foncteurs $F$ de $\F$
  tels que $\mathbf{Ext}^*(A^{\otimes n},F)=0$ pour $n$ assez grand
  vérifie l'hypothèse~\ref{semep}, il suffit de démontrer
  l'assertion lorsque $F=\omega(X)$, où $X$ est un objet fini de $\F_\Gr$.

Comme $A(0)=0$, le foncteur $A$ admet une résolution projective dont
les termes sont du type $\bar{P}_\kk\otimes B_i$, donc $A^{\otimes n}$  admet une résolution projective dont
les termes sont du type $\bar{P}_\kk^{\otimes n}\otimes C_i$. Pour tout entier $n>\deg(X)$ et tout $E\in {\rm Ob}\,\E^f_\Gr$,
$${\rm hom}_\Gr(\iota(\bar{P}_\kk^{\otimes n}\otimes C_i)\otimes
P^\Gr_E,X)\simeq {\rm hom}_\Gr(\iota(C_i)\otimes
P^\Gr_E,(\Delta^\Gr)^n X)=0$$ (par le corollaire~\ref{cretgr}), donc
${\rm Ext}^*_\Gr(\iota(A^{\otimes n})\otimes
P^\Gr_E,X)=0$,
puisque $(\iota(\bar{P}_\kk^{\otimes n}\otimes
C_i)\otimes P^\Gr_E)_i$ est une
résolution projective de $\iota(A^{\otimes n})\otimes
P^\Gr_E$. L'isomorphisme (\ref{eq-hi}) de la proposition~\ref{prdf-hid}
fournit alors $\mathbf{Ext}^*_\Gr(\iota(A^{\otimes n}),X)=0$ pour $n>\deg(X)$.

La conclusion résulte donc de la proposition~\ref{hl1}.
\end{proof}

\begin{rem}\begin{enumerate}\item L'argument  d'annulation
    cohomologique utilisé dans cette démonstration, analogue à la
    proposition~$1.5.1$ de~\cite{Franjou}, est dû originellement à
    Pirashvili (cf. remarque~\ref{rq-pira}.\,\ref{rqitpira}).
\item Cette proposition est surtout significative en degré $0$,
  i.e. pour les foncteurs hom internes. En particulier, elle illustre
  le comportement radicalement différent du foncteur différence
  $\Delta\simeq\mathbf{Hom}_\F(\bar{P}_\kk,\cdot)$ et de
  $\mathbf{Hom}_\F(F,\cdot)$ pour $F$ fini.
\item Pour $\kk=\FF$, le cas le plus intéressant est celui où
  $A=\Lambda^1$ (cf.~\cite{art1}).
\end{enumerate}
\end{rem}

\begin{pr} Si $F$ est un foncteur oméga-adapté de $\F$,
  pour tout $i\in\mathbb{N}$, l'ensemble des classes d'isomorphisme
  d'objets simples $S$ de $\F$ tels que ${\rm Ext}^i_\F (S,F)\neq 0$ est fini.
\end{pr}

\begin{proof} Il suffit de montrer la propriété pour $F$ de la
  forme $\omega_n(X)$, où $n\in\mathbb{N}$ et $X\in {\rm
    Ob}\,\F^f_{\Gr,n}$. Pour $n=0$, cela vient de ce que les objets
  finis de $\F$ sont co-pf$_\infty$ (cf.~\cite{these} pour les
  détails) ; pour $n>0$, tous les groupes
${\rm Ext}^i_\F (S,\omega_n(X))$ sont nuls par le
  corollaire~\ref{fondcr}, d'où la proposition.
\end{proof}

Cet énoncé est à comparer à la forme de la conjecture artinienne
donnée par l'assertion~\ref{ascae} de la proposition~\ref{cael1}.

\paragraph*{Le groupe de Grothendieck $G_0^{tf}(\F)$} Le foncteur
exact $\omega$ préserve les objets de type fini, il induit donc un
morphisme de groupes $\omega_* : G^f_0(\F_\Gr)\to G_0^{tf}(\F)$, qui est un
morphisme d'anneaux lorsqu'on munit la source de la structure
multiplicative induite par le produit tensoriel total
(cf. proposition~\ref{comom}). Si $F$ est un foncteur oméga-adapté, la
classe de $F$ dans $G_0^{tf}(\F)$ appartient à l'image de $\omega_*$,
donc la conjecture artinienne extrêmement forte implique que ce
morphisme est {\em surjectif}. On peut compléter cette remarque comme
suit.

\begin{nota} Soient $\mathcal{S}$ un système complet de représentants
  des objets simples de $\F$ et $\widehat{G}^f_0(\F)$ le groupe
  produit $\mathbb{Z}^\mathcal{S}$. Le morphisme de groupes canonique
  $j : G_0^{tf}(\F)\to\widehat{G}^f_0(\F)$ s'obtient en
  associant à un objet de type fini $F$ de $\F$ la famille des
  multiplicités d'un élément $S$ de $\mathcal{S}$ dans~$F$.
\end{nota}

Dans~\cite{art3}, nous établirons le résultat suivant.

\begin{theo} Le morphisme
  $G^f_0(\F_\Gr(\FF))\xrightarrow{\omega_*}G_0^{tf}(\F(\FF))\xrightarrow{j}\widehat{G}^f_0(\F(\FF))$ est injectif.
\end{theo}

\begin{cor} Si la conjecture artinienne extrêmement forte pour $\FF$
  est vérifiée, alors le morphisme $\omega_* : G^f_0(\F_\Gr(\FF))\to
  G_0^{tf}(\F(\FF))$ est un isomorphisme, et le morphisme canonique
  $j : G_0^{tf}(\F(\FF))\to\widehat{G}^f_0(\F(\FF))$ est un monomorphisme.
\end{cor}

\section{Résultats d'annulation cohomologique dans
  $\F_{inj}$}\label{parfij}

La comparaison entre les groupes d'extensions dans les catégories $\F$
et des modules sur un groupe linéaire s'opère naturellement par
l'intermédiaire de la catégorie $\F_{inj}$ (ou $\F_{surj}$) ; on peut l'illustrer par
le diagramme commutatif
$$\xymatrix{\F\ar[r]^{o_{inj}}\ar[d]_{{\rm ev}_{E_n}} & \F_{inj}\ar[d]^{{\rm ev}_n} \\
_{\kk[\mathcal{M}_n(\kk)]}\mathbf{Mod}\ar[r] & _{\kk[GL_n(\kk)]}\mathbf{Mod}
}$$
dans lequel la flèche inférieure est le foncteur de restriction : il
s'agit d'étudier le comportement cohomologique du foncteur composé
$\F\to\,_{GL_n(\kk)}\mathbf{Mod}$ de ce diagramme.

Ce principe est implicite dans la démonstration de Suslin, donnée dans
l'appendice de \cite{FFSS}, du théorème selon lequel ce foncteur induit un
isomorphisme entre les groupes d'extensions entre deux foncteurs {\em
  finis} de la catégorie $\F$, pourvu que $n$ soit assez grand.

Le résultat principal de cette section constitue une généralisation de
ce théorème. Il est démontré dans le paragraphe~\ref{par-inj-f}, à
partir d'une propriété d'annulation cohomologique déduite du
théorème~\ref{prfondo} relatif au foncteur $\omega : \F_\Gr\to\F$ et
de la proposition~\ref{comdom2} relative à l'auto-dualité du foncteur $\wt{\omega}$.

Le premier paragraphe établit, de façon directe, une autre propriété
d'annulation cohomologique du foncteur d'oubli $o_{inj} : \F\to\F_{inj}$.
Dans le paragraphe~\ref{p-krfi}, nous formulons une conjecture sur la
filtration de Krull de la catégorie $\F_{inj}$ analogue à la
conjecture artinienne extrêmement forte. Nous en démontrons un cas
particulier à l'aide des résultats des deux paragraphes précédents.

\subsection{Une propriété élémentaire} Le résultat que nous
établissons dans ce paragraphe (proposition~\ref{prcoh-fi1}) repose
uniquement sur les résultats de la section~\ref{sct-surj}. Pour les
appliquer plus commodément, pour traitons d'abord de la
catégorie~$\F_{surj}$.

\begin{lm}\label{lmf-acfsj} Il existe dans $\F_{surj}$ un épimorphisme
  scindé $o(F)\ptt
P^{surj}_V\twoheadrightarrow o(F)$ naturel en les objets $V$ de $\E^f_{surj}$ et $F$ de $\F$.
\end{lm}

\begin{proof} On définit un morphisme naturel $o(F)\ptt
P^{surj}_V\to o(F)$ comme l'adjoint de la projection canonique
$$\varpi(o(F)\ptt
P^{surj}_V)\simeq\varpi o(F)\otimes\varpi(P^{surj}_V)\simeq\varpi
o(F)\otimes P_V\twoheadrightarrow\varpi
o(F)\to F$$
déduit de l'épimorphisme $P_V\twoheadrightarrow\kk$ et de la coünité
de l'adjonction, via la proposition~\ref{pre-vpo}. 

D'autre part, l'épimorphisme $V\twoheadrightarrow 0$ induit un
morphisme ${\rm Is}_0\simeq P^{surj}_0\to P^{surj}_V$ dans
$\F_{surj}$, d'où un morphisme naturel $o(F)\simeq o(F)\ptt {\rm
  Is}_0\to o(F)\ptt P^{surj}_V$.

On vérifie aussitôt que le morphisme $o(F)\to o(F)\ptt P^{surj}_V\to
o(F)$ composé des deux flèches précédentes est l'identité.
\end{proof}

\begin{pr}\label{prcoh-fs1} Soient $X$ un objet fini de $\F_{surj}$ et
  $F$ un objet de $\F$. On a $${\rm Ext}^*_{\F_{surj}} (o(F),X)=0.$$
\end{pr}

\begin{proof} Par la proposition~\ref{ftffs}, il existe un espace vectoriel
  de dimension finie $V$ tel que $\delta_V(X)=0$. La proposition~\ref{agdscr} et le corollaire~\ref{adj-ah} montrent alors que ${\rm Ext}^*_{\F_{surj}} (o(F)\ptt
P^{surj}_V,X)\simeq {\rm
    Ext}^*_{\F_{surj}} (o(F),\delta_V(X))=0$. Mais  ${\rm Ext}^*_{\F_{surj}} (o(F),X)$ est facteur
direct de ${\rm Ext}^*_{\F_{surj}} (o(F)\ptt
P^{surj}_V,X)$ par le lemme~\ref{lmf-acfsj}, d'où la proposition.
\end{proof}

Revenons maintenant à la catégorie $\F_{inj}$, plus
naturelle pour les considérations suivantes. 

\begin{pr}\label{prcoh-fi1} Soient $X$ un objet localement fini de $\F_{inj}$ et
  $F$ un objet de $\F$. On a $${\rm Ext}^*_{\F_{inj}}
  (X,o_{inj}(F))=0.$$
\end{pr}

\begin{proof} Le cas où $X$ est fini est dual du précédent. Le cas
  général s'en déduit par passage à la colimite.
\end{proof}

\subsection{Propriétés utilisant le foncteur $\omega$}\label{par-inj-f} Nous présentons maintenant d'autres résultats d'annulation
cohomologique dans la catégorie~$\F_{inj}$, reposant sur le
théorème~\ref{prfondo} et la proposition suivante.

\begin{pr}\label{lmdt-inj} Les endofoncteurs $\varpi_{inj}o_{inj}$ et
  $\omega\kappa$ de $\F$ sont isomorphes.
\end{pr}

\begin{proof} Soit $\widetilde{\kappa} : \F\to\widetilde{\F}_\Gr$ le foncteur de précomposition par le foncteur de réduction
  $\widetilde{\mathfrak{K}} : \widetilde{\E}^f_\Gr\to\E^f$. On a des
  isomorphismes canoniques $\varpi_{inj}o_{inj}\simeq\widetilde{\omega}'\circ\widetilde{\kappa}$ et
  $\omega\kappa\simeq\widetilde{\omega}\circ\widetilde{\kappa}$, où
  $\wt{\omega}$ (resp. $\wt{\omega}'$) est le foncteur introduit dans la notation~\ref{ndw}
  (resp. \ref{nwtp}), de sorte que la proposition~\ref{comdom2} donne la conclusion.\end{proof}

On introduit à présent un foncteur très analogue à~$\I$.

\begin{defi}\label{j-def} On définit un endofoncteur $\J$ de $\F_\Gr$
  par
$$\J(X)(V,B)=\bigoplus_{W\in\Gr(B)} X(V/W,B/W)\,,$$
l'application linéaire $$\J(X)(f) : \J(X)(V,B)\to\J(X)(V',B')$$ (où $f :
(V,B)\to (V',B')$ est une flèche de $\E^f_\Gr$) ayant pour composante
$$X(V/W,B/W)\to X(V'/W',B'/W')$$ l'application linéaire induite par le
morphisme $(V/W,B/W)\to (V'/W',B'/W')$ induit par $f$ si $f(W)=W'$,
$0$ sinon.
\end{defi}

La propriété suivante justifie l'apparition du foncteur $\J$ dans ce paragraphe.

\begin{pr}\label{oko} Il existe un isomorphisme de foncteurs $\omega\kappa\omega\simeq\omega\J$.
\end{pr}

\begin{proof} On l'obtient par la suite d'isomorphismes naturels
$$\omega\kappa\omega(X)(V)=\bigoplus_{W\in\Gr(V)}\omega(X)(V/W)\simeq\bigoplus_{W\subset
  B\subset V} X(V/W,B/W)$$
$$\simeq\bigoplus_{B\in\Gr(V)}\J(X)(V,B)=\omega\J(X)(V).$$
\end{proof}

La proposition suivante est la plus importante de cette section. Elle
joue dans $\F_{inj}$ un rôle analogue à celui qu'occupe le
théorème~\ref{prfondo} dans la catégorie $\F_\Gr$.

\begin{pr}\label{gl-suslin} Soient $X$ un  objet localement fini de
$\F_\Gr$ et $Y$ un objet de $\F_\Gr$. Il existe un
isomorphisme naturel gradué ${\rm Ext}^*_{\F_{inj}}
(o_{inj}\omega(X),o_{inj}\omega(Y))\simeq {\rm Ext}^*_{\F_\Gr} (X,\I\J(Y))$.
\end{pr}

\begin{proof} Le foncteur $o_{inj}$ étant adjoint à gauche à
  $\varpi_{inj}$ (proposition \ref{pdfd-inj}), le corollaire
  \ref{adj-ah}, la proposition~\ref{lmdt-inj}, la
  proposition~\ref{oko} puis le théorème~\ref{prfondo} procurent des
  isomorphismes gradués naturels 
$${\rm Ext}^*_{\F_{inj}}
(o_{inj}\omega(X),o_{inj}\omega(Y))\simeq {\rm
  Ext}^*_\F(\omega(X),\varpi_{inj}o_{inj}\omega(Y))$$
$$\simeq {\rm
  Ext}^*_\F(\omega(X),\omega\kappa\omega(Y))\simeq {\rm
  Ext}^*_\F(\omega(X),\omega\J(Y))\simeq {\rm Ext}^*_\Gr (X,\I\J(Y)).$$
\end{proof}

\begin{rem}\label{cidoit} Ainsi, pour obtenir un énoncé relatif à la
  catégorie $\F_{inj}$, nous avons transité par la catégorie $\F_\Gr$,
  dont l'étude requiert celle de $\F_{surj}$. Cela explique pourquoi
  nous avons dû introduire les deux catégories $\F_{inj}$ et $\F_{surj}$.
\end{rem}

Le corollaire suivant permet d'utiliser la proposition~\ref{gl-suslin}
tout en s'affranchissant des foncteurs $\I$ et $\J$.

\begin{cor}\label{cohkfi} Soient $k,n\in\mathbb{N}$, $X$ un objet
  localement fini de $\F_{\Gr,k}$ et $Y$ un objet de $\F_{\Gr,n}$. 
\begin{enumerate}\item Si $k<n$, alors  ${\rm Ext}^*_{\F_{inj}}
(o_{inj}\omega_k(X),o_{inj}\omega_n(Y))=0$.
\item Si $k=n$, le foncteur exact $o_{inj}\omega_n$ induit un
  isomorphisme naturel gradué $${\rm Ext}^*_{\Gr,n}
(X,Y)\xrightarrow{\simeq} {\rm Ext}^*_{\F_{inj}}
(o_{inj}\omega_n(X),o_{inj}\omega_n(Y)).$$
\end{enumerate}
\end{cor}

\begin{proof} On déduit ce résultat de la proposition \ref{gl-suslin}
  de la même façon que l'on a déduit le corollaire \ref{fondcr} du théorème \ref{prfondo}. En effet, si
  $A$ est un objet de $\F_\Gr$ de coniveau au moins égal à $n$,
  l'inclusion canonique $A\hookrightarrow\J(A)$ induit après
  application du foncteur $\mathcal{R}_{\mathbb{N},\leq n}$ un isomorphisme.
\end{proof}

\paragraph*{Application à la $K$-théorie stable} La démonstration par
Suslin (cf. appendice de~\cite{FFSS}) de l'isomorphisme entre la $K$-théorie stable de $\kk$ et
l'homologie dans $\F(\kk)$ pour des foncteurs finis se décompose en
deux étapes\,\footnote{Scorichenko a étendu le résultat de Suslin au
  cas des modules sur anneau arbitraire --- cf. le dernier chapitre de
\cite{FFPS}.}.

Rappelons que l'on dispose d'un foncteur exact ${\rm e} :
\F_{inj}\to\,_{\kk[GL(\kk)]}\mathbf{Mod}$, conformément
au~§\,\ref{par-dwyer}. Nous noterons $F\mapsto F(\kk^\infty)$ le foncteur
$\F\to\,_{\kk[GL(\kk)]}\mathbf{Mod}$ composé de $o_{inj} :
\F\to\F_{inj}$ et du précédent. Il induit un morphisme gradué naturel
$${\rm Ext}^*_{\F_{inj}(\kk)}(o_{inj}(F),o_{inj}(G))\to {\rm
  Ext}^*_{\kk[GL(\kk)]}(F(\kk^\infty),G(\kk^\infty))$$
pour $F, G\in {\rm Ob}\,\F$. La première étape de la démonstration de
Suslin consiste à démontrer que ce morphisme est un isomorphisme
lorsque $F$ et $G$ sont finis ; c'est vrai en fait lorsque $F$ est
pf$_\infty$ et $G$ fini (cf. \cite{Pira}, proposition~4.3). Le lien
avec la $K$-théorie stable vient de ce que, à dualisation près, les
groupes d'extensions ${\rm
  Ext}^*_{\kk[GL(\kk)]}(F(\kk^\infty),G(\kk^\infty))$ sont des groupes de
$K$-théorie stable de $\kk$ (voir encore \cite{Pira}).

\begin{cor} Soient $X$ un objet fini de $\F_\Gr$ et $F$ un objet fini
  de $\F$. Il existe un isomorphisme naturel gradué
$${\rm  Ext}^*_{\kk[GL(\kk)]}(\omega(X)(\kk^\infty),F(\kk^\infty))\simeq {\rm
  Ext}^*_\Gr (X,\mathcal{I}\kappa(F)).$$
\end{cor}

\begin{proof} Comme $X$ est pf$_\infty$ et que le foncteur exact
  $\omega$ préserve les projectifs de type fini, $\omega(X)$ est
  pf$_\infty$, de sorte que ${\rm Ext}^*_{\F_{inj}(\kk)}(o_{inj}\omega(X),o_{inj}(F))\simeq {\rm
  Ext}^*_{\kk[GL(\kk)]}(\omega(X)(\kk^\infty),F(\kk^\infty))$. On conclut en
appliquant la proposition~\ref{gl-suslin}, puisque $\mathcal{J}\mathcal{P}_{0,\mathbb{N}}\simeq\kappa$.
\end{proof}

Lorsque $X$ est de niveau $0$, i.e. que  $\omega(X)$ est un objet fini
de $\F$, on retrouve le résultat de Betley-Suslin affirmant que le
morphisme naturel ${\rm Ext}_\F^*(F,G)\to {\rm
  Ext}^*_{\kk[GL(\kk)]}(F(\kk^\infty),G(\kk^\infty))$ est un isomorphisme si
$F$ et $G$ sont finis.

\subsection{La filtration de Krull de la catégorie $\F_{inj}$}\label{p-krfi} Le corollaire~\ref{cohkfi} et la proposition~\ref{prcoh-fi1} nous amènent à formuler la conjecture suivante.

\begin{conj}[Conjecture artinienne extrêmement forte pour $\F_{inj}$]\label{cakinj2} Pour tout $n\in\mathbb{N}$, le foncteur
  $o_{inj}\omega_n : \F_{\Gr,n}\to\F_{inj}$ induit une équivalence
  entre les catégories $\F_{\Gr,n}^{lf}$ et $\K_{n+1}(\F_{inj})/\K_n(\F_{inj})$.
\end{conj}

Si la conjecture artinienne extrêmement forte est
satisfaite, cet énoncé équivaut au suivant.

\begin{conj}\label{cakinj} Pour tout $n\in\mathbb{N}$, le foncteur $o_{inj}$ induit une équivalence entre les catégories $\K_{n}(\F)/\K_{n-1}(\F)$ et $\K_{n+1}(\F_{inj})/\K_n(\F_{inj})$.
\end{conj}

Le résultat suivant donne une
réponse positive partielle à la conjecture~\ref{cakinj2} pour $n=0$.

\begin{pr}\label{prfki-i} Le foncteur $o_{inj} : \F\to\F_{inj}$ induit une
  équivalence entre la sous-catégorie $\F^{f}$ de $\F$ et une
  sous-catégorie épaisse de $\F_{inj}/\F_{inj}^{lf}$.
\end{pr}

\begin{proof} Si $S$ est un objet simple de $\F$, le
 corollaire~\ref{cr-strffse} montre, par dualité, que $o_{inj}(F)$ est
 un objet simple noethérien de type $1$ de $\F_{inj}$, donc un objet
 simple de $\F_{inj}/\F_{inj}^{lf}$. Le corollaire~\ref{cohkfi} (appliqué
 avec $n=0$) et la proposition~\ref{prcoh-fi1} établissent par ailleurs que le
 foncteur $o_{inj}$ induit un foncteur pleinement fidèle entre $\F^f$ et une
 sous-catégorie de $\F_{inj}/\F_{inj}^{lf}$ stable par extensions. Cela
 prouve la proposition.
\end{proof}

\begin{rem}\begin{enumerate}\item Dans \cite{these}, nous établissons ce résultat de façon
  différente. La méthode suivie, moins élémentaire que celle de la proposition~\ref{prfki-i}, est une variante
 moins technique de la démonstration du théorème~\ref{tha31} que nous
 donnerons dans \cite{art3}, théorème que l'on ne peut pas en revanche
 montrer de manière aussi directe que la proposition~\ref{prfki-i}, le
 corollaire~\ref{cr-strffse} n'ayant pas d'équivalent dans la catégorie~$\F$.
\item La proposition~\ref{finjmlf} esquisse une autre approche de la
  catégorie $\F_{inj}/\F_{inj}^{lf}$. Malheureusement, celle-ci ne
  semble pas suffisante pour établir simplement que tous les objets
  simples de la catégorie  $\F_{inj}/\F_{inj}^{lf}$ sont dans l'image
  du foncteur $\F\xrightarrow{o_{inj}}\F_{inj}\twoheadrightarrow\F_{inj}/\F_{inj}^{lf}$.
\end{enumerate}
\end{rem}

\appendix

\section{Adjonctions}\label{apa}

\subsection{Algèbre homologique} La propriété immédiate suivante
est d'un usage très courant dans les catégories de foncteurs.

\begin{pr}\label{adj-ah} Supposons que $\A$ et
  $\B$ sont deux catégories abéliennes possédant  suffisamment d'objets injectifs. Si $F : \B\to\A$ et $G : \A\to\B$ sont des foncteurs {\em
  exacts} tels que $F$ est
adjoint à droite à $G$, alors l'isomorphisme naturel
  ${\rm hom}_\A(X,F(Y))\simeq{\rm hom}_\B(G(X),Y)$
  s'étend en un isomorphisme naturel gradué 
$${\rm Ext}^*_\A(X,F(Y))\simeq{\rm Ext}^*_\B(G(X),Y).$$
\end{pr}

Nous utiliserons parfois cette proposition via le corollaire suivant.

\begin{cor}\label{aux-adah} Supposons que : 
\begin{enumerate}\item $\A$ et $\B$ sont deux catégories abéliennes,
  $\A$ possédant assez d'injectifs ;
\item $F$ et $G$ sont deux endofoncteurs exacts de $\A$, avec $G$
  adjoint à gauche à $F$ ; 
\item $H : \B\to\A$ est un foncteur exact ; 
\item $K$ est un endofonteur de $\B$ tel que $F\circ H\simeq H\circ K$ ;
\item il existe une transformation naturelle injective
  $j : id_\B\hookrightarrow K$. 
\end{enumerate}

Soient $X$ un objet de $\A$ tel qu'il existe $n\in\mathbb{N}$ tel que
$G^n(X)=0$ et $Y$ un objet de~$\B$. On a ${\rm Ext}^*_\A (X,H(Y))=0$.
\end{cor}

\begin{proof} Quitte à remplacer $F$ par $F^n$, $G$ par $G^n$ et $K$
  par $K^n$, on  peut supposer $n=1$.

Pour tout objet $Y$ de $\B$, on a une suite exacte naturelle
$$0\to Y\xrightarrow{f_1} K(Y_1)\xrightarrow{f_2} K(Y_2)\xrightarrow{f_3}\dots\xrightarrow{f_n} K(Y_n)\xrightarrow{f_{n+1}}\cdots$$
où $Y_1=Y$ et  $f_1=j_Y$, et pour $n>1$, $Y_n=coker\,f_{n-1}$ et $f_n$
est la composée $K(Y_{n-1})\twoheadrightarrow Y_n\xrightarrow{j_{Y_n}}
K(Y_n)$. En appliquant le foncteur exact $H$ à cette suite exacte et
en utilisant l'isomorphisme $F\circ H\simeq H\circ K$, on en déduit
une résolution de $H(Y)$ par les objets $F(H(Y_i))$. Comme ${\rm Ext}^*_\A
(X,F(A))=0$ pour tout $A\in {\rm Ob}\,\A$ par la
proposition~\ref{adj-ah}, cela donne la conclusion. 
\end{proof}

\subsection{Monades et comonades} Pour les démonstrations des
propriétés rappelées ci-dessous, nous renvoyons à~\cite{Mac},
chapitre VI.

\begin{defi} Soit $\C$ une catégorie.
\begin{enumerate}\item Une {\em monade} sur $\C$ est
  un triplet $(T,\eta,\mu)$ formé d'un endofoncteur $T$ de $\C$ et de
  transformations naturelles $\eta : id_\C\to T$ et $\mu : T^2\to T$
  telles que $\mu\circ\eta_T=id_T=\mu\circ T(\eta)$ et $\mu\circ T(\mu)=\mu\circ\mu_T$. La transformation naturelle $\mu$ est appelée {\em
  multiplication} de la
monade.
\item Une {\em comonade} sur $\C$ est une monade
  sur $\C^{op}$, i.e. un triplet $(T,\gamma,\delta)$ formé d'un
  endofoncteur $T$ de $\C$ et de transformations naturelles $\gamma :
  T\to id_\C$ et $\delta : T\to T^2$ telles que
  $\gamma_T\circ\delta=id_T=T(\gamma)\circ\delta$ et $\delta_T\circ\delta=T(\delta)\circ\delta$.
\end{enumerate}
\end{defi}

\begin{conv}
Dans la suite de cet appendice, nous considérons des catégories $\A$ et $\B$
et des foncteurs $F : \B\to\A$ et $G : \A\to\B$ tels que $F$ est
adjoint à droite à $G$.
\end{conv}

\begin{pr}\label{mon-adj} L'adjonction entre $F$ et $G$ détermine :
\begin{enumerate}\item une monade $(FG,\eta,F(\gamma_G))$ sur $\A$ ;
\item une comonade $(GF,\gamma,G(\eta_F))$ sur $\B$.
\end{enumerate}
\end{pr}

\begin{defi}\label{dfalgm}
\begin{enumerate}\item Soit $(T,\eta,\mu)$ une
    monade sur une catégorie $\C$. Un {\em
      module} sur cette monade
    est la donnée d'un objet $X$ de $\C$ et d'un morphisme $m :
    T(X)\to X$ tel que $m\circ \eta_X=id_X$ et $m\circ
    T(m)=m\circ\mu_X : T^2(X)\to X$.
\item Soit $(T,\gamma,\delta)$ une comonade sur une catégorie $\C$. Un
  {\em comodule} sur cette comonade est un module sur la monade de
  $\C^{op}$ associée, i.e. un objet $X$ de $\C$ muni d'un morphisme $c
  : X\to T(X)$ tel que $\gamma_X\circ c=id_X$ et $T(c)\circ
  c=\delta_X\circ c : X\to T^2(X)$.
\end{enumerate}
\end{defi}

On définit de façon évidente la notion de morphisme de modules sur une
monade, ou de comodules sur une comonade. On prendra garde au fait
que nous nommons {\em module} (resp. {\em comodule}) ce que Mac Lane
appelle {\em algèbre} (resp. {\em coalgèbre}).

\begin{pr}\label{prmon} Il existe un foncteur de $\A$ vers la
  catégorie des comodules sur la comonade $(GF,\gamma,G(\eta_F))$ ; il
  s'obtient sur les objets en munissant $G(X)$ (où $X\in {\rm
    Ob}\,\A$) de la structure de comodule donnée par le morphisme
  $G(\eta_X) : G(X)\to GFG(X)$.
\end{pr}

Cette proposition possède une variante duale en terme de monade.

Le résultat suivant, qui se déduit de \cite{Mac}, chapitre~VI, §\,7, théorème~1, est un cas particulier du {\em théorème de Beck}.

\begin{pr}\label{monadique} Faisons les hypothèses suivantes :
\begin{enumerate}\item les catégories $\A$ et $\B$ sont abéliennes ; 
\item le foncteur $G$ est exact et fidèle.
\end{enumerate}

Alors le foncteur de la proposition précédente de $\A$ vers la catégorie des comodules
sur la comonade déterminée par l'adjonction entre $F$ et $G$ est une
équivalence de catégories.
\end{pr}

Dualement, on a le résultat suivant.

\begin{pr}\label{dualmonadique} Faisons les hypothèses suivantes :
\begin{enumerate}\item les catégories $\A$ et $\B$ sont abéliennes ; 
\item le foncteur $F$ est exact et fidèle.
\end{enumerate}

Alors la catégorie $\B$ est équivalente à la catégorie des modules
sur la monade déterminée par l'adjonction entre $F$ et $G$.
\end{pr}

\section{Propriétés de finitude dans les catégories abéliennes}\label{apfca}

Cet appendice donne les définitions, notations et propriétés relatives
aux notions de finitude utilisées dans cet article. La plupart des ces
notions se trouvent dans \cite{Pop}, \cite{Gab} ou \cite{CR} ; d'autres
références et les démonstrations sont données dans~\cite{these}.

\begin{conv} Dans cet appendice, $\A$ désigne une catégorie de Grothendieck.
\end{conv}

On rappelle que, par définition, une {\em catégorie de Grothendieck} est une
  catégorie abélienne possédant des colimites, un ensemble de
  générateurs et dans laquelle les colimites filtrantes sont
  exactes. Une telle catégorie possède toujours des limites, des enveloppes injectives et un cogénérateur
injectif (cf.~\cite{Pop}, §\,3.7 et \cite{Gab}, chapitre~II, §\,6).

\subsection{Définitions}\label{par-anfd} Soit $A$ un objet de $\A$. 

\begin{defi}\label{defi-anf} On dit que $A$ est :
\begin{enumerate}\item {\em noethérien} si toute suite croissante de
  sous-objets de $A$ stationne ; 
\item {\em artinien}  si toute suite décroissante de
  sous-objets de $A$ stationne ;
\item {\em de type fini}  si toute suite croissante de
  sous-objets de $A$ de colimite $A$ stationne ;
\item {\em de co-type fini}  si toute suite décroissante de
  sous-objets de $A$ de limite nulle stationne.
\end{enumerate}
\end{defi}

\begin{nota}\label{not-prfi} Nous désignerons par $\A^{tf}$ la sous-catégorie pleine des objets de type fini de $\A$.
\end{nota}

\begin{defi}\label{cola} On dit que la catégorie $\A$ est {\rm
    localement noethérienne} (resp. {\rm co-localement artinienne}) si
  elle possède un ensemble de générateurs (resp. cogénérateurs)
  noethériens (resp. artiniens).
\end{defi}

\begin{hyp}\label{hypf1} Il existe dans $\A$ un ensemble de
  générateurs projectifs de type fini.
\end{hyp}

Le premier point de la définition suivante est valable sous
l'hypothèse~\ref{hypf1} : la \go bonne\gf définition d'un objet pf$_n$ est différente lorsque
cette hypothèse n'est pas vérifiée.

\begin{defi} Soit $n\in\mathbb{N}^*$. On dit que $A$ est : 
\begin{enumerate}\item de {\em $n$-présentation finie}, en abrégé
  pf$_n$, s'il existe une suite exacte $P_n\to P_{n-1}\to\cdots\to
  P_0\to A\to 0$, où les $P_i$ sont des objets projectifs de type fini
  de $\A$.
\item de {\em $n$-co-présentation finie}, en abrégé co-pf$_n$, s'il
  existe une suite exacte $0\to I_0\to\cdots\to I_n$, où les $I_i$
  sont des objets injectifs de co-type fini de $\A$.
\end{enumerate}
\end{defi}

On dit que $A$ est de
présentation finie s'il est de $1$-présentation finie, pf$_\infty$
s'il est pf$_n$ pour tout $n$. On adopte des simplifications
terminologiques analogues pour les objets co-pf$_i$ (ainsi, on écrira
co-pf pour co-pf$_1$).

La manipulation des objets pf$_\infty$ est facilitée par la
proposition suivante :

\begin{pr}\label{pr:pfi} Soit $0\to A\to B\to C\to 0$ une suite exacte
  de $\A$. Si deux des objets $A$, $B$ et $C$ sont pf$_\infty$, il en
  est de même du troisième.
\end{pr}

\begin{defi}\label{dff-tr} On dit que $A$ est :
\begin{enumerate}\item {\em simple} s'il est non nul mais que tous ses
  sous-objets stricts sont nuls ;
\item {\em fini} s'il admet une filtration finie de sous-quotients simples ; 
\item {\em localement fini} s'il est colimite de ses sous-objets finis
  ;
\item {\em co-localement fini} s'il est limite de ses quotients finis.
\end{enumerate}
\end{defi}

\begin{nota}\label{ndf-tr} On note $\A^f$ (respectivement $\A^{lf}$) la sous-catégorie pleine
  des objets finis  (resp. localement finis) de $\A$.
\end{nota}

\begin{defi}\label{dfserre} Une {\em sous-catégorie de Serre} de $\A$
  est une sous-catégorie pleine $\C$ de $\A$ stable par
  sommes directes finies et par sous-quotients. 

Une {\em sous-catégorie épaisse} de $\A$ est une sous-catégorie de Serre de
  $\A$ stable par extensions. 
\end{defi}

\begin{pr}\label{prevff}\begin{enumerate}\item La sous-catégorie
    $\A^f$ de $\A$ est épaisse.
\item La catégorie $\A^{lf}$ est une sous-catégorie de Serre de $\A$. Si les objets simples de $\A$ sont de présentation finie, elle est épaisse.
\end{enumerate}
\end{pr}

La dernière assertion de cette proposition est un cas particulier du
résultat suivant.

\begin{pr}\label{pr:pf}  Soit $\C$ une sous-catégorie épaisse de
  $\A$. Notons $\overline{\C}$ la sous-catégorie pleine de $\A$ dont
  les objets sont les colimites d'objets de $\C$.

Si les objets de $\C$ sont de présentation finie, alors
$\overline{\C}$ est une sous-catégorie épaisse de~$\A$.
\end{pr}

\begin{defi}\label{defsocs}\begin{enumerate}\item On appelle {\em
      socle} de $A$ la somme ${\rm soc} (A)$ des sous-objets simples
    de $A$.
\item On nomme {\em radical} de $A$ l'intersection ${\rm rad}(A)$ des
  sous-objets stricts maximaux de $A$. On appelle {\em cosocle} de $A$ le quotient $A/{\rm rad}(A)$,
  noté ${\rm cosoc}(A)$.
\end{enumerate}
\end{defi}

\begin{pr}\label{cotfs} Si $A$ est localement fini et que ${\rm soc} (A)$ est
  fini, alors $A$ est de co-type fini.
\end{pr}

\subsection{Effet de foncteurs exacts} Les propositions énoncées dans
ce paragraphe s'établissent par des méthodes standard ; on en trouvera
une démonstration complète dans \cite{these}.
\begin{pr}\label{preltf} Soient $\B$ une catégorie de Grothendieck et $F : \A\to\B$ un foncteur exact et fidèle.
\begin{enumerate}
\item Si $X$ est un objet de $\A$ tel que $F(X)$ est noethérien
  (resp. artinien) dans
  $\B$, alors $X$ est noethérien (resp. artinien) dans $\A$.
\item Si $F$ commute aux colimites (resp. limites) filtrantes et si
  $X$ est un objet de $\A$ tel que $F(X)$ est de type fini (resp. de
  co-type fini), alors $X$ est de type fini (resp. de co-type fini).
\item Supposons que $F$ est plein et que son image est une
  sous-catégorie de Serre de $\B$. Si $X$ est un objet noethérien
  (resp. de type fini, artinien, de co-type fini) de $\A$, alors $F(X)$ est un
  objet noethérien (resp. de type fini, artinien, de co-type fini) de $\B$.
\end{enumerate}
\end{pr}

\begin{pr}\label{crpf2} Soient $i\in\mathbb{N}^*$ et $F$ un foncteur exact et fidèle de $\A$ dans une catégorie de Grothendieck, commutant aux colimites filtrantes et
  préservant les objets projectifs de type fini. Si $X$
  est un objet de $\A$ tel que $F(X)$ est pf$_i$, alors $X$ est~pf$_i$.
\end{pr}

\begin{pr}\label{pr-ff} Soient $\B$ une catégorie de Grothendieck et $F : \A\to\B$ un foncteur exact.
\begin{enumerate}
\item Supposons $F$ fidèle. Si $X$ est un objet de $\A$ tel que $F(X)$ est fini dans $\B$, alors $X$ est fini dans $\A$.
\item Supposons que $F$ est plein et que son image est une
  sous-catégorie de Serre de $\B$. Si $X$ est un objet fini de $\A$, alors $F(X)$ est un
  objet fini de $\B$.
\end{enumerate}
\end{pr}

\subsection{Groupes de Grothendieck}
\begin{defi} Soit $\C$ une sous-catégorie
  pleine et {\em petite} de $\A$, contenant $0$. On appelle {\em groupe de Grothendieck} de $\C$ relativement à $\A$ le
  groupe abélien noté $G_0(\C;\A)$ défini par générateurs et relations
  de la manière suivante.
\begin{itemize}\item \textbf{Générateurs} : un générateur $[A]$ pour chaque objet $A$ de $\C$.
\item \textbf{Relations} : $[A]=[B]+[C]$ pour toute suite exacte
  courte $0\to
  B\to A\to C\to 0$ de $\A$ dont tous les objets $A$, $B$, $C$ sont
  dans $\C$.
\end{itemize}
On a $[A]=[B]$ dans $G_0(\C;\A)$ si $A$ et $B$ sont
deux objets isomorphes de $\C$, ce qui permet de définir ce groupe
lorsque $\C$ est seulement  {\em essentiellement petite}.
\end{defi}

\begin{lm} La sous-catégorie $\A^{tf}$ de $\A$ est essentiellement petite.
\end{lm}

Cela permet de donner la notation suivante.

\begin{nota}\label{not-grot} Nous noterons respectivement $G_0^f(\A)$, $G_0^{tf}(\A)$ et $K_0(\A)$ les groupes de Grothendieck $G_0(\A^f
  ; \A)$, $G_0(\A^{tf},\A)$ et $G_0(\C ; \A)$, où $\C$ est la
  sous-catégorie pleine des objets projectifs de type fini de $\A$.
\end{nota}

Le théorème de Jordan-Hölder dit que le groupe abélien $G_0^f(\A)$ est
libre, les classes des objets simples de $\A$ en formant une base.

Nous rappelons maintenant la notion de {\em recollement} de
catégories abéliennes.

\begin{defi}\label{def-rec}
Un {\em diagramme de recollement} est un diagramme du type
$$\xymatrix{\C\ar[r]|-i & \A\ar[r]|-e\ar@/_/[l]_-q\ar@/^/[l]^-p & \B\ar@/_/[l]_-l\ar@/^/[l]^-r
}$$
dans lequel :
\begin{itemize}
\item $\A$, $\B$ et $\C$ sont des catégories abéliennes.  
\item Le foncteur $l$ est adjoint à gauche à $e$ et $e$ est adjoint à gauche à $r$ (en
  particulier, $e$ est exact).
\item L'unité $id_\B\to el$ et la coünité $er\to id_\B$ sont des isomorphismes.
\item Le foncteur $q$ est adjoint à gauche à $i$ et $i$ est adjoint à gauche à $p$ (en
  particulier, $i$ est exact).
\item L'unité $id_\C\to pi$ et la coünité $qi\to id_\C$ sont des isomorphismes.
\item  Le foncteur $i$ est un plongement pleinement fidèle d'image $ker\,e$ (en
  particulier, $i$ identifie $\C$ à une sous-catégorie épaisse de $\A$).
\end{itemize}
\end{defi}

Dans cette situation, le foncteur $e$ induit
une équivalence $\A/\C\xrightarrow{\simeq}\B$ (nous renvoyons le
lecteur à \cite{Gab} pour ce qui concerne la notion de quotient d'une catégorie
abélienne par une sous-catégorie épaisse).

Nous renvoyons à \cite{K2} à ce sujet, ainsi que pour la démonstration
de la proposition suivante (et la description explicite de l'isomorphisme).

\begin{pr}\label{crgrrec} On a un isomorphisme de groupes
  $G_0^f(\A)\simeq G_0^f(\B)\oplus G_0^f(\C)$.
\end{pr}

\subsection{Filtration de Krull}\label{par-krull}
\begin{defi} La {\em filtration de
    Krull} de la catégorie $\A$ est
  la suite croissante de sous-catégories épaisses stables par
  colimites $\mathcal{K}_n(\A)$ de $\A$ définie
  inductivement comme suit.
\begin{itemize}\item La catégorie $\mathcal{K}_n(\A)$ est réduite à
  l'objet nul pour $n<0$.
\item Pour $n\geq 0$,  $\mathcal{K}_n(\A)$ est l'image réciproque par
  le foncteur canonique $\A\to\A/\mathcal{K}_{n-1}(\A)$ de la
plus petite  sous-catégorie épaisse et stable par colimites de
  $\A/\mathcal{K}_{n-1}(\A)$ contenant tous les objets simples de cette
  catégorie.
\end{itemize}
\end{defi}

\begin{rem} On peut étendre de manière claire la définition de la
  filtration de Krull à  tout ordinal (cf. \cite{Gab}) ; nous n'avons
  introduit que les termes indicés par $\mathbb{N}$ car eux seuls
  interviendront dans notre contexte.
\end{rem}

\begin{defi}\label{dsnt} On définit par
  récurrence sur $n\in\mathbb{N}\cup\{-1\}$ les notions d'objet {\em simple
    noethérien de type $n$} et {\em  noethérien de type
    $n$} de $\A$ de la façon suivante.
\begin{itemize}\item Un objet est simple noethérien de type $-1$ s'il
  est nul.
\item Un objet $X$ est noethérien de type $n$ s'il possède une
  filtration finie $0=F_0\subset F_1\subset\dots\subset F_k=X$ telle
  que, pour tout $i\in\{1,\dots,k\}$, le quotient $F_i/F_{i-1}$ est
  simple noethérien de type $a(i)$, pour un certain entier $a(i)\leq n$.
\item Un objet est simple noethérien de type $n$ pour $n\in\mathbb{N}$
  si et seulement s'il n'est pas noethérien de type $n-1$ et que tous
  ses quotients stricts sont noethériens de type $n-1$.
\end{itemize}
\end{defi}

\begin{nota}\label{notnt} Nous désignerons par $\A_{\mathbf{NT}(n)}$ la sous-catégorie pleine de
  $\A$ formée des objets noethériens de type $n$. Elle est épaisse.
\end{nota}

\begin{pr}\label{prkrull} Soient $n\in\mathbb{N}\cup\{-1\}$ et $X$ un objet noethérien de type $n$ de $\A$. Alors $X$ est noethérien, appartient à $\mathcal{K}_{n}(\A)$ et son
image dans $\A/\mathcal{K}_{n-1}(\A)$ (si $n\geq 0$)
est finie. Elle est même simple si $X$ est simple noethérien de type $n$.

Réciproquement, un objet noethérien de $\A$ qui appartient à $\mathcal{K}_{n}(\A)$ est noethérien de type $n$.
\end{pr}
%

\section{Catégories de foncteurs}\label{apfct}

Nous rappelons ici quelques propriétés d'usage courant des catégories
de foncteurs. Bien que standard, elles ne sont pas toujours facilement
accessibles dans la littérature --- voir~\cite{these} pour plus de détails.

\begin{conv} Dans cet appendice, $\I$ et $\J$ désignent des catégories
essentiellement petites, $\A$ et $\B$ des catégories abéliennes et $A$
un anneau.
\end{conv}

\subsection{Généralités}

\begin{nota}\label{postpre}\begin{enumerate}\item Si $F$ est un foncteur de $\J$ vers $\I$, nous désignerons
    par $F^* : \mathbf{Fct}(\I,\A)\to\mathbf{Fct}(\J,\A)$ le {\em foncteur de
    précomposition} par $F$.
\item Pour tout foncteur $G : \A\to\B$, nous désignerons par $G_* :
  \mathbf{Fct}(\I,\A)\to\mathbf{Fct}(\I,\B)$\index{Nota}{$G_*$} le {\em foncteur de
  postcomposition} par $G$.
\item Soit $E$ un objet de $\I$. Nous noterons $ev_E
  : \mathbf{Fct}(\I,\A)\to\A$ et ${\rm ev}_E : \mathbf{Fct}(\I,_A\mathbf{Mod})\to
_{A[M]}\negmedspace\mathbf{Mod}$
les {\em foncteurs d'évaluation} en $E$, donnés par la précomposition par le
foncteur $*\to\I$ d'image $E$ et par le foncteur pleinement fidèle  $\underline{{\rm
End}_\I(E)}\to\I$ d'image $E$ respectivement. \end{enumerate}
\end{nota}

\begin{pr}\label{crcfab}\begin{enumerate}\item La catégorie
  $\mathbf{Fct}(\I,\A)$ est abélienne. 
\item L'exactitude se teste \go argument par
argument\gf : une suite $X\to Y\to Z$ de $\mathbf{Fct}(\I,\A)$ est
exacte si et seulement si la suite $X(E)\to Y(E)\to Z(E)$ est
exacte dans $\A$ pour tout objet $E$ de $\I$.
\item Si $\A$ est une catégorie {\em $A$-linéaire}, $\mathbf{Fct}(\I,\A)$
hérite d'une structure de catégorie $A$-linéaire.
\item Supposons que la catégorie abélienne $\A$ est monoïdale symétrique. Le bifoncteur
 $$\mathbf{Fct}(\I,\A)\times\mathbf{Fct}(\I,\A)\simeq
  \mathbf{Fct}(\I,\A\times\A)\xrightarrow{\otimes_*}\mathbf{Fct}(\I,\A)$$
  définit une structure monoïdale symétrique
  sur $\mathbf{Fct}(\I,\A)$.
\end{enumerate}
\end{pr}

On rappelle que les notions de sous-catégorie de Serre et
sous-catégorie épaisse sont introduites dans la définition~\ref{dfserre}.

\begin{pr}\label{lmfora} Soit $ F: \A\to \B$ un foncteur.
\begin{enumerate}\item Si $F$ est exact (resp. additif, exact à
  gauche, exact à droite), il en est de même pour $F_* : \mathbf{Fct}(\I,\A)\to\mathbf{Fct}(\I,\B)$.
\item Si $F$ est fidèle, alors $F_*$ est fidèle.
\item Si $F$ est pleinement fidèle, alors $F_*$ est pleinement fidèle. Si de
  plus l'image de $F$ est une sous-catégorie de Serre (resp. épaisse)
  de $\B$, alors l'image de $F_*$ est une sous-catégorie de Serre (resp. épaisse)
  de $\mathbf{Fct}(\I,\B)$.
\item Si $F$ possède un adjoint à gauche $G$, alors $G_*$ est adjoint
  à gauche à $F_*$. 
\end{enumerate}
\end{pr}

\begin{pr}\label{lm-form} Soit $F: \J\to\I$ un foncteur. 
\begin{enumerate}
\item Le foncteur $F^* : \mathbf{Fct}(\I,\A)\to\mathbf{Fct}(\J,\A)$ commute aux limites et aux colimites ;
  il est en particulier exact.
\item Si $\mathcal{A}$ est une catégorie monoïdale symétrique, $F^*$ commute au produit tensoriel.
\item Si $F$ est essentiellement surjectif, $F^*$ est fidèle.
\item Si $F$ est plein et essentiellement surjectif, $F^*$ est
  pleinement fidèle, et son image est une sous-catégorie de Serre de
  $\mathbf{Fct}(\J,\A)$ stable par limites et colimites.
\item Si $F$ possède un adjoint à gauche $G$, alors $G^*$ est adjoint
  à droite à $F^*$.  
\end{enumerate}
\end{pr}

 Dans la suite de cette annexe, on ne s'intéresse plus qu'au cas où la
 catégorie but est la catégorie d'espaces vectoriels $\E_\kk$.

\subsection{Générateurs projectifs}\label{sgs-p} \'Etant donné un objet $E$ de $\I$, nous noterons $P_E^\I$
  l'objet $\kk[{\rm hom}_\I(E,\cdot)]$ de
  $\mathbf{Fct}(\I,\E_\kk)$. On rappelle que  $\kk[.]$ désigne le
  foncteur de $\kk$-linéarisation ; ainsi, $\kk[F]$ désigne la
$\kk$-linéarisation $\kk[.]_*F$ d'un foncteur $F : \I\to\mathbf{Ens}$.

La bifonctoralité de ${\rm hom}_\I$ permet de considérer $E\mapsto
P_E^\I$ comme un foncteur $\I^{op}\to\mathbf{Fct}(\I,\E_\kk)$.

\begin{pr}[Lemme de Yoneda linéaire]\label{yoneda} Il existe un isomorphisme
$${\rm hom}_{\mathbf{Fct}(\I,\E_\kk)}(P_E^\I,F)\simeq F(E)$$
naturel en les objets $E$ de $\I$ et $F$ de $\mathbf{Fct}(\I,\E_\kk)$.
\end{pr}

\begin{prdef}\label{prgpcf} Les objets $P_E^\I$ forment un ensemble de
  générateurs projectifs de type fini de
  $\mathbf{Fct}(\I,\E_\kk)$ lorsque $E$ parcourt un squelette
  de $\I$. On les appelle {\em générateurs projectifs
    standard} de  $\mathbf{Fct}(\I,\E_\kk)$.
\end{prdef}

\begin{cor}\label{corgrct} La catégorie de foncteurs
  $\mathbf{Fct}(\I,\E_\kk)$ est une catégorie de
  Grothendieck vérifiant l'hypothèse~\ref{hypf1}.
\end{cor}

\begin{cor}\label{cr-adjp} Pour tout foncteur $F : \J\to\I$, le
  foncteur de précomposition $F^* : \mathbf{Fct}(\I,\E_\kk)\to\mathbf{Fct}(\J,\E_\kk)$ admet un adjoint à droite et un adjoint à gauche.
\end{cor}

 Les {\em extensions de Kan} (cf. \cite{Mac}, ch. X) permettent de donner une construction des
  adjoints.

\begin{rem}\label{rqprcf} Si deux objets $V$ et
  $W$ de $\I$ possèdent une somme, on a un isomorphisme canonique
  $P_V^\I\otimes P^\I_W\simeq P^\I_{V\amalg W}$.
\end{rem}

\begin{cor}\label{crptfcf} Supposons que   $\I$ possède des sommes finies. Alors le produit tensoriel de $\mathbf{Fct}(\I,\E_\kk)$ préserve  les objets projectifs et les objets de type fini.
\end{cor}

\begin{prdef}[Dualité entre catégories de foncteurs]\label{prdcff} Notons
$D_{\I,\kk} : \mathbf{Fct}(\I,\E_\kk)^{op}\simeq
\mathbf{Fct}(\I^{op},\E_\kk^{op})\to\mathbf{Fct}(\I^{op},\E_\kk)$
le foncteur  de postcomposition par le foncteur de dualité $(.)^*={\rm
  hom}_\E(.,\kk) : \E^{op}_\kk\to\E_\kk$, et $D'_{\I,\kk}=D_{\I^{op},\kk} :
\mathbf{Fct}(\I^{op},\E_\kk)^{op}\to\mathbf{Fct}(\I,\E_\kk)$.
\begin{enumerate}\item Les foncteurs $D_{\I,\kk}$ et $D'_{\I,\kk}$ sont
  exacts et fidèles.
\item Le foncteur $D_{\I,\kk}$ est adjoint à droite à
  $(D'_{\I,\kk})^{op}$.
\item Ces foncteurs induisent des équivalences de catégories
  réciproques l'une de l'autre entre $\mathbf{Fct}(\I,\E_\kk^f)^{op}$
  et $\mathbf{Fct}(\I^{op},\E_\kk^f)$.
\end{enumerate}  
Ces foncteurs seront appelés {\rm foncteurs de
  dualité}. 
\end{prdef}

Ce résultat entraîne formellement, via la proposition \ref{yoneda}, la propriété suivante.

\begin{prdef}[Cogénérateurs injectifs
  standard]\label{yonedual} Pour tout
  objet $E$ de $\I$, nous noterons $I_E^\I$ 
  l'objet $D'_{\I,\kk}(P_E^{\I^{op}})$ de $\mathbf{Fct}(\I,\E_\kk)$
  (ainsi, on a $I_E^\I(V)=\kk^{{\rm hom}_\I(V,E)}$). Cette construction est
  fonctorielle contravariante en $E$.
\begin{enumerate}\item Il existe un isomorphisme ${\rm
    hom}_{\mathbf{Fct}(\I,\E_\kk)}(F,I_E^\I)\simeq F(E)^*$ naturel en les objets $F$ de $\mathbf{Fct}(\I,\E_\kk)$ et $E$ de $\I$.  
\item Les $I_E^\I$ forment un ensemble de cogénérateurs injectifs lorsque $E$ décrit un squelette de $\I$. On les appelle
  {\em cogénérateurs injectifs standard} de $\mathbf{Fct}(\I,\E_\kk)$.
\end{enumerate}
\end{prdef}

Dans tous les cas que nous considérerons, l'hypothèse
suivante sera vérifiée\,\footnote{Les catégories de foncteurs dont
  la source ne vérifie pas l'hypothèse~\ref{hypf3} ont un comportement
  profondément différent de celles que nous étudierons.}.

\begin{hyp}\label{hypf3} Les ensembles ${\rm hom}_\I (V,W)$ sont finis
  pour tous $V, W\in {\rm Ob}\,\I$.
\end{hyp}

Cette hypothèse assure que les foncteurs projectifs standard, donc
tous les foncteurs de type fini, prennent des valeurs de dimension
finie, de sorte que la dernière assertion de la
proposition~\ref{prdcff} donne un lien très rigide entre
$\mathbf{Fct}(\I,\E_\kk)^{op}$ et $\mathbf{Fct}(\I^{op},\E_\kk)$.

\begin{cor}\label{crctfcf} Lorsque l'hypothèse \ref{hypf3} est
  satisfaite, les injectifs standard de la catégorie
  $\mathbf{Fct}(\I,\E_\kk)$ sont de co-type fini.
\end{cor}

\subsection{Foncteurs hom internes et foncteurs de division}

\begin{prdef}\label{prdf-hid}
\begin{enumerate}\item Pour
    tout objet $X$ de $\mathbf{Fct}(\I,\E_\kk)$, l'endofoncteur
    $\,\cdot\,\otimes X$ de $\mathbf{Fct}(\I,\E_\kk)$ admet un adjoint à
    droite, noté $\mathbf{Hom}_{\mathbf{Fct}(\I,\E_\kk)}(X,.)$ ; on dispose ainsi d'un
    {\em foncteur hom interne}
$$\mathbf{Hom}_{\mathbf{Fct}(\I,\E_\kk)}
:\mathbf{Fct}(\I,\E_\kk)^{op}\times\mathbf{Fct}(\I,\E_\kk)\to\mathbf{Fct}(\I,\E_\kk).$$
On note $\mathbf{Ext}^*_{\mathbf{Fct}(\I,\E_\kk)}$ les foncteurs dérivés droits de ce bifoncteur exact à gauche.
%
\item On a un isomorphisme naturel gradué
\begin{equation}\label{eq-hi}\mathbf{Ext}^*_{\mathbf{Fct}(\I,\E_\kk)}(X,Y)(E)\simeq{\rm
  Ext}^*_{\mathbf{Fct}(\I,\E_\kk)}(P^\I_E\otimes X,Y)\,,\end{equation}
où $E$ est un objet de $\I$ et $X$ et $Y$ sont des objets de
$\mathbf{Fct}(\I,\E_\kk)$.
\item Pour tout objet $A$ de $\mathbf{Fct}(\I,\E_\kk^f)$, l'endofoncteur
    $\,\cdot\,\otimes A$ de $\mathbf{Fct}(\I,\E_\kk)$ admet un adjoint à
    gauche, noté
    $(\,\cdot:A)_{\mathbf{Fct}(\I,\E_\kk)}$ ; on dispose
    ainsi d'un {\em foncteur de division}
$$(\,\cdot : \cdot\,)_{\mathbf{Fct}(\I,\E_\kk)}
:\mathbf{Fct}(\I,\E_\kk)\times\mathbf{Fct}(\I,\E_\kk^f)^{op}\to\mathbf{Fct}(\I,\E_\kk).$$
\end{enumerate} 
Les indices seront omis dans ces notations quand il n'y a pas d'ambiguïté.
\end{prdef}

La terminologie de foncteur {\em hom interne} est standard (cf. \cite{Mac}, ch. VII, §\,7). Le terme de {\em
  foncteur de division} a quant à lui été introduit par Lannes
(cf.~\cite{Lannes}) dans le
cadre des modules instables sur l'algèbre de Steenrod, voisin de celui
des catégories de foncteurs (cf.~\cite{GP4}, §\,3).


\begin{pr}\label{prdffdec} Soient $F$ un endofoncteur de $\I$ et
  $T$ un objet de $\mathbf{Fct}(\I,\mathbf{Ens})$ tels qu'il existe une bijection
$${\rm hom}_\I(F(V),W)\simeq {\rm hom}_\I(V,W)\times T(W)$$
naturelle en les objets $V$ et $W$ de $\I$.

Alors le foncteur $\mathbf{Hom}_{\mathbf{Fct}(\I,\E_\kk)}(\kk[T],.)$ est
isomorphe au foncteur de précomposition $F^*$ ; autrement dit, $F^*$
est adjoint à droite à $\cdot\otimes \kk[T]$.
\end{pr}

\begin{pr}\label{prdffdec2} Soient $F$ un endofoncteur de $\I$ et
  $T$ un objet de $\mathbf{Fct}(\I^{op},\mathbf{Ens}^f)$ tels qu'il existe une bijection
$${\rm hom}_\I(V,F(W))\simeq {\rm hom}_\I(V,W)\times T(V)$$
naturelle en les objets $V$ et $W$ de $\I$.

Alors le foncteur $(\,\cdot : \kk^T)_{\mathbf{Fct}(\I,\E_\kk)}$ est
isomorphe au foncteur de précomposition $F^*$ ; autrement dit, $F^*$
est adjoint à gauche à $\cdot\otimes\kk^T$.
\end{pr}

\subsection{Décomposition scalaire}  On suppose ici que $\I$ est une catégorie {\em $\kk$-linéaire}. 

\begin{nota}\label{not-ds} \'Etant donné un entier naturel $i$, on désigne par
  $\mathbf{Fct}(\I,\E_\kk)_i$ la sous-catégorie pleine de
  $\mathbf{Fct}(\I,\E_\kk)$ formée des foncteurs $F$ tels que pour
  tout $\lambda\in\kk$ et tout objet $E$ de $\I$, on a
  $F(\lambda.id_E)=\lambda^i .id_{F(E)}$ (on convient ici que $0^0=1$).
\end{nota}

Le produit tensoriel induit des foncteurs $\mathbf{Fct}(\I,\E_\kk)_i\times\mathbf{Fct}(\I,\E_\kk)_j\to\mathbf{Fct}(\I,\E_\kk)_{i+j}$.

Comme le groupe cyclique fini $\kk^\times$ est d'ordre premier à la
caractéristique de $\kk$, il existe un isomorphisme de $\kk$-algèbres
$\kk[\kk^\times]\simeq\kk^{q-1}$, et la $\kk$-algèbre du monoïde
multiplicatif sous-jacent à $\kk$ est isomorphe à $\kk^q$. On en déduit la proposition suivante
(cf.~\cite{K1}, §\,3.3).

\begin{prdef}\label{pra-ds} Les inclusions
  $\mathbf{Fct}(\I,\E_\kk)_i\hookrightarrow\mathbf{Fct}(\I,\E_\kk)$
  induisent une équivalence de catégories
$$\mathbf{Fct}(\I,\E_\kk)\simeq\prod_{i=0}^{q-1}\mathbf{Fct}(\I,\E_\kk)_i.$$

Nous noterons $F\simeq\bigoplus_{i=0}^{q-1} F_i$ la décomposition
canonique d'un foncteur $F$ de $\mathbf{Fct}(\I,\E_\kk)$ qu'on en
déduit, où $F_i$ appartient à $\mathbf{Fct}(\I,\E_\kk)_i$. On
l'appelle {\em décomposition scalaire} de $F$.
\end{prdef}

%
%

\subsection{Produit tensoriel extérieur}\label{pap-pte} C'est le foncteur
$$\boxtimes :
\mathbf{Fct}(\I,\E_\kk)\times\mathbf{Fct}(\J,\E_\kk)\to\mathbf{Fct}(\I\times\J,\E_\kk)$$
défini par la composition
$$\mathbf{Fct}(\I,\E_k)\times\mathbf{Fct}(\J,\E_k)\xrightarrow{\pi_\I^*\times\pi_\J^*}\mathbf{Fct}(\I\times\J,\E_k)\times\mathbf{Fct}(\I\times\J,\E_k)\xrightarrow{\otimes}\mathbf{Fct}(\I\times\J,\E_k),$$
où l'on note  $\pi_\I : \I\times\J\to\I$ et $\pi_\J : \I\times\J\to\J$
les foncteurs de projection. Autrement dit, $(F\boxtimes
G)(A,B)=F(A)\otimes G(B)$.

 On a ainsi des isomorphismes canoniques
$P^{\I\times\J}_{(A,B)}\simeq P^\I_A\boxtimes P^\J_B$ ; cette
observation et ses conséquences, ainsi que les propriétés que nous
rappelons ensuite, justifient la convention de notation
suivante.

\begin{nota}\label{ptcat} Nous désignerons  par
  $\mathbf{Fct}(\I,\E_\kk)\otimes\mathbf{Fct}(\J,\E_\kk)$, par abus, la
  catégorie $\mathbf{Fct}(\I\times\J,\E_\kk)$.
\end{nota}

\begin{pr}[Simples de $\mathbf{Fct}(\I\times\J,\E_\kk)$]\label{simple-pte}\begin{enumerate}\item Si $\I$ et $\J$
    vérifient l'hypothèse \ref{hypf3}, pour tout objet simple $S$ de
    $\mathbf{Fct}(\I\times\J,\E_\kk)$, il existe un objet simple $S_1$
    de  $\mathbf{Fct}(\I,\E_\kk)$, un objet simple $S_2$ de
    $\mathbf{Fct}(\I,\E_\kk)$ et un épimorphisme $S_1\boxtimes
    S_2\twoheadrightarrow S$.
\item Soient $S$ un objet simple de $\mathbf{Fct}(\I,\E_\kk)$ et $S'$ un
  objet simple de $\mathbf{Fct}(\J,\E_\kk)$ tel que le corps ${\rm
    End}_{\mathbf{Fct}(\J,\E_\kk)}(S')$ est réduit à $\kk$. Alors
  $S\boxtimes S'$ est un objet simple de $\mathbf{Fct}(\I\times\J,\E_\kk)$.
 \item Soient $S_{1}$, $S_{2}$ deux objets simples de $\mathbf{Fct}(\I,\E_\kk)$ et $S'_{1}$, $S'_{2}$ deux objets simples de $\mathbf{Fct}(\J,\E_\kk)$. Si
 $S_{1}\boxtimes S'_{1}\simeq S_{2}\boxtimes S'_{2}$, alors $S_{1}\simeq S'_{1}$ et $S_{2}\simeq S'_{2}$.
\end{enumerate}
\end{pr}

\begin{cor}\label{gro-pte} Supposons que $\I$ et $\J$ vérifient
  l'hypothèse \ref{hypf3} et que les corps d'endomorphismes des objets
  simples de $\mathbf{Fct}(\J,\E_\kk)$ sont réduits à $\kk$. 

Le produit tensoriel extérieur induit des isomorphismes de groupes abéliens
$$G_0^f(\mathbf{Fct}(\I\times\J,\E_\kk))\simeq
G_0^f(\mathbf{Fct}(\I,\E_\kk))\otimes G_0^f(\mathbf{Fct}(\J,\E_\kk)),$$
$$K_0(\mathbf{Fct}(\I\times\J,\E_\kk))\simeq
K_0(\mathbf{Fct}(\I,\E_\kk))\otimes K_0(\mathbf{Fct}(\J,\E_\kk)).$$
\end{cor}

Ces propriétés se démontrent de façon analogue à celles du produit
tensoriel extérieur usuel en théorie des représentations ---
cf. \cite{CR}, chapitre~$1$, §\,$10$~E. Une démonstration des énoncés
fonctoriels précédents est donnée dans \cite{these}.

\subsection{Recollements}\label{p-prec} Nous terminons cet appendice
avec quelques résultats communs destinés à définir et utiliser commodément le prolongement par zéro
dans un cadre assez général. On les trouvera établis dans \cite{these}.

\begin{defi}\label{df-scc} Soit $\C$ une sous-catégorie pleine de
  $\I$. Nous dirons que $\C$ est~:
\begin{itemize}\item une {\em sous-catégorie relativement connexe} de $\I$
  si pour tous objets $A$, $B$ et $X$ de $\I$ tels que $A, B\in {\rm
    Ob}\,\C$, ${\rm hom} (A,X)\neq\varnothing$ et ${\rm hom}
  (X,B)\neq\varnothing$, on a $X\in {\rm
    Ob}\,\C$ ;
\item une {\em sous-catégorie complète à
    gauche} de $\I$ si pour tout objet $E$ de $\I$, $E$ est objet de
  $\C$ dès que ${\rm hom}_\I (E,X)\neq\varnothing$ pour un objet $X$
  de $\C$ ;
\item une {\em sous-catégorie complète à droite} de $\I$ si $\C^{op}$ est une sous-catégorie complète à
    gauche de $\I^{op}$.
\end{itemize}
\end{defi}

\begin{prdef}[Prolongement par zéro]\label{prfrec} Soit $\C$
  une sous-catégorie pleine relativement connexe de $\I$. On note $\D$ la sous-catégorie pleine de $\I$ dont la classe d'objets est
  le complémentaire de celle de $\C$, et  $\mathcal{R} : \mathbf{Fct}(\I,\E_\kk)\to\mathbf{Fct}(\C,\E_\kk)$ le foncteur de restriction.
\begin{enumerate}\item On définit un foncteur $\mathcal{P}_{\I,\C} : \mathbf{Fct}(\C,\E_\kk)\to\mathbf{Fct}(\I,\E_\kk)$
(noté simplement $\mathcal{P}$ lorsqu'il
n'y a pas d'ambiguïté) appelé {\em prolongement par
  zéro} en posant $\mathcal{P}(F)(E)=F(E)$ ($F\in {\rm Ob}\,\mathbf{Fct}(\C,\E_\kk)$) si $E\in {\rm Ob}\,\C$, $0$ si $E\in
{\rm Ob}\,\D$ ;
$\mathcal{P}(F)(t)=F(t)$ si $t$ est une flèche de $\C$, $0$ si c'est
une autre flèche de $\I$ ; et $\mathcal{P}(T)_E=T_E$ si $T$ est une
flèche de $\mathbf{Fct}(\I,\E_\kk)$ et $E$ un objet de $\C$, $0$ sinon.
\item Supposons $\C$ complète à gauche ; $\D$ est donc
  complète à droite.
\begin{enumerate}\item Le foncteur $\mathcal{P}$ est adjoint à droite
  à $\mathcal{R}$.
\item Le foncteur $\mathcal{P}$ est adjoint à gauche au foncteur
  $\mathcal{N} : \mathbf{Fct}(\I,\E_\kk)\to\mathbf{Fct}(\C,\E_\kk)$ défini
  comme suit :
\begin{enumerate}\item pour tout foncteur $F : \I\to\E_\kk$ et tout objet
  $X$ de $\C$,
$$\N(F)(X)=ker\,\Big(F(X)\xrightarrow{\prod F(f)}\underset{Y\in {\rm Ob}\,\D}{\prod_{X\xrightarrow{f}Y}}F(Y)\Big)\,;$$ 
\item si $F : \I\to\E_\kk$ est un foncteur, pour toute flèche
  $X\xrightarrow{u}X'$ de $\C$, $\N(F)(u)$ est induite par $F(u)$ ;
\item si $T : F\to G$ est une flèche de $\mathbf{Fct}(\C,\E_\kk)$, $\N(T)_X$
  est induite par $T_X$ pour tout objet $X$ de $\C$.
\end{enumerate}
\end{enumerate}
\end{enumerate}
\end{prdef}

On a un énoncé analogue dans le cas d'une sous-catégorie complète à droite.

\begin{cor}\label{crf-rec} On conserve les notations précédentes. Si
  la sous-catégorie $\C$ de $\I$ est complète à gauche, on a un
  diagramme de recollement
$$\xymatrix{\mathbf{Fct}(\C,\E_\kk)\ar[r]|-{\mathcal{P}} &
  \mathbf{Fct}(\I,\E_\kk)\ar[r]|-{\mathcal{R}}\ar@/_/[l]_-{\mathcal{R}}\ar@/^/[l]^-{\mathcal{N}} &
  \mathbf{Fct}(\D,\E_\kk)\ar@/_/[l]_-{\mathcal{P}}\ar@/^/[l].
}$$

Dans le cas où $\C$ est complète à droite, 
on a un diagramme de recollement
$$\xymatrix{\mathbf{Fct}(\C,\E_\kk)\ar[r]|-{\mathcal{P}} &
  \mathbf{Fct}(\I,\E_\kk)\ar[r]|-{\mathcal{R}}\ar@/_/[l]\ar@/^/[l]^-{\mathcal{R}} &
  \mathbf{Fct}(\D,\E_\kk)\ar@/_/[l]\ar@/^/[l]^-{\mathcal{P}}.}$$
\end{cor}

                           
\nocite{*}
\bibliographystyle{smfalpha}
\bibliography{bibart2}

\end{document}